\documentclass[12pt,reqno]{amsart}
\usepackage{mathtools,amssymb,tensor,stackrel,amsmath}
\usepackage{import}
\usepackage[T1]{fontenc}
\usepackage[dvipsnames]{xcolor}
\usepackage{enumitem}
\setitemize{leftmargin=*}
\setenumerate{leftmargin=*}
\setlist[enumerate,1]{label=\textup{(\arabic*)}, 
  ref=\arabic*} 
\setlist[enumerate,2]{label = \textup{(\roman*)},
  ref = \theenumi.\roman*}
\usepackage{setspace}
\usepackage{tikz}
\usetikzlibrary{knots,calc,positioning}
\usetikzlibrary{decorations.markings}
\usepackage{etoolbox, textcomp}
\usepackage{tcolorbox}
\usepackage[colorlinks=true,citecolor=red,linkcolor=blue]{hyperref}
\usepackage{cancel}
\usepackage[normalem]{ulem}

\usepackage[capitalise,noabbrev]{cleveref}
\crefname{enumi}{Part}{Parts}
\Crefname{enumi}{Part}{Parts}

\usepackage{mathtools} 
\usepackage{stmaryrd} 

\numberwithin{equation}{section}

\theoremstyle{plain} 
\newtheorem{theorem}{Theorem}[section]
\newtheorem{lemma}[theorem]{Lemma}
\newtheorem{prop}[theorem]{Proposition}
\newtheorem{cor}[theorem]{Corollary}
\newtheorem{conj}[theorem]{Conjecture}

\newtheorem*{mtheorem*}{Main Theorem}

\theoremstyle{definition}
\newtheorem{defn}[theorem]{Definition}

\newtheorem{remark}[theorem]{Remark}

\newcommand{\vac}{\mathbf{1}}

\Crefname{prop}{Proposition}{Propositions}
\Crefname{cor}{Corollary}{Corollaries}
\Crefname{defn}{Definition}{Definitions}

\usepackage[margin=20mm]{geometry}

\def\CC{\mathbb{C}}

\def\HH{\mathbb{H}}

\def\NN{\mathbb{N}}

\def\RR{\mathbb{R}}

\def\ZZ{\mathbb{Z}}

\def\calJ{\mathcal{J}}

\DeclareMathOperator{\id}{id}

\DeclareMathOperator{\lf}{lf}
\DeclareMathOperator{\opp}{opp}
\DeclareMathOperator{\spn}{span}
\newcommand{\cspn}[1]{\spn_\CC\set*{#1}}

\DeclareMathOperator{\supp}{supp}
\let\hom\relax
\DeclareMathOperator{\hom}{Hom}

\DeclareMathOperator{\ext}{Ext}
\DeclareMathOperator{\aut}{Aut}
\DeclareMathOperator{\res}{Res}
\DeclareMathOperator{\rep}{Rep}
\DeclareMathOperator{\ind}{Ind}
\DeclareMathOperator{\re}{Re}
\DeclareMathOperator{\com}{Com}

\DeclareMathOperator{\im}{im}
\newcommand{\cten}{\otimes} 
\newcommand{\fuse}{\boxtimes} 
\newcommand{\delprod}{\widehat{\boxtimes}} 

\newcommand{\alg}[1]{\mathfrak{#1}} 

\newcommand{\nlie}{\alg{N}} 
\newcommand{\bclie}{\alg{bc}} 
\newcommand{\slt}{\alg{sl}\brac*{2}} 
\newcommand{\slta}{\widehat{\alg{sl}}\brac*{2}}
\newcommand{\Nt}{\mathcal{N}}

\newcommand{\gfin}{\alg{g}} 
\newcommand{\gaff}{\hat{\alg{g}}} 
\DeclareMathOperator{\aff}{aff}

\newcommand{\UEA}[1]{\mathcal{U}\brac*{#1}} 

\let\originalleft\left     
\let\originalright\right
\renewcommand{\left}{\mathopen{}\mathclose\bgroup\originalleft}
\renewcommand{\right}{\aftergroup\egroup\originalright}

\DeclarePairedDelimiter{\brac}{\lparen}{\rparen} 
\DeclarePairedDelimiter{\sqbrac}{\lbrack}{\rbrack} 
\DeclarePairedDelimiter{\set}{\lbrace}{\rbrace}
\newcommand{\st}{\mspace{5mu} \vert \mspace{5mu}} 
\DeclarePairedDelimiter{\abs}{\lvert}{\rvert}

\DeclarePairedDelimiter{\ang}{\langle}{\rangle}
\DeclarePairedDelimiter{\normord}{{} :}{: {}} 
\DeclarePairedDelimiter{\powser}{\llbracket}{\rrbracket} 

\DeclarePairedDelimiterX{\comm}[2]{\lbrack}{\rbrack}{#1 , #2}  
\DeclarePairedDelimiterX{\acomm}[2]{\lbrace}{\rbrace}{#1 , #2} 
\DeclarePairedDelimiterX{\super}[2]{\lparen}{\rparen}{#1 \delimsize\vert \mathopen{} #2} 

\newcommand{\killing}[2]{\kappa \brac[\big]{#1 , #2}} 
\newcommand{\pair}[2]{\ang*{#1,#2}} 

\newcommand{\usl}[1]{\mathsf{A}_1\brac*{#1}} 
\newcommand{\mmsl}[2]{\mathsf{A}_{1}\brac*{#1,#2}} 
\newcommand{\slE}[3]{E_{#1 ; #2,#3}} 
\newcommand{\slD}[3]{D_{#1,#2}^{#3}} 
\newcommand{\slEr}[3]{E_{#1,#2}^{#3}} 
\newcommand{\heis}[1]{\mathsf{H}\brac*{#1}} 
\newcommand{\lat}[1]{\mathsf{L}\brac*{#1}} 
\newcommand{\bc}{\mathsf{BC}} 
\newcommand{\uN}[1]{\mathsf{N}\brac*{#1}} 
\newcommand{\mmN}[2]{\mathsf{N}\brac*{#1,#2}} 

\newcommand{\NM}[3]{\mathcal{M}_{#1, #2,#3}} 
\newcommand{\NL}[3]{\mathcal{L}_{#1 , #2,#3}} 
\newcommand{\NLs}[3]{\mathcal{L}_{#1 ; #2,#3}} 
\newcommand{\smod}[1]{#1\text{-}\mathsf{mod}^{\text{smth}}}
\newcommand{\wtmod}[1]{#1\text{-}\mathsf{mod}^{\text{wt}}} 
\newcommand{\pwtmod}[1]{#1\text{-}\mathsf{mod}_{\ge0}^{\text{wt}}} 
\newcommand{\svec}{\mathsf{sVec}}
\newcommand{\mmV}[2]{\mathsf{M}\brac*{#1,#2}} 
\newcommand{\mmVM}[2]{\mathsf{S}_{#1,#2}} 

\newcommand{\scr}{\mathcal{Q}} 

\newcommand{\cft}{conformal field theory}

\newcommand{\va}{vertex algebra}

\newcommand{\voa}{vertex operator algebra}

\newcommand{\vosa}{vertex operator superalgebra}
\newcommand{\ope}{operator product expansion}

\newcommand{\rhs}{right-hand side}

\newcommand{\dd}{\mathrm{d}}
\newcommand{\ii}{\mathsf{i}}
\newcommand{\ee}{\mathsf{e}}
\newcommand{\bee}{\mathbf{\ee}}

\newcommand\quash[1]{}

\makeatletter
\newlength{\@pxlwd} \newlength{\@rulewd} \newlength{\@pxlht}
\catcode`.=\active \catcode`B=\active \catcode`:=\active
\catcode`|=\active
\def\sprite#1(#2,#3)[#4,#5]{
   \edef\@sprbox{\expandafter\@cdr\string#1\@nil @box}
   \expandafter\newsavebox\csname\@sprbox\endcsname
   \edef#1{\expandafter\usebox\csname\@sprbox\endcsname}
   \expandafter\setbox\csname\@sprbox\endcsname =\hbox\bgroup
   \vbox\bgroup
  \catcode`.=\active\catcode`B=\active\catcode`:=\active\catcode`|=\active
      \@pxlwd=#4 \divide\@pxlwd by #3 \@rulewd=\@pxlwd
      \@pxlht=#5 \divide\@pxlht by #2
      \def .{\hskip \@pxlwd \ignorespaces}
      \def B{\@ifnextchar B{\advance\@rulewd by \@pxlwd}{\vrule
         height \@pxlht width \@rulewd depth 0 pt \@rulewd=\@pxlwd}}
      \def :{\hbox\bgroup\vrule height \@pxlht width 0pt depth
0pt\ignorespaces}
      \def |{\vrule height \@pxlht width 0pt depth 0pt\egroup
         \prevdepth= -1000 pt}
   }
\def\endsprite{\egroup\egroup}
\catcode`.=12 \catcode`B=11 \catcode`:=12 \catcode`|=12\relax
\makeatother

\sprite{\FormOfHboxtr}(25,25)[0.5 em, 1.2 ex]

:BBBBBBBBBBBBBBBBBBBBBBBBB | :BB......................B |
:B.B.....................B | :B..B....................B |
:B...B...................B | :B....B..................B |
:B.....B.................B | :B......B................B |
:B.......B...............B | :B........B..............B |
:B.........B.............B | :B..........B............B |
:B...........B...........B | :B............B..........B |
:B.............B.........B | :B..............B........B |
:B...............B.......B | :B................B......B |
:B.................B.....B | :B..................B....B |
:B...................B...B | :B....................B..B |
:B.....................B.B | :B......................BB |
:BBBBBBBBBBBBBBBBBBBBBBBBB |

\endsprite
\sprite{\FormOfShboxtr}(25,25)[0.3 em, 0.72 ex]

:BBBBBBBBBBBBBBBBBBBBBBBBB | :BB......................B |
:B.B.....................B | :B..B....................B |
:B...B...................B | :B....B..................B |
:B.....B.................B | :B......B................B |
:B.......B...............B | :B........B..............B |
:B.........B.............B | :B..........B............B |
:B...........B...........B | :B............B..........B |
:B.............B.........B | :B..............B........B |
:B...............B.......B | :B................B......B |
:B.................B.....B | :B..................B....B |
:B...................B...B | :B....................B..B |
:B.....................B.B | :B......................BB |
:BBBBBBBBBBBBBBBBBBBBBBBBB |

\endsprite

\newcommand{\ket}[1]{\left|#1\right\rangle}
\newcommand{\kket}[1]{\left\Vert #1\right\rangle}
\newcommand{\bra}[1]{\left\langle #1\right|}
\newcommand{\bbra}[1]{\left\langle #1\right\Vert}

\makeatletter
\renewcommand\author@andify{%
  \nxandlist {\unskip ,\penalty-1 \space\ignorespaces}%
    {\unskip {} \@@and~}%
    {\unskip \penalty-2 \space \@@and~}%
}
\makeatother

\allowdisplaybreaks

\begin{document}

\title[Fusion and rigidity for admissible \(\slt\) and \(\Nt=2\)]{Fusion rules and rigidity for weight modules over the simple admissible affine
  \(\slt\) and \(\Nt=2\) superconformal vertex operator superalgebras}

\author[H.~Nakano]{Hiromu Nakano}
\address[Hiromu Nakano]{
  Advanced Mathematical Institute\\
  Osaka Metropolitan University\\
  Osaka, 558-8585, Japan.
 }
\email{d23155u@omu.ac.jp}

\author[F.~Orosz Hunziker]{Florencia Orosz Hunziker}
\address[Florencia Orosz Hunziker]{
  Department of Mathematics\\
  University of Colorado Boulder\\
  Boulder, CO 80309-0395, United States of America.}
\email{florencia.orosz@colorado.edu}

\author[A.~Ros Camacho]{Ana Ros Camacho}
\address[Ana Ros Camacho]{
  Facultat de Ciències Matemàtiques, Universitat de València \\
  Avinguda Vicent Andrés Estellés 19 \\
  46100 Burjassot, València, Spain.}
\email{ana.ros.camacho@uv.es}

\author[S.~Wood]{Simon Wood}
\address[Simon Wood]{
  School of Mathematics \\
  Cardiff University \\
  Cardiff, CF24 4AG, United Kingdom\\
  and\\
  Research Institute for Mathematical Sciences\\
  Kyoto University\\
  Kyoto, 606-8502, Japan.
}
\email{woodsi@cardiff.ac.uk}


\begin{abstract}
  We prove that the categories of weight modules over the simple
  \(\slt\) and \(\Nt=2\) superconformal \vosa{s} at fractional
  admissible levels and central charges are rigid (and hence the
    categories of weight modules are braided ribbon categories) and that
  the decomposition formulae of fusion products of simple projective modules
  conjectured by Thomas Creutzig, David Ridout and collaborators hold
  (including when the decomposition involves summands that are
    indecomposable yet not simple). In
  addition to solving this old open problem, we develop new techniques for the
  construction of intertwining operators by means of integrating screening
  currents over certain cycles, which are expected to be of independent
  interest, due to their applicability to many other algebras. In the example
  of \(\slt\) these new 
  techniques allow us to give explicit formulae for a logarithmic intertwining
  operator from a pair of simple projective modules to the projective cover of
  the tensor unit, namely,  the \voa{} as a module over itself.
\end{abstract}

\maketitle

\section{Introduction}

Within the body of literature on \cft{} and \voa{s} the class of rational
theories (satisfying a number of technical niceness conditions including
\(C_2\)-cofiniteness and the category of admissible modules being
semisimple) stand out as being exceptionally intensively studied and well
understood. One particularly compelling aspect of these rational theories is
the abundance of rich mathematical structures they exhibit. For example, the fact that their
categories 
of modules are modular tensor
categories and that the categorical action of
the modular group via Hopf links and twists matches the modular
transformations of characters \cite{HuaVer08}. Apart from their intrinsic
mathematical beauty these modularity results also have practical implications: they
allow the efficient computation of tensor products (also
called fusion) of \voa{} modules  in terms of the modular transformation
formulae of characters. In practice, this provides an enormous reduction in
the computational effort required to understand fusion products.

There is ample evidence to suggest that the modularity properties enjoyed by
rational theories generalise to large classes of non-rational ones. The
purpose of this paper is to prove that these conjectured modularity properties
do indeed hold for so-called admissible \(\slt\) and \(\Nt=2\)
superconformal theories. In order to  better understand these conjectures, we begin with a historical overview.

A first hint of modularity beyond rationality was discovered by Kac and Wakimoto \cite{KM1,KM2,KM3} when they computed modular
transformation formulae for characters of simple
highest weight modules over affine simple Lie algebras at admissible levels and
weights. Since the non-negative integral levels (which are the only ones
giving rise to rational theories) are a proper subset of all admissible
levels, it is somewhat surprising for such formulae to exist at all. Shortly
afterwards Koh and Sorba \cite{KoSo} plugged these modular transformation properties
for \(\slt\) into the Verlinde formula, which, rather startlingly, predicted negative
multiplicities for some summands appearing in certain fusion products. This in
turn led to concern within the academic community that these
non-integral admissible (also called fractional admissible) theories may
suffer from some intrinsic ``sickness'' \cite{CFTbook}.

New light was later shed on this
riddle by Ridout after a careful analysis of \(\slt\) at level
\(k=-\frac{1}{2}\) \cite{Ridk12} led him to note that the characters of
highest weight modules at this level were only convergent in certain domains
and that the modular transformation formulae of Kac and Wakimoto hold only for
the analytic continuations of characters rather than their series
expansions. This is because these modular transformations do not preserve
domains of convergence. He also noted that in the category of weight
  modules the analytic continuation of characters of certain highest weight
  modules was the negative of the analytic continuation of characters of
  certain other weight modules. This gave a first hint as to why signs were
  appearing in the Verlinde formula of Koh and Sorba.

Therefore, the larger category of weight modules, which contains all highest weight modules but also contains relaxed highest weight modules,  was proposed as the appropriate category to study modularity in the non-rational setting. Creutzig and Ridout started demonstrating the importance of this larger category by first studying weight modules over
  affine \(\alg{gl}\brac*{1|1}\) and the modular properties of its characters \cite{glCR}.
In that case study, they showed that while category \(\mathcal{O}\) is semisimple
and finite, the category of weight modules is neither. 
They checked, however, that the span of
characters (taken as specific series expansions rather than their analytic
continuations) 
of a distinguished class of weight modules (called \emph{standard modules})
carries an action of the modular group. Further, they demonstrated that
evaluating the Verlinde formula (in a generalised version conjectured to hold for infinite categories
of modules) using this action predicts non-negative fusion
multiplicities. This new generalised Verlinde formula shed light on
the riddle of negative fusion multiplicities of Koh and Sorba: there are linear
relations between the analytic continuations of characters of simple weight
modules. So a negative multiplicity in the Verlinde formula for category
\(\mathcal{O}\) can be interpreted as a positive multiplicity of a different
weight module outside of category \(\mathcal{O}\).
This work was then generalised to \(\slt\) at all admissible levels
  \cite{CR1,CR2}, where the same pattern holds; namely, the category of weight modules is neither semisimple nor finite and there is a generalised Verlinde formula in terms of standard modules. Creutzig and Ridout together with other authors later showed that these properties hold for the category of weight modules for many other families of algebras
\cite{CMtheta,RWbghosts,RWtrip,CRRn11,CRRn12,CMSqdim,ACRpar,RSWosp,KRWsl3,FRbp,FRRsl3}.

Two key conjectures or hopes that crystallised from Creutzig and Ridout's work
in \cite{glCR,CR1,CR2} were that the category of weight modules is rigid
(this implies that the fusion product is exact, a necessary condition for
anything akin to a Verlinde formula being well-defined) and that the fusion
product decomposition 
formulae predicted by the generalised Verlinde formulae of Creutzig and Ridout
are true. These conjectures were proved for levels \(k=-\frac12,-\frac43\) in \cite{CMY2} and our work here proves them for all  admissible levels in general.

This paper makes use of a coset
realisation of the \(\Nt=2\) superconformal algebra in terms of
\(\slt\) and a fermionic ghost (or \(bc\)) system \cite{KZ1,KZ2,CL2,EG}. Because of this the
representation theories of \(\Nt=2\) and \(\slt\) are closely 
intertwined and any statement about the representation theory for one algebra
has a corresponding version for the other.
In particular, rigidity holds for \(\slt\) if and only if it holds for \(\mathcal{N}=2\). 
  Our proofs rely crucially on this connection as certain statements are easier to prove for \(\slt\) and others for \(\mathcal{N}=2\). 

The description of the logarithmic tensor structures for categories of modules
over Virasoro \va{s} and their $\mathcal{N}=1$ and $\mathcal{N}=2$ super
extensions is an important problem that has been intensively studied in the
literature in the past years \cite{CJOHRY, MYc=25, Myc=1p, CMOHYN1, MS,
  CLRW}. Since the monoidal structure on categories of \voa{} modules
automatically comes with a braiding and a twist that is balanced with respect
to this braiding, the rigidity of these categories implies that they are
braided ribbon categories.
The fusion rules for the  strongly rational theories arising from unitary
$\mathcal{N}=2$ minimal models were studied by Adamovi\'c in \cite{A2} while
the tensor category structure for their module categories was established by
Huang and Milas in \cite{HMII} using the theory of Lepowsky and Huang
developed in \cite{HL}. In the Virasoro and $\mathcal{N}=1$ cases, the central
charges in which the universal vertex superalgebras are non-simple are exactly
the central charges in which the simple quotients are strongly rational vertex
operator superalgebras. Interestingly, the $\mathcal{N}=2$ vertex
superalgebras have richer behavior since the central charges in which the
universal $\mathcal{N}=2$ vertex superalgebra admits a strongly rational
simple quotient vertex superalgebra form a proper subset of the central
charges in which the universal vertex superalgebra is  non-simple. While the
logarithmic tensor product theory developed by Huang, Lepowsky and Zhang in
\cite{HLZ} has recently been used to establish the tensor category structure
on natural module categories for the universal Virasoro \cite{CJOHRY} and
$\mathcal{N}=1$ vertex superalgebras \cite{CMOHYN1}, in the case of the
$\mathcal{N}=2$ vertex superalgebras the logarithmic tensor structure has not
been applied to the universal algebra, but only to the simple fractional
minimal model central charges, namely, to the intermediate family of central
charges in which the universal $\mathcal{N}=2$ vertex superalgebra admits a
non-trivial ideal such that its simple quotient yields an irrational
representation theory \cite{Csl2ten}. 

This paper is organised as follows. In \cref{sec:algebras} we review all the
\vosa{s} that will be needed as well as their corresponding categories of
weight modules. In particular, we recall some \(\slt\) fusion product
decompositions, which will be needed later, that have already been proved and
state the conjectured ones (\cref{thm:EEfusionconj}) which we shall prove in later sections.

In \cref{sec:coset} we review the coset construction of the \(\Nt=2\)
superconformal algebra in terms of affine \(\slt\) and how the fusion products
of these two algebras are interrelated. We give a bijection between spaces of intertwining operators in \cref{thm:homspaceeq}, we give 
a sufficient condition for transporting rigidity results from the
\(\slt\) side to the \(\Nt=2\) side in
\cref{thm:1121rules}, and we give a sufficient condition
for the semisimple fusion products of \cref{thm:EEfusionconj} to hold in terms
of the dimensions of certain spaces of intertwining operators.

\Cref{sec:upperbounds,sec:lowerbounds} are then dedicated to proving
\cref{thm:1121rules} by proving upper and lower bounds, and observing that
these bounds are equal. The upper bounds are computed in
\cref{sec:upperbounds} by use of the Zhu algebra formalism on the
\(\Nt=2\) side. The lower bounds are computed in
\cref{sec:lowerbounds} on the \(\slt\) side by means of a free field
realisation and screening operators. The appearance of screening operators and
their associated integrals over intricate cycles necessitates the use of
\(P(w)\) rather than logarithmic intertwining operators, that is, intertwining
operators where the variable is considered as a complex number rather than a
formal variable.

The non-semisimple fusion products of \cref{thm:EEfusionconj} are tackled in
\cref{sec:intoprigid,sec:nonssifusion}. Specifically, the projective cover of
the tensor unit for \(\slt\) is constructed in \cref{sec:pcover} using the free field realisation
of \cref{sec:lowerbounds}. In \cref{sec:logintops} this
projective module is used to analyse
one of the fusion products appearing in \cref{thm:EEfusionconj} and to construct the
evaluation morphisms for a certain family of simple projective modules,
\(\slE{\mu}{1}{1}\), that are needed to prove their rigidity. In
\cref{sec:logintops} we prove that the \(\slE{\mu}{1}{1}\) modules are rigid.
Finally, in \cref{sec:nonssifusion} the results of \cref{sec:intoprigid} are
combined to prove that the categories of weight modules over simple admissible
affine 
\(\slt\) and \(\Nt=2\) are rigid, as well as that the non-semisimple fusion
product decomposition of \cref{thm:EEfusionconj} hold.

\subsection*{Acknowledgements}
FOH's work is supported by the US National Science Foundation
under grant No. DMS-2102786. 
ARC's work is supported by the London Mathematical
Society through the Emmy Noether Fellowship REF EN-2324-02. 
SW's work is supported by the Engineering
and Physical Sciences Research Council (EPSRC) grant EP/V053787/1 and by the
Alexander von Humboldt Foundation. 
The authors would like to thank T. Creutzig for helpful discussions.

As this manuscript was nearing completion, the authors became aware of other
work
\cite{CMYrigid} also proving the rigidity of simple admissible affine \(\slt\)
and
\cite{Cfusion} proving generalised Verlinde formulae with the
Verlinde formula for $\slt$ at admissible levels as an example.
Since the methods in these other works are completely different from those
used here, it was agreed to complete all manuscripts independently
and submit them to arXiv simultaneously.

\section{Algebras and categories of modules}
\label{sec:algebras}
We review and fix notation for the four families of vertex algebras that will be considered
repeatedly below.

\subsection{Affine vertex algebras}
Let $\gfin$ be a simple (possibly abelian) finite-dimensional complex Lie algebra with invariant
symmetric non-degenerate bilinear form $\kappa$ and choice of Cartan
subalgebra \(\alg{h}\). Consider the affinisation of
\(\gfin\),
\begin{equation}
  \gaff=\gfin\cten\CC[t,t^{-1}]\oplus \CC K,
\end{equation}
where \(K\) is central and \(x\cten f(t), y\cten g(t)\in
\gfin\cten\CC[t,t^{-1}]\) satisfy the commutation relation
\begin{equation}
  \comm{x\cten f(t)}{y\cten g(t)}=[x,y]\cten f(t)g(t)+
  \killing{x}{y}\res_t f(t)dg(t) K.
\end{equation}
We do not include the degree operator in the construction of \(\gaff\) because
we will always identify it with the negative for the Virasoro \(L_0\) operator
obtained from the Sugawara construction.
Denote \(x_n=x\cten t^n\) for \(x\in \gfin\), \(n\in \ZZ\) and define
\begin{equation}
  \gaff_{\pm}=\cspn{x_{\pm n}\st x\in \gfin,\ n\ge1},\quad
  \gaff_0=\gfin\oplus \CC K,\quad
  \gaff_{\ge0}=\gaff_0\oplus \gaff_+,
\end{equation}
to obtain the triangular decomposition \(\gaff=\gaff_{-}\oplus \gaff_{0}\oplus
\gaff_{+}\).
\begin{defn}
  Let \(M\) be a \(\gaff\)-module on which the central element \(K\) acts as
  \(k\cdot\id\) for some \(k\in\CC\).
  \begin{enumerate}
  \item The module \(M\) is called \emph{smooth}, if for every \(m\in M\)
    \begin{equation}
      x_n m =0 \ \ \forall x\in \gfin,\ n\gg0.
    \end{equation}
    Denote the full subcategory of the category of all \(\gaff\)-modules whose
    objects are all smooth modules by \(\smod{\gaff_k}\).
  \item The module \(M\) is called \emph{weight}, if it is smooth, finitely
    generated, simultaneously graded by generalized $L_0$ eigenvalues and
    \(\alg{h}\) eigenvalues, and all homogeneous spaces are finite-dimensional. That is, \(M\) 
    decomposes into a direct sum
    \begin{equation}
      M=\bigoplus_{\substack{\ \lambda\in \alg{h}^\ast\\ h\in\CC}}
      M_{\lambda, h},\qquad M_{\lambda,h}=\set{m\in M\st
        (x_0-\lambda(x))m=0=(L_0-h)^Nm,\ N\gg0,  \ x \in \alg{h}},
    \end{equation}
    where \(\dim M_{\lambda,h}<\infty\) for all \(\lambda\in \alg{h}^\ast\), \(h\in\CC\) and for fixed \(\lambda\in
    \alg{h}^\ast\), \(M_{\lambda,h}=0\) for \(\re(h)\ll 0\).

    Denote the full subcategory of \(\smod{\gaff_k}\) whose objects are all
    weight modules by \(\wtmod{\gaff_k}\).
  \item A weight module \(M\) is called \emph{positive energy}, if  there exists a real number \(h_{\min}\) such that for all \(\lambda\in
    \alg{h}^\ast\) and all \(h\in \CC\) satisfying \(\re(h)<h_{\min}\) we have
    \(M_{\lambda,h}=0\).
    Denote the full subcategory of \(\wtmod{\gaff_k}\) whose objects are all
    positive energy modules by \(\pwtmod{\gaff_k}\).   
  \item Let \(M\) be weight. A non-zero homogeneous vector \(v\in M\) is
    called \emph{relaxed  highest weight}, if it is
    homogeneous and \(\gaff_+ v=0\). If \(v\in M\) is
    relaxed  highest weight and generates \(M\), then
    \(M\) is called a \emph{relaxed highest weight module}. 
  \end{enumerate}
\end{defn}

We can now easily construct examples of smooth \(\gaff\)-modules. For example,
generalised Verma modules and their simple quotients are smooth. These are constructed
as follows. Let \(\overline{M}\) be a \(\gfin\) weight module, which becomes a
\(\gaff_{\geq0}\)-module by defining \(K\) to act as \(k\cdot\id\), \(k\in\CC\)
and \(\gaff_+\) to act trivially.
The generalised Verma module $\mathfrak{V}(k,\overline{M})$, is the induced module
\begin{align}
  \mathfrak{V}(k,\overline{M})= \ind_{\gaff_{\ge0}}^{\gaff}\overline{M}.
  \label{eq:affverma}
\end{align}
If \(\overline{M}\) is simple, then \(\mathfrak{V}(k,\overline{M})\) has a unique maximal submodule hence a unique
simple quotient
\begin{align}
  \mathfrak{L}(k,\overline{M})=\frac{\mathfrak{V}(k,\overline{M})}{\ang{\text{maximal submodule}}}.    
\end{align}
\begin{prop}[Frenkel, Zhu \cite{FZ}]
  Let \(\gfin\) be simple (possibly abelian), \(k\in \CC\) and let \(\CC\) be the trivial \(\gfin\)-module. Then the
  parabolic Verma module
  \(\mathfrak{V}(k,\CC)\) admits the
  structure of a \voa{} by defining the field map on \(x_{-1}\vac\), \(x\in
  \gfin\) to be
  \begin{align}
    Y(x_{-1}\vac,z)=x(z)=\sum_{n\in\ZZ} x_n z^{-n-1}
  \end{align}
  and continuing by derivatives and normal ordering. This \voa{} is called
  the \emph{universal affine \(\gfin\) \voa{} at level \(k\)}. Any choice of
  basis of \(\gfin\) corresponds to a set of strong generators and any two
  such generators \(x,y\in\gfin\)
  satisfy the \ope{} relations
  \begin{align}
    x(z)y(w)\sim \frac{\killing{x}{y}k}{(z-w)^2}+\frac{\comm{x}{y}(w)}{z-w}.
  \end{align}
  There is a distinguished choice of conformal vector, called the
  \emph{Sugawara vector} 
  \begin{equation}
    \omega_k=\frac{1}{2(k+h^\vee)}\sum_i (x^i)_{-1}(y^i)_{-1}\vac,
  \end{equation}
  of central charge \(c=\frac{k \dim
    \gfin}{k+h^\vee}\), where \(h^\vee\) is the dual Coxeter number (defined to be \(0\), if
  \(\gfin\) is abelian) \(\set{x^i}\) is any choice of basis of \(\gfin\) and \(\set{y^i}\) is
  its dual with respect to \(\kappa\). 
  Finally, \(\smod{\gaff_k}\) is the category of all \(\mathfrak{V}(k,\CC)\)-modules.
  \label{thm:affinevoaconstr}
\end{prop}

\subsubsection{The Heisenberg vertex algebra}
Let \(\CC\) be the trivial Lie algebra with non-degenerate invariant bilinear
form characterised by \(\killing{1}{1}=2\).
The rank one Heisenberg vertex algebra $\heis{t}$ 
is the affine \voa{} constructed by applying \cref{thm:affinevoaconstr} at
level \(t\in \CC^\times\) to the trivial Lie algebra \(\CC\). Higher rank
Heisenberg \va{s} are constructed by taking tensor products of the 
rank one Heisenberg \va{s}. Note also that any two Heisenberg \va{s} of equal rank are
isomorphic regardless of the choice of non-zero level.

The (parabolic) Verma modules for the Heisenberg algebra, constructed from weight modules over the abelian Lie algebra \(\CC\), are just Fock
spaces. In particular, we denote by $\mathcal{F}_p$ the Fock space of highest
weight $p \in \CC$ (as a module of the Heisenberg Lie algebra
$\heis{t}=\mathcal{F}_0$).
\begin{prop}
  Let \(t\in\RR^\times\). As a linear abelian category \(\wtmod{\hat\CC_t}\) is semisimple with Fock spaces 
  forming a complete set of representatives for all simple isomorphism
  classes. Let \(\wtmod{\hat\CC_t}_\RR\) be the full subcategory generated by Fock spaces with real weights. Then \(\heis{t}\) furnishes \(\wtmod{\hat\CC_t}_\RR\)
  with the structure of a rigid braided monoidal category, which is braided
  equivalent to the category of finite-dimensional \(\RR\)-graded vector spaces, with \(q(\alpha)=\ee^{\ii \alpha^2},\ \alpha\in\RR\) as choice of quadratic form.
\end{prop}

The linear category structure of \(\wtmod{\hat\CC_t}\) (in particular the
semisimplicity and classification of simples) is due to \cite[Prop
3.6]{LepRam82}, the fusion rules are due to \cite{DoLGVA93} and the remainder of the monoidal structure to \cite{CKLR}.
See \cite[Prop 3.11, Thm 3.12]{ALSW} for a discussion of the equivalence between \(\wtmod{\hat\CC_t}_{\RR}\) and \(\RR\)-graded vector spaces.

\subsubsection{The $\slt$ vertex operator algebras}

Consider the smallest non-abelian simple complex Lie algebra
\(\slt=\cspn{e,h,f}\) with choice of Cartan subalgebra \(\alg{h}=\CC h\). We
spell out the well known relations and
normalisations of \(\slt\) to fix notation. The non-vanishing commutation relations are
\begin{equation}
  \comm{e}{f}=h,\qquad \comm{h}{e}=2e,\qquad \comm{h}{f}=-2f.
\end{equation}
We normalise the Killing form such that its non-vanishing pairings are
\begin{equation}
  \killing{h}{h}=2,\qquad \killing{e}{f}=\killing{f}{e}=1
\end{equation}
Recall that the dual Coxeter number of \(\slt\) is \(h^\vee=2\).
The vertex operator algebra $\usl{k}$ constructed from \(\slt\) by applying
\cref{thm:affinevoaconstr} at level \(k\in \CC\setminus\set{-2}\) is called the
\emph{universal \(\slt\) \voa{} at level \(k\)}.

Recall that \(\slta\) admits a family of automorphisms called
\emph{spectral flow}.
\begin{defn}\leavevmode
  \begin{enumerate}
  \item The {\it spectral flow} automorphisms $\sigma^\ell, \ell\in \ZZ$, of \(\slta\)
    are determined on basis vectors and the Virasoro operators constructed by
    the Sugawara conformal vector by
    \begin{align}
      \sigma^\ell(e_n)&=e_{n-\ell},&
      \sigma^\ell(h_n)&=h_n-\ell\delta_{n,0}K,&
      \sigma^\ell(f_n)&=f_{n+\ell},\nonumber\\
      \sigma^\ell(L_n)&=L_n-\frac12\ell h_n+\delta_{n,0}\frac14\ell^2K,&
      \sigma^\ell(K)&=K.
    \end{align}
  \item Let \(g\in \aut\brac*{\slta}\) and \(M\) an \(\slta\)-module, then
    the \(g\)-twist of \(M\) is an \(\slta\)-module whose underlying vectors
    space is that of \(M\) be with the action
    \begin{equation}
      x\cdot_g m=g^{-1}(x)\cdot m,\qquad \forall x\in \slta,\ \forall m\in M,
    \end{equation}
    where \(\cdot\) is the original action of \(\slta\) on \(M\).
  \end{enumerate} 
  \label{def:sl2spectralflow}
\end{defn}
Note that the categories \(\smod{{\slta}_k}\) and \(\wtmod{{\slta}_k}\) are both
closed under spectral flow automorphisms.

\begin{prop}[Bernard-Felder   \cite{BerFel90}] 
  The \voa{} $\usl{k}$ admits a non-trivial proper ideal
  if and only if there exist positive integers \(u,v\in \ZZ_{>0}\), \(\gcd\brac*{u,v}=1\), \(u\ge2\)
  such that \(k+2=\frac{u}{v}\). The ideals at these levels are unique and
  hence maximal. The simple quotient of $\usl{\frac{u}{v}-2}$ will be denoted
  \(\mmsl{u}{v}\) and is called the \(\slt\) \((u,v)\)-minimal model.
\end{prop}

For \(u\ge2\), \(\mmsl{u}{1}\) is the much studied \(\alg{su}_2\)
Wess-Zumino-Witten model at non-negative integral level. Here we shall
primarily be interested in \(\mmsl{u}{v}\) with \(v\ge2\).
Modules over \(\mmsl{u}{v}\) are of course naturally identified with modules over
\(\usl{\frac{u}{v}-2}\) on which the maximal ideal acts trivially, which leads
to the following definition.
\begin{defn}
  Let \(\smod{\mmsl{u}{v}}\) be the full subcategory of
  \(\smod{\usl{\frac{u}{v}-2}}\) whose objects are all smooth modules on which
  the maximal ideal of \(\usl{\frac{u}{v}-2}\) acts trivially. 
  
  The category of weight modules \(\wtmod{\mmsl{u}{v}}\) is defined to be
    the full subcategory of all weight modules in
    \(\smod{\mmsl{u}{v}}\) for which the eigenvalues of the Cartan generator
    \(h_0\) are real. The category of positive energy modules
    \(\pwtmod{\mmsl{u}{v}}\)  is the full subcategory of
    all positive energy modules in \(\wtmod{\mmsl{u}{v}}\).
\end{defn}

\begin{prop}
  For \(r,s\in\ZZ\) and \(t\in \CC^\ast\), let
  \begin{equation}
    \lambda_{r,s,t}=r-1-ts,\qquad \Delta^{\aff}_{r,s,t}=\frac{\brac*{r-st}^2-1}{4t}.
  \end{equation}\label{eq:affineDelta}
    \begin{enumerate}
    \item \emph{(Block \cite{Blosl279}).} Every simple weight module over \(\slt\)
      (modules on which $h$ acts semisimply and each weight space is finite-dimensional) is isomorphic to one of the 
      following mutually inequivalent modules. 
    \begin{enumerate}
    \item \(\overline{L(\mu)}\), \(\mu\in \ZZ_{\ge0}\), the simple finite-dimensional module
      of highest weight \(\mu\) and dimension \(\mu+1\).
    \item \(\overline{D^+(\mu)}\), \(\mu \in \CC\setminus \ZZ_{\ge0}\), the Verma module of
      highest weight \(\mu\).
    \item \(\overline{D^-(\mu)}\cong \overline{D^+(-\mu)}^\ast\), \(\mu \in \CC\setminus\ZZ_{\le 0}\),
      the lowest weight modules which are the duals of the Verma modules above
      (they are also Verma modules for the choice of Borel \(\cspn{h,f}\)).
    \item \(\overline{E(\mu,q)}\), \(\mu\in \CC/2\ZZ\), \(q\in
      \CC\setminus\set{\frac{\nu}{2}\brac*{\nu+2}\st \nu\in \mu}\),
      the dense module with weight support \(\mu\) on which the quadratic Casimir
      \( \frac12 h^2 +ef +fe\)
      acts as \(q\cdot\id\). 
    \end{enumerate}
  \item \emph{(Adamovi\'c, Milas \cite{AM}).} For coprime \(u,v\in \ZZ_{>0}\), \(u\ge 2\), every simple relaxed highest weight module over \(\mmsl{u}{v}\)
    (equivalently every simple object in \(\pwtmod{\mmsl{u}{v}}\)) is
    isomorphic to one of the following mutually inequivalent simple quotients
    of an induction of a simple \(\slt\) module.
    \begin{enumerate}
    \item \(L_{r}=\mathfrak{L}\left(\frac{u}{v}-2,\overline{L(r-1)}\right)\), \(r\in \set{1,\dots,u-1}\).
    \item
      \(D^+_{r,s}=\mathfrak{L}\left(\frac{u}{v}-2,\overline{D^+(\lambda_{r,s,\frac{u}{v}})}\right)\),
      \(1\le r\le u-1\), \(1\le s\le v-1\).
    \item \(D^-_{r,s}=D^{+\,\prime}_{r,s}\), \(1\le r\le u-1\), \(1\le s\le v-1\).
    \item
      \(\slE{\mu}{r}{s}=\mathfrak{L}\left(\frac{u}{v}-2,\overline{E(\mu,2\frac{u}{v}\Delta^{\aff}_{r,s,\frac{u}{v}})}\right)\),  \(1\le r\le u-1\), \(1\le s\le v-1\), \(vr+us<uv\),
      \(\mu\in \RR/2\ZZ\), \(\lambda_{r,s,\frac{u}{v}}, \lambda_{u-r,v-s,\frac{u}{v}}\notin \mu\). The
      inequality involving both \(r\) and \(s\) is required due to the
      isomorphism \(\slE{\mu}{r}{s}\cong \slE{\mu}{u-r}{v-s}\).
    \end{enumerate}
    Further, there are the spectral flow relations \(\sigma^{\pm1} L_r\cong
    D^{\pm}_{u-r,v-1}\) and \(\sigma D^-_{r,s}\cong D^+_{u-r,v-s-1}\), where
    \(1\le r\le u-1\) and \(1\le s\le v-2\).
  \item \emph{(Adamovi\'{c}, Kawasetsu, Ridout \cite{bpAKR}).} Every simple module in
    \(\wtmod{\mmsl{u}{v}}\) is isomorphic to one
    of the following mutually inequivalent spectral flow twists of simple
    modules in \(\pwtmod{\mmsl{u}{v}}\).
    \begin{enumerate}
    \item \(\sigma^\ell D^+_{r,s}\), \(1\le r\le u-1\), \(1\le s\le v-1\), \(\ell\in\ZZ\).
    \item \(\sigma^\ell \slE{\mu}{r}{s}\), \(1\le r\le u-1\), \(1\le s\le v-1\), \(vr+us<uv\),
      \(\mu\in \RR/2\ZZ\), \(\lambda_{r,s,\frac{u}{v}}, \lambda_{u-r,v-s,\frac{u}{v}}\notin \mu\).
    \end{enumerate}
  \item \emph{(Adamovi\'{c}, Kawasetsu, Ridout \cite{bpAKR}).} At each of the two
    disallowed weights \(\mu\) above there are two  
    mutually inequivalent reducible yet indecomposable relaxed highest weight modules
    uniquely characterised by the non-split short exact sequences
    \begin{align}
      0\rightarrow D^{\pm}_{r,s}\rightarrow &\slEr{r}{s}{\pm}\rightarrow
      D^{\mp}_{u-r,v-s}\rightarrow 0,\nonumber \\
       0\rightarrow D^{\pm}_{u-r,v-s}\rightarrow &\slEr{u-r}{v-s}{\pm}\rightarrow
      D^{\mp}_{r,s}\rightarrow 0, 
    \end{align}
   for \(1\le r\le u-1\), \(1\le s\le v-1\), \(vr+us<uv\), where \(\slEr{r}{s}{+}\) and
    \(\slEr{u-r}{v-s}{-}\) correspond to \(\mu=\sqbrac{\lambda_{r,s,\frac{u}{v}}}\), and \(\slEr{r}{s}{-}\) and
    \(\slEr{u-r}{v-s}{+}\) correspond to \(\mu=\sqbrac{\lambda_{u-r,v-s,\frac{u}{v}}}\).
  \item \emph{(Arakawa, Creutzig, Kawasetsu \cite{ACK}).} The simple modules \(\sigma^\ell \slE{\mu}{r}{s}\), \(1\le r\le u-1\), \(1\le s\le v-1\), \(vr+us<uv\), \(\ell\in\ZZ\), \(\mu\in \RR/2\ZZ\),
      \(\lambda_{r,s,\frac{u}{v}},\lambda_{u-r,v-s,\frac{u}{v}}\notin \mu\), are
      projective and injective in \(\wtmod{\mmsl{u}{v}}\).
  \item \emph{(Adamovi{\'c} \cite{A3}; Arakawa, Creutzig, Kawasetsu \cite{ACK}).} For \(1\le r\le u-1\), \(1\le s\le v-1\), \(\ell\in \ZZ\) the
    projective cover and injective hull of \(D^+_{r,s}\) are isomorphic and
    denoted
    \(P_{r,s}\). These modules are uniquely characterised by the non-split
    exact sequences
    \begin{align}
      0 \rightarrow \sigma^\ell \slEr{r}{s}{+} \rightarrow &\sigma^\ell P_{r,s}
      \rightarrow \sigma^{\ell+1} \slEr{r}{s+1}{+} \rightarrow 0,\qquad 1\le s\le v-2,\nonumber\\
      0 \rightarrow \sigma^\ell \slEr{r}{v-1}{+} \rightarrow &\sigma^\ell P_{r,v-1}\rightarrow
      \sigma^{\ell+2} \slEr{u-r}{1}{+} \rightarrow 0,
      \label{eq:projpchar}
    \end{align}
    and they are also uniquely characterised by the non-split exact sequences
    \begin{align}
      0 \rightarrow \sigma^\ell \slEr{u-r}{v-s-1}{-}
      \rightarrow &\sigma^{\ell-1} P_{r,s}\rightarrow
      \sigma^{\ell-1} \slEr{u-r}{v-s}{-}
      \rightarrow 0,\qquad 1\le s\le v-2,\nonumber\\
      0 \rightarrow \sigma^\ell \slEr{r}{v-1}{-} \rightarrow &\sigma^{\ell-2} P_{r,v-1}\rightarrow
    \sigma^{\ell-2} \slEr{u-r}{1}{-}\rightarrow 0.
    \label{eq:projmchar}
  \end{align}
    Since the \(L_r\) and \(D_{r,s}^-\) modules are related to the \(D^+_{r,s}\)
    by spectral flow their projective covers and injective hulls satisfy the same
    spectral flow identifications. In particular the projective cover of
    \(L_1\cong \sigma^{-1}(D_{u-1, v-1}^+),\) (namely, \(\mmsl{u}{v}\) as a module over itself and hence the tensor
    unit) is \(\sigma^{-1}P_{u-1,v-1}\).
  \end{enumerate}
  \label{thm:sl2modclass}
\end{prop}

We recall the following result of Creutzig that shows that the category of
weight modules for $\mmsl{u}{v}$ admits a braided tensor category structure,
with the framework of Lepowsky, Huang and Zhang \cite{HLZ}.
\begin{theorem}[Creutzig \cite{Csl2ten}; Creutzig, Huang, Yang \cite{CHY}]
    For coprime \(u\ge2\), \(v\ge1\) the category \(\wtmod{\mmsl{u}{v}}\) admits a
  tensor product induced from intertwining operators and thereby is a braided monoidal 
  category.
  Further, for \(a,b,c,d\in \NN_0\), \( b,c,d\le a-1\), define the integers
   \begin{equation} 
     N^{(a)\,b}_{c,d}=
     \begin{cases}
       1&\text{ if } \abs{c-d}+1\le b\le \min\set{c+d-1,2b-c-d-1},\ b+c+d\equiv
          1\pmod{2},\\
       0&\text{ otherwise}.
     \end{cases} \label{eq:sl2fusionrules}
   \end{equation}
   The simple \(\mmsl{u}{v}\) weight modules \(L_r\) satisfy the fusion product
  decompositions
  \begin{align}
    L_r\fuse L_{r'}&\cong \bigoplus^{u-1}_{r''=1}N^{(u)\,r''}_{r,r'}
                         L_{r''},\quad 1\le r,r'\le u-1,\nonumber\\
    L_r\fuse
    \slE{\mu'}{r'}{s'}&\bigoplus^{u-1}_{r''=1}N^{(u)\,r''}_{r,r'}\slE{\mu'+r-1}{r''}{ s'},
    \quad \mu' \in \RR/2\ZZ,\ \lambda_{r,s,\frac{u}{v}}, \lambda_{u-r,v-s,\frac{u}{v}}\notin
                    \mu',\nonumber\\ 
                       &\hspace{4.4em} 1\le r,r'\le u-1,\ 1\le s'\le v-1,\ vr+us<uv.
                         \label{eq:LEfusion}
  \end{align}
  Further, the simple modules \(L_r\), \(1\le r\le u-1\) are  rigid. 
  \label{thm:sl2monoidal}
\end{theorem}

\begin{conj}[Creutzig, Ridout  \cite{CR2}; Creutzig, Kanade, Liu, Ridout \cite{ospCKLR}]
  For  \(1\le r,r'\le u-1\), \(1\le s,s'\le v-1\), \(vr+us<uv\), \(\mu,\mu'\in \RR/2\ZZ\),
      \(\lambda_{r,s,\frac{u}{v}},\lambda_{u-r,v-s,\frac{u}{v}}\notin \mu\),
      \(\lambda_{r',s',\frac{u}{v}},\lambda_{u-r',v-s',\frac{u}{v}}\notin \mu'\),
  \begin{align}
  \slE{\mu}{r}{s}\fuse
  \slE{\mu'}{r'}{s'}&\cong\bigoplus_{r'',s''}N^{(u)\, r''}_{r,r'}\brac*{N^{(v)\,
  s''}_{s,s'-1}+N^{(v)\, s''}_{s,s'+1}} \slE{\mu+\mu'}{r''}{s''}\nonumber\\
  &\quad\oplus \bigoplus_{r'',s''}N^{(u)\, r''}_{r,r'}N^{(v)\, s''}_{s,s'}\brac*{\sigma^{-1} \slE{\mu+\mu'+\frac{u}{v}}{r''}{s''}\oplus
    \sigma \slE{\mu+\mu'- \frac{u}{v}}{r''}{s''}},
    \label{eq:EEfusion}
  \end{align}
  where additionally \(\mu+\mu'\) must be such that all of the relaxed highest
  weight modules appearing in the direct sum decomposition on the \rhs{} are
  simple (this only excludes a finite number of values for the sum \(\mu+\mu'\)).
  For \(1\le r\le u-1\), \(\mu,\mu'\in \RR/2\ZZ\), \(\pm \lambda_{1,1,\frac{u}{v}}\notin \mu\), \(\pm
  \lambda_{r,1,\frac{u}{v}}\notin \mu'\) and \(\lambda_{r,0,\frac{u}{v}}=r-1\in \mu+\mu'\)
  \begin{align}
    \slE{\mu}{1}{1}\fuse \slE{\mu'}{r}{1}\cong \sigma^{-1}P_{u-r,v-1}\oplus
    \brac*{1-\delta_{v,2}}\slE{\mu+\mu'}{r}{2},
    \label{eq:r1logfusion}
  \end{align}
  while for \(1\le r\le u-1\), \(2\le s\le v-2\), \(\mu,\mu'\in \RR/2\ZZ\), \(\pm \lambda_{1,1,\frac{u}{v}}\notin \mu\), \(\pm
  \lambda_{r,s,\frac{u}{v}}\notin \mu'\),
  \begin{align}
    \slE{\mu}{1}{1}\fuse \slE{\mu'}{r}{s}\cong
    \begin{cases}
      P_{r,s-1}\oplus \sigma^{-1} \slE{\mu+\mu'+\frac{u}{v}}{r}{s}\oplus
      \slE{\mu+\mu'}{r}{s+1},&\text{if }\ \lambda_{r,s-1,\frac{u}{v}}\in \mu+\mu',\\
      P_{u-r,v-s-1}\oplus\sigma^{-1}\slE{\mu+\mu'+\frac{u}{v}}{r}{s}\oplus
      \slE{\mu+\mu'}{r}{s-1}&\text{if }\ \lambda_{u-r,v-s-1,\frac{u}{v}}\in \mu+\mu',\\
      \sigma^{-1}P_{r,s}\oplus \sigma \slE{\mu+\mu'-\frac{u}{v}}{r}{s}\oplus
      \slE{\mu+\mu'}{r}{s-1},&\text{if }\ \lambda_{r,s+1,\frac{u}{v}}\in \mu+\mu',\\
      \sigma^{-1}P_{u-r,v-s}\oplus\sigma \slE{\mu+\mu'-\frac{u}{v}}{r}{s}\oplus
      \slE{\mu+\mu'}{r}{s+ 1}&\text{if }\ \lambda_{u-r,v-s+1,\frac{u}{v}}\in \mu+\mu'.
    \end{cases}
  \end{align}
  \label{thm:EEfusionconj}
\end{conj}

Note that from \eqref{eq:EEfusion} we can see that all simple relaxed highest weight modules of the form
\(\slE{\mu}{r}{s}\) are conjecturally generated as direct summands by repeated application of
\(\slE{\mu'}{1}{1}\) and \(\slE{\mu'}{2}{1}\) to
\(\slE{\mu''}{1}{1}\), and further that all other simple modules appear
as subquotients of these repeated tensor products, due to the appearance
  of the indecomposable reducible projective modules \(P_{r,s}\). 
Since the \(\slE{\mu}{r}{s}\) are projective, if they are
simple, a sufficient condition for their rigidity is that \(\slE{\mu'}{1}{1}\) and
\(\slE{\mu'}{2}{1}\) are rigid. Further, from \eqref{eq:LEfusion}, we see that
\(\slE{\mu'}{2}{1}\cong L_2\fuse \slE{\mu'+1}{1}{1}\) and hence by the rigidity of
\(L_2\), the rigidity of \(\slE{\mu'}{1}{1}\) implies the rigidity of 
\(\slE{\mu'}{2}{1}\). Thus sufficient conditions for the rigidity of
\(\wtmod{\mmsl{u}{v}}\) are given by \cref{thm:EEfusionconj} being true and
\(\slE{\mu'}{1}{1}\) being rigid.

\subsection{The $\Nt=2$ superconformal Lie algebra}

In order to define the $\Nt=2$ superconformal vertex algebra we first introduce the superconformal Lie algebra. It is the infinite-dimensional Lie
superalgebra with even and odd components of the basis respectively given by \(\set{L_n, J_n, C \st n\in \ZZ}\) and \(\set{G^{\pm}_r\st r\in \frac{1}{2}+ \ZZ}\) subject to the relations
\begin{align}
  \comm{L_n}{L_m}&=(n-m)L_{n+m}+\frac{n^3-n}{12}\delta_{n,-m}C, 
  &\comm{L_n}{G^\pm_r}&=\left( \frac{1}{2}n-r\right)G^\pm_{n+r},\nonumber\\
  \acomm{G^+_{r}}{G^-_s}&=2L_{r+s}+(r-s)J_{r+s}+\frac{4r^2-1}{12}\delta_{r,-s}C,
  &\comm{L_n}{J_m}&=-mJ_{n+m},\nonumber\\
  \comm{J_n}{J_m}&= \frac{n}{3}\delta_{n,-m}C,\qquad
  \comm{J_n}{G^\pm_r}=\pm G^\pm_{n+r},
  &\acomm{G^\pm_{r}}{G^\pm_s}&=0, \label{eq:comm}
\end{align}
for all $n,m \in \ZZ$, $r,s\in \frac{1}{2}+\ZZ$ and where \(C\) is central.
We denote
\begin{align}
  \nlie_\pm&= \cspn{L_{\pm(n+1)}, J_{\pm(n+1)}, G^+_{\pm(n+\frac{1}{2})}, G^-_{\pm(n+\frac{1}{2})} \st
    n\in \NN_0},\nonumber\\
  \nlie_0&=\cspn{L_0, J_0, C},\nonumber\\
  \nlie_{\geq 0}&= \nlie_+\oplus \nlie_0,
\end{align}
so that we have a triangular decomposition $\nlie=\nlie_-\oplus \nlie_0 \oplus
\nlie_+$.

\begin{defn}
  Let \(M\) be a \(\ZZ_2\)-graded \(\nlie\)-module on which \(C\) acts as \(c\cdot\id\) for
  some \(c\in \CC\).
  \begin{enumerate}
  \item \(M\) is called \emph{smooth}, if for every \(m\in M\)
    \begin{equation}
    E_n m=0=O_{n+\frac{1}{2}} m,\qquad n\gg 0,
  \end{equation}
    where \(E\in \set{L,J},\ O\in \set{G^+,G^-}\). Denote the
    full subcategory of all \(\nlie\)-modules whose objects are all smooth
    modules by \(\smod{\nlie_c}\).
  \item \(M\) is called \emph{weight}, if it is smooth, finitely
    generated, \( \nlie_{>0} \) acts locally nilpotently (for every  \(m\in
    M\), \(\dim\brac*{\UEA{\nlie_{>0}}m}<\infty\))  and \(J_0\) acts semisimply.
    This implies that \(M\) is graded by \(J_0\)-eigenvalues and generalised
    \(L_0\)-eigenvalues, that is, it decomposes into a direct sum of
    homogeneous spaces
    \begin{equation}
      M=\bigoplus_{q, h \in\CC} M_{q,h},\qquad M_{q,h}=\set{m\in M\st
        (J_0-q)m=0=(L_0-h)^Nm,\ N\gg 0},
    \end{equation}
    where \(\dim M_{q,h}<\infty\), the weight support
    \(\supp(M)=\set{(q,h),\in \CC^2\st M_{q,h}\neq 0}\subset \CC^2\) is
    contained within a finite number of \(\ZZ^2\)-cosets of \(\CC^2\) and
    \(M_{q,h}=0\) for \(\re\brac*{h}\ll 0\).
      Denote the full subcategory of \(\smod{\nlie_c}\) whose objects are all weight
  modules by \(\wtmod{\nlie_c}\).
\item Let \(M\) be weight. A non-zero homogeneous vector \(v\in M\) is called
  \emph{singular} if \(\nlie_+ v=0\).
  Note that \(G^\pm_{\frac12}v=0=J_1v\) is a necessary and sufficient
  condition for 
  singularity. If \(v\) is singular and generates \(M\), then \(v\) is called
  \emph{highest weight} and \(M\) is called a \emph{highest weight module}.
  \end{enumerate}
\end{defn}

Recall that \(\nlie\) admits a family of automorphisms  called
\emph{spectral flow}.
\begin{defn}\leavevmode
  \begin{enumerate}
  \item The {\it spectral flow} automorphisms $\sigma^\ell, \ell\in \ZZ$ of \(\nlie\),
    are determined on basis vectors by
    \begin{equation}
      \sigma^\ell(L_n)=L_n-\ell J_n+\frac{1}{6}\ell^2\delta_{n,0}C,\quad
      \sigma^\ell(J_n)=J_n-\frac{1}{3}\ell\delta_{n,0}C,\quad
      \sigma^\ell(G_s^\pm)=G_{s+\ell}^\pm,\quad
      \sigma^\ell(C)=C.
    \end{equation}
  \item Let \(g\in \aut\brac*{\nlie}\) and \(M\) a \(\nlie\)-module, then
    the \(g\)-twist of \(M\) is a \(\nlie\)-module whose underlying vectors
    space is that of \(M\) but  with the action
    \begin{equation}
      x\cdot_g m=g^{-1}(x)\cdot m,\qquad \forall x\in \nlie,\ \forall m\in M,
    \end{equation}
    where \(\cdot\) is the original action of \(\nlie\) on \(M\).
  \end{enumerate} 
  \label{def:spectralflowconjugation}
\end{defn}
Note that the categories \(\smod{\nlie_c}\) and \(\wtmod{\nlie_c}\) are both
closed under spectral flow automorphisms. Note that in principle the spectral
flow parameter could be a half odd integer, however, the resulting twisted
module would then no longer satisfy the locality axiom. That is, the
  series expansion of fields can contain half odd integers. Such twisted
  modules are called \emph{Ramond} (as opposed to \emph{Neveu-Schwarz}) modules. Half odd
integral spectral flow interchanges the Neveu-Schwarz and Ramond sectors.

We can now easily construct examples of smooth \(\nlie\)-modules. For example,
the Verma modules and their simple quotients are smooth. These are constructed
as follows. For \(q,h,c\in\CC\), let $\CC\vac_{q,h,c}$ denote the even \(1\)-dimensional $\nlie_{\ge0}$-module
characterised by
\begin{align}
 J_0 \vac_{q,h,c}&=q\vac_{q,h,c},
  &C \vac_{q,h,c}&=c\vac_{q,h,c}\nonumber\\
  L_0 \vac_{q,h,c}&=h\vac_{q,h,c},,
  &\nlie_+\vac_{q,h,c}&=0.
\end{align}
The Verma module $\NM{q}{h}{c}$, is the induced module
\begin{align}
  \NM{q}{h}{c}=\UEA{\nlie}\cten_{\UEA{\nlie_{\geq 0}}}\CC\vac_{q,h,c},
  \label{eq:verma}
\end{align}
while the simple quotient of \(\NM{q}{h}{c}\) is 
\begin{align}\label{defL}
\NL{q}{h}{c}=\frac{\NM{q}{h}{c}}{\ang{\text{maximal submodule}}}.    
\end{align}

\begin{prop} \label{thm:univn=2}
  Let \(c\in \CC\). Then the \(\nlie\)-module
  \(\uN{c}=\NM{0}{0}{c}/\ang{G^{\pm}_{-\frac{1}{2}}\vac_{0,0,c}}\) admits the
  structure of a \vosa{} by defining the field map on the strong generators
  \begin{align}
    Y(L_{-2}\vac,z)&=L(z)=\sum_{n\in \ZZ} L_{n}z^{-n-2},&
    Y(J_{-1}\vac,z)&=J(z)=\sum_{n\in \ZZ} J_{n}z^{-n-1},\nonumber\\
    Y(G^\pm_{-\frac{3}{2}}\vac,z)&=G^\pm(z)=\sum_{n\in \ZZ} G^\pm_{n+\frac{1}{2}}z^{-n-2},
  \end{align}
and continuing by derivatives and normal ordering. This \vosa{} is called
  the \emph{universal \(\Nt=2\) superconformal \va{}}. The strong generators above
  satisfy the \ope{} relations
  \begin{align}
    L(z)L(w)&\sim \frac{\frac{c}{2}}{(z-w)^4}+\frac{2L(w)}{(z-w)^2}+\frac{\partial_w L(w)}{z-w},
    &L(z)J(w)&\sim \frac{J(w)}{(z-w)^2}+\frac{\partial_wJ(w)}{z-w},\nonumber\\
    L(z)G^\pm(w)&\sim \frac{\frac32 G^\pm(w)}{(z-w)^2}+\frac{\partial_wG^\pm(w)}{z-w},
    &J(z)J(w)&\sim \frac{\frac{c}{3}}{(z-w)^2},\nonumber\\
    J(z)G^{\pm}(w)&\sim\frac{\pm G^{\pm}(w)}{(z-w)},
    &G^\pm(z)G^\pm(w)&\sim 0,\nonumber\\
    G^{\pm}(z)G^{\mp}(z)&\sim \frac{\frac{2c}{3}}{(z-w)^3}\pm\frac{2J(w)}{(z-w)^2}+ \frac{2L(w)\pm \partial J(w)}{z-w}.
  \end{align}
  Finally, \(\smod{\nlie_c}\) is the category of all \(\uN{c}\)-modules.
\end{prop}

\begin{prop}[Gorelik, Kac \cite{GK1}] The universal $\Nt=2$ vertex algebra $\uN{c}$ is not simple if and only  if there exist $u, v \in\ZZ$ with $u\geq 2$, $v\geq 1$ and $gcd(u,v)=1$ such that the central charge satisfies 
\begin{align}
  c=3-\frac{6v}{u}.
  \label{eq:cuv}
\end{align}
For these central charges \(\uN{c}\) admits a non-trivial maximal ideal
generated by a singular vector of conformal weight \((u-1)v\) and \(J\)-weight
\(0\).
\label{thm:minn=2}
\end{prop}

When $\uN{c}$ is not simple, so that $c=3-\frac{6v}{u}$ for some $u,v$ as in
Proposition \ref{thm:minn=2}, we denote its simple quotient by $\mmN{u}{v}$
and refer to it as the \(\Nt=2\) minimal model at central charge
$3-\frac{6v}{u}$.
\begin{defn}
For coprime \(u,v\in\ZZ_{>0}\), \(u\ge2\), let \(\wtmod{\mmN{u}{v}}\) be the
  full subcategory of \(\wtmod{\uN{3-\frac{6v}{u}}}\) of modules
  annihilated by the maximal ideal of \(\uN{3-\frac{6v}{u}}\) and for which
  the eigenvalues of \(J_0\) are real.
\end{defn}
\begin{prop}[Adamovi\'c \cite{A1}]
  For coprime \(u,v\in\ZZ_{>0}\), \(u\ge2\), denote
  \begin{equation} 
    h_{r,s;q}=\frac{\brac*{r-\frac{u}{v}s}^2-1}{4\frac{u}{v}}-\frac{u}{4v}q^2,\qquad
    r,s\in\ZZ,\ q\in\RR. \label{eq:N=2confweights}
  \end{equation}
  Then the simple highest weight modules
  \(\NL{q}{h}{3-\frac{6v}{u}}\) and their parity reversals \(\Pi \NL{q}{h}{3-\frac{6v}{u}}\)
  whose conformal and \(J\)-weight \((q,h)\) lie in one of the following sets
  \begin{enumerate}
  \item \(\set*{\brac*{\frac{pv}{u},h_{r,0;\frac{pv}{u}}}\st
      1\le r\le u-1,\ 1-r\le p\le r-1, p+r \equiv 1 \pmod{2}}\),
  \item \(\set*{\brac*{q,h_{r,s;q}}\st
      1\le r\le u-1,\ 1\le s\le v-1,\ vr+us<uv,  q\in \mathbb{R}}  \),
  \end{enumerate}
  form a complete set of representatives of simple isomorphism classes in
  \(\wtmod{\mmN{u}{v}}\). We denote \(\NL{q}{h_{r,s;q}}{3-\frac{6v}{u}}=\NL{q}{r}{s}\).
  \label{thm:n2simpleclass}
\end{prop}
We provide a new proof of the theorem above in \cref{thm:newclass}.

\begin{prop} [Creutzig, {\cite[Thm 8.2]{Csl2ten}}]
  The category \(\wtmod{\mmN{u}{v}}\)
    admits a tensor product induced from intertwining operators and
      thereby is a braided monoidal category.
\end{prop}

\subsection{The fermionic ghost vertex algebra}

Ghost systems date back to the early days in string theory in the physics
literature. For a formal mathematical description of them, we refer to Kac's book \cite{KacBeg}.
To define the fermionic ghost vertex algebra \(\bc\) we first introduce its underlying
Lie superalgebra \(\bclie\) with even and odd components of the basis respectively given by
\(\set{\vac}\) and \(\set{b_r,c_r\st r\in \frac{1}{2}+\ZZ}\). These basis vectors
satisfy the commutation relations
\begin{equation}
  \acomm{b_r}{c_s}=\delta_{r+s,0}\vac,\qquad r,s\in \frac{1}{2}+\ZZ,
\end{equation}
 where \(\vac\) is central (and will always be taken to act as the identity
on modules). We then have the triangular decomposition
\begin{align}
  \bclie_\pm&=\cspn{b_{\pm r},c_{\pm r}\st r\in \frac{1}{2}+\NN_0},&
                                                                     \bclie_0&=\CC\vac,&
  \bclie_{\ge}&=\bclie_+\oplus \bclie_0.
\end{align}
Since we require \(\vac\) to act as the identity there is only one Verma
module (up to isomorphism and parity reversal) for this decomposition.
\begin{equation}
  \bc = \ind_{\bclie_{\ge}}^{\bclie}\CC \Omega,\qquad b_r\Omega=c_r\Omega=0,\quad
  \forall r\in \frac{1}{2}+\NN_0.
\end{equation}

\begin{prop}
  The \(\bclie\)-module \(\bc\) admits the structure of a vertex operator
  superalgebra by defining the field map on the odd strong generators
  \begin{equation}
    Y(b_{-\frac12}\Omega,z)=b(z)=\sum_{r\in\ZZ+\frac12}b_{r}z^{-r-\frac12},\qquad
    Y(c_{-\frac12}\Omega,z)=c(z)=\sum_{r\in\ZZ+\frac12}c_{r}z^{-r-\frac12},
  \end{equation}
  and continuing by derivatives and normal ordering. This \va{} has a number
  of names in the literature including \(\bc\) \emph{vertex algebra},
  \emph{fermionic ghost system} and \emph{charged free fermions}. The strong
  generators above satisfy the \ope{} relations
  \begin{equation}
    b(z)c(w)\sim\frac{1}{z-w},\qquad b(z)b(w)\sim c(z)c(w)\sim 0.
  \end{equation}
  The \(\bc\) \emph{vertex algebra} admits a conformal vector
  \begin{equation}
    \omega_{\bc}=\frac{1}{2}\brac*{b_{-\frac32}c_{-\frac12}+c_{-\frac32}b_{-\frac12}}\Omega,\quad
    Y(\omega_{\bc},z)=T^{\bc}(z)=\frac{1}{2}\brac*{\normord{(\partial
        b(z))c(z)}+\normord{(\partial c(z))b(z)}}.
  \end{equation}
  of central charge \(1\) and an additional distinguished vector
  \begin{equation}
    Q=b_{-\frac12}c_{-\frac12}\Omega,\qquad Y(Q,z)=Q(z)=\normord{b(z)c(z)},
  \end{equation}
  such that \(b(z),\ c(z)\) are conformal weight \(\frac12\) primaries,
  \(Q(z)\) is a conformal weight \(1\) primary and \(Q(z)\) satisfies the
  \ope{} relations
  \begin{equation}
    Q(z)b(w)\sim \frac{b(w)}{z-w},\qquad
    Q(z)c(w)\sim \frac{-c(w)}{z-w},\quad
    Q(z)Q(w)\sim \frac{1}{(z-w)^2}.
  \end{equation}
  Namely, \(Q(z)\) generates a subalgebra isomorphic to the Heisenberg \va{}
  \(\heis{\frac12}\) and its zero mode \(Q_0\) (the coefficient of \(z^{-1}\)) gives
  a \(\ZZ\)-grading to \(\bc\), \emph{called ghost weight}, which assigns
  weight \(0\) to \(\Omega\), weight \(1\) to \(b(z)\) modes and weight \(-1\)
  to \(c(z)\) modes.
\end{prop}

\begin{defn}
Let \(M\) be a \(\ZZ_2\)-graded \(\bclie\)-module on which \(\vac\) acts as
the identity.
  \begin{enumerate}
  \item \(M\) is called \emph{smooth}, if for every \(m\in M\)
    \begin{equation}
    b_r m=c_r m=0,\qquad r\gg 0.
  \end{equation}
  Denote the full subcategory of all \(\bclie\)-modules whose objects are all smooth
  modules by \(\smod{\bclie}\). This category is also the category of all \(\ZZ_2\)-graded \(\bc\)-modules.
  \item \(M\) is called \emph{weight}, if it is smooth, finitely
    generated, \(\bclie_{\ge0}\) acts locally nilpotently (for every \(m\in
    M\), \(\dim\brac*{\UEA{\bclie_{\ge0}}m}<\infty\)) and \(Q_0\) acts semisimply.
    This implies that \(M\) is graded by \(Q_0\)-eigenvalues and generalised
    \(L_0\)-eigenvalues, that is, it decomposes into a direct sum of
    homogeneous spaces
    \begin{equation}
      M=\bigoplus_{q,h\in\CC} M_{q,h},\qquad M_{q,h}=\set{m\in M\st
        (Q_0-q)m=0=(L_0-h)^Nm,\ N\gg 0},
    \end{equation}
    where \(\dim M_{q,h}<\infty\), the weight support
    \(\supp(M)=\set{(q,h),\in \CC^2\st M_{q,h}\neq 0}\subset \CC^2\) is
    contained within a finite number of \(\ZZ^2\)-cosets of \(\CC^2\) and
    \(M_{q,h}=0\) for \(\re\brac*{h}\ll 0\).
      Denote the full subcategory of \(\smod{\bclie}\) whose objects are all weight
  modules by \(\wtmod{\bclie}\).
  \end{enumerate}
\end{defn}

We conclude this section by recalling the known equivalence between the category of $\bclie$ weight modules and super vector spaces.
\begin{prop}
  The category \(\wtmod{\bclie}\) is semisimple with two simple isomorphism
  classes \(\bc\) and its parity reversal \(\Pi\bc\). As a braided monoidal
  category (with the tensor product constructed from intertwining operators) \(\wtmod{\bclie}\) is equivalent to \(\svec\), the symmetric
  monoidal category of finite-dimensional super vector spaces. In particular,
  \(\bc\) corresponds to \(\CC^{1|0}\), the even \(1\)-dimensional super vector
  space and \(\Pi\bc\) corresponds to \(\CC^{0|1}\) the odd \(1\)-dimensional
  super vector space.
\end{prop}

\section{The coset realisation of the \(\Nt=2\) 
superconformal algebra}
\label{sec:coset}

In this section we recall a well known coset construction of the
universal \(\Nt=2\) superconformal \vosa{s} \(\uN{c}\) and their simple
quotients \(\mmN{u}{v}\)  which appeared in the physics literature in  \cite{KZ1, KZ2, EG} and was more recently studied in detail in \cite{CL2}.

\begin{prop}
  Consider the \vosa{} \(\usl{t-2}\cten \bc\), \(t\in \CC^\times\),
  then there exists an embedding \(\phi_h:\heis{2t}\hookrightarrow\usl{t-2}\cten \bc
  \) characterised on the strong generator \(a_{-1}\ket{0}\) and continued to
  the conformal vector \(\frac{1}{4t}a_{-1}^2\ket{0}\) by
  \begin{align}
    a=a_{-1}\ket{0}&\stackrel{\phi_h}{\longmapsto} h\cten\Omega+2\vac\cten Q=
                     \brac*{h_{-1}+2b_{-\frac12}c_{-\frac12}}\vac\cten\Omega,\nonumber\\
    \frac{1}{4t}a_{-1}^2\ket{0}&\stackrel{\phi_h}{\longmapsto}\frac{1}{4t}\brac*{h_{-1}^2+4\brac*{b_{-\frac32}c_{-\frac12}+c_{-\frac32}b_{-\frac12}}+4h_{-1}b_{-\frac12}c_{-\frac12}}\vac\cten\Omega.
    \label{eq:heisembed}
  \end{align}
  There exists an additional embedding \(\phi_n:\uN{3-\frac{6}{t}}\hookrightarrow\usl{t-2}\cten \bc \) characterised  by
  \begin{align}
    J=J_{-1}\vac&\stackrel{\phi_n}{\longmapsto}
    \frac{h\cten\Omega-(t-2)\vac\cten Q}{t}=\frac{h_{-1}-(t-2)b_{-\frac12}c_{-\frac12}}{t}\vac\cten\Omega,\nonumber\\
    G^+=G^+_{-\frac32}\vac&\stackrel{\phi_n}{\longmapsto} \sqrt{\frac{2}{t}}e\cten
                            c=\sqrt{\frac{2}{t}}e_{-1}c_{-\frac12}\vac\cten \Omega,\nonumber\\
    G^-=G^-_{-\frac32}\vac&\stackrel{\phi_n}{\longmapsto} \sqrt{\frac{2}{t}}f\cten
    b=\sqrt{\frac{2}{t}}f_{-1}b_{-\frac12}\vac\cten \Omega,\nonumber\\
    \omega=L_{-2}\vac &\stackrel{\phi_n}{\longmapsto}
                        \frac{1}{2t}\brac*{e_{-1}f_{-1}+f_{-1}e_{-1}+\brac*{t-2}\brac*{b_{-\frac32}c_{-\frac12}+c_{-\frac32}b_{-\frac12}}-2h_{-1}b_{-\frac12}c_{-\frac12}}\vac\cten\Omega.
                        \label{eq:N2embed}
  \end{align}
  Further, the image of \(\phi_n\) is the commutant (also known as a coset in
  the physics literature) of the image of \(\phi_h\), that is
  \begin{equation}
    \phi_n\brac*{\uN{3-\frac{6}{t}}}=\com\brac*{\phi_h(\heis{2t}), \usl{t-2}\cten \bc}.
  \end{equation}
  Further, if \(t=\frac{u}{v}\), for coprime \(u\ge 2\), \(v\ge1\), then
  \(\uN{3-\frac{6}{t}}\) and \(\usl{t-2}\) admit non-trivial ideals, and the
  embedding \(\phi_n\) and commutant factor through the respective simple
  quotients, that is
  \begin{equation}
     \phi_n\brac*{\mmN{u}{v}}=\com\brac*{\phi_h(\heis{2t}),\mmsl{u}{v}\cten \bc}.
   \end{equation}
   \label{thm:coset}
\end{prop}

\begin{prop}
  For \(t\in \CC^\times\) consider the embedding 
  \( \phi =\phi_h\cten \phi_n:
  \heis{2t}\cten \uN{3-\frac{6}{t}}\hookrightarrow \usl{t-2}\cten \bc \),
  where \(\phi_h,\phi_n\) were characterised in \cref{thm:coset}.
   Let \(M\) be an indecomposable \(\usl{t-2}\) weight module with \(\slt\)-weight
   support \(\supp(M)=\lambda+2\ZZ, \lambda\in \CC\). Then the following hold.
   \begin{enumerate}
   \item \(M\cten \Pi^{i}\bc\), \(i=0,1\) is an
     \(\heis{2t}\cten \uN{3-\frac{6}{t}}\)-module by restriction (that is,
     pulling back along \(\phi\)) and decomposes as
     \begin{equation}
       \res\brac*{M\cten \Pi^{i}\bc} \cong \bigoplus_{p\in \supp(M)} \mathcal{F}_p\cten
       C^{\sqbrac*{i}}_p(M),
       \label{eq:moduledecomp}
     \end{equation}
     where the \(C^{\sqbrac*{i}}_p(M)\) are indecomposable
     \(\uN{3-\frac{6}{t}}\)-modules.
     \label{itm:moduledecomp}
   \item A vector \(m\in M\cten \Pi^i\bc\) is homogeneous for
     \(\usl{t-2}\cten \bc\) if and only if it is homogeneous for
     \(\heis{2t}\cten \uN{3-\frac{6}{t}}\). The grading operators (and hence
     weights) for the
     four algebras are interrelated by
     \begin{equation}
       L_{0}^{\slt}+L_{0}^\bc = L_{0}^{\mathsf{H}}+
       L_{0}^{\Nt=2},\qquad
       a_0=h_0+2Q_0,\qquad J_0=\frac{h_0-(t-2)Q_0}{t}.
     \end{equation}
     \label{itm:wtrel}
   \item For every \(\slta\) relaxed highest
     weight vector \(v\in M\) of \(\slt\) weight \(\mu\in \supp(M)\) and
     conformal weight \(h\) there exists a \(\nlie\) highest weight vector
     \(\chi\in C^{\sqbrac*{i}}_{\mu}(M)\) whose conformal weight is
     \(h-\frac{\mu^2}{4t}\) and whose \(J\)-weight is \(\frac{\mu}{t}\). In
     particular, in the decomposition \eqref{eq:moduledecomp} above
     \(v\cten \Omega=\ket{\mu}\cten \chi\).
     \label{itm:singrel} 
   \item  If \(t=\frac{u}{v}\) for
     coprime \(u\ge2\), \(v\ge1\) and \(M\) is one of the simple relaxed
     highest weight modules \(\slE{\mu}{r}{s}\), \(\mu\in \CC/2\ZZ\),
     \(\pm\lambda_{r,s}\notin\mu\) of 
     \cref{thm:sl2modclass}, then
     \begin{equation}
       C^{\sqbrac*{i}}_p(\sigma^\ell \slE{\mu}{r}{s})=\Pi^{i+\ell}\NLs{\frac{pv}{u}-\ell}{r}{s},\qquad
       p\in \mu+\ell\tfrac{u}{v},
     \end{equation}
     with \(\NL{q}{r}{s}\) as in \cref{thm:n2simpleclass}.
     \label{itm:resid}
   \end{enumerate}
   \label{thm:coset-and-sing-vectors}
 \end{prop}
 \begin{proof}
   With the exception of Part \ref{itm:singrel} above
   \cref{thm:coset-and-sing-vectors} is well known. 
   Most of the proof will
   therefore focus on Part \ref{itm:singrel}.

   \cref{itm:moduledecomp}:
   The form of the restriction \eqref{eq:moduledecomp} and the range of the
   Fock space weights \(p\) follows by applying \cite[Thm 3.8]{CKLR} to
   the embedding \(\phi=\phi_h\cten \phi_n\). The indecomposability of
   \(C_p(M)\) then follows from the indecomposability of \(M\) also by using
   \cite[Thm 3.8]{CKLR}.

   \cref{itm:wtrel}:
   The relations between grading operators follow directly from the embedding
   formulae \eqref{eq:heisembed} and \eqref{eq:N2embed}.

   \cref{itm:singrel}:
   Let \(v\in M\) be a relaxed highest weight vector of \(\slt\) weight
   \(\mu\) and conformal weight \(h\).
   The formula for the embedding \eqref{eq:heisembed} implies that
   \begin{equation}
     a_n v\cten \Omega=\brac*{h_n+2Q_n}v\cten
     \Omega=\mu\delta_{n,0}v\cten \Omega,\qquad n\ge0.
   \end{equation}
   Note that we have suppressed the embedding \(\phi\). Thus \(v\cten
   \Omega\) is a Heisenberg highest weight vector of Heisenberg weight \(\mu\), that
   is, \(v\cten \Omega\in \ket{\mu}\cten\chi\), \(\chi\in C_\mu(M)\).
   Further, the formulae for the embedding \eqref{eq:N2embed} imply that
   \begin{align}
     G^+_r v\cten \Omega&= \sqrt{\frac{2}{t}}
     \sum_{n\in\NN_0}e_{-n}v\cten c_{r+n}\Omega=0,\quad
                            r\ge\frac{1}{2},\nonumber\\
     G^-_r v\cten \Omega&= \sqrt{\frac{2}{t}}
                            \sum_{n\in\NN_0}f_{-n}v\cten b_{r+n}\Omega=0,\quad r\ge\frac{1}{2},\nonumber\\
     J_n v\cten \Omega&=\frac{h_n-kQ_n}{t}v\cten
                          \Omega=\delta_{n,0}\frac{\mu}{t}v\cten
                          \Omega,\quad n\ge 0.
   \end{align}
   Thus \(\chi\) is an \(\nlie\) singular vector with \(J\)-weight
   \(\frac{\mu}{t}\). Finally, since the image of the sum of the conformal
   vectors of \(\heis{2t}\) and \(\uN{3-\frac{6}{t}}\) is the conformal vector
   of \(\usl{t-2}\cten \bc\) the conformal weight \(h\) of \(v\) must be the
   sum of the conformal weights of \(\ket{\mu}\) and \(\chi\). Hence the
   conformal weight of \(\chi\) is \(h-\frac{\mu^2}{4t}\).

   \cref{itm:resid}:
   The relevant formulae for \(C^{\sqbrac*{i}}_p(M)\) are given in
   \cite[Eq 4.16]{CLRW} (note that \(i\) here corresponds to \(2i\) in  \cite{CLRW}).
   
 \end{proof}

For coprime \(u\ge2,\ v\ge1\), the categories \(\wtmod{\heis{2\frac{u}{v}}}\),
\(\wtmod{\mmN{u}{v}}\), \(\wtmod{\mmsl{u}{v}}\) and \(\wtmod{\bc}\) are all
locally finite abelian categories whose objects all have finite Jordan-H\"older
length (see \cite[Thm 1.2]{ACK} and \cite[Cor 5.1]{Csl2ten}), thus by
\cite[Thm 4.11]{MCDelPRod} we have the following
equivalences of braided monoidal categories.
\begin{align}
\wtmod{\brac*{\heis{2\frac{u}{v}}\cten \mmN{u}{v}}}&\cong
\wtmod{\heis{2\frac{u}{v}}}\delprod \wtmod{\mmN{u}{v}},\nonumber\\
\wtmod{\brac*{\mmsl{u}{v}\cten \bc}}&\cong
\wtmod{\mmsl{u}{v}}\delprod \wtmod{\bc},
\end{align}
where \(\delprod\) denotes the Deligne tensor product of abelian categories.
Further, \(\mmsl{u}{v}\cten \bc\) is a commutative algebra
 object in (a direct limit completion of) \(\wtmod{\heis{2\frac{u}{v}}\cten
   \mmN{u}{v}}\) (see \cite{DLcomp} for details on direct limit completions in
 the context of \vosa{s}).
We can therefore consider the category \(\rep{\brac*{\mmsl{u}{v}\cten \bc}}\), the
category of not necessarily local \(\mmsl{u}{v}\cten \bc\)-modules in
\(\wtmod{\brac*{\heis{2\frac{u}{v}}\cten \mmN{u}{v}}}\). These will
generally be twisted modules, the full subcategory of non-twisted (or local)
modules is denoted \(\rep^0 \brac*{\mmsl{u}{v}\cten \bc}\cong \wtmod{\brac*{\mmsl{u}{v}\cten \bc}}\). Further, the restriction
functor
\(\res: \rep \brac*{\mmsl{u}{v}\cten \bc} \to \wtmod{\brac*{\heis{2\frac{u}{v}}\cten
    \mmN{u}{v}}}\) has a left adjoint, the induction functor
\(\ind: \wtmod{\brac*{\heis{2\frac{u}{v}}\cten \mmN{u}{v}}}\to \rep\brac*{
\mmsl{u}{v}\cten \bc}\). This functor admits 
a monoidal structure (see \cite[Sec 2\&3]{CKM} for details).

Since the categories \(\wtmod{\heis{2\frac{u}{v}}}\) and \(\wtmod{\bc}\)
are semisimple and pointed (that is, tensoring by any simple object is
invertible) many hom spaces are \(1\)-dimensional which leads to convenient
factorisations in the Deligne tensor product categories. For example,
\begin{align}
  &\hom_{\mmN{u}{v}\cten\heis{2\frac{u}{v}}}\brac*{\brac*{M\cten \mathcal{F}_{p}}\fuse
  \brac*{N\cten \mathcal{F}_q},\brac*{P\cten \mathcal{F}_s}}\nonumber\\
  &\qquad\cong
  \hom_{\mmN{u}{v}}\brac*{M\fuse
    N,P}\cten\hom_{\heis{2\frac{u}{v}}}\brac*{F_p\fuse F_q,F_s}
    \cong \delta_{p+q,s}\hom_{\mmN{u}{v}}\brac*{M\fuse N,  P},\nonumber\\
  &\hom_{\mmsl{u}{v}\cten\bc}\brac*{\brac*{R\cten \Pi^i\bc}\fuse
    \brac*{S\cten \pi^j \bc},\brac*{T\cten \Pi^k\bc}}\nonumber\\
  &\qquad \cong\hom_{\mmsl{u}{v}}\brac*{R\fuse S,T}\cten
    \hom_{\bc}\brac*{\Pi^i\bc\fuse \Pi^j\bc,\Pi^k\bc}\cong
    \delta_{i+j,k}\hom_{\mmsl{u}{v}}\brac*{R\fuse S,T},
\end{align}
where \(M,N,P\in \wtmod{\mmN{u}{v}}\), \(p,q,s\in \RR\),
\(R,S,T\in\wtmod{\mmsl{u}{v}}\) and \(i,j,k\in\set{0,1}\).
 
\begin{theorem}
  For any \(M\in \wtmod{\mmsl{u}{v}}\), let \(C^{[i]}_p(M)\), \(p\in\supp{M}\), \(i=0,1\) denote the \(\mmN{u}{v}\)-module appearing as the counter part to the Fock space \(\mathcal{F}_p\) in the formula \eqref{eq:moduledecomp} for the restriction of \(M\otimes \Pi^i \bc\).
  Then there exists a functorial linear isomorphism
  \begin{equation}
    \binom{O}{M,\ N}_{\mmsl{u}{v}}\cong
    \binom{C^{\sqbrac*{i+j}}_{p+q}(O)}{C^{\sqbrac*{i}}_p(M),\ C^{\sqbrac*{j}}_q(N)}_{\mmN{u}{v}},
  \end{equation}
  where \(M,N,O\in \wtmod{\mmsl{u}{v}}\), \(p\in\supp(M)\), \(q\in\supp(N)\), \(i,j\in\set{0,1}\) and where the subscript of the intertwiner spaces indicates their corresponding \va{.}
  Further for \(s\in \supp{O}\)
  \begin{equation}
    \binom{C^{\sqbrac*{n}}_{s}(O)}{C^{\sqbrac*{i}}_p(M),\
      C^{\sqbrac*{j}}_q(N)}_{\mmN{u}{v}}=0,\qquad \text{if}\ s\neq p+q \text{
      or }\ n \not\equiv i+j \pmod{2}.
  \end{equation}
  \label{thm:homspaceeq}
\end{theorem}
The above theorem was implicit in the presentation of the conjectured fusion
rules in \cite[Sec 6]{CLRW} and has also been generalised for larger families of
algebras in \cite{CGNR}. Nevertheless, we provide a specialised proof for clarity.
\begin{proof}
  The theorem follows from induction and restriction being adjoint functors
  and induction being a monoidal functor \cite[Thm 2.59]{CKM}).
  More specifically,
  \begin{align}
    \binom{O}{M,\ N}_{\mmsl{u}{v}}&\cong \hom_{\mmsl{u}{v}}\brac*{M\fuse
      N,O} \nonumber\\
    &\cong \hom_{\mmsl{u}{v}\cten \bc}\brac*{(M\cten \Pi^i\bc)\fuse
      (N\cten \Pi^j\bc),O\cten\Pi^{i+j}\bc}\nonumber\\
    &\cong \hom_{\mmsl{u}{v}\cten
      \bc}\brac*{\ind\brac*{\brac*{\mathcal{F}_p\cten C^{\sqbrac*{i}}_p(M)}\fuse
        \brac*{\mathcal{F}_q\cten C^{\sqbrac*{j}}_q(N)}},O\cten\Pi^{i+j}\bc}\nonumber\\
    &\cong\hom_{\heis{2\frac{u}{v}}\cten\mmN{u}{v} }\brac*{\brac*{\mathcal{F}_p\cten C^{\sqbrac*{i}}_p(M)}\fuse
      \brac*{\mathcal{F}_q\cten C^{\sqbrac*{j}}_q(N)},\res\brac*{O\cten\Pi^{i+j}\bc}}\nonumber\\
    &\cong\bigoplus_{s\in \supp{O}}\hom_{{\heis{2\frac{u}{v}}\cten\mmN{u}{v} }}\brac*{\mathcal{F}_{p+q}\cten \brac*{C^{\sqbrac*{i}}_p(M)\fuse
        C^{\sqbrac*{j}}_q(N)},\mathcal{F}_{s}\cten C^{\sqbrac*{i+j}}_s(O)}\nonumber\\
    &\cong\hom_{\mmN{u}{v}}\brac*{C^{\sqbrac*{i}}_p(M)\fuse
      C^{\sqbrac*{j}}_q(N), C^{\sqbrac*{i+j}}_{p+q}(O)}\cong \binom{C^{\sqbrac*{i+j}}_{p+q}(O)}{C^{\sqbrac*{i}}_p(M),\ C^{\sqbrac*{j}}_q(N)}_{\mmN{u}{v}},
  \end{align}
  where we have used that \(\Pi^i\bc\fuse \Pi^{j}\bc\cong \Pi^{i+j}\bc\) and \(\mathcal{F}_p\fuse \mathcal{F}_q\cong \mathcal{F}_{p+q}\).
\end{proof}

While it is known that monoidal functors such as induction preserve duals (see
\cite[Exe 2.10.6]{EGNO}), sufficient conditions for the converse to hold
will also prove helpful below.
\begin{lemma}
  Let \(M\in\wtmod{\mmsl{u}{v}}\) be simple and let \(C_p^{[i]}(M)\cten
  \mathcal{F}_p\), \(p\in \supp\brac*{M}\), \(i\in \set{0,1}\) be a direct summand
  of the restriction \(\res M\cten\bc\). If \(M\) is rigid with dual
  \(M^\vee\), then \(C_p^{[i]}(M)\) is too,
  with rigid dual \(C_p^{[i]}\brac*{M}^\vee= C_{-p}^{[i]}\brac*{M^\vee}\).
  \label{thm:rigiditytransport}
\end{lemma}

\begin{proof}
  Since all tensor categories involved are braided, it is sufficient to only
  consider left duals.
  Further, since all the \vosa{s} involved are isomorphic to their
  contragredients, it is sufficient to show that only one of the rigidity zig
  zag relations is non-zero by
  \cite[Lem 4.2.1 and Cor 4.2.2]{CMcY}.
  The \voa{} \(\mmsl{u}{v}\) is self contragredient, thus, the rigid dual \(M^\vee\) of
  \(M\), which exists by assumption, must be isomorphic to the contragredient
  \(M'\). Hence the weight support of \(M^\vee\) must be the negative of the
  weight support of \(M\), that is,
  \(\supp\brac*{M^\vee}=-\supp\brac*{M}\). Therefore the restriction of
  \(M^\vee\cten \bc\) must admit a non-trivial direct summand of the form
  \(C_{-p}^{[i]}\brac*{M^\vee}\cten \mathcal{F}_{-p}\). Since \(\bc\) is rigid
  and self dual \(M\cten \bc\) is also rigid with the first of its zig zag
  relations given by
  \begin{align}
    \id=&M\cten\bc \xrightarrow{l^{-1}} \brac*{\mmsl{u}{v}\cten \bc}\fuse M\cten\bc
    \xrightarrow{i_{M\cten\bc}\fuse\id_{M\cten\bc}} (\brac*{M\cten\bc}\fuse
         \brac*{M^\vee\cten\bc})\fuse \brac*{M\cten\bc}\nonumber\\
    &\xrightarrow{\mathcal{A}^{-1}}\brac*{M\cten\bc}\fuse\brac*{\brac*{M^\vee\cten\bc}\fuse
    \brac*{M\cten\bc}}\nonumber\\
        &\xrightarrow{{\rm 1}\fuse e_{M\cten\bc}}\brac*{M\cten\bc}\fuse \brac*{\mmsl{u}{v}\cten\bc}\xrightarrow{r} \brac*{M\cten\bc}.
  \end{align}
  Applying the restriction functor to this composition of maps, using that
  induction is monoidal and discarding all but the \(C_p^{[i]}\brac*{M}\)
  summand then yields the following composition.
    \begin{align}
    \id=&C_{p}^{[i]}\brac*{M}\cten\mathcal{F}_{p} \xrightarrow{l^{-1}} \brac*{\mmN{u}{v}\cten \mathcal{F}_{0}}\fuse C_{p}^{[i]}\brac*{M}\cten\mathcal{F}_{p}\nonumber\\
    &\xrightarrow{?} (\brac*{C_{p}^{[i]}\brac*{M}\cten\mathcal{F}_{p}}\fuse \brac*{C_{-p}^{[i]}\brac*{M^\vee}\cten\mathcal{F}_{-p}})\fuse \brac*{C_{p}^{[i]}\brac*{M}\cten\mathcal{F}_{p}}\nonumber\\
      &\xrightarrow{\mathcal{A}^{-1}}\brac*{C_{p}^{[i]}\brac*{M}\cten\mathcal{F}_{p}}\fuse\brac*{\brac*{C_{-p}^{[i]}\brac*{M^\vee}\cten\mathcal{F}_{-p}}\fuse
        \brac*{C_{p}^{[i]}\brac*{M}\cten\mathcal{F}_{p}}}\nonumber\\
      &\xrightarrow{?}\brac*{C_{p}^{[i]}\brac*{M}\cten\mathcal{F}_{p}}\fuse
        \mmN{u}{v}\cten\mathcal{F}_{0}\xrightarrow{r}
        \brac*{C_{p}^{[i]}\brac*{M}\cten\mathcal{F}_{p}},
        \label{eq:Fnotfactored}
    \end{align}
    where the arrows marked by a question mark are to be determined.
    Next we use that \(\wtmod{\mmN{u}{v}\cten\heis{2\frac{u}{v}}}\) is braided
    equivalent to the Deligne tensor product
    \(\wtmod{\mmN{u}{v}}\delprod\wtmod{\heis{2\frac{u}{v}}}\) to factor out
    all the Fock spaces to obtain the composition
    \begin{align}
      \id=&C_{p}^{[i]}\brac*{M}\delprod\mathcal{F}_{p} \xrightarrow{l^{-1}}
            \brac*{\mmN{u}{v}\fuse
            C_{p}^{[i]}\brac*{M}}\delprod\brac*{\mathcal{F}_{0}\fuse \mathcal{F}_{p}}\nonumber\\
          &\xrightarrow{?} \brac*{\brac*{C_{p}^{[i]}\brac*{M}\fuse
            C_{-p}^{[i]}\brac*{M^\vee}}\fuse C_{p}^{[i]}\brac*{M}}\delprod
            \brac*{\brac*{\mathcal{F}_{p}\fuse\mathcal{F}_{-p}}\fuse \mathcal{F}_{p}}\nonumber\\
          &\xrightarrow{\mathcal{A}^{-1}}\brac*{C_{p}^{[i]}\brac*{M}\fuse
            \brac*{C_{-p}^{[i]}\brac*{M^\vee}\fuse C_{p}^{[i]}\brac*{M}}}\delprod
            \brac*{\mathcal{F}_{p}\brac*{\fuse\mathcal{F}_{-p}\fuse \mathcal{F}_{p}}}\nonumber\\
          &\xrightarrow{?}\brac*{C_{p}^{[i]}\brac*{M}\fuse \mmN{u}{v}}\delprod\brac*{
            \mathcal{F}_p\fuse \mathcal{F}_{0}}\xrightarrow{r}
            C_{p}^{[i]}\brac*{M}\delprod\mathcal{F}_{p}.
    \end{align}
    Note that each of the products of Fock spaces above is simple (and
    isomorphic to \(\mathcal{F}_p\)) and so each of the arrows in the
    composition above lies in a Deligne tensor product hom space of the form
    \(\hom_{\mmN{u}{v}}\brac*{A,B}\cten\hom_{\heis{2\frac{u}{v}}}\brac*{\mathcal{F}_p,\mathcal{F}_p}\). Since
    the Heisenberg hom space is \(1\)-dimensional, every \(f\in
    \hom_{\mmN{u}{v}}\brac*{A,B}\cten\hom_{\heis{2\frac{u}{v}}}\brac*{\mathcal{F}_p,\mathcal{F}_p}\)
    can be written as a tensor product \(f_N\cten f_H\),
    \(f_N\in\hom_{\mmN{u}{v}}\brac*{A,B}\), \(f_H\in
    \hom_{\heis{2\frac{u}{v}}}\brac*{\mathcal{F}_p,\mathcal{F}_p}\) as opposed
    to some linear combination of such tensor products. We can therefore
    factor out all the Heisenberg morphisms to obtain the composition
    \begin{align}
        \id=&C_{p}^{[i]}\brac*{M} \xrightarrow{l^{-1}} \mmN{u}{v}\fuse C_{p}^{[i]}\brac*{M}
              \xrightarrow{?} (\brac*{C_{p}^{[i]}\brac*{M}}\fuse C_{-p}^{[i]}\brac*{M^\vee})\fuse C_{p}^{[i]}\brac*{M}\nonumber\\
            &\xrightarrow{\mathcal{A}^{-1}}C_{p}^{[i]}\brac*{M}\fuse\brac*{C_{-p}^{[i]}\brac*{M^\vee}\fuse
              C_{p}^{[i]}\brac*{M}}
              \xrightarrow{?}C_{p}^{[i]}\brac*{M}\fuse \mmN{u}{v}\xrightarrow{r} C_{p}^{[i]}\brac*{M}.
      \end{align}
      The two arrows marked by a question mark lie in 
      \(\hom_{\mmN{u}{v}}\brac*{\mmN{u}{v},C_p^{[i]}(M)\fuse
        C_{-p}^{[i]}(M^\vee)}\cten\id\) and
      \(\id\cten\hom_{\mmN{u}{v}}\brac*{C_p^{[i]}(M^\vee)\fuse
    C_{-p}^{[i]}(M),\mmN{u}{v}}\), respectively. Thus there exists a pair of morphisms from these hom spaces and
  such that the composition yields
  the identity morphism \(C_p^{[i]}\brac*{M}\to C_p^{[i]}\brac*{M}\).
\end{proof}

\begin{theorem}
   For coprime \(u\ge2,\ v\ge1\), \(1\le r\le u-1\) and \(1\le s\le v-1\) and $q\in \RR$, let \(\NLs{q}{r}{s}\) be the simple highest weight
   \(\mmN{u}{v}\)-module \(\NL{q}{h_{r,s;q}}{3-\frac{6v}{u}}\) as given in the second
   part of \cref{thm:n2simpleclass}. Additionally consider the simple relaxed highest
   weight \(\mmsl{u}{v}\)-modules \(\slE{[q\frac{u}{v}]}{r}{s}\), where
   \([q\frac{u}{v}]=q\frac{u}{v}+2\ZZ\). Then for all
   \(1\le r',r''\le u-1\), \(1\le s',s''\le v-1\), \(q'\in \RR\),
   \(\frac{u}{v}q\notin \pm \lambda_{r,s,\frac{u}{v}}+2\ZZ\)
   the non-vanishing simple fusion rules between \(\NLs{q}{r}{s}, \NLs{q'}{r'}{s'}\)
   and between \(\slE{[q\frac{u}{v}]}{r}{s}, \slE{[q'\frac{u}{v}]}{r'}{s'}\) are
   \begin{align}
     \dim \binom{\NLs{q+q'}{r''}{s''}}{\NLs{q}{r}{s},\ \NLs{q'}{r'}{s'}}&
     = \dim \binom{\slE{[\brac*{q+q}'\frac{u}{v}]}{r''}{s''}}{\slE{[q\frac{u}{v}]}{r}{s},\ \slE{[q'\frac{u}{v}]}{r'}{s'}}
     =N^{(u)\,
       r''}_{r,r'}\brac*{N^{(v)\, s''}_{s,s'-1}+N^{(v)\, s''}_{s,s'+1}},\nonumber\\
     \dim \binom{\Pi \NLs{q+q'\pm 1}{r''}{s''}}{\NLs{q}{r}{s},\ \NLs{q'}{r'}{s'}}
     &=\dim \binom{\sigma^{\mp 1} \slE{[\brac*{q+q'\pm 1}\frac{u}{v}]}{r''}{s''}}{\slE{[q\frac{u}{v}]}{r}{s},\ \slE{[q'\frac{u}{v}]}{r'}{s'}}
       =N^{(u)\, r''}_{r,r'}N^{(v)\, s''}_{s,s'},
       \label{eq:intertwinerdims}
     \end{align}
   where the $N$-coefficients are defined in \eqref{eq:sl2fusionrules}.
   \label{thm:intertwinerdims}
 \end{theorem}
The equality of the dimensions of \(\slt\) and \(\Nt=2\) intertwiner
spaces \eqref{eq:intertwinerdims} is just
\cref{thm:homspaceeq} applied to the relaxed highest weight modules
\(\slE{[q\frac{u}{v}]}{r}{s}\). The computation of these dimensions will form the focus of the next two sections. Note that when \(q,q'\) in \eqref{eq:intertwinerdims}
are such that \(\frac{u}{v}(q+q'+\epsilon),\notin \pm \lambda_{r,s,\frac{u}{v}}+2\ZZ\),
\(\epsilon=-1,0,1\), then \eqref{eq:intertwinerdims} together with the
projectivity of the simple modules $\sigma^\ell \slE{\mu}{r}{s}$ established
in \cite{ACK} implies that the tensor
products of \(\slE{[q\frac{u}{v}]}{r}{s},\slE{[q'\frac{u}{v}]}{r'}{s'}\) and of
\(\NLs{q}{r}{s}, \NLs{q'}{r'}{s'}\) are semisimple 
with direct sum decomposition given
by
\begin{align}
  \slE{[q\frac{u}{v}]}{r}{s}\fuse
  \slE{[q'\frac{u}{v}]}{r'}{s'}&=\bigoplus_{r'',s''} N^{(u)\, r''}_{r,r'} \brac*{N^{(v)\,
  s''}_{s,s'-1}+N^{(v)\, s''}_{s,s'+1}} \slE{[(q+q')\frac{u}{v}]}{r''}{s''}\nonumber\\
  &\quad\oplus \bigoplus_{r'',s''}N^{(u)\, r''}_{r,r'}N^{(v)\, s''}_{s,s'}\brac*{\sigma^{-1} \slE{[\brac*{q+q'+1}\frac{u}{v}]}{r''}{s''}\oplus
    \sigma \slE{[\brac*{q+q'- 1}\frac{u}{v}]}{r''}{s''}},\nonumber\\
  \NLs{q}{r}{s}\fuse \NLs{q'}{r'}{s'}&=\bigoplus_{r'',s''} N^{(u)\, r''}_{r,r'}\brac*{N^{(v)\,
                                   s''}_{s,s'-1}+N^{(v)\,
                                   s''}_{s,s'+1}}\NLs{q+q'}{r'}{s'}\nonumber\\
                           &\quad \oplus\bigoplus_{r'',s''}N^{(u)\, r''}_{r,r'}N^{(v)\, s''}_{s,s'}
                             \brac*{\Pi \NLs{q+q'+ 1}{r''}{s''}\oplus \Pi \NLs{q+q'- 1}{r''}{s''}}.
\end{align}
By associativity these tensor products are generated by
\(\slE{[q\frac{u}{v}]}{1}{1},\slE{[q\frac{u}{v}]}{2}{1}\) and by
\(\NLs{q}{1}{1},\NLs{q}{2}{1}\). Moreover, by specialising \eqref{eq:LEfusion}, we
observe $\slE{[q\frac{u}{v}]}{2}{1}\cong L_2\fuse \slE{[q\frac{u}{v}]-1}{1}{1}$. 
Hence all of the relaxed \(\slt\) fusion rules can be obtained from the
$\slE{[q\frac{u}{v}]}{ 1}{1}$ fusion rules and by tensoring those with \(L_2\).  That is, \cref{thm:intertwinerdims} holds if and only if the
following lemma holds.
\begin{lemma}
  For coprime \(u\ge2,\ v\ge1\), and \(1\le r,r'\le u-1\), \(1\le s,s'\le v-1\), \(q,q'\in \RR\),
  \(\frac{u}{v}q'\notin \pm \lambda_{r,s,\frac{u}{v}}+2\ZZ\). If \(q\notin
  \pm\lambda_{1,1,\frac{u}{v}}+2\ZZ\),
  \begin{align}
    \dim \binom{\NLs{q+q'}{r'}{s'}}{\NLs{q}{1}{1},\ \NLs{q'}{r}{s}}&
    = \dim \binom{\slE{[\brac*{q+q'}\frac{u}{v}]}{r'}{s'}}{\slE{[q\frac{u}{v}]}{1}{1},\ \slE{[q'\frac{u}{v}]}{r}{s}}
    =\delta_{r,r'} \brac*{\delta_{s-1,s'}+\delta_{s+1,s'}},\nonumber\\
    \dim \binom{\Pi \NLs{q+q'\pm 1}{r'}{s'}}{\NLs{q}{1}{1},\ \NLs{q'}{r}{s}}
    &=\dim \binom{\sigma^{\mp 1} \slE{[\brac*{q+q'\pm 1}\frac{u}{v}]}{r'}{s'}}{\slE{[q\frac{u}{v}]}{1}{1},\ \slE{[q'\frac{u}{v}]}{r}{s}}
    =\delta_{r,r'}\delta_{s,s'},\nonumber\\
    \dim \binom{X}{\NLs{q}{1}{1},\ \NLs{q'}{r}{s}}&
    = \dim \binom{Y}{\slE{[q\frac{u}{v}]}{1}{1},\ \slE{[q'\frac{u}{v}]}{r}{s}}
    =0,
    \label{eq:11rule}
  \end{align}
  where \(X\) and \(Y\) are any simple modules other than those appearing
    in the first two lines of \eqref{eq:11rule}.
  \label{thm:1121rules}
\end{lemma}
We will prove \cref{thm:1121rules} by computing upper bounds for the
dimensions in \eqref{eq:11rule} on the \(\Nt=2\)
side using Zhu's algebra and its bimodules in \cref{sec:upperbounds} (see
\cref{thm:upperbounds}) and lower bounds on the \(\slt\) side in
\cref{sec:lowerbounds} (see Propositions \ref{thm:lowe1stline} and
\ref{thm:lowerb2}) by explicitly constructing certain intertwining operators.

\section{Upper bounds \texorpdfstring{for \cref{thm:1121rules}}{} via Zhu algebras}
\label{sec:upperbounds}

We recall some definitions and results from Zhu and Frenkel's work \cite{FZ, Z1} and the generalisation by Kac and Wang \cite{KW} to vertex
  superalgebras where all odd vectors 
  with respect to the \(\ZZ_2\) grading have half odd conformal weights. We then apply these results to \(\mmN{u}{v}\) to compute upper
bounds for \cref{thm:1121rules} in \cref{thm:upperbounds}.

\subsection{Zhu algebras, Zhu modules and intertwining operators}
\label{sec:Zhuthy}
For any super vector space \(V\), we denote the subspaces of even and odd vectors
  by \(V^{\bar0}\) and \(V^{\bar1}\) respectively.
For a vertex operator superalgebra $V$ in which the space of odd vectors
coincides with the space of vectors of half odd integer conformal weight, the
Zhu algebra is constructed from the two bilinear operations on $V$
characterised for homogenous $a,b \in V$ by
\begin{align}
a*b&=\begin{cases} {\rm Res}_z \brac*{ Y(a,z) \frac{(1+z)^{{\rm wt }a }}{z} b},
  & a,b\in V^{\bar0}, \\ 
  0, & a\in V^{\bar1}\ \text{or}\ b\in V^{\bar1},
\end{cases}  \nonumber\\
a\circ b&=
\begin{cases} {\rm Res}_z \brac*{ Y(a,z) \frac{(1+z)^{{\rm wt }a }}{z^2} b},
  & a\in V^{\bar0}, \\ 
  {\rm Res}_z \brac*{ Y(a,z) \frac{(1+z)^{{\rm wt }a -\frac{1}{2} }}{z} b},
  & a\in V^{\bar1}. 
\end{cases}
\end{align} 

\begin{prop}[Zhu \cite{Z1}; Kac, Wang \cite{ KW}]
  Let \(V\) be a \vosa{} for which the space of odd vectors and the space of
  vectors of half odd conformal weight coincide and define \(O(V)=\cspn{a\circ
    b\st a,b\in V}\). Then \(A(V)=V/O(V)\) is an associative unital algebra
  with multiplication given by \(\ast\). The identity element is
  \([\vac]=\vac+O(V)\), the equivalence class of the vacuum vector, and the
  class \([\omega]\) of the conformal vector is central.
  \label{thm:AV}
\end{prop}

Note that for odd \(a\in V^{\bar{1}}\), \(a\circ \vac=a\) and hence
  \(V^{\bar1}\subset O(V)\).
We next recall some standard Zhu algebra results. 

\begin{lemma}[Zhu \cite{Z1}; Kac, Wang \cite{KW}]
  Let $V$ be a vertex operator superalgebra. For all
  homogeneous elements
  $a,b\in V$ 
  and  $m\geq n\geq 0$ we have that 
  \begin{enumerate}
  \item $\left( L_{-1}+L_0 \right) a \in O \left( V \right)$, 
     $a\in V^{\bar{0}}$,
  \item ${\rm Res}_x\left(Y \left( a,x \right) \frac{\left( 1+x
        \right)^{{\rm wt}a+n}}{x^{2+m}}b\right) \in O \left( V
    \right)$, \(a\in V^{\bar0}\), \label{eq:ru}
  \item ${\rm Res}_x\left(Y \left( a,x \right) \frac{(1+x)^{{\rm
            wt}a+n-\frac{1}{2}}}{x^{1+m}}b\right) \in O \left( V
    \right)$, \(a\in V^{\bar1}\), 
  \item $a*b=Res_x\left( Y \left( b,x \right) \frac{(x+1)^{{\rm
              wt}b}-1}{x}a\right) + O(V)$,
    \(a,b\in V^{\bar0}\).
  \end{enumerate}
  \label{thm:powers}
\end{lemma}

\begin{theorem}[Zhu \cite{Z1}; Kac, Wang \cite{KW}]\label{M0}
  Let $V$ be a vertex operator superalgebra and let $M$ be a 
  $V$-module. Define the \emph{top space} of \(M\) to be
  \(M^{top}=\set{m\in M\st v_n m=0\ \forall v\in V, \forall n>0}\).
  \begin{enumerate}
  \item Any element \([a]\in A(V)\) acts on \(M^{\mathrm{top}}\) by the zero mode
    \(a_0\) where $Y(a,z)=\sum_{n\in \ZZ}a_n z^{-wt(a)-n}$,
    hence giving \(M^{\mathrm{top}}\) the structure of a left
    \(A(V)\)-module. Further, if \(M\) is simple as a module over \(V\), then
    \(M^{\mathrm{top}}\) is simple over \(A(V)\).
  \item Any left \(A(V)\)-module \(M^{\mathrm{top}}\) can be induced to a \(V\) module
    with conformal weights bounded below and with \(M^{\mathrm{top}}\) contained
    in the top space. If \(M^{\mathrm{top}}\) is simple over \(A(V)\), then the
    induction admits a unique simple quotient.
  \item Simple \(A(V)\) modules are in bijection with simple \(V\)-modules
    with conformal weights bounded below.
    \end{enumerate}
  \end{theorem}

  Ideals in a \vosa{} and in its
  corresponding Zhu algebra are interrelated as follows.

\begin{lemma}[Zhu \cite{Z1}; Kac, Wang \cite{KW}]\label{thm:qu}
  Let $V$ be a vertex operator superalgebra and $I$ a $\mathbb{Z}_2$-graded
  ideal of $V$ (with grading consistent with that of $V$). 
  Assume that $\vac\notin I, \omega \notin I$. Then the Zhu algebra of the quotient
  $A \left( V/I \right)$ is isomorphic to the quotient
  $A \left( V \right)/\left[I\right]$ where $\left[I\right]$ is the image of $I$ in $A \left( V \right)$. 
\end{lemma}

If $M$ is a $V$-module, we define an $A \left( V \right)$-bimodule $A \left( M \right)$ as follows.
First, we set the left and right actions of $V$ on $M$ by
\begin{align}
  &a*m= \begin{cases} \res_z \brac*{Y(a,z)\frac{(1+z)^{{\rm wt }a}}{z} m },
    &\text{if}\ a\in V^{\bar{0}}, \\
    0 &\text{if}\ a\in V^{\bar1},
  \end{cases}  \nonumber\\
  &m*a=\begin{cases} \res_z \brac*{Y(a,z)\frac{(1-z)^{{\rm wt }a-1}}{z} m}, &
    \text{if}\ a\in V^{\bar{0}}, \\
    0 &\text{if}\ a\in V^{\bar1}.
  \end{cases}
  \label{eq:Zhuaction}
\end{align}
where \(m\in M\).
Next, we let $O \left( M \right)$ be the subspace of $M$ linearly spanned by
elements of the form
\begin{align}
  &\res_z \brac*{Y \left( a,z \right)\frac{\left( 1+z \right)^{{\rm wt
        }a}}{z^2} m },\qquad a\in V^{\bar0},
  \nonumber\\
  &\res_z \brac*{Y(a,z)\frac{\left( 1+z \right)^{{\rm wt }a-\frac{1}{2}}}{z}
    m},\qquad
  a\in V^{\bar{1}},
\end{align}
where \(m\in M\).

\begin{prop}[Zhu \cite{Z1}; Kac, Wang \cite{KW}] 
  Let $V$ be a vertex operator superalgebra and $M$ a $V$-module. Then
  \(A(M)=M/O(M)\) is an \(A(V)\)-bimodule with left and right actions given by
  \eqref{eq:Zhuaction}. For any homogeneous elements
  \(a\in V\), \(m\in M\)
  and $p\geq q\geq 0$ we have that:
  \begin{enumerate}
  \item $\res_x\left(Y(a,x)\frac{(1+x)^{{\rm wt}a+q
        }}{x^{2+p}}m\right) \in O \left( M \right)$,
    \(a\in V^{\bar0}\), and
  \item $\res_x\left(Y(a,x)\frac{(1+x)^{{\rm wt}a+q-\frac{1}{2}
        }}{x^{1+p}}m\right) \in O \left( M \right)$,
    \(a\in V^{\bar1}\).
  \end{enumerate}
  \label{thm:Omod}
\end{prop}

As in \cref{thm:qu} the Zhu bimodule structure is compatible with quotients by
submodules as described in the following result. 

\begin{lemma}[Zhu \cite{Z1}; Kac, Wang \cite{KW}]
  Let $V$ be a vertex operator superalgebra, $M$ a $V$-module and \(\widetilde{M}\subset{M}\) a submodule.
  \begin{enumerate}
  \item Then $A\brac{M/\widetilde{M}}\cong A(M)/A\brac{\widetilde{M}}$, where $A\brac{\widetilde{M}}$ denotes the image of $\widetilde{M}$ under the projection $M\twoheadrightarrow A(M).$
  \item If $I$ is an ideal of $V$, $\vac\notin I, \omega \notin I$ and $I.M\subset \widetilde{M}$, then the $A(V/I)$-module $A(M)/A\brac{\widetilde{M}}$ is isomorphic to $A\brac{M/\widetilde{M}}.$
  \end{enumerate}
  \label{thm:muz}
\end{lemma}

We will denote the vector space of intertwining operators of type
$\binom{W_3}{W_1 \ \ W_2}$ by the same symbol $\binom{W_3}{W_1 \ \ W_2}$ and
for the moment postpone the definition of intertwining operators until they
are explicitly used in \cref{sec:lowerbounds}. One of the main benefits of Zhu
bimodules is that they provide a means for computing dimensions of spaces of
intertwining operators.

\begin{lemma}[Huang, Yang {\cite[Prop 5.8]{HYzhudims}}] 
  Let $V$ be a vertex operator superalgebra and let $M_1, M_2$ and $M_3$ be 
  $V$-modules with \(M_2\) and \(M_3'\) (the contragrediant of \(M_3\)) generated from their top spaces. 
  Then 
  \begin{align}
    \dim \binom{M_3}{M_1 \ \ M_2}\leq \dim\hom_{A(V)}\left( A(M_1)\cten_{A(V)} {M_2^\mathrm{top}}, {M_3^\mathrm{top}}\right). 
  \end{align}
  \label{lemma:Zhu}
\end{lemma}
The relationship between spaces of intertwining operators and hom spaces was, of
  course, also considered in Zhu \cite{Z1}; Kac, Wang \cite{KW}, but not with
  the assumptions that we require here.

\subsection{Determining the Zhu algebra for $\mmN{u}{v}$}
Let $c$ be any complex number and let $\uN{c}$ denote the universal $\Nt=2$
vertex operator superalgebra of central charge $c$ as defined in Proposition
\ref{thm:univn=2}.

\begin{theorem}[Eholzer, Gaberdiel \cite{EG}] \leavevmode
  \begin{enumerate}
  \item For any \(c\in\CC\), the Zhu algebra $A(\uN{c})$ is isomorphic to the
    polynomial algebra $\mathbb{C} [\Delta, \eta]$, where \(\Delta\) is the
    image of the conformal vector \([\omega]\) and \(\eta\) is the image of
    \([J_{-1}\vac]\).
  \item  For \(u,v\) coprime \(u\ge2\), \(v\ge1\), and
    \(c=3-\frac{6v}{u}\), the image in \(A(\uN{c})\) of the non-trivial
    maximal ideal of $\uN{c}$ generated by a singular vector
    \(\chi\in \uN{c}\) of conformal weight \((u-1)v\) and \(J\)-weight 0 is
    \(\ang{p_1(\Delta,\eta),p_2(\Delta,\eta)}\), where \(p_1\) is the image of
    \(\chi\) and \(p_2\) is the image of \(G^+_{-\frac12}G^-_{-\frac12}\chi\).
  \item The Zhu algebra of \(\mmN{u}{v}\) admits the presentation
    \begin{equation}
      A\brac*{\mmN{u}{v}}\cong \frac{\CC[\Delta.\eta]}{\ang{ p_1(\Delta,\eta),p_2(\Delta,\eta)}}.
    \end{equation}
  \end{enumerate}
  \label{thm:A0}
\end{theorem}

The goal of this section is to compute the polynomials \(p_1,p_2\). 
\begin{theorem} \label{thm:Zhumm}
  Fix $c=3-\frac{6v}{u}$ for coprime $u\geq 2$, $v\geq 1$  and let
  $\calJ$ denote
  the maximal ideal $\calJ \subset \uN{c}$. Then,  as an ideal
  of $A(\uN{c})\cong \mathbb{C}[\Delta,\eta]$ under the isomorphism given in
  \cref{thm:A0}, $[\calJ]$ is generated by the following two polynomials 
  \begin{align}
    p_1(\eta,\Delta)&=f_u\brac*{\eta,\Delta,\frac{u}{v}}
    \prod_{(r,s)\in K(u,v)} \left(\Delta- h_{r,s;\eta} \right), 
    \nonumber\\
    p_{2}(\eta,\Delta)&=\left(f_u\brac*{\eta+\frac{2v}{u}, \Delta -\eta
        -\frac{v}{u},\frac{u}{v}}
      -f_u\brac*{\eta,\Delta,\frac{u}{v}}\right)\left(
      2\Delta-\eta\right)\prod_{(r,s)\in K(u,v)} \left(\Delta-
      h_{r,s,\eta}\right),
  \end{align}
  where the $h_{r,s;\eta}$ are the conforma weights defined in \eqref{eq:N=2confweights}, $K(u,v)$ is the set of pairs $(r,s) \in \{1, \dots,
  u-1\}\times \{ 1, \dots v-1\}$ satisfying the additional constraint \(
  vr+us< uv\), and the polynomials \(f_n\brac*{x,y,z},\ n\ge2\)
  are defined recursively by
  \begin{align}
    f_{n+2}\brac*{x,y,z}
    =\frac{(2n+1)xz}{(n+1)^2} f_{ n+1}\brac*{x,y,z} 
    -\frac{4yz+x^2z^2-(n-1)(n+1)}{(n+1)^2}f_{n}\brac*{x,y,z},
    \label{eq:gtu+2v}
  \end{align}   
  with 
  \begin{align}
    f_2(x,y,z) &=xz,&
    f_3(x,y,z) &=\frac{x^2z^2}{2}-yz.
    \label{eq:gt2v}
  \end{align}
  In particular, the Zhu algebra for the minimal model superconformal algebra at central charge $c=3-\frac{6v}{u}$ admits the presentation
  \begin{align}
    A(\mmN{u}{v})\cong \frac{\mathbb{C}[\Delta,\eta]}{\ang*{ p_1(\Delta,
        \eta), p_2(\Delta, \eta)}}.
    \label{eq:Zhualgpres}
  \end{align}
  \label{Thm:ZhuJgen}
\end{theorem}

\begin{proof}
  Eholzer and Gaberdiel established in \cite{EG} that, under the isomorphism
  $\varphi$ in \cref{thm:A0}, the ideal $[\calJ]\subset \CC[\Delta, \eta]$ is generated by the images of the singular vector \(\chi\), which
    generates \(\calJ\) and its descendant
    \(G^+_{-\frac12}G^-_{-\frac12}\chi\). Thus all that is left to show is that
  these images are equal to the polynomials \(p_1,p_2\).
For this we use the coset realisation in \cref{thm:coset} and the relations
  between singular vectors in \cref{thm:coset-and-sing-vectors}(3) for the
  vacuum module $\usl{\frac{u}{v}-2}$. Let $w_{u,v}$ denote the singular
  vector generating the maximal $\slta$ submodule in
  $\usl{\frac{u}{v}-2}$. Then $f_0^{u-1}w_{u,v}$ is a relaxed singular vector
  in $\usl{\frac{u}{v}-2}$ whose image in the Zhu algebra of
  $\usl{\frac{u}{v}-2}$, up to a non-zero multiple, was shown in \cite{RW} to
  be given by
  \begin{align}
    I\left([h], [T_{\slt}] \right)=\prod_{(r,s)\in K(u,v)}\left([T_{ \slt}]-\Delta^{\aff}_{r,s,t}\right)g_{u,v} \left( [h], [T_{\slt}] \right), \label{eq:Isl2data}
  \end{align}
 where $[h]$ and  $[T_{ \slt}]$ denote the images of $h_{-1}\vac$ and and the
 conformal vector $T_{ \slt}=\omega_{\frac{u}{v}-2}=\frac{1}{2\frac{u}{v}} \left(\frac{1}{2}h_{-1}^2-e_{-1}f_{-1}- f_{-1}e_{-1}\right)\vac$ in $A(\usl{\frac{u}{v}-2})\cong \UEA{\slt}$ and $g_{u,v}(x,y)$ are polynomials recursively defined in \cite{RW} by
 \begin{align}
&g_{2,v}(x,y)=x,\nonumber\\
&g_{3,v}(x,y)=\frac{3}{4}x^2-ty,
\end{align}
and 
\begin{align}
g_{u+2,v}(x,y)=\frac{(2u+1)x}{(u+1)^2}g_{u+1,v}(x,y)-\frac{4ty-(u-1)(u+1)}{(u+1)^2}g_{u,v}(x,y),
\end{align}
where $t$ is a parameter that will be later specialised to $t=\frac{u}{v}$.
By \cref{thm:coset-and-sing-vectors}.(3) we have that
\begin{align}
  \phi(\ket{0} \cten \chi)=f_0^{u-1}w_{u,v}\cten \Omega,  \label{eq:rhwidentification}
\end{align}
where \(\phi\) is the coset embedding of \cref{thm:coset}, $\ket{0}$ is the vacuum vector in $\heis{\frac{u}{v}}$ and $\Omega$ is the vacuum vector in $\bc$.

In particular, the zero modes $(\ket{0}\cten\chi)_0=\id\cten\chi_0$ and
$(f^{u-1}w_{u,v}\cten \Omega)_0=(f^{u-1}w_{u,v})_0\cten \id$ must act
equally on any \(\usl{\frac{u}{v}-2}\cten \bc\) module. Thus, for an
\(\usl{\frac{u}{v}-2}\) relaxed highest weight vector \(m(\lambda,\Delta^{\aff})\) of
\(\slt\) weight \(\lambda\) and conformal weight \(\Delta^{\aff}\) we have
\begin{align}
  \chi_0 \brac*{m(\lambda,\Delta^{\aff})\cten \Omega} = \brac*{(f^{u-1}w_{u,v})_0
  m(\lambda,\Delta^{\aff})}\cten \Omega = I(\lambda,\Delta^{\aff}) m(\lambda,\Delta^{\aff})\cten \Omega,
\end{align}
where in the second step we have used that Zhu's algebra models the algebra of
zero modes acting on relaxed highest weight vectors.
Recall that the \(\slt\) data \((\lambda,\Delta^{\aff})\) corresponding to the \(J_0\)
eigenvalue \(\eta\)
and \(\Nt=2\) conformal weight \(\Delta\) are interrelated via \(\lambda=t \eta\) and
\(\Delta^{\aff}=\Delta+\frac{t}{4}\eta^2\). Thus
\(p_1(\eta,\Delta)=I(t\eta,\Delta+\frac{t}{4}\eta^2)\), which matches the
formula given in the Theorem after identifying \(f_u(\eta,\Delta,t)=g_{u,v}(t\eta,\Delta+\frac{t}{4}\eta^2)\)  and $t=\frac{u}{v}$.

Next we compute the image of \(G^+_{-\frac12}G^-_{-\frac12}\chi\) in the
  Zhu algebra, again via the coset realisation. Following on from the
  identification \eqref{eq:rhwidentification} we have
  \begin{align}
    \phi(\ket{0} \cten
    G^+_{-\frac12}G^-_{-\frac12}\chi)&=\phi(G^+_{-\frac12}G^-_{-\frac12})f_0^{u-1}w_{u,v}\cten
                                           \Omega\nonumber\\
    &=\frac{2}{t}\brac*{e_0\cten c_{-\frac12}+e_{-1}\cten
      c_{\frac12}}f_0f_0^{u-1}w_{u,v}\cten b_{-\frac12}\Omega\nonumber\\
    &=-\frac{2}{t}e_0f_0f_0^{u-1}w_{u,v}\cten Q+\frac{2}{t}e_{-1}f_0f_0^{u-1}w_{u,v}\cten\Omega,
  \end{align}
  where \(Q\) is the Heisenberg vector within \(\bc\). Applying both
    sides of the above identity to the same test vector
    \(m(\lambda,\Delta^{\aff})\cten \Omega\) as before yields
  \begin{multline}
    \brac*{G^+_{-\frac12}G^-_{-\frac12}\chi}_0 m(\lambda,\Delta^{\aff})\cten \Omega
    \\=
    -\frac{2}{t}\brac*{e_0f_0f_0^{u-1}w_{u,v}}_0\cten Q_0
    \brac*{m(\lambda,\Delta^{\aff})\cten
    \Omega}+\frac{2}{t}\brac*{e_{-1}f_0f_0^{u-1}w_{u,v}}_0 \brac*{m(\lambda,\Delta^{\aff})\cten
    \Omega}.
\end{multline}
Note that the first summand vanishes due to \(Q_0\Omega=0\), while the second
can be simplified using the Zhu relations \eqref{eq:Zhuaction} to give
\begin{align}
  \brac*{G^+_{-\frac12}G^-_{-\frac12}\chi}_0 &m(\lambda,\Delta^{\aff})\cten \Omega
  =
  \frac{2}{t}\brac*{f\ast\brac*{f_0^{u-1}w_{u,v}}\ast
  e-\brac*{f_0^{u-1}w_{u,v}}\ast f \ast e}_0 m(\lambda,\Delta^{\aff})\cten
    \Omega\nonumber\\
  &=\frac{2}{t}\brac*{f_0 \brac*{f_0^{u-1}w_{u,v}}_0
  e_0-\brac*{f_0^{u-1}w_{u,v}}_{0}f_0 e_0} m(\lambda,\Delta^{\aff})\cten
    \Omega\nonumber\\
  &=\frac{2}{t}\brac*{I\brac*{\lambda+2,\Delta^{\aff}}-I\brac*{\lambda,\Delta^{\aff}}}\brac*{2t\Delta^{\aff}-\frac{\lambda^2}{2}-\lambda}
  m(\lambda,\Delta^{\aff})\cten \Omega,
\end{align}
where for the final equality we used that \(m(\lambda,\Delta^{\aff})\) is a
relaxed highest weight vector and hence \(f_0
e_0m(\lambda,\Delta^{\aff})=\brac*{2t
  L_0-\frac{1}{2}h_0^2-h_0}m(\lambda,\Delta^{\aff})\). The \(\slt\) and
\(\Nt=2\) weights are again interrelated by \(\lambda=t \eta\) and
\(\Delta^{\aff}=\Delta+\frac{t}{4}\eta^2\) and thus the eigenvalue of
\(\brac*{G^+_{-\frac12}G^-_{-\frac12}\chi}_0\) is the polynomial \(p_2(\eta,\Delta)\).

\end{proof}
 Using \cref{thm:Zhumm} we can give a new proof of the classification of $\mmN{u}{v}$-irreducible modules obtained by Adamovi{\'c} in \cite{A1}.
\begin{theorem}
  A complete set of inequivalent simple modules over the Zhu algebra
  \eqref{eq:Zhualgpres} is given by the 1-dimensional vector spaces on which
  \((\eta,\Delta)\) act as any of the 
  pairs of scalars within the following sets.
  \begin{enumerate}
  \item \(\set*{\brac*{\frac{pv}{u},h_{r,0;\frac{pv}{u}}}\st
      1\le r\le u-1,\ 1-r\le p\le r-1, p+r \equiv 1 \pmod{2}}\),
  \item \(\set*{\brac*{q,h_{r,s;q}}\st 1\le r\le u-1,\ 1\le s\le v-1,\ vr+us<uv,\ q\in \CC}\).
  \end{enumerate}
  where \(h_{r,s;q}\) is given in \eqref{eq:N=2confweights}.
  \label{thm:newclass}
 \end{theorem}

  \begin{proof}
    The theorem follows by showing that the pairs of scalars above constitute
    all the simultaneous zeros of the polynomials \(p_1(\eta,\Delta)\),
    \(p_2(\eta,\Delta)\) of \cref{Thm:ZhuJgen}. Note first that both
    polynomials contain the divisors \((\Delta-h_{r,s;\eta})\). These
    correspond to the second set of pairs \(\brac*{q,h_{r,s;q}}\) above.

    If a simultaneous zero of \(p_1,\ p_2\) is not of the form
    \(\brac*{q,h_{r,s;q}}\), then it must necessarily be a simultaneous zero of
    \(f_u\brac*{\eta,\Delta,\frac{u}{v}}\) and
    \(f_u\brac*{\eta,\Delta,\frac{u}{v}}-f_u\brac*{\eta+\frac{2v}{u},\Delta-\eta-\frac{v}{u},\frac{u}{v}}\),
    or of \(f_u\brac*{\eta,\Delta,\frac{u}{v}}\) and \((2\Delta-\eta)\). Using
    B\'ezout's Theorem, 
    which in its simplest
    form states that two bivariate
    polynomials can have at most as many simultaneous zeros as the product of
    their degrees we will show that there are at most \(\frac{u(u-1)}{2}\)
    such simultaneous zeros. This is precisely the number of pairs of scalars
    in the first set above. After that has been established, all that will
    remain to be shown is that the above
    pairs are indeed simultaneous zeros, which we will show inductively.

    We start by giving on upper bound on the simultaneous zeros of \(f_u\brac*{\eta,\Delta,\frac{u}{v}}\) and
    \(f_u\brac*{\eta,\Delta,\frac{u}{v}}-f_u\brac*{\eta+\frac{2v}{u},\Delta-\eta-\frac{v}{u},\frac{u}{v}}\).
    From the given formulae for \(f_2\) and \(f_3\), and the recursion
    relation for \(f_u\), we see that \(f_u\brac*{\eta,\Delta,\frac{u}{v}}\)
    has degree \(u-1\) in \(\eta\) and degree at most \(\frac{u-1}{2}\) in
    \(\Delta\), in particular the monomial term with the largest power of
    \(\eta\) has no factor of \(\Delta\) and the monomial term with the
    largest power of \(\Delta\)
    has no factor of \(\eta\). Further the shifted version
    \(f_u\brac*{\eta+\frac{2v}{u},\Delta-\eta-\frac{v}{u},\frac{u}{v}}\)
    preserves the power of \(\eta\) in its first argument and \(\Delta\) in
    its second, so it too will have degree \(u-1\) in \(\eta\) and degree at most \(\frac{u-1}{2}\) in
    \(\Delta\), with the monomial term with the largest power of
    \(\eta\) having no factor of \(\Delta\) and the monomial term with the
    largest power of 
    \(\Delta\) having no factor of \(\eta\). Crucially, these respective
    greatest 
    monomial terms of \(f_u\brac*{\eta,\Delta,\frac{u}{v}}\) and
    \(f_u\brac*{\eta+\frac{2v}{u},\Delta-\eta-\frac{v}{u},\frac{u}{v}}\) have
    identical coefficients so the difference
    \(\tilde{f}_u(\eta,\Delta)=f_u\brac*{\eta,\Delta,\frac{u}{v}} -
    f_u\brac*{\eta+\frac{2v}{u},\Delta-\eta-\frac{v}{u},\frac{u}{v}}\) has
    degree at most \(u-2\) in \(\eta\) and at most \(\frac{u-1}{2}-1\) in
    \(\Delta\). So if we introduce an additional variable \(\gamma\), which is
    to take on the role of a square root of \(\Delta\), then
    \(f_u\brac*{\eta,\gamma^2,\frac{u}{v}}\) is a polynomial in
    \(\eta,\gamma\) of degree \(u-1\) and \(\tilde{f}_u\brac*{\eta,\gamma^2}\)
    is a polynomial of degree \(u-2\) in \(\eta,\gamma\), that is, replacing
    all the \(\Delta\) by squares does not increase the total degrees of the polynomials. Thus by B\'ezout's
    Theorem they have at most \((u-1)(u-2)\) simultaneous zeros. However,
    these zeros come in pairs, as if \((\eta_0,\gamma_0)\) is a zero, then
    \((\eta_0,-\gamma_0)\) is too, yet both of these zeros correspond to the
    same zero in \(\eta,\Delta\) given by \((\eta_0,\gamma_0^2)\). So \(f_u\)
    and \(\tilde{f}_u\) have at most \(\frac{(u-1)(u-2)}{2}\) simultaneous
    zeros in \(\eta,\Delta\). B\'ezout's Theorem also immediately bounds the number of
    simultaneous zeros of \(f_u\) and \(\brac*{2\Delta-\eta}\) from above by
    \(u-1\). Thus there are at most
    \(\frac{(u-1)(u-2)}{2}+u-1=\frac{u(u-1)}{2}\) simultaneous zeros of \(f_u\brac*{\eta,\Delta,\frac{u}{v}}\) and
    \(f_u\brac*{\eta+\frac{2v}{u},\Delta-\eta-\frac{v}{u},\frac{u}{v}}(2\Delta-\eta)\).

    Next we show by induction that all pairs in the first set listed in the theorem are indeed
    simultaneous zeros of \(f_u\brac*{\eta,\Delta,\frac{u}{v}}\) and
    \(f_u\brac*{\eta+\frac{2v}{u},\Delta-\eta-\frac{v}{u},\frac{u}{v}}(2\Delta-\eta)\).
    For our base case we need to consider \(u=2\) and \(u=3\). First \(u=2\)
    and \(v\) odd. Then
    \begin{equation}
      f_2\brac*{\eta,\Delta,\frac{2}{v}}=\eta\frac{2}{v}=0\quad\text{and}\quad
      f_2\brac*{\eta+v,\Delta-\eta-\frac{v}{2},\frac{2}{v}}(2\Delta-\eta)=\brac*{\eta+v}\frac{2}{v}\brac*{2\Delta-\eta}=0.
    \end{equation}
    The first equation requires \(\eta=0\) and then the second can only hold
    if \(\Delta=0\), which reproduces the first set of zeros for \(u=2\).
    Similarly for \(u=3\) and \(v\) not a multiple of \(3\), one can easily
    verify that the three candidate simultaneous zeros for \(u=3\) solve the equations
    \begin{align}
      f_3\brac*{\eta,\Delta,\frac{3}{v}}&=\eta^2\frac{9}{2v^2}-\frac{3}{v}\Delta=0,\nonumber\\
      f_3\brac*{\eta+\frac{2v}{3},\Delta-\eta-\frac{v}{3},\frac{3}{v}}(2\Delta-\eta)&=
        \brac*{\brac*{\eta+\frac{2v}{3}}^2\frac{9}{2v^2}-\frac{3}{v}\brac*{\Delta-\eta-\frac{v}{3}}}\brac*{2\Delta-\eta}=0.
    \end{align}

    Next for \(u\ge 3\) and \(v\) coprime to \(u\) assume that the pairs
    \((\frac{pv}{u},h_{r,0;\frac{pv}{u}})\), \(1\le r\le k-1\), \(1-r\le p\le
    r-1\), \(p+r\equiv 1\pmod{2}\) are zeros of
    \(f_k(\eta,\Delta,\frac{u}{v})\) for \(2\le k<u\). Then the recursion
    relation \eqref{eq:gtu+2v} implies that
    \((\frac{pv}{u},h_{r,0;\frac{pv}{u}})\) with \(r\le u-3\) is a zero,
    because it is a zero of \(f_{u-1}\) and \(f_{u-2}\). For \(r=u-2\), note
    that plugging \((\frac{pv}{u},h_{u-2,0;\frac{pv}{u}})\) into the recursion
    relation \eqref{eq:gtu+2v} leads to the coefficient in front of
    \(f_{u-2}\) vanishing and hence \((\frac{pv}{u},h_{u-2,0;\frac{pv}{u}})\)
    is a zero of \(f_u\) because it is a zero of \(f_{u-1}\). Finally,
    \(f_u\brac*{\frac{pv}{u},h_{u-1,0;\frac{pv}{u}},\frac{u}{v}}\) was
    computed in \cite[Eq (4.30)]{RW} in terms of \(\slt\)
    data. Expressed in terms of \(\Nt=2\) data this becomes
    \begin{equation}
      f_u\brac*{\frac{pv}{u},h_{u-1,0;\frac{pv}{u}},\frac{u}{v}}=\binom{2(u-1)}{u-1}\frac{(p+u-2)(p+u-4)\cdots (p-u+2)}{2^{u-1}(u-1)!},
    \end{equation}
    which is zero for \(2-u\le p\le u-2\), \(p+u\equiv 0\pmod{2}\).

    Finally we consider the zeros of
    \(f_u\brac*{\eta+\frac{2v}{u},\Delta-\eta-\frac{v}{u},\frac{u}{v}}(2\Delta-\eta)\). The
    base case of \(u=2\) and \(u=3\) has already been established. For the
    induction assume again that the pairs 
    \((\frac{pv}{u},h_{r,0;\frac{pv}{u}})\), \(1\le r\le k-1\), \(1-r\le p\le
    r-1\), \(p+r\equiv 1\pmod{2}\) are zeros of
    \(f_k(\eta+\frac{2v}{u},\Delta-\eta-\frac{u}{v},\frac{u}{v})(2\Delta-\eta)\)
    for \(2\le k<u\). Then the 
    pairs \((\frac{pv}{u},h_{r,0;\frac{pv}{u}})\) are for zeros of
    \(f_u(\eta+\frac{2v}{u},\Delta-\eta-\frac{u}{v},\frac{u}{v})(2\Delta-\eta)\)
    for \(r\le u-3\) by the recursion relation \eqref{eq:gtu+2v}. For
    \(r=u-2\), plugging \((\frac{pv}{u},h_{u-2,0;\frac{pv}{u}})\) into the
    recursion relation \eqref{eq:gtu+2v} for
    \(f_u(\eta+\frac{2v}{u},\Delta-\eta-\frac{u}{v},\frac{u}{v})(2\Delta-\eta)\)
    leads to the coefficient in front of
    \(f_{u-2}(\eta+\frac{2v}{u},\Delta-\eta-\frac{u}{v},\frac{u}{v})(2\Delta-\eta)\)
    vanishing and hence \((\frac{pv}{u},h_{u-2,0;\frac{pv}{u}})\) is a zero
    for
    \(f_{u}(\eta+\frac{2v}{u},\Delta-\eta-\frac{u}{v},\frac{u}{v})(2\Delta-\eta)\)
    because it is a zero for
    \(f_{u-1}(\eta+\frac{2v}{u},\Delta-\eta-\frac{u}{v},\frac{u}{v})(2\Delta-\eta)\). Finally
    we can again use \cite[Eq (4.30)]{RW} to compute
    \begin{equation}
      f_u\brac*{\frac{pv}{u}+\frac{2v}{u},h_{u-1,0;\frac{pv}{u}}-\frac{pv}{u}-\frac{v}{u},\frac{u}{v}}=\binom{2(u-1)}{u-1}\frac{(p+u)(p+u-2)\cdots (p-u+4)}{2^{u-1}(u-1)!},
    \end{equation}
    and note that
    \begin{equation}
      2h_{u-1,0;\frac{vp}{u}}-\frac{vp}{u}=\frac{v}{2u}\brac*{u^2-2u-p^2-2p}=-\frac{vp}{2u}\brac*{p-u+2}\brac*{p+u}.
      \end{equation}
      Thus
      \(f_u\brac*{\frac{pv}{u}+\frac{2v}{u},h_{u-1,0;\frac{pv}{u}}-\frac{pv}{u}-\frac{v}{u},\frac{u}{v}}\brac*{2h_{u-1,0;\frac{vp}{u}-\frac{vp}{u}}}\)
        vanishes for \(2-u\le p\le u-2\), \(u\equiv p\pmod{2}\).
        
\end{proof}

\subsection{Fusion rules from Zhu bimodules}
Using the theory on Zhu algebras outlined in \cref{sec:Zhuthy} 
we are now in a position to describe the Zhu bimodules associated to $\uN{c}$-Verma  modules $\NM{q}{h}{c}$.

 We first introduce the following auxiliary spaces:
\begin{align}
\overline{W}_{q,h}&= {\rm span}_{\mathbb{C}}\{\vac_{q,h}, G^+_{-\frac{1}{2}}\vac_{q,h},  G^-_{-\frac{1}{2}}\vac_{q,h},  G^-_{-\frac{1}{2}}G^+_{-\frac{1}{2}}\vac_{q,h} \}\subset \NM{q}{h}{c},\nonumber\\
    W_{q,h}&=\mathbb{C}[x,y,z] \cten \overline{W}_{q,h}.
\end{align} 
Next, we endow $W_{q,h}$ with a $\mathbb{C}[\Delta,\eta]$-bimodule structure as follows.
\begin{align}
    \Delta.(f(x,y,z)\cten v) &= x\cdot f(x,y,z)\cten v, \nonumber\\
    (f(x,y,z)\cten v)  .\Delta&= y\cdot f(x,y,z)\cten v, \nonumber\\
    \eta.(f(x,y,z)\cten v )&= z\cdot f(x,y,z)\cten v+f(x,y,z)\cten J_0 v, \nonumber \\
    (f(x,y,z)\cten v). \eta &= z\cdot f(x,y,z)\cten v.\label{eq:deltaetaaction}
\end{align}
Using standard arguments in Zhu theory (see \cite{Z1, FZ, KW}) we obtain the following description for the Zhu bimodules $A(\NM{q}{h}{c})$.
\begin{prop}
As an $A(\uN{c})\cong \mathbb{C}[\Delta,\eta]$-bimodule 
\begin{align}
A(\NM{q}{h}{c})\cong W_{q,h},
\end{align}
where the isomorphism is given by 
\begin{align}
 \mathbb{C}[x,y,z]\cten\overline{W}_{q,h}  &\longrightarrow A(\NM{q}{h}{ c}) \nonumber\\
 f(x,y,z)\cten w &\longmapsto    [ f(L_{-2}+2L_{-1}+L_0, L_{-2}+L_{-1}, J_{-1})w]. \label{eq:isoZV}
\end{align}
\label{thm:zhubimodrel}
\end{prop}
Denote the left \(A(\uN{c})\)-module formed by the space of least conformal
weight of a Verma module \(\NM{q}{h}{c}\) by \(N_{q,h}\). That is $N_{q,h}=\CC
  n_{q,h}$, with \(\Delta. n_{q,h}\)=\(h n_{q,h}\) and
  \(\eta. n_{q,h}=q n_{q,h}\). As a \(\ZZ_2\)-graded vector space \(N_{q,h}\) is even because we defined the highest weight vectors of Verma modules to be even. We denote the parity reversal of \(N_{q,h}\) by \(\Pi N_{q,h}\). We describe
  $A(\NM{q}{h}{c})\cten_{A(\uN{c})}N_{q_2,h_2}$ as a left $A(\uN{c})$-module.
  \begin{lemma}
Let \(h_1,h_2,q_1,q_2\in\CC\).
\begin{enumerate}
\item Then, the set
  \begin{equation}
    \left\{ x^a\cten {G^+_{-\frac12}}^i{G^-_{-\frac12}}^j\vac_{q_1, h_1}\cten_{A(\uN{c})} n_{q_2,h_2}| \ a\in \ZZ_{\ge0}, i,j\in\{0,1\}\right\}
    \label{eq:modbasis}
  \end{equation}
  is a \(\CC\)-basis of \(W_{q_1, h_1}\cten_{A(\uN{c})} N_{q_2,h_2}\).
\item The action of \(\Delta\) and \(\eta\) on the \(\CC\)-basis above is
  \begin{align}
    \Delta. x^a\cten w\cten_{A(\uN{c})} n_{q_2,h_2}&=x^{a+1}\cten w\cten_{A(\uN{c})} n_{q_2,h_2},\nonumber\\
    \eta. x^a\cten w\cten_{A(\uN{c})} n_{q_2,h_2}&=
    \brac*{q_2+\tilde{q}}x^a\cten w\cten_{A(\uN{c})}
    n_{q_2,h_2}, 
  \end{align}
  where \(w={G^+_{-\frac12}}^i{G^-_{-\frac12}}^j\vac_{q_1, h_1}\) and \(\tilde{q}=q_1+i-j\)
  is the \(J_0\)-weight of \(w\).
\end{enumerate}
\label{thm:basisZhumodVer}   
\end{lemma}

\begin{theorem} 
  Let \(h_i,q_i\in \CC\), \(i=1,2,3\), then
  \begin{equation}
    \dim \hom \left( W_{q_1, h_1} \cten_{A(\uN{c})} N_{q_2,h_2}, \ \Pi^j N_{q_3,h_3}
    \right)=
    \begin{cases}
      1&\text{if}\quad \left|q_1+q_2-q_3\right|=1,\ j=1,\\
      2&\text{if}\quad q_1+q_2=q_3,\ j=0,\\
      0&\text{else}.
    \end{cases}
  \end{equation}
  \label{thm:fusionVermas}
\end{theorem}

\begin{proof}
  Since the action of \(\Delta\) and \(\eta\) on \(W_{q_1, h_1}
  \cten_{A(\uN{c})} N_{q_2,h_2}\) does not mix the basis vectors
  \({G^+_{-\frac12}}^i{G^-_{-\frac12}}^j\vac_{q_1, h_1}\) of
  \(\overline{W}_{q_1, h_1}\), 
  we have the direct sum decomposition
  \begin{equation}
    W_{q_1, h_1}
    \cten_{A(\uN{c})} N_{q_2,h_2}=\bigoplus_{i,j=0}^{ 1} \CC[x]\cten\CC{G^+_{-\frac12}}^i{G^-_{-\frac12}}^j\vac_{q_1, h_1}\cten_{A(\uN{c})} \CC n_{q_2,h_2},
  \end{equation}
  where each summand is the \(\eta\)-eigenspace of respective eigenvalue
  \(q_1+q_2+i-j\). This immediately implies the ``else'' line of the theorem.

  Next assume there exists a homomorphism \(\rho:  W_{q_1, h_1} \cten_{A(\uN{c})}
  N_{q_2,h_2} \to \ N_{q_3,h_3}\), then
  \begin{equation}
    0=(\Delta-h_3)\rho(-)=\rho((\Delta-h_3) - ),
  \end{equation}
  that is, the image of \(\Delta-h_3\) in \(W_{q_1, h_1} \cten_{A(\uN{c})}
  N_{q_2,h_2}\) lies in the kernel of \(\rho\). The image of \((\Delta-h_3)\)
  is just the submodule
  \begin{equation}
    S_{h_3}=\ang*{ (x-h_3)\cten
      w\cten_{A(\uN{c})} n_{q_2,h_2}| w\in \overline{W}_{q_1, h_1}}.
  \end{equation}
  Consider the quotient by \(S_{h_3}\) and we immediately obtain
  \begin{align}
    \frac{W_{q_1, h_1}\cten_{A(\uN{c})}N_{q_2,h_2}}{S_{h_3}}&\cong
    \bigoplus_{i,j=0}^{ 1}
    \frac{\CC[x]}{\ang*{x-h_3}}\cten\CC{G^+_{-\frac12}}^i{G^-_{-\frac12}}^j\vac_{q_1, h_1}\cten_{A(\uN{c})} \CC n_{q_2,h_2}\nonumber\\
    &\cong
    \Pi N_{h_3,q_1+q_2-1}\oplus 2N_{h_3,q_1+q_2}\oplus \Pi N_{h_3,q_1+q_2+1},
  \end{align}
  which implies the theorem.
\end{proof}

\begin{prop}
  Let \(t,q\in \CC\) and consider the Verma module
  \(\NM{q}{h_{1,1;q}(t)}{c(t)}\), where
  \(h_{1,1;q}(t)=\frac{(1-t)^2-1-q^2}{4t}\). Then the vector
  \begin{equation}
    w_{1,1} \left( t,q \right)=\left( \left( q-1 \right) L_{-1}+\frac{t}{2} \left( q^2-1 \right) J_{-1}+G_{-\frac{1}{2}}^+ G_{- \frac{1}{2}}^- \right)  \vac_{ q, h_{1,1;q}(t),c(t)}
  \end{equation}
  is singular in $\NM{q}{h_{1,1;q}(t)}{c(t)}$. Furthermore, under the isomorphism in \eqref{eq:isoZV}, the images
  of the singular
  vector \(\sqbrac*{w_{1,1}}\in A(\NM{q}{h_{1,1;q}(t)}{c(t)})\) and its \(G^\pm_{-\frac12}\) descendants are
  \begin{align}
    \sqbrac*{w_{1,1}}&\mapsto f^{1}(x,y,z,q) \cten \sqbrac*{\vac_{q,h_{1,1;q}(t)}}+1
    \cten \sqbrac*{G^+_{-\frac12}G^-_{-\frac12}\vac_{q,h_{1,1;q}(t)}},\nonumber\\
    \sqbrac*{G^+_{-\frac12}w_{1,1}}&\mapsto f^+(x,y,z,q)\cten
                                     \sqbrac*{G^+_{-\frac12}\vac_{q,h_{1,1;q}(t)}}, \nonumber\\
    \sqbrac*{G^{-}_{-\frac12}w_{1,1}}&\mapsto f^-(x,y,z,q) \cten
                                       \sqbrac*{G^-_{-\frac12}\vac_{q,h_{1,1;q}(t)}},\nonumber\\ 
    \sqbrac*{G^+_{-\frac12}G^{-}_{-\frac12}w_{1,1}}&\mapsto
    \frac{t}{2} (q^2-1)(2y-z)\cten\sqbrac*{\vac_{q,h_{1,1;q}(t)}}, \nonumber\\
    &\qquad +f^{G}(x,y,z,q) \cten \sqbrac*{G^+_{-\frac12}G^-_{-\frac12}\vac_{q,h_{1,1;q}(t)}},
      \label{eq:homspacevec}
  \end{align}
  where 
  \begin{align}
    f^1(x,y,z,q)&=\brac*{q-1}\brac*{x-y+\frac{1}{2}+\frac{t}{4}\brac*{q+1}\brac*{2z+q-1}},\nonumber\\
    f^+(x,y,z,q)&=\brac*{q-1}\brac*{x-y+\frac{t}{4}\brac*{q+1}\brac*{2z+q+1}},\nonumber\\
    f^-(x,y,z,q)&=\brac*{q+1}\brac*{x-y+\frac{t}{4}\brac*{q-1}\brac*{2z+q-1}},\nonumber\\
    f^G(x,y,z,q)&=\brac*{q+1}\brac*{x-y-\frac{1}{2}+\frac{t}{4}\brac*{q-1}\brac*{2z+q+1}}.
  \end{align}
\end{prop}
\begin{proof}
  By direct computation we see that \(G_{\frac12}^\pm w_{1,1}=0=J_1 w_{1,1}\)
  and hence \(w_{1,1}\) is singular. The simplifications within the Zhu
  bimodule follow from the relations in \cref{thm:zhubimodrel}, which allows us
  to make the following replacements
  \begin{align}
    L_{-2}&\mapsto 2y-x+L_0,\quad L_{-1}\mapsto x-y-L_0,\quad J_{-1}\mapsto z,\nonumber\\
    \sqbrac*{J_{-2}v}&= -\sqbrac*{J_{-1}v},\qquad \sqbrac*{G^\pm_{-\frac{3}{2}}v}=-\sqbrac*{G^\pm_{-\frac{1}{2}}v}.
  \end{align}
\end{proof}

\begin{prop} 
  For \(q,t\in \CC\), \(q\neq\pm1\), let $c(t)=3-\frac{6}{t}$.
  Let $r,s$  be positive integers, $q_2 \in
  \CC$  and $h_{r,s; q_2}(t)=\frac{(r-ts)^2-1-(tq_2)^2}{4t}$. Finally, let  $q_3,h_3\in \CC$. Then, 
  \begin{equation}
    \dim \mathrm{Hom} \left(
      A\brac*{\frac{\NM{q}{h_{1,1;q}(t)}{c(t)}}{\ang{w_{1,1}}}}
      \cten_A 
      N_{q_2,h_{r,s}}, \ \Pi^jN_{q_3,h_3}
    \right)\leq 
    \begin{cases}
      1,&\hspace{-3mm}\text{if}\ q_3=q+q_2+1,\  h_3=h_{r,s, q_3},\ j=1,\\
      1,&\hspace{-3mm}\text{if}\ q_3=q+q_2-1,\  h_3=h_{r,s, q_3},\ j=1,\\
      1,&\hspace{-3mm}\text{if}\ q_3=q+q_2,\  h_3=h_{r, s\pm 1, q_3},\ j=0,\\  
        0,&\hspace{-3mm}\text{otherwise},
    \end{cases}
    \label{eq:homspacebounds}
  \end{equation}
where \(\cten_A\) denotes the tensor product over \(A(\uN{c})\).
  \label{thm:upperbounds}
\end{prop}
\begin{proof}
  We use the canonical identification of vector spaces
  \begin{equation}
    \mathrm{Hom} \left( A({\NM{q}{h_{1,1;q}(t)}{c(t)}}) \cten_{A} N_{q_2, h_{r,s;q_2}(t)}, \ N_{q_3,h_3}
    \right)\cong N_{q_3,h_3}^\ast\cten_{A}A({\NM{q}{h_{1,1;q}(t)}{c(t)}}) \cten_{A} N_{q_2,h_{r,s;q}(t)}
  \end{equation}
  to compute dimensions, where \(N_{q_3,h_3}^\ast\) is the dual of
  \(N_{q_3,h_3}\). By \cref{thm:muz}, the dimension of the hom space in \eqref{eq:homspacebounds}
  is therefore the codimension of
  \(N_{q_3,h_3}^\ast\cten_{A}A(\ang{w_{1,1}}) \cten_{A}
  N_{q_2, h_{r,s;q}(t)}\) in \(N_{q_3,h_3}^\ast\cten_{A}A({\NM_{q,h_{1,1;q}(t)}{c(t)}}) \cten_{A} N_{q_2, h_{r,s;q}(t)}\).
  If \(q_3-q-q_2\neq 0,\pm1\), then the dimension of the hom space
  \eqref{eq:homspacebounds} must be zero by \cref{thm:fusionVermas}. So we
  first consider the case \(q_3=q+q_2+1\). After tensoring the vectors in
  \eqref{eq:homspacevec} with the basis vector \(n_{q+q_2+1, h_3}\) of
  \(N_{q_3,h_3}^\ast\) we see that only the second vector in
  \eqref{eq:homspacevec} does not vanish and yields
  \begin{align}
    n_{q+q_2+1,h_3}\cten_AG^+_{-\frac{1}{2}}w_{1,1}\cten_{A}n_{q_2,h_2}=f^+(h_3,h_{r,s;q_2}(t),q_2,q)n_{q+q_2+1,h_3}\cten_{A}  G^+_{-\frac12}\vac_{q, h_{1,1}}
  \cten_{A}n_{q_2,h_2},
  \end{align}
  where
  \begin{align}
    f^+(h_3,h_{r,s;q_2}(t),q_2,q)&=(q-1)\brac*{h_3-h_{r,s;q_2}(t)-\frac{t}{4}(1-q^2)+\frac{t}{2}(q+1)(q_2+1)}\nonumber\\
    &=(q-1)\brac*{h_3-h_{r,s;q+q_2+1}(t)},
  \end{align}
  that is, the variables \(x,y,z\) evaluate to \(h_3,h_{r,s;q_2}(t),q_2\), respectively.
  Hence \(N_{q_3,h_3}^\ast\cten_{A}A(\ang{w_{1,1}})
  \cten_{A} N_{q_2, h_{r,s;q}(t)}\) has codimension \(0\) in
  \(N_{q_3,h_3}^\ast\cten_{A}A(\NM{q}{h_{1,1;q}(t)}{c(t)})
  \cten_{A} N_{q_2, h_{r,s;q}(t)}\) unless
  \(h_3=h_{r,s;q+q_2+1}(t)\), which proves the first line of
  \eqref{eq:homspacebounds}.
  Similarly, for \(q_3=q+q_2-1\) only the third vector in \eqref{eq:homspacevec} does not
  vanish when tensoring with \(n_{q+q_2-1, h_3}\) and yields
  \begin{align}
    n_{q+q_2-1, h_3}\cten_AG^-_{-\frac{1}{2}}w_{1,1}\cten_{A}n_{q_2,h_2}=
    f^-(h_3,h_{r,s;q_2}(t),q_2,q)n_{q+q_2-1, h_3}\cten_{A}  G^+_{-\frac12}\vac_{q, h_{1,1}}
  \cten_{A}n_{q_2,h_2},
  \end{align}
  where
  \begin{align}
    f^-(h_3,h_{r,s;q_2}(t),q_2,q)&=(q+1)\brac*{h_3-h_{r,s;q_2}(t)-\frac{t}{4}(1-q^2)+\frac{t}{2}
                                   (q-1) (q_2-1)}\nonumber\\
    &=(q+1)\brac*{h_3-h_{r,s;q+q_2-1}(t)}.
  \end{align}
  For \(q_3+q+q_2\) the second and third vectors in \eqref{eq:homspacevec}
  vanish when tensoring with \(n_{q+q_2, h_3}\), while the others yield
  \begin{align}
    n_{q+q_2, h_3}\cten_Aw_{1,1}\cten_{A}n_{q_2,h_2}=
    &f^{1}(h_3,h_{r,s;q_2}(t),q_2,q) n_{q+q_2, h_3}\cten_{A}
    \vac_{q, h_{1,1}}\cten_{A}n_{q_2,h_2}\nonumber\\
    &\ + n_{q+q_2, h_3}  \cten_{A}
      G^+_{-\frac12}G^-_{-\frac12}\vac_{q, h_{1,1}}\cten_{A}
      n_{q_2,h_2},\nonumber\\
    n_{q+q_2, h_3}\cten_AG^+_{-\frac{1}{2}}G^-_{-\frac{1}{2}}w_{1,1}\cten_{A}n_{q_2,h_2}&=
    \frac{t}{2}(q^2-1)
      (2h_{r,s;q_2}(t)-q_2)\cten_A\vac_{q, h_{1,1}}\cten_{A}n_{q_2,h_2}
      \nonumber\\
    &\qquad+f^{G}(h_3,h_{r,s;q_2}(t),q_2,q) \cten
      G^+_{-\frac12}G^-_{-\frac12}\vac_{q, h_{1,1}}\cten_{A} n_{q_2,h_2}.
      \label{eq:homvects}
  \end{align}
  The constant coefficient in front of one of the summands above means that
  \(N_{q+q_2, h_3}^\ast\cten_{A}A(\ang{w_{1,1}}) \cten_{A}
  N_{q_2, h_{r,s;q}(t)}\) must always have codimension at most \(1\) in
  \(N_{q+q_2, h_3}^\ast\cten_{A}A(\NM{q}{h_{1,1;q}(t)}{c(t)}) \cten_{A}
  N_{q_2, h_{r,s;q}(t)}\). To determine when the codimension can be greater
  than \(0\) we compute the determinant of the coefficients
  \eqref{eq:homvects} to obtain a polynomial in \(h_3\). The roots of this
  polynomial are precisely where the codimension can be greater than \(0\).
  \begin{align}
    &\det
    \begin{pmatrix}
       f^1(h_3,q,h_{r,s;q_2}(t),q_2)&1\\
      \frac{t}{2}(q^2-1)(2h_{r,s;q_2}(t)-q_2)&f^G(h_3,q,h_{r,s;q_2}(t),q_2)
    \end{pmatrix}\nonumber\\
    &\qquad=\brac*{q^2-1}\brac*{h_3-h_{r,s+1;q+q_2}(t)}\brac*{h_3-h_{r,s-1;q+q_2}(t)}.
  \end{align}
  Therefore \(h_3=h_{r,s\pm1;q+q_2}(t)\) are the only values for which the
  codimension need not be \(0\).
  
\end{proof}

Note that these bounds on hom space dimensions in particular provide an upper bound
  for the fusion rules \eqref{eq:11rule} in \cref{thm:1121rules} after
  setting \(t=\frac{u}{v}\).

\section{Lower bounds \texorpdfstring{for \cref{thm:1121rules}}{} via free field realisations}
\label{sec:lowerbounds}

In this section we prove lower bounds for the dimensions appearing in
\cref{thm:1121rules} by explicitly constructing suitable intertwining
operators in a free field realisation of \(\mmsl{u}{v}\). To this end
we recall some facts regarding  $P(w)$-intertwining operators following
\cite{CKM,HLZ3} before discussing the
free field realisation
of \(\mmsl{u}{v}\) and its screening operators.

Given a generalised
$V$-module $M=\bigoplus\limits_{h \in \CC} M_{[h]},$ its
\textit{algebraic completion} $\overline{M}$ is defined as the superspace
$\overline{M}=\prod\limits_{h \in \CC} M_{[h]}$, where
$\overline{M}^i=\overline{M^i}$ for $i=\overline{0},\overline{1}$ while its {\it contragredient} or {\it graded dual} module $M'$ is defined as the vector space 
\begin{align}
M'=\bigoplus_{h\in \mathbb{C}} M^*_{[h]}, \ \ \ M^*_{[h]}=\hom(M_{[h]}, \mathbb{C}),
\end{align}
together with the action \(Y_{M'}\) characterised by
\begin{align}
  \pair{ Y_{M'}(v,z)m'}{ m}=\pair{  m'}{ Y_M^{\opp}(v,z)m}, 
\end{align}
where $m'\in M', m\in M$, $v\in V$, and where
\begin{equation}
  Y_M^{\opp}(v,z)=Y_M(e^{z L_1}(-z^{-2})^{L_0}v, z^{-1}),
\end{equation}
is the \emph{opposed action}  (note that the module \(M\) needs to satisfy some
lower boundedness conditions, which will never be an issue here, in order for
the action on \(M'\) to
be well-defined). Similarly, for any triple of modules \(M_1,M_2,M_3\) and an
intertwining operator \(\mathcal{Y}\in \binom{M_3}{M_1,\ M_2}\) the opposed
intertwining operator
\begin{equation}
  \mathcal{Y}^{\opp}(m_1,z)=\mathcal{Y}(e^{z L_1}(-z^{-2})^{L_0}m_1, z^{-1})
  \label{eq:oppformula}
\end{equation}
characterises an intertwining operator \(\widetilde{\mathcal{Y}}\in
\binom{M_2'}{M_1,\ M_3'}\) via
\begin{equation}
  \pair{m_2}{\widetilde{\mathcal{Y}}(m_1,z)m_3'}=\pair{m_3'}{\mathcal{Y}^{\opp}(m_1,z)m_2}.
\end{equation}

\begin{defn}\cite[Definition 4.1]{KaRi} 
  Let $w\in\CC^\times$ and let $M_1$, $M_2$, $M_3$ be modules over a
  \vosa{} \(V\). A \textit{$P \left( w \right)$-intertwining map of type $\binom{M_3}{M_1 \quad M_2}$} is a bilinear map $I \colon M_1 \cten M_2 \to \overline{M_3}$, for $\overline{M_3}$ the algebraic completion of $M_3$,  that satisfies the following properties:
  \begin{enumerate}
  \item \emph{Lower truncation}: For any $\psi_1 \in M_1$, $\psi_2 \in M_2$, $\pi_h \brac*{ I \brac*{\psi_1 \cten \psi_2}}=0$ for all $\re(h)\ll0$, where $\pi_h$ denotes the projection onto the generalised eigenspace $\left( M_3 \right)_{[h]}$ of $L \left( 0 \right)$-eigenvalue $h$;
  \item \emph{Convergence}:  For any $\psi_1 \in M_1$, $\psi_2 \in M_2$, $\psi_3 \in
    M_3^\prime$, the contragredient of
    \(M_3\), 
    the series defined by
    \begin{align}
      & \pair{ \psi_3}{ Y_3\left( v,z \right) 
        I\brac{\psi_1,w}{\psi}_2}, \nonumber\\
      & \pair{ \psi_3}{I \brac{Y_1 \left( v,z-w \right) \psi_1,w}\psi_2}, \quad \mathrm{and} \nonumber\\ 
      & \pair{ \psi_3}{I \brac{\psi_1,w}Y_2 \left( v,z \right) \psi_2}
    \end{align}
    are absolutely convergent in the regions $\vert z \vert > \vert w \vert
    >0$, $\vert w \vert > \vert z-w \vert >0$ and $\vert w \vert > \vert z
    \vert >0$, respectively, where the subscript under each $Y$ indicates the module being acted on.
  \item \emph{Cauchy-Jacobi identity}: Given any $f \left( t \right) \in
    R_{P(w)}=\mathbb{C} \left[ t,t^{-1},\left( t-w \right)^{-1}
    \right]$, the space of rational functions whose poles lie in some
    subset of $\lbrace 0,w,\infty \rbrace$, we have the following identity.
    \begin{align}
      \oint_{0,w} f \left( z \right) \pair{ \psi_3}{Y_3 \left( v,z
        \right)I \brac{\psi_1,w} \psi_2 
      } \dd z &=
      \oint_w f\left(z\right)\pair{\psi_3}{I \brac{Y_1\left(v,z-w\right)\psi_1,w}\psi_2
      }\dd z \nonumber\\
      & +\oint_0 f \left( z \right) \pair{ \psi_3}{ 
        I \brac{\psi_1,w}Y_2 \left( v,z \right) \psi_2 
      } \dd z,
      \label{eq:CauchyJacobi}
    \end{align}
    where the subscript on each integral indicates the points in $\lbrace
    0,w,\infty \rbrace$ that must be enclosed by simple positively oriented
    contours. 
  \end{enumerate}
\end{defn}
  
Next we recall a free field construction of \(\mmsl{u}{v}\) due to Feigin-Semikhatov-Tipunin \cite{FST}, and Adamovi{\'c} \cite{A3}.
Consider the 2-dimensional trivial Lie algebra with basis denoted
\(\set{a,b}\) and with invariant non-degenerate symmetric form normalised such
that
\begin{equation}
  \pair{a}{a}=-\pair{b}{b}=\frac{t}{2}-1,
  \qquad
  \pair{a}{b}=0,\qquad t\in\RR\setminus{2}.
  \label{eq:latpairing}
\end{equation}
We denote the rank 2 Heisenberg \va{} constructed from the affinisation of this Lie algebra by
\(\heis{2,t}\). This \va{} is isomorphic to the tensor product of two rank 1 Heisenberg
algebras of respective levels \(\frac{t-2}{4}\) and \(\frac{2-t}{4}\), that is
\(\heis{2,t}\cong \heis{\frac{t-2}{4}}\cten \heis{\frac{2-t}{4}}\). Let \(\set{a(z),b(z)}\)
be the generating currents corresponding to the basis \(\set{a,b}\) whose \ope{s}
therefore satisfy
\begin{equation}
  a(z)a(w)\sim -b(z)b(w)\sim \frac{\frac{t}{2}-1}{\brac*{z-w}^2},\qquad
  a(z) b(w)\sim 0.
\end{equation}
Further, consider the rank 1 lattice
$\mathbb{L}=\ZZ\frac{2}{t-2}(a-b)$ (note that the pairing \eqref{eq:latpairing} restricted
to this lattice vanishes) and
corresponding lattice \va{}
\begin{equation}
  \lat{t}=\bigoplus_{p\in \mathbb{L}}\mathcal{F}_p,
\end{equation}
where \(\mathcal{F}_p\) denotes the rank 2 Fock space over \(\heis{2,t}\) whose
highest weight vector \(\ket{p}\) satisfies
\begin{equation}
  a_0\ket{p}=\pair{a}{p}\ket{p},\qquad b_0\ket{p}=\pair{b}{p}\ket{p}.
\end{equation}
Note that even though \(\heis{2,t}=\mathcal{F}_0\) has two independent
generating currents (that is, it is rank 2) the lattice \(\mathbb{L}\) is only
rank 1. The hence dual of \(\mathbb{L}\)
\begin{equation}
  \mathbb{L}^\ast=\set*{v\in \RR a\oplus\RR b
    \st
    \frac{2}{t-2}\pair{v}{a-b}\in \ZZ}=\RR(a-b)\oplus \ZZ b
\end{equation}
is not discrete. Finally, for any \(\mu\in \mathbb{L}^\ast/\mathbb{L}\), the sum of Fock spaces
\begin{equation}
  \mathbb{F}_\mu=\bigoplus_{p\in \mu}\mathcal{F}_p
\end{equation}
is naturally an \(\lat{t}\) module.

\begin{theorem}[Conjectured by Feigin, Semikhatov, Tipunin \cite{semi}; proved
  by Adamovi{\'c} {\cite[Thm 5.5, Prop 9.2]{A3}} ]
For coprime \(u,v\ge2\), let \(\mmV{u}{v}\) denote the Virasoro minimal model
\voa{} at central charge \(c_{u,v}=1-6\frac{\brac*{u-v}^2}{uv}\). Denote the
conformal vector of \(\mmV{u}{v}\) by
\(\omega_{u,v}=L_{-2}^{\brac*{u,v}}\vac\).
\begin{enumerate}
\item There exists an embedding of \voa{s}
  \begin{align}
    \psi:\mmsl{u}{v}\hookrightarrow \mmV{u}{v} \cten \lat{\frac{u}{v}}
  \end{align}
  characterised on generators by
  \begin{align}
    &e_{-1}\vac \mapsto \vac\cten \ket{\frac{2v}{u-2v}(a-b)},\\
    &h_{-1}\vac \mapsto 2a_{-1}\vac\cten \ket{0},\\
    &f_{-1}\vac\mapsto \brac*{\frac{u}{v}L^{\brac*{u,v}}_{-2}-b_{-1}^2-\brac*{\frac{u}{v}-1}b_{-2}}\vac\cten\ket{-\frac{2v}{u-2v}(a-b)}.
  \end{align}
  The image of the conformal vector in \(\mmsl{u}{v}\) under this embedding is
  \begin{align}
    \frac{1}{2\frac{u}{v}}\brac*{\frac12
    h_{-1}^2+e_{-1}f_{-1}+f_{-1}e_{-1}}\vac \mapsto \brac*{L_{-2}^{\brac*{u,v}}+\frac{v}{u-2v}\brac*{a_{-1}^2-b_{-1}^2}-b_{-2}}\vac\cten\ket{0}.
  \end{align}
  \label{itm:ffrformulae}
\item Let \(\mmVM{r}{s}\), \(1\le r\le u-1\), \(1\le s\le v-1\) denote the
  simple \(\mmV{u}{v}\)-module of highest conformal weight
  \begin{equation}
    h_{r,s}=\frac{\brac*{us-vr}^2-\brac*{u-v}^2}{4uv}.
  \end{equation}
  Recall that \(\mmsl{u}{v}\) admits an automorphism \(\gamma\) called
  \emph{conjugation}, which is characterised by
  \begin{equation}
    e_{-1}\vac\mapsto f_{-1}\vac,\qquad h_{-1}\vac \mapsto -h_{-1}\vac,\qquad
    f_{-1}\vac\mapsto e_{-1}\vac.
  \end{equation}
  Thus \(\psi\circ \gamma\) is also an embedding \(\mmsl{u}{v}\hookrightarrow
  \mmV{u}{v} \cten \lat{\frac{u}{v}}\). For any \(\mmV{u}{v} \cten
  \lat{\frac{u}{v}}\)-module we can therefore pull back along \(\psi\) and
  \(\psi\circ \gamma\) to obtain \(\mmsl{u}{v}\) modules. Under these
  pullbacks, for \(\ell\in\ZZ\), \(1\le r\le u-1\), \(1\le s\le v-1\) we have the identifications
  \begin{align}
    \psi^\ast \mmVM{r}{s}\cten \mathbb{F}_{\mu_{r,s;l}}&\cong \sigma^{\ell}\slEr{u-r}{v-s}{-},&\mu_{r,s;l}&=\sqbrac*{\brac*{\tfrac{v\lambda_{r,s}}{u-2v}+\ell}\brac*{a-b}+\brac*{l-1}b},\nonumber\\
    \brac*{\psi\circ\gamma}^\ast \mmVM{r}{s}\cten
    \mathbb{F}_{\tilde\mu_{r,s;l}}&\cong
    \sigma^{\ell}\slEr{r}{s}{+},&\tilde{\mu}_{r,s;l}&=\sqbrac*{\brac*{-\tfrac{v\lambda_{r,s}}{u-2v}-\ell}\brac*{a-b}-\brac*{l+1}b},\nonumber\\
    \psi^\ast \mmVM{r}{s}\cten \mathbb{F}_{\mu}&\cong \sigma^{\pair{\frac{2v}{u-2v}(a-b)}{\mu+b}}\slE{\pair{2b}{\mu+b}}{r}{s},
    &&\mu\in \mathbb{L}^\ast/\mathbb{L},\
    \pair{2b}{\mu+b}\neq\sqbrac*{\pm
      \lambda_{r,s,\frac{u}{v}}},\nonumber\\
    (\psi\circ\gamma)^\ast \mmVM{r}{s}\cten \mathbb{F}_{\mu}&\cong
    \sigma^{\pair{-\frac{2v}{u-2v}(a-b)}{\mu+b}}\slE{-\pair{2b}{\mu+b}}{r}{s}.
  \end{align}
  \label{itm:moduleidentifications}
\item Consider the simple \(\mmV{u}{v}\) module \(\mmVM{1}{2}\) of highest
  conformal weight \(h_{1,2}=\frac{3u}{4v}-\frac12\) generated by a highest weight vector $v_{1,2}$. 
  Tensor products with \(\mmVM{1}{2}\) decompose, for $1\le r\le u-1$, as
  \begin{equation}
    \mmVM{1}{2}\fuse \mmVM{r}{s}\cong
    \begin{cases}
      \mmVM{r}{2}&\text{if } s=1,\\
      \mmVM{r}{s-1}\oplus\mmVM{r}{s+1}&\text{if }2\le  s\le v-2,\\
      \mmVM{r}{v-2}&\text{if } s=v-1.
    \end{cases}
  \end{equation}
  Further, let \(I^+_{r,s}\in \binom{\mmVM{r}{s+1}}{\mmVM{1}{2},\
    \mmVM{r}{s}}\), \(1\le s\le v-2\) and \(I^-_{r,s}\in \binom{\mmVM{r}{s-1}}{\mmVM{1}{2},\
    \mmVM{r}{s}}\), \(2\le s\le v-1\)
  be surjective \(\mmV{u}{v}\)-intertwining operators and let
  \(I_\mu\) a surjective \(\lat{\frac{u}{v}}\)-intertwining operator of type
  \(\binom{\mathbb{F}_{\mu+b}}{\mathbb{F}_{b+\mathbb{L}},\ \mathbb{F}_{\mu}}\)
  such that their tensor product \(\mathcal{Y}^\pm_{\mu;r,s}=I^\pm_{r,s}\cten
  I_{\mu}\) is an \(\mmV{u}{v} \cten \lat{\frac{u}{v}}\)-intertwining operator of type
  \(\binom{\mmVM{r}{s\pm 1}\cten
    \mathbb{F}_{b+\mu}}{\mmVM{1}{2}\cten\mathbb{F}_{b+\mathbb{L}},\
    \mmVM{r}{s}\cten \mathbb{F}_{\mu}}\). The currents
  \begin{equation}
      \scr^{\pm}(z)=\mathcal{Y}^\pm_{\mu;r,s}(v_{1,2}\cten\ket{b},z),
      \label{eq:scrcur}
  \end{equation}
corresponding to the vector \(v_{1,2}\cten \ket{b}\in
  \mmVM{1}{2}\cten \mathbb{F}_{b+\mathbb{L}}\) and  where the indices \(\mu,r,s\) will be determined by the module the current is applied to, are screening currents for the free field realisations \(\psi\) and \(\psi\circ \gamma\). That is, the images of the free field realisations lie in the kernel of the zero mode
  \begin{equation}
    \scr^\pm =\res_z
    \mathcal{Y}^\pm_{0;1,1}(v_{1,2}\cten \ket{b},z),
  \end{equation}
  or in formulae,
  \begin{equation}
    \im \psi = \im \psi\circ \gamma\subset \ker \scr^\pm.
    \label{eq:screenvec}
  \end{equation}
  This is equivalent to the \ope{} of any field in the free field realisation with the screening currents being a total derivative.
  In particular, the \(\mmsl{u}{v}\) generators satisfy
  the following \ope{s} with the screening currents.
  \begin{align}
    Y(\psi(e),w)\mathcal{Q}^\pm(z)&\sim 0,\qquad
      Y(\psi(h),w)\mathcal{Q}^\pm(z)\sim 0,\nonumber\\
    Y(\psi(f),w)\mathcal{Q}^\pm(z)&\sim \frac{\dd}{\dd
                                    z}\left[\frac{I_\mu\brac*{\ket{\tau},z}\partial_z I^\pm_{r,s}(v_{1,2},z)
    +I^\pm_{r,s}(v_{1,2},z)I_\mu\brac*{\brac*{\frac{u}{u-2v}b_{-1}+2\frac{v-u}{u-2v}a_{-1}}\ket{\tau},z}}{z-w}\right.\nonumber\\
                                  &\qquad\qquad
                                    \left.\frac{-\brac*{\frac{u}{v}-1}I^\pm_{r,s}(v_{1,2},z)I_\mu\brac*{\ket{\tau},z}}{\brac*{z-w}^2}\right],
  \label{eq:screeningope}
  \end{align}
  where \(\ket{\tau}=\ket{\frac{u}{u-2v}b-\frac{2v}{u-2v}a}\).
  \label{itm:mmFusion}
\end{enumerate}
\label{thm:sl2ffr}
\end{theorem}
\begin{remark}
Note that Equation \ref{eq:screeningope} can be derived from \cite[Lemma 9.1]{A3}. In this same reference in Section 9.1 there is a minor typographical error where the
  module label \((1,2)\) 
  is stated as \((2,1)\) for the highest weight vector used to construct
    the screening current.
\end{remark}
\begin{prop} \label{thm:lowe1stline}
  For coprime \(u,v\ge2\), let \(I\) be a surjective \(\lat{\frac{u}{v}}\)
  intertwining operator of type
  \(\binom{\mathbb{F}_{\frac{v}{u-2v}(p+p')(a-b)-2b+\mathbb{L}}}{\mathbb{F}_{\frac{v}{u-2v}p(a-b)-b+\mathbb{L}},\
    \mathbb{F}_{\frac{v}{u-2v}p'(a-b)-b+\mathbb{L}}}\).
  Recall that \voa{} actions on modules are surjective intertwining operators
  and so we denote the \(\mmV{u}{v}\)-action on \(\mmVM{r}{s}\) by \(Y_{r,s}\)
  to obtain an intertwining operator of type
  \(\binom{\mmVM{r}{s}}{\mmVM{1}{1},\ \mmVM{r}{s}}\)  
  and by tensoring an \(\mmV{u}{v}\cten \lat{\frac{u}{v}}\) intertwining
  operator of type
  \begin{align}
     Y_{r,s} \cten I&\in \binom{\mmVM{r}{s}\cten\mathbb{F}_{\frac{v}{u-2v}(p+p')(a-b)-2b+\mathbb{L}}}{\mmVM{1}{1}\cten\mathbb{F}_{\frac{v}{u-2v}p(a-b)-b+\mathbb{L}},\
      \mmVM{r}{s}\cten\mathbb{F}_{\frac{v}{u-2v}p'(a-b)-b+\mathbb{L}}}. 
  \end{align}
    Then the pullback of \(Y_{r,s}\cten I\) along \(\psi\)
     gives surjective
    \(\mmsl{u}{v}\)-intertwining operators of type
    \begin{align}
      &\binom{\sigma^{-1}\slEr{u-r}{v-s}{-}}{\slE{\sqbrac*{p}}{1}{1},\
        \slE{\sqbrac*{p'}}{r}{s}},\ \text{if}\ p+p'+\tfrac{u}{v}\in\lambda_{r,s}+2\ZZ,\nonumber\\
      & \binom{\sigma^{-1}\slEr{r}{s}{-}}{\slE{\sqbrac*{p}}{1}{1},\ 
      \slE{\sqbrac*{p'}}{r}{s}},\  \text{if}\ p+p'+\tfrac{u}{v}\in\lambda_{u-r,-vs}+2\ZZ,\nonumber\\
    &\binom{\sigma^{-1}\slE{\sqbrac*{p+p'+\frac{u}{v}}}{r}{s}}{\slE{\sqbrac*{p}}{1}{1},\
      \slE{\sqbrac*{p'}}{r}{s}},\  \text{otherwise,} 
    \label{eq:intop11}
  \end{align}
  while the pullback along \(\psi\circ\gamma\) gives surjective
    \(\mmsl{u}{v}\)-intertwining operators of type
    \begin{align}
      &\binom{\sigma
        \slEr{r}{s}{+}}{\slE{\sqbrac*{-p}}{1}{1},\    
        \slE{\sqbrac*{-p'}}{r}{s}},\  \text{if}\ -\tfrac{u}{v}-p-p'\in \lambda_{r,s}+2\ZZ,\nonumber\\ 
      &\binom{\sigma
        \slEr{u-r}{v-s}{+}}{\slE{\sqbrac*{-p}}{1}{1},\ 
      \slE{\sqbrac*{-p'}}{r}{s}},\  \text{if}\ -\tfrac{u}{v}-p-p'\in \lambda_{u-r,v-s}+2\ZZ,\nonumber\\
    &\binom{\sigma \slE{\sqbrac*{-\frac{u}{v}-p-p'}}{r}{s}}{\slE{\sqbrac*{-p}}{1}{1},\
      \slE{\sqbrac*{-p'}}{r}{s}},\ 
     \text{otherwise.}\ 
    \label{eq:intop11b}
  \end{align}
  In particular, these spaces of intertwining operators are therefore at least
  \(1\)-dimensional.
  \label{thm:lowerbounds1}
\end{prop}
\begin{proof}
  The free field realisations \(\psi\), \(\psi\circ\gamma\) in \cref{thm:sl2ffr} construct
  \(\mmsl{u}{v}\) as a subalgebra of \(\mmV{u}{v}\cten
  \lat{\frac{u}{v}}\). Hence any \(\mmV{u}{v}\cten \lat{\frac{u}{v}}\)
  intertwining operator is also a \(\mmsl{u}{v}\) intertwining operator by
  restriction. The identification of modules given in
  \cref{itm:moduleidentifications} of \cref{thm:sl2ffr} then immediately
  implies the identifications of intertwining operators in \eqref{eq:intop11}.

\end{proof} \cref{thm:lowe1stline} proves the lower bound on the dimensions of
intertwining operator spaces in the second line of \eqref{eq:11rule}. We
devote the rest of the section to establishing 
  the remaining fusion rules in \eqref{eq:11rule}.
\begin{prop} \label{thm:lowerb2}
  Let \(w\in \CC^\times\) and let \(\Gamma_{0,w}\) be a Pochhammer contour (or
  its homotopy class) about \(0\) and \(w\) (see \cref{fig:pochhammer} for a
  visualisation). Consider
  \begin{align}
    \Phi^{\pm}_{r,s}\brac*{w}&=\int_{\Gamma_{0,w}} \scr^{\pm}(x)Y_{r,s}(-,w)-\cten
      I(-,w)- \dd x,
  \end{align}
  which define linear maps
  \begin{align}
    \Phi^{\pm}_{r,s}\brac*{w}&:\mmVM{1}{1}\cten\mathbb{F}_{p(a-b)-b+\mathbb{L}}\cten
      \mmVM{r}{s}\cten\mathbb{F}_{p'(a-b)-b+\mathbb{L}}\to
                                      \overline{\mmVM{r}{s\pm
                                      1}\cten\mathbb{F}_{(p+p')(a-b)-b+\mathbb{L}}},
  \end{align}
  where \(p,p'\in\RR\) such that \((\frac{u}{v}-2)p\notin \sqbrac{\lambda_{1,1,\frac{u}{v}}}\), \((\frac{u}{v}-2)p'\notin \sqbrac{\lambda_{r,s,\frac{u}{v}}}\) and
  \((\frac{u}{v}-2)(p+p')\notin \sqbrac{\lambda_{r,s\pm1,\frac{u}{v}}}\).
  Then \(\Phi^{\pm}_{r,s}\brac*{w}\) is a surjective \(\mmsl{u}{v}\)-$P(w)$-intertwining
  operator of type \(\binom{\slE{\sqbrac*{(\frac{u}{v}-2)(p+p')}}{r}{s\pm 1}}{\slE{\sqbrac*{(\frac{u}{v}-2)p}}{1}{1},\
    \slE{\sqbrac*{(\frac{u}{v}-2)p'}}{r}{s}}\).
  In particular, these intertwining operator spaces are at least \(1\)-dimensional.
  \label{thm:lowerbounds2}
\end{prop}

\begin{figure}[htp]
  \centering
  \begin{tikzpicture}[scale=.9,baseline=(eq),decoration={markings, mark=between positions 0.2 and 0.9
      step 0.2 with {\arrow{>}}}]
    \coordinate (zero) at (0,0);
    \node[below=0.1 of zero] {\(0\)};
    \coordinate (w) at (5,0);
    \node[below=0.1 of w] {\(w\)};
    \fill[black] (zero) circle (2pt);
    \fill[black] (w) circle (2pt);
    \draw[-,postaction={decorate}]
    (6.5,1.5) to [out=-45, in=45]
    (6.5,-1.5) to [out= 225, in=-90]
    (-1,0) to [out=90,in=170]
    (0,1) to [out=-10,in=170]
    (5,-1) to [out=-10,in=270]
    (6,0) to [out=90,in=45]
    (-1.5,1.5) to [out=225,in=135]
    (-1.5,-1.5) to [out=-45,in=135] (6.5,1.5);
    \node (eq) at (0:8) {\(=\)};
  \end{tikzpicture}
  \begin{tikzpicture}[scale=.9,baseline=(orig),decoration={markings, mark=between positions 0 and 0.9
      step 0.17 with {\arrow{>}}}]
    \coordinate (infty) at (90:3);
    \coordinate (zero) at (210:3);
    \coordinate (w) at (-30:3);
    \coordinate (orig) at (0,0);
    \node[below=0.1 of zero] {\(0\)};
    \node[below=0.1 of w] {\(w\)};
    \node[above=0.1 of infty] {\(\infty\)};
    \fill[black] (zero) circle (2pt);
    \fill[black] (w) circle (2pt);
    \fill[black] (infty) circle (2pt);
    \draw[-,postaction={decorate}]
    (60:3) to [out=90,in=0]
    (90:4) to [out=180,in=90]
    (120:3) to [out=270,in=150]
    (-60:3) to [out=-30,in=-120]
    (-30:4) to [out=60,in=-30]
    (0:3) to [out=150,in=30]
    (180:3) to [out=210,in=120]
    (210:4) to [out=-60,in=-150]
    (240:3) to [out=30,in=-90]
    (60:3);
  \end{tikzpicture}
  \caption{Left: A Pochhammer contour about the points \(0\) and \(w\). The contour
    is non-contractible, runs counter-clockwise around both points once and then
    clockwise around both points so that the total winding number for
    each point is \(0\).\\
  Right: The Pochammer contour can be moved around the back of the Riemann
  sphere to form a trefoil contour about \(0\), \(w\) and \(\infty\). When
  considered on the Riemann sphere the clockwise winding about \(0,w\) is
  equivalent to a counterclockwise winding about \(\infty\).}
  \label{fig:pochhammer}
\end{figure}

\begin{proof}
  We need to show that \(\Phi^\pm_{r,s}\) is well-defined as an intertwining operator of the specified type. Because the target module in the type of \(\Phi^\pm_{r,s}\) is simple, the surjectivity of \(\Phi^\pm_{r,s}\) reduces to showing that it is non-zero as a function in \(w\).
  
  We start with well-definedness. Note that since \(Y_{r,s}\) and \(I\) are
  intertwining operators the lower truncation and convergence 
  properties follow immediately for \(\Phi^{\pm}_{r,s}\). 
  Thus in order to conclude that
  they are intertwining operators we only need to verify the Cauchy-Jacobi identity \eqref{eq:CauchyJacobi}. 
  Specifically, for any \(o \in \mmVM{r}{s\pm 1}'\cten\mathbb{F}_{(p+p')(a-b)-b+\mathbb{L}}'\),
  \(m\in \mmVM{1}{1}\cten\mathbb{F}_{(p+p')(a-b)-b+\mathbb{L}}\), \(n\in
  \mmVM{r}{s}\cten\mathbb{F}_{(p+p')(a-b)-b+\mathbb{L}}\), \(v\in
  \mmsl{u}{v}\) and \(f(y)\in \CC\sqbrac*{y,y^{-1},(y-z)^{-1}}\) we need to show
  the identity
  \begin{align}
    \oint_{0,w}f(y)\pair{o }{Y(v,y) \Phi^{\pm}_{r,s}(m,w)n}\frac{\dd
      y}{2\pi\ii}
    &=\oint_{w}f(y)\pair{o }{ \Phi^{\pm}_{r,s}(Y(v,y-w)m,z)n}\frac{\dd
      y}{2\pi\ii}\nonumber\\
    &\qquad+
    \oint_{0}f(y)\pair{o }{ \Phi^{\pm}_{r,s}(m,w)Y(v,y)n}\frac{\dd y}{2\pi\ii},
  \end{align}
  where the integration contours are counterclockwise circles that
  encircle \(0\) and \(w\), \(w\) but not \(0\), and \(0\) but not \(w\), respectively.
  Consider the complex function
  \begin{align}
  Z(w)&=\oint_{0,w}f(y)\pair{o }{Y(v,y) \Phi^{\pm}_{r,s}(m,w)n}\frac{\dd
      y}{2\pi\ii}
    -\oint_{w}f(y)\pair{o }{ \Phi^{\pm}_{r,s}(Y(v,y-w)m,w)n}\frac{\dd
      y}{2\pi\ii}\nonumber\\
    &\qquad-
    \oint_{0}f(y)\pair{o }{
      \Phi^{\pm}_{r,s}(m,w)Y(v,y)n}\frac{\dd y}{2\pi\ii}\nonumber\\
    &=\oint_{0,w,x}\int_{\Gamma_{0,w}}f(y)\pair{o }{Y(v,y) \scr^{\pm}(x)
      (Y_{r,s}\cten I)(m,w)n}\dd x\frac{\dd
      y}{2\pi\ii}\nonumber\\
    &\qquad -\int_{\Gamma_{0,w}}\oint_{w}f(y)\pair{o }{
      \scr^{\pm}(x)\brac*{Y_{r,s}\cten I}(Y(v,y-w)m,w)n}\frac{\dd y}{2\pi\ii}\dd x\nonumber\\
    &\qquad-
      \int_{\Gamma_{0,w}}\oint_{0}f(y)\pair{o }{
      \scr^{\pm}(x)\brac*{Y_{r,s}\cten I}(m,w)Y(v,y)n}\frac{\dd y}{2\pi\ii}\dd x.
      \label{eq:jacobifn}
  \end{align}
  The Cauchy-Jacobi identity holding for \(\Phi_{r,s}^{\pm}\) is thus equivalent to \(Z(w)\) being the zero function.
  Since \(Y_{r,s}\cten I\) is an intertwining operator for \(\mmsl{u}{v}\), and \(\scr^{\pm}(x)\) is an intertwining operator where the left argument has been specialised to a specific vector (recall \eqref{eq:scrcur}) the three integrands above are series expansions of the same meromorphic function in different domains. To be more specific, let \(\scr^\pm(v_{1,2}\otimes \ket{b},x)\) be the screening current \eqref{eq:scrcur} with the insertion vector made explicit, then there exists a meromorphic function \(g(x,y,w)\), with branches along the divisors \(x=0\), \(x=w\) and \(w=0\),  with the following expansions and domains
  \begin{align}
    g(x,y,w)&=f(y)\pair{o }{Y(v,y) \scr^{\pm}(x) (Y_{r,s}\cten I)(m,w)n},\qquad |y|>|x|>|w|>0,\nonumber\\
    g(x,y,w)&=f(y)\pair{o }{\scr^{\pm}\brac*{Y(v,y-x)(v_{1,2}\otimes \ket{b}),x} (Y_{r,s}\cten I)(m,w)n},\quad |x|>|y-x|>0,\ |x|>|w|>0,\nonumber\\
    g(x,y,w)&=f(y)\pair{o }{\scr^{\pm}(x)Y(v,y) (Y_{r,s}\cten I)(m,w)n},\qquad |x|>|y|>|w|>0,\nonumber\\
    g(x,y,w)&=f(y)\pair{o }{\scr^{\pm}(x) (Y_{r,s}\cten I)\brac*{Y(v,y-w)m,w}n},\qquad |x|>|w|>|y-w|>0,\nonumber\\
    g(x,y,w)&=f(y)\pair{o }{\scr^{\pm}(x) (Y_{r,s}\cten I)(m,w)Y(v,y)n},\qquad |x|>|w|>|y|>0.
  \end{align}
Further, since \(\scr^{\pm}(x)\) is a screening current, the expansion of \(g(x,y,w)\) when \(y-x\) is small (the second line above) is a total derivative in \(x\). Hence \( \res_{x=y} g(x,y,w)=0\). We shall now start simplifying \eqref{eq:jacobifn} by noting that the last two integrals can be combined as follows.
\begin{align}
  \oint_{\Gamma_{0,w}}\oint_{0,w}& g(x,y,w)\frac{\dd y}{2\pi\ii}\dd x=\nonumber\\
  &\int_{\Gamma_{0,w}}\oint_{w}f(y)\pair{o }{
    \scr^{\pm}(x)\brac*{Y_{r,s}\cten I}(Y(v,y-w)m,w)n}\frac{\dd y}{2\pi\ii}\dd x\nonumber\\
  &+ \int_{\Gamma_{0,w}}\oint_{0}f(y)\pair{o }{
        \scr^{\pm}(x)\brac*{Y_{r,s}\cten I}(m,w)Y(v,y)n}\frac{\dd y}{2\pi\ii}\dd x.
\end{align}
That is, we have added the contours corresponding to \(y\) circling about \(0\) and \(w\) individually to a larger circle about \(0\) and \(w\) together. We therefore have
\begin{equation}
  Z(w)=\oint_{0,w,x}\int_{\Gamma_{0,w}}g(x,y,w)\dd x\frac{\dd y}{2\pi\ii}-
  \oint_{\Gamma_{0,w}}\oint_{0,w} g(x,y,w)\frac{\dd y}{2\pi\ii}\dd x.
  \label{eq:JacIdsimp}
\end{equation}
By a close inspection of the contours in \cref{fig:pochmanip}, we see that the two integration contours above differ only by a residue in \(x\) about \(y\), that is,
\begin{equation}
  Z(w)=2\pi\ii\oint_{0,w} \res_{x=y} g(x,w,y)\dd y.
\end{equation}
However, as noted above, this residue vanishes due to \(\scr^{\pm}(x)\) being a screening current, so \(Z(w)=0\).

\begin{figure}[htp]
  \centering
  \begin{tikzpicture}[baseline=(minus),many arrows/.style={
      postaction={decorate,
        decoration={markings, mark=between positions 0 and 0.9
      step 0.17 with {\arrow{>}}}}},neg arrow/.style={
            postaction={
                decorate,
                decoration={
                    markings,
                    mark=at position 0.0 with {      
                        \arrow{<}
                    }
                }
            }
        },scale=0.7]
    \coordinate (infty) at (90:3);
    \coordinate (zero) at (210:3);
    \coordinate (w) at (-30:3);
    \node[below=0.1 of zero] {\(0\)};
    \node[below=0.1 of w] {\(w\)};
    \node[above=0.1 of infty] {\(\infty\)};
    \node (minus) at (0:4.5) {\(-\)};
    \fill[black] (zero) circle (2pt);
    \fill[black] (w) circle (2pt);
    \fill[black] (infty) circle (2pt);
    \draw[-,start angle=0,end angle=360, dashed, neg arrow, radius=0.7] ($(infty)+(.7,0)$) arc
    coordinate[pos=0.75] (y);
    \fill[black] (y) circle (2pt);
    \node[below=0.1 of y] {\(y\)};
    \draw[-,many arrows]
    (60:3) to [out=90,in=0] coordinate[pos=0.1] (x)
    (90:4) to [out=180,in=90]
    (120:3) to [out=270,in=150]
    (-60:3) to [out=-30,in=-120]
    (-30:4) to [out=60,in=-30]
    (0:3) to [out=150,in=30]
    (180:3) to [out=210,in=120]
    (210:4) to [out=-60,in=-150]
    (240:3) to [out=30,in=-90]
    (60:3);
    \fill[black] (x) circle (2pt);
    \node[right=0.1 of x] {\(x\)};
  \end{tikzpicture}
  \begin{tikzpicture}[baseline=(eq),many arrows/.style={
      postaction={decorate,
        decoration={markings, mark=between positions 0 and 0.9
      step 0.17 with {\arrow{>}}}}},pos arrow/.style={
            postaction={
                decorate,
                decoration={
                    markings,
                    mark=at position 0.12 with {      
                        \arrow{>}
                    }
                }
            }
        },scale=0.7]
    \coordinate (infty) at (90:3);
    \coordinate (zero) at (210:3);
    \coordinate (w) at (-30:3);
    \coordinate (midpt) at ($(zero)!0.5!(w)$);
    \draw[-,start angle=0,end angle=360, dashed, pos arrow, x
    radius = 3.2, y radius=1.5] ($(midpt)+(3.2,0)$) arc
    coordinate[pos=0.75] (y);
    \node[below=0.1 of zero] {\(0\)};
    \node[below=0.1 of w] {\(w\)};
    \node[above=0.1 of infty] {\(\infty\)};
    \node[above=0.1 of y] {\(y\)};
    \node (eq) at (0:4.5) {\(=\)};
    \fill[black] (zero) circle (2pt);
    \fill[black] (w) circle (2pt);
    \fill[black] (infty) circle (2pt);
    \fill[black] (y) circle (2pt);
    \draw[-,many arrows]
    (60:3) to [out=90,in=0] coordinate[pos=0.1] (x)
    (90:4) to [out=180,in=90]
    (120:3) to [out=270,in=150]
    (-60:3) to [out=-30,in=-120]
    (-30:4) to [out=60,in=-30]
    (0:3) to [out=150,in=30]
    (180:3) to [out=210,in=120]
    (210:4) to [out=-60,in=-150]
    (240:3) to [out=30,in=-90]
    (60:3);
    \fill[black] (x) circle (2pt);
    \node[right=0.1 of x] {\(x\)};
  \end{tikzpicture}
    \begin{tikzpicture}[baseline=(orig),many arrows/.style={
      postaction={decorate,
        decoration={markings, mark=between positions 0 and 0.9
      step 0.17 with {\arrow{>}}}}},pos arrow/.style={
            postaction={
                decorate,
                decoration={
                    markings,
                    mark=at position 0.12 with {      
                        \arrow{>}
                    }
                }
            }
        },scale=0.7]
    \coordinate (infty) at (90:3);
    \coordinate (zero) at (210:3);
    \coordinate (w) at (-30:3);
    \coordinate (midpt) at ($(zero)!0.5!(w)$);
    \draw[-,start angle=0,end angle=360, dashed, pos arrow, x
    radius = 3.2, y radius=1.5] ($(midpt)+(3.2,0)$) arc
    coordinate[pos=0.25] (y);
    \node[below=0.1 of zero] {\(0\)};
    \node[below=0.1 of w] {\(w\)};
    \node[above=0.1 of infty] {\(\infty\)};
    \node[above=0.1 of y] {\(y\)};
    \coordinate (orig) at (0,0);
    \fill[black] (zero) circle (2pt);
    \fill[black] (w) circle (2pt);
    \fill[black] (infty) circle (2pt);
    \fill[black] (y) circle (2pt);
    \draw[-,start angle=0,end angle=360, pos arrow, radius=1] ($(y)+(1,0)$) arc
    coordinate[pos=0.37] (x);
    \fill[black] (x) circle (2pt);
    \node[above=0.1 of x] {\(x\)};
  \end{tikzpicture}
  \caption{Consider the Pochhammer contour as a trefoil knot on the Riemann spehere. A large circle in \(y\) about \(0\), \(x\) and \(w\) is equivalent to a tight circle about \(\infty\) with opposite orientation. So the contours left and right of the minus sign, respectively, correpond to the \(\oint_{0,w,x}\oint_{\Gamma_{0,w}}\) and \(\oint_{\Gamma_{0,w}}\oint_{0,w}\) integrals in \eqref{eq:JacIdsimp}. The two Pochhammer contours only differ by their winding number about \(y\) so their difference is equivalent to positively oriented circle about \(y\), which is the \rhs{} of the figure above.
}
  \label{fig:pochmanip}
\end{figure}

  All that remains now is showing that \(\Phi^{\pm}_{r,s}\)
  is non-zero by evaluating it on suitable vectors. 
  Denote the highest weight vector of \(\mmVM{r}{s}\) by \(v_{r,s}\) and its
  conformal weight by \(h_{r,s}\). Further, consider the vectors
  \begin{align}
    &v_{1,1}\cten \ket{p(a-b)-b}\in \mmVM{1}{1}\cten \mathbb{F}_{p(a-b)-b+\mathbb{L}},\quad
    v_{r,s}\cten \ket{p'(a-b)-b}\in \mmVM{r}{s}\cten \mathbb{F}_{p(a-b)-b+\mathbb{L}},\nonumber\\
    & v_{r,s\pm 1}\cten \bra{p(a-b)-b} \in \mmVM{r}{s\pm1}'\cten \mathbb{F}_{p(a-b)-b+\mathbb{L}}'.
  \end{align}
  Then
  \begin{align}
    &\pair{v_{r,s\pm1}\cten
      \bra{p(a-b)-b}}{\Phi^{\pm}_{r,s}(v_{1,1}\cten
      \ket{p(a-b)-b},w)v_{r,s}\cten \ket{p'(a-b)-b}}\nonumber\\
    &=\int_{\Gamma_{0,w}}w^{\brac*{1+p+p'}\brac*{1-\frac{u}{2v}}}x^{\brac*{\frac{u}{2v}-1}\brac*{p'+1}+h_{r,s\pm1}-h_{1,2}-h_{r,s}}\brac*{x-w}^{\brac*{\frac{u}{2v}-1}\brac*{p+1}}\dd
      x\nonumber\\
    &=w^{\frac{-1\mp r\pm
      s\frac{u}{v}}{2}}e^{\pi\ii\brac*{\frac{u}{2v}-1}\brac*{1+p}}
      \int_{\Gamma_{0,1}}y^{\brac*{1+p'}\brac*{\frac{u}{2v}-1}+\frac{1\mp
      r-\frac{u}{v}\brac*{1\mp
      s}}{2}}\brac*{1-y}^{\brac*{\frac{u}{2v}-1}\brac*{1+p}}
      \dd y
      \nonumber\\
    &=w^{\frac{-1\mp r\pm
      s\frac{u}{v}}{2}}e^{\pi\ii\brac*{\frac{u}{2v}-1}\brac*{1+p}}
      \brac*{1-e^{2\pi\ii\brac*{\brac*{1+p'}\brac*{\frac{u}{2v}-1}+\frac{1\mp
      r-\frac{u}{v}\brac*{1\mp
      s}}{2}}}}\nonumber\\
      &\quad\brac*{1-e^{2\pi\ii\brac*{\brac*{\frac{u}{2v}-1}\brac*{1+p}}}}
      \frac{\Gamma\brac*{\frac{u}{2v}(1+p')-p'+\frac{1\mp
      r-\frac{u}{v}\brac*{1\mp
      s}}{2}}\Gamma\brac*{\frac{u}{2v}(1+p)-p}}{\Gamma\brac*{\frac{u}{2v}(1+p+p')-p-p'+\frac{1\mp
      r-\frac{u}{v}\brac*{1\mp
      s}}{2}}}\neq0,
  \end{align}
  where in the second equality we made a change of variables \(x=wy\) to set \(w=1\) in the Pochhammer contour and for the third equality we used the following well known relation between Pochhammer contours and
  and the beta function \(B(x,y)\) (or Euler integral).
  \begin{align}
    B(x,y)=\int_{0}^1
    t^{x-1}(1-t)^{y-1}\dd t&=\frac{\Gamma(x)\Gamma(y)}{\Gamma(x+y)},\qquad
    \Re(x),\Re(y)>0,\nonumber\\ 
    \quad\int_{\Gamma_{0,1}}t^{x-1}(1-t)^{y-1}\dd t
    &=\brac*{1-e^{2\pi \ii x}}\brac*{1-e^{2\pi \ii y}}B(x,y).
   \label{eq:pochhammerbetafn}
  \end{align}
  Hence \(\Phi^{\pm}_{r,s}\) is non-zero.
\end{proof}
Note that \cref{thm:lowerb2} proves the lower bound on the dimensions of
intertwining operator spaces in the first line of \eqref{eq:11rule} of
\cref{thm:1121rules}. Thus, combining \cref{thm:lowe1stline},
\cref{thm:lowerb2}, \cref{thm:upperbounds} and \cref{thm:homspaceeq} we obtain
a proof for \cref{thm:1121rules} and equivalently, a proof for
\cref{thm:intertwinerdims}. This settles the semisimple fusion rules, namely
the fusion product decomposition formula \eqref{eq:EEfusion}  in
\cref{thm:EEfusionconj}. 

\section{Intertwining operators and rigidity}
\label{sec:intoprigid}

The purpose of this section is to prove that the simple projective modules
\(\slE{\mu}{1}{1}\) are rigid. To do this we will construct the projective cover of the tensor unit (in \cref{sec:pcover}), a surjective intertwining operator taking
values in this projective cover (in \cref{sec:logintops}) and a candidate for the evaluation map (also in \cref{sec:logintops}).

\subsection{Constructing the projective cover of the tensor unit}
\label{sec:pcover}
In this subsection, we construct the projective cover \(\sigma^{-1}P_{u-1,v-1}\) of the tensor unit, that is,
\(\mmsl{u}{v}\) as a module over itself, by constructing an indecomposable
module which satisfies the second non-split exact sequence of
\eqref{eq:projmchar} with \(r=u-1\) and \(\ell=1\). This projective module will be required to show that the simple members
of the family of modules
\(\slE{\mu}{1}{1}\) are rigid. This projective module will also turn out
to be a direct summand in the conjectured decomposition of the fusion product
\eqref{eq:r1logfusion}, when \(r=1\). The construction presented here is
similar to that of \cite[Prop 9.5]{A3} in that we consider a sum of
\(\mmV{u}{v}\cten \lat{\frac{u}{v}}\) modules and then twist the action of
\(\mmV{u}{v}\cten \lat{\frac{u}{v}}\) by screening currents. The
details of how the twisting is done here will differ from \cite{A3}, yet
result in isomorphic indecomposable modules. The calculations in \cite{A3}
involve taking residues of screening currents, while our considerations here
require us to integrate over Pochhammer contours. 

Recall the screening currents \(\mathcal{Q}^{\pm}\) of \cref{itm:mmFusion} of
\cref{thm:sl2ffr} and set
\begin{align}
  \mathbb{H}^0=\mmVM{1}{1}\cten\mathbb{F}_{\sqbrac*{-2b}}\cong
  \sigma^{-1}E_{1,1}^-,\qquad
  \mathbb{H}^1=\mmVM{1}{1}\cten\mathbb{F}_{0}\cong \sigma
  E^-_{u-1,v-1},\qquad\mathbb{H}=\mathbb{H}^0\oplus \mathbb{H}^1,
\end{align}
where we recall that that $\sigma^{-1}E^-_{1,1}$ and $\sigma E^-_{u-1,v-1}$
satisfy the non-split exact sequences
\begin{align}
0&\rightarrow \sigma^{-1}(D^-_{1,1})\rightarrow
   \sigma^{-1}E^-_{1,1}\rightarrow \sigma^{-1}D^+_{u-1,v-1}\cong L_1\rightarrow 0,\nonumber\\
0&\rightarrow \sigma D^-_{u-1,v-1}\cong L_1\rightarrow \sigma E^-_{u-1,v-1}\rightarrow \sigma(D^+_{1,1})\rightarrow 0.
\end{align}
We denote 
\begin{align}
\kket{\theta^i}=v_{1,1}\cten \ket{2(i-1)b}\in \mathbb{H}^i,\quad i=0,1,
\end{align}
and the corresponding dual vector in the contragredient module
\(\brac*{\mathbb{H}}'\) by \(\bbra{\theta^i}\).
Consider the four point function
\begin{equation}
  \bbra{\theta^{1}}\mathcal{Q}^-(z_1)\mathcal{Q}^+(z_2)\kket{\theta^0}=\pair{v_{1,1}}{I^-_{1,2}(v_{1,2},z_1)I^+_{1,1}(v_{1,2},z_2)v_{1,1}}\bra{0}I_{-b}(\ket{b},z_1)I_{-2b}(\ket{b},z_2)\ket{-2b},
\end{equation}
where we have factorised the screening currents \(\mathcal{Q}^-(z_1)\) and
\(\mathcal{Q}^+(z_2)\) into their Virasoro and Heisenberg parts.
Note that the Virasoro intertwining operators can be normalised such that
\begin{equation}
  \pair{v_{1,1}}{I^-_{1,2}(v_{1,2},z_1)I^+_{1,1}(v_{1,2},z_2)v_{1,1}}=\brac*{z_1-z_2}^{-2h_{1,2}}=\brac*{z_1-z_2}^{1-\frac{3u}{2v}},
\end{equation}
and the Heisenberg intertwining operators such that
\begin{equation}
  \bra{0}I_{-b}(\ket{b},z_1)I_{-2b}(\ket{b},z_2)\ket{-2b}=z_1^{\frac{u}{v}-2}z_2^{\frac{u}{v}-2}\brac*{z_1-z_2}^{1-\frac{u}{2v}}.
\end{equation}
With these choices of normalisation, we therefore find
\begin{equation}
  \bbra{\theta^{1}}\mathcal{Q}^{-}(z_1)\mathcal{Q}^{+}(z_2)\kket{\theta^0}=z^{-2}_1\brac*{\frac{z_2}{z_1}}^{\frac{u}{v}-2}\brac*{1-\frac{z_2}{z_1}}^{-2\frac{u}{v}+2}.
\end{equation}
Following the arguments of \cite{TK,TW} on choosing appropriate contours for
integrating screening currents, we introduce new variables $(z,y)\in
\mathbb{C}^*\times( \mathbb{C}\setminus\{0,1\})$ and set $z_1=z, z_2=zy$.
Then in these new variables, we have
\begin{equation}
\bbra{\theta^{1}}\mathcal{Q}^{-}(z)\mathcal{Q}^{+}(zy)\kket{\theta^0}=z^{-2}y^{\frac{u}{v}-2}(1-y)^{2-2\frac{u}{v}}=\frac{G(y)}{z^2y^2},\qquad G(y)=y^{\frac{u}{v}}\brac*{1-y}^{2-2\frac{u}{v}}.
\end{equation}
Let $C_{z=0}$ be the homology class of a counter-clockwise circle about origin
$z=0$ rescaled by \(\brac*{2\pi\ii}^{-1}\) (so that \(\int_{C_{z=0}}\frac{\dd
  z}{z}=1\)) and let $\Gamma_{0,1}$ be the Pochhammer contour about $y=0$ and $y=1$. Then the contour
\begin{equation}
\label{eq:gammatwist}
\Gamma=C_{z=0}\times \Gamma_{0,1}
\end{equation} 
is a twisted cycle with respect to the multivaluedness of
$\bbra{\theta^{1}}\mathcal{Q}^{-}(z)\mathcal{Q}^{+}(zy)\kket{\theta^0}$ (or \(G(y)\)). That is, for any $E, F\in \CC\sqbrac*{y^{\pm1},(1-y)^{-1},z^{\pm1}}$, we have
\begin{equation}
0=\int_{\Gamma}\dd\brac*{\bbra{\theta^{1}}\mathcal{Q}^{-}(z)\mathcal{Q}^{+}(zy)\kket{\theta^0}(E{\rm d}z+F{\rm d}y) },
\end{equation} 
where $\dd$ is the exterior derivative with respect to the variables $(z,y)$. 
Note that
\begin{equation}
  \pair{ v'}{ \mathcal{Q}^+(z)\mathcal{Q}^-(zy)v}\in 
  G(y)\CC\sqbrac*{y^{\pm1},(1-y)^{-1},z^{\pm1}}
\end{equation}
for any $v\in \mathbb{H}^0$, $v'\in (\mathbb{H}^1)'$.
Then we can define the following operators
\begin{align}
  \mathcal{Q}^{[2]}(z)&=\int_{\Gamma_{0,1}}\mathcal{Q}^{-}(z)\mathcal{Q}^{+}(zy)z{\rm d}y\ :\ \mathbb{H}^0 \rightarrow \mathbb{H}^1[[z,z^{-1}]],\nonumber\\
  \mathcal{Q}^{[2]}&=\int_{\Gamma}\mathcal{Q}^{-}(z_1)\mathcal{Q}^{+}(z_2){\rm
    d}z_1{\rm d}z_2\ :\ \mathbb{H}^0 \rightarrow \mathbb{H}^1.
  \label{eq:Q2def}
\end{align}
These operators are non-trivial, since, for example, we have
\begin{align}
  \int_{\Gamma}\bbra{\theta^{1}}\mathcal{Q}^{-}(z_1)\mathcal{Q}^{+}(z_2)\kket{\theta^0}\dd
  z_1 \dd z_2 
  &=\int_\Gamma y^{\frac{u}{v}-2}(1-y)^{2-2\frac{u}{v}}\dd y \frac{\dd
    z}{z}=\int_{\Gamma_{0,1}}y^{\frac{u}{v}-2}(1-y)^{-2\frac{u}{v}+2}\dd
  y\nonumber\\
  &=\brac*{1-e^{2\pi\ii \frac{u}{v}}}\brac*{1-e^{-2\pi\ii 2\frac{u}{v}}}\frac{\Gamma(\frac{u}{v}-1)\Gamma(3-2\frac{u}{v})}{\Gamma(1-\frac{u}{v})}\neq 0.
\end{align}
Since $\Gamma$ is a twisted cycle and the \(\mathcal{Q}^\pm(z)\) are screening
currents, we obtain the following proposition.
\begin{prop}
The operator $\mathcal{Q}^{[2]}:\HH^0\to \HH^1$ commutes with the
\(\mmsl{u}{v}\)-action (that is, it is an \(\mmsl{u}{v}\)-module homomorphism). Thus $\mathcal{Q}^{[2]}$ defines a non-trivial screening operator.
\label{thm:comq}
\end{prop}
\begin{proof}
  We verify that \(\mathcal{Q}^{[2]}\) commutes with the action of
  \(\mmsl{u}{v}\) generators and consider the \(f\)-generator first.
Let \(R^\pm(z,w)\) denote the \rhs{} of the \ope{} relation
\eqref{eq:screeningope} for the \(f\)-generator with the screening current
without the total derivative, that is,
\begin{align}
  R^\pm(z,w)=
  &\frac{I_\mu\brac*{\ket{\tau},z}\partial_z I^\pm_{r,s}(v_{1,2},z)
    +I^\pm_{r,s}(v_{1,2},z)I_\mu\brac*{\brac*{\frac{u}{u-2v}b_{-1}+2\frac{v-u}{u-2v}a_{-1}}\ket{\tau},z}}{z-w}\nonumber\\
  &\qquad\qquad
  \frac{-\brac*{\frac{u}{v}-1}I^\pm_{r,s}(v_{1,2},z)I_\mu\brac*{\ket{\tau},z}}{\brac*{z-w}^2}.
\end{align}
Then from the definition \eqref{eq:Q2def} of $\mathcal{Q}^{[2]}(z)$, we have
\begin{multline}
  \comm{Y(\psi(f),w)}{\mathcal{Q}^{[2]}}=Y(\psi(f),w)\mathcal{Q}^{[2]}(z)\sim
  \int_{\Gamma_{0,1}}({\partial}_{z}R^-(z,w))\mathcal{Q}^{+}(zy)z{\rm d}y\\
  +\int_{\Gamma_{0,1}}\mathcal{Q}^{-}(z)({\partial}_{\tilde{z}}R^+(\tilde{z},w))|_{\tilde{z}=zy}z{\rm d}y.
  \label{eq:com1}
\end{multline}
Note that
${\rm d}z\land (z{\rm d}y)={\rm d}z\land {\rm d}(zy)$. 
Then, from (\ref{eq:com1}), we have
\begin{align}
\comm{\mathcal{Q}^{[2]}}{Y(\psi(f),w)}&=\res_{z=w}Y(\psi(f),w)\mathcal{Q}^{[2]}(z)\nonumber\\
&= \int_{C_{z=w}}\int_{\Gamma_{0,1}} \dd
\Bigl(R^-(z,w)\mathcal{Q}^{+}(\tilde{z}){\rm d}\tilde{z}
+\mathcal{Q}^{-}(z)R^+(\tilde{z},w)\dd z\Bigr),
\label{eq:com2}
\end{align}
where $C_{z=w}$ is a counter clock-wise circle (scaled by \((2\pi\ii)^{-1}\)
about $w$, and $\dd $ is the
exterior derivative with respect to the variables $(z,\tilde{z})$. Since
$C_{z=w}\times\Gamma_{0,1}$ is a twisted cycle and the integrand is exact \eqref{eq:com2} must
vanish. Thus $\mathcal{Q}^{[2]}$ commutes with the \(f\)-generator. For the
the \(e\) and \(h\)-generators, note that the \ope{s} \eqref{eq:screeningope}
are regular and hence the integral calculations analogous to those above
vanish trivially. Thus $\mathcal{Q}^{[2]}$ commutes with the action of \(\mmsl{u}{v}\).
\end{proof}

Consider the operator
\begin{align}
  \Delta_\epsilon(-;z):\mmsl{u}{v}&\to
  \hom_{\CC}(\HH^0,\HH^1)\powser{z^{\pm1},\epsilon^{\pm1}}\nonumber\\
  g&\mapsto \frac{1}{\epsilon}\oint_{z}\mathcal{Q}^{[2]}(w)
  \bee^{-\epsilon}Y(g,z) \bee^{\epsilon}\dd w,
  \label{eq:wtshift}
\end{align}
where $\bee^{\pm \epsilon}$ is the Heisenberg weight
shifting operator characterised by 
\begin{align}
  [a_0,\bee^{\pm \epsilon}]&=\pm\frac{\epsilon}{2}\bee^{\pm \epsilon},
  \qquad[b_0,\bee^{\pm \epsilon}]=\pm\frac{\epsilon}{2}\bee^{\pm \epsilon},\nonumber\\
  [a_n,\bee^{\pm \epsilon}]&=[b_n,\bee^{\pm \epsilon}]=0,\ n\neq 0,\qquad
                             \bee^{\pm \epsilon}\ket{\mu}=\ket{\mu\pm \epsilon\frac{ v}{u-2v}\brac*{a-b}},
  \label{eq:defrels}
\end{align}
and \(\epsilon\) is a deformation parameter.
Note that the integral in \eqref{eq:wtshift} can also be expressed as the
commutator
\begin{align}
  \Delta_\epsilon(g;z)=\frac{1}{\epsilon} \comm{\mathcal{Q}^{[2]}}
  {\bee^{-\epsilon}Y(g,z) \bee^{\epsilon}}.
\end{align}

\begin{prop}
  Consider the operator \(\Delta_\epsilon(-;z)\) defined above.
  \begin{enumerate}
  \item For all \(g\in \mmsl{u}{v}\),
    \(\Delta_\epsilon(g;z)\in
    \hom_{\CC}(\HH^0,\HH^1)\powser{z^{\pm1},\epsilon}\), that is, no negative
    powers of \(\epsilon\) appear.
    \label{itm:Dwelldef}
  \item For all \(g\in \mmsl{u}{v}\), set
    \(\Delta(g;z)=\left.\Delta_\epsilon(g;z)\right|_{\epsilon=0}\). Then
    on the generators of \(\mmsl{u}{v}\) and the conformal vector we have
    \begin{align}
      \Delta(e;z)&=0,\qquad
      \Delta(h;z)=0,
      \nonumber\\
      \Delta(f;z)&= -\frac{1}{z}\comm*{\mathcal{Q}^{[2]}}{Y\brac*{b_{-1}\ket{\frac{-2v}{u-2v}},z}}
      +\frac{1}{z^2}\frac{u-v}{2v}\comm*{\mathcal{Q}^{[2]}}{Y\brac*{\ket{\frac{-2v}{u-2v}},z}},\nonumber\\
      \Delta(\omega_{\frac{u-2v}{v}};z)&=\frac{1}{z}\frac{v}{u-2v}\comm{\mathcal{Q}^{[2]}}{Y\brac*{\brac*{a_{-1}-b_{-1}}\ket{0},z}}.
    \end{align}
    \label{itm:Dformulae}
  \end{enumerate}
  \label{thm:Dprops}
\end{prop}
\begin{proof}
  For any $v\in \mathbb{H}^0$, $u\in (\mathbb{H}^1)'$ and $g\in \mmsl{u}{v}$,
  the matrix element
  \begin{align}
    \res_{z=1}\pair{u}{\mathcal{Q}^{[2]}(z)\mathbf{e}^{-\epsilon}Y(g,1)\mathbf{e}^{\epsilon}v} \in
    \CC\powser*{\epsilon}
    \label{eq:matelser}
  \end{align}
  is a power series in non-negative powers of \(\epsilon\), because the
  Heisenberg weight shift operators can only introduce non-negative powers by
  the relations \eqref{eq:defrels}. Next, 
  recall from the proof of \cref{thm:comq} that
  $\res_{z=w}\mathcal{Q}^{[2]}(z)Y(g,w)=0$ $(g\in \mmsl{u}{v})$, hence the power
  series \eqref{eq:matelser} has vanishing constant term, that is,
  \begin{equation}
\res_{z=1}\pair{ u}{\mathcal{Q}^{[2]}(z)\mathbf{e}^{-\epsilon}Y(g,1)\mathbf{e}^{\epsilon}v}\in \epsilon\CC\powser*{\epsilon}.
\end{equation}
Thus we have \(\Delta_\epsilon(g;z)\in
\hom_{\CC}(\HH^0,\HH^1)\powser{z^{\pm1},\epsilon}\) and hence \(\Delta(g;z)\)
is well-defined.

We evaluate \(\Delta(-;z)\) on generators and the conformal vector by direct computation.
First note that
\begin{align}
  \bee^{-\epsilon}Y\brac*{a_{-1}\ket{0},z}\bee^{\epsilon}&=Y\brac*{a_{-1}\ket{0},z}+\frac{\epsilon/2}{z},\qquad
  \bee^{-\epsilon}Y\brac*{b_{-1}\ket{0},z}\bee^{\epsilon}=Y\brac*{b_{-1}\ket{0},z}+\frac{\epsilon/2}{z},\nonumber\\
  \bee^{-\epsilon}Y\brac*{\ket{\frac{\pm 2v}{u-2v}(a-b)},z}\bee^{\epsilon}&=Y\brac*{\ket{\frac{\pm 2v}{u-2v}(a-b)},z}.
\end{align}
For images under the free field realisation \(\psi\) of \cref{itm:ffrformulae}
of \cref{thm:sl2ffr} we therefore have
\begin{align}
  \bee^{-\epsilon}Y\brac*{\psi(e_{-1}\vac),z}\bee^{\epsilon}&=Y\brac*{\psi(e_{-1}\vac),z},\qquad
  \bee^{-\epsilon}Y\brac*{\psi(h_{-1}\vac),z}\bee^{\epsilon}=Y\brac*{\psi(h_{-1}\vac),z}+\frac{\epsilon}{z},\nonumber\\
  \bee^{-\epsilon}Y\brac*{\psi(f_{-1}\vac),z}\bee^{\epsilon}&=Y\brac*{\psi(f_{-1}\vac),z}-\frac{\epsilon}{z}Y\brac*{b_{-1}\ket{\frac{-2v}{u-2v}},z}
                                                              -\frac{\epsilon^2/4}{z^2}Y\brac*{\ket{\frac{-2v}{u-2v}},z}\nonumber\\
                                                            &\quad+\frac{u-v}{2v}\frac{\epsilon}{z^2}Y\brac*{\ket{\frac{-2v}{u-2v}},z},\nonumber\\
  \bee^{-\epsilon}Y\brac*{\psi(\omega_{\frac{u-2v}{v}}),z}\bee^{\epsilon}&=Y\brac*{\psi(\omega_{\frac{u-2v}{v}}),z}+\frac{\epsilon}{z}\frac{v}{u-2v}Y\brac*{\brac*{a_{-1}-b_{-1}}\ket{0},z}+\frac{\epsilon/2}{z^2}.
\end{align}
The evaluations of \(\Delta(-,z)\) are therefore just the commutator of
\(\mathcal{Q}^{[2]}\) with the terms linear in \(\epsilon\) above.
\end{proof}

\begin{prop}
  Consider the \(\mmV{u}{v}\cten \lat{\frac{u}{v}}\) module \(\HH\) and
  denote the \(\mmV{u}{v}\cten \lat{\frac{u}{v}}\)-action by \(Y_{\HH}\) (which is
  also an \(\mmsl{u}{v}\) action by restriction). For \(d\in \CC^\times\) define the linear operator
  \(\widetilde{Y}_d(-,z):\mmsl{u}{v}\to \hom_{\mathbb{C}}(\HH,\HH)\powser{z^{\pm1}}\)
    by
    \begin{align}
      {\left.\widetilde{Y}_{d}(-,z)\right|}_{\HH^0}=Y_{\HH}(-,z)+d\cdot\Delta(-;z),\qquad
      {\left.\widetilde{Y}_{d}(-,z)\right|}_{\HH^1}=Y_{\HH}(-,z).
    \end{align}
    Then \(\widetilde{Y}_d\) is an action of \(\mmsl{u}{v}\) on
    \(\HH\). Further, \((\HH,\widetilde{Y}_d)\) is isomorphic to the projective
    cover \(\sigma^{-1}P_{u-1,v-1}\) of \(L_1\cong \mmsl{u}{v}\).
    \label{thm:projmod}
  \end{prop}
  The construction above is inspired by the construction of logarithmic modules in \cite[Sec 2]{AM09}, but with some differences in the choice of integration cycle which
    will turn out to be more convenient for the intertwining operators to be considered below.
  Note that the extension group \(\ext^1\brac*{\sigma
    \slEr{u-1}{v-1}{-},\sigma^{-1}\slEr{1}{1}{-}}\) is \(1\)-dimensional. The
  parameter \(d\) above hence reflects this one degree of freedom.
  \begin{proof}
    To show that \(\widetilde{Y}_d\) defines an \(\mmsl{u}{v}\)-action, we need
    to show that \(\widetilde{Y}_d(\vac,z)=\id_{\HH}\), that
    \(\widetilde{Y}_d(w,z)h\) has at most a finite order pole for any
    \(w\in\mmsl{u}{v}\) and any \(h\in\HH\), and that \(\widetilde{Y}_d\)
    satisfies the Jacobi identity. The first two conditions hold by
    construction and so we only need to consider the Jacobi identity. Note
    that the factor \(d\) can be absorbed by rescaling the entire space
    \(\HH^1\) by \(d\), hence without loss of generality it is sufficient to
    verify the Jacobi identity for \(d=1\). We thus set
    \(\widetilde{Y}=\widetilde{Y}_1\). 
    Let \(f(z)\in\CC\sqbrac*{z^{\pm1},(z-w)^{-1}}\), \(\phi\in \HH\),
    \(\psi\in \HH'\) and \(a,b\in \mmsl{u}{v}\), then we need to show that the
    sum of integrals
    \begin{align}
      &\oint_{0,w}f(z)\pair{\psi}{\widetilde{Y}(a,z)\widetilde{Y}(b,w)\phi}\dd z
      - \oint_{w}f(z)\pair{\psi}{\widetilde{Y}\brac*{Y(a,z-w)b,w}\phi}\dd z\nonumber\\
      &\qquad - \oint_{0}f(z)\pair{\psi}{\widetilde{Y}(b,w)\widetilde{Y}(a,z)\phi}\dd z
    \end{align}
    vanishes.
    Note that for \(\phi\in \HH^1, \psi\in \HH^{1\prime}\) or \(\phi\in \HH^0,
    \psi\in \HH^{0\prime}\) the above expression reduces to the Jacobi
    identity for the action \(Y_\HH\) and hence vanishes and that it is
    trivially zero for \(\phi\in \HH^1, \psi\in \HH^{0\prime}\). So we only
    need to consider the case \(\phi\in \HH^0, \psi\in \HH^{1\prime}\), where
    the above expression specialises to
    \begin{align}
      &\oint_{0,w}f(z)\pair{\psi}{\brac*{\Delta(a,z)Y_\HH(b,w)+Y_\HH(a,z)\Delta(b,w)}\phi}\dd z\nonumber\\
      &\qquad- \oint_{w}f(z)\pair{\psi}{\Delta\brac*{Y(a,z-w)b,w}\phi}\dd z\nonumber\\
      &\qquad -
        \oint_{0}f(z)\pair{\psi}{\brac*{Y_\HH(b,w)\Delta(a,z)+\Delta(b,w)Y_\HH(a,z)}\phi}\dd
        z.
        \label{eq:Jacspec}
    \end{align}
    The first and third integrands above can be simplified by noting
    \begin{align}
      \Delta(a,z)Y_\HH(b,w)+Y_\HH(a,z)\Delta(b,w)=\left.\frac{1}{\epsilon}\oint_{z,w}\int_{\Gamma_{0,1}}\hspace{-4mm}\mathcal{Q}^-(x)\mathcal{Q}^+(xy)\bee^{-\epsilon}Y_\HH(a,z)Y_\HH(b,w)\bee^{\epsilon}\dd
      x\dd y\right|_{\epsilon=0}\hspace{-2mm}.
    \end{align}
    Let \(\set{\gamma_i}\) be a basis of \(\HH^1\) with dual basis
    \(\set{\gamma^i}\), then \eqref{eq:Jacspec} is equal to
    \begin{align}
      \oint_{z,w}\int_{\Gamma_{0,1}}\sum_{i}\pair{\psi}{\mathcal{Q}^-(x)\mathcal{Q}^+(xy)\gamma_i}
      &\left(\oint_{0,w}f(z)\pair{\gamma^i}{Y_\HH(a,z)Y_\HH(b,w)\phi}\dd
        z\right.\nonumber\\
      &\quad - \oint_{w}f(z)\pair{\gamma^i}{Y_{\HH}\brac*{Y(a,z-w)b,w}\phi}\dd
        z\nonumber\\
      &\quad  \left.- \oint_{0}f(z)\pair{\gamma^i}{Y_\HH(b,w) Y_\HH(a,z)\phi}\dd z\right),
    \end{align}
    which vanishes, because \(Y_\HH\) satisfies the Jacobi identity.
    Further, from the formulae in \cref{itm:Dformulae} of \cref{thm:Dprops} we see that the
    \(h_0\)-grading remains unchanged, while the \(L_0\) operator acquires a
    nilpotent part but its generalised eigenvalues do not change, hence
    \((\HH,\widetilde{Y})\) is a weight module.
    To  conclude that \((\HH,\widetilde{Y})\) is isomorphic to
    \(\sigma^{-1}P_{u-1,v-1}\) note that \(\HH^1\) is a submodule isomorphic to
    \(\sigma E^-_{u-1,v-1}\), that the quotient of \((\HH,\widetilde{Y})\) by
    \(\HH^1\) is isomorphic to \(\HH^0\) and thus \(\sigma^{-1}E^-_{1,1}\) and
    that \((\HH,\widetilde{Y})\) is indecomposable due to \(\Delta\) mapping
    between \(\HH^0\) and \(\HH^1\). Hence \((\HH,\widetilde{Y})\) satisfies
    the characterising non-split exact sequence \eqref{eq:projmchar} of \(\sigma^{-1}P_{u-1,v-1}\).
    
  \end{proof}

\subsection{Constructing logarithmic intertwining operators}
\label{sec:logintops}
This section assumes some
familiarity with twisted cycles and twisted De Rham cohomology. We refer
readers unfamiliar with these notions to Aomoto and Kita's book on
hypergeometric functions \cite{hgeom} for an comprehensive account or to
\cite[Sec 3.2]{TW} for a short summary.

Let $w\in\mathbb{R}_{>0}$, and consider the complex manifold
$Z_{w}=\set{(z_1,z_2)\in \mathbb{C}^2\ \st\ z_1\neq z_2,\ z_1,z_2\neq w,\ z_1,z_2\neq 0}$. 
For $\alpha,\beta,\gamma\in \mathbb{C}$, let 
\begin{equation}
\mathcal{U}_w(\alpha,\beta,\gamma,z)=\prod_{i=1}^2z^\alpha_i\prod_{j=1}^2(z_j-w)^\beta(z_1-z_2)^{\gamma}(z_2-z_1)^{\gamma}.
\end{equation}
be a multivalued function on $Z_{w}$.
The logarithmic derivative of $\mathcal{U}_w(\alpha,\beta,\gamma,z)$
\begin{equation}
\omega_w(\alpha,\beta,\gamma,z)={\rm d}\ {\rm
  log}\mathcal{U}_{w}(\alpha,\beta,\gamma,z)=\sum_{i=1,2}\Bigl(\frac{\alpha}{z_i}+\frac{\beta}{z_i-w}\Bigr){\rm
  d}z_i+2\gamma\frac{{\rm d}z_1-{\rm d}z_2}{z_1-z_2}
\label{eq:dU}
\end{equation}
defines the twisted differential $\nabla_{\omega_w}={\rm d}+\omega_w(\alpha,\beta,\gamma,z)\wedge$.
Let $\mathcal{L}_w(\alpha,\beta,\gamma)$ be the local system defined by the local solutions of $\nabla_{\omega_w}g(z_1,z_2)=0$, and let $\mathcal{L}^\vee_w(\alpha,\beta,\gamma)={\rm Hom}_{\mathbb{C}}(\mathcal{L}_w(\alpha,\beta,\gamma),\mathbb{C})$ be the dual local system.
The twisted homology groups with coefficients in $\mathcal{L}^\vee_w(\alpha,\beta,\gamma)$ are denoted by $H_p(Z_{w},\mathcal{L}^\vee_w(\alpha,\beta,\gamma))$ and the twisted cohomology groups by
\begin{equation} 
H^p(Z_{w},\mathcal{L}_w(\alpha,\beta,\gamma))={\rm Hom}_{\mathbb{C}}(H_p(Z_{w},\mathcal{L}^\vee_w(\alpha,\beta,\gamma)),\mathbb{C}).
\end{equation}
It is known that the twisted cohomology groups are isomorphic to the twisted
de Rham cohomology groups (see \cite[Sec 2]{hgeom})
\begin{equation}
H^p(Z_{w},\mathcal{L}_w(\alpha,\beta,\gamma))\simeq H^p(Z_{w},\nabla_{\omega_w}).
\end{equation}
By permuting the variables \(z_1,z_2\) there is a natural action of the
symmetric group. We denote the subspace that is invariant this action
by \(H^2(Z_{w},\nabla_{\omega_w})^\mathfrak{S}\)
and the subspace on which the transposition \((1,2)\) acts by \(-1\) by
\(H^2(Z_{w},\nabla_{\omega_w})^\mathfrak{A}\).

\begin{lemma}
  Consider the hyperplane arrangements
  \begin{align}
     \mathcal{H}&=\bigcup_{i=1}^2\set{(\alpha,\beta,\gamma)\in
   \CC^3\st i(\alpha+(i-1)\gamma)\in \mathbb{Z}}\nonumber\\
 &\quad\cup\bigcup_{j=1}^2\set{(\alpha,\beta,\gamma)\in
 \CC^3\st j(\beta+(j-1)\gamma)\in \mathbb{Z}}\cup \set{(\alpha,\beta,\gamma)\in
 \CC^3\st 2\gamma\in \mathbb{Z}},\nonumber\\
    \mathcal{G}&=\mathcal{H}\cup\set{(\alpha,\beta,\gamma)\in \CC^3\st
                 \alpha+\beta+2\gamma\notin \ZZ},
                 \label{eq:hplanes}
  \end{align}
  and let $(\alpha,\beta,\gamma)\in \mathbb{C}^3\setminus \mathcal{G}$.
  The skew symmetric \(2\)-form \(\mathcal{U}_w(\alpha,\beta,\gamma,z)\)
  represents a non-zero cohomology class within the cohomology group
  \(H^2(Z_{w},\nabla_{\omega_w})^{\mathfrak{A}}\)  and \(\dim
   H^2(Z_{w},\nabla_{\omega_w})^{\mathfrak{A}}=1\).
  Further, for any symmetric Laurent polynomial in
  $F(z_1,z_2)\in \mathbb{C}[z^{\pm 1}_1,z^{\pm
    1}_2,(z_1-w)^{-1},(z_2-w)^{-1}]^{\mathfrak{S}}\sqbrac*{\brac*{z_1-z_2}^{-2}}$ there exist a rational
  function $c(\alpha,\beta,\gamma)\in \mathbb{C}(\alpha,\beta,\gamma)$ such
  that
  \begin{equation}
    \mathcal{U}_w(\alpha,\beta,\gamma,z_1,z_2)F(z_1,z_2)\dd z_1\wedge\dd z_2=c(\alpha,\beta,\gamma)\mathcal{U}_w(\alpha,\beta,\gamma,z_1,z_2)\dd z_1\wedge\dd z_2
    +\dd(\cdots),
    \label{eq:lem-one}
  \end{equation}
  where \(\dd(\cdots)\) denotes a total derivative, that is, the two form
  \(\mathcal{U}_w(\alpha,\beta,\gamma,z_1,z_2)F(z_1,z_2)\dd z_1\wedge{\rm
    d}z_2\) is cohomologous to a scalar multiple of 
  \(\mathcal{U}_w(\alpha,\beta,\gamma,z_1,z_2)\dd z_1\wedge{\rm d}z_2\).
  Finally, the rational function $c(\alpha,\beta,\gamma)$ is holomorphic on
  $\mathbb{C}^3\setminus\mathcal{H}$,
  that is, all of its singularities lie in
  \(\mathcal{H}\).
  \label{thm:one-dim}
\end{lemma}

By \cite[Sec 2]{hgeom} \cref{thm:one-dim} is known for generic values of
\(\alpha,\beta,\gamma\), however, as we shall see below, here we are
primarily interested in the case \(\alpha+\beta+\gamma=-1\), which is not
generic.
\begin{proof}
  By \cite[Eq. (2.53) and Thm 2.5]{hgeom} local systems cohomology can be identified with
  twisted de Rham cohomology of rational forms.  Hence we restrict our attention to rational forms.
  Fix $(\alpha,\beta,\gamma)\in \mathbb{C}^3\setminus\mathcal{G}$ and denote
  $\mathcal{U}_w=\mathcal{U}_w(\alpha,\beta,\gamma,z_1,z_2)$. 
  Consider the \(1\)-form 
  \(\mathcal{U}_w((z_1-w)(z^{k}_1z^l_2)\dd z_2-(z_2-w)(z^{k}_2z^l_1)\dd z_1)\), 
  $k,l\in \mathbb{Z}$, \(k\ge l\) and compute its exterior derivative.
  \begin{align}
    \dd \brac*{\mathcal{U}_w((z_1-w)(z^{k}_1z^l_2)\dd z_2-(z_2-w)(z^{k}_2z^l_1)\dd z_1)}\qquad&\nonumber\\
    =(\alpha+\beta+2\gamma+k+1)(z^k_1z^l_2+z^k_2z^l_1)\mathcal{U}_w\dd z_1\wedge\dd z_2&\nonumber\\
    \qquad-w(\alpha+2\gamma+k)(z^{k-1}_1z^l_2+z^{k-1}_2z^l_1)\mathcal{U}_w\dd z_1\wedge\dd z_2&\nonumber\\
    \qquad+2\gamma \brac*{\sum_{i=1}^{k-l-1}z^{i+l}_1z^{k-i}_2
      -w\sum_{j=1}^{k-l-2}z^{j+l}_1z^{k-1-j}_2}\mathcal{U}_w\dd
    z_1\wedge\dd z_2, &&k\ge l+2.\nonumber\\
    \dd \brac*{\mathcal{U}_w((z_1-w)(z_1^{l+1}z^l_2)\dd z_2-(z_2-w)(z_2^{l+1}z^l_1)\dd z_1)}\qquad&\nonumber\\
    =(\alpha+\beta+2\gamma+l+2)(z^{l+1}_1z^l_2+z^{l+1}_2z^l_1)\mathcal{U}_w\dd z_1\wedge\dd z_2&\nonumber\\
    \qquad-w2(\alpha+\gamma+l+1)z_!^{l}z_2^l\mathcal{U}_w\dd z_1\wedge\dd z_2,&
    &k = l+1.\nonumber\\
    \dd \brac*{\mathcal{U}_w((z_1-w)(z^{l}_1z^l_2)\dd z_2-(z_2-w)(z^{l}_2z^l_1)\dd z_1)}&\nonumber\\
    =2(\alpha+\beta+\gamma+l+1)z^l_1z^l_2\mathcal{U}_w\dd z_1\wedge\dd z_2&\nonumber\\
    \qquad-w(\alpha+l)(z^{l-1}_1z^l_2+z^{l-1}_2z^l_1)\mathcal{U}_w\dd z_1\wedge\dd z_2,&
    &k=1.
    \label{eq:rec1}
  \end{align}
  Note that with the exception of the coefficient on the second to last line, the coefficients in the above identities cannot vanish because of the
  assumptions on \(\alpha,\beta,\gamma\).
  The above exterior derivative implies that a skew symmetric \(2\)-form of the
  form
  \begin{equation}
    \brac*{z^n_1z^m_2+z^n_2z^m_1}\mathcal{U}_w\dd
    z_1\wedge\dd z_2, \qquad m,n\in \ZZ
    \label{eq:monomial2form1}
  \end{equation}
  is cohomologous to linear combinations of such two forms,
  where the difference \(|m-n|\) of the exponents \(m,n\) has been reduced by
  \(1\) or more. This
  procedure can be iterated until the exponents are equal. Once the exponents
  are equal we can apply the relations generated from \eqref{eq:rec1} twice to
  shift both exponents by \(1\), that is,
  \begin{multline}
    \brac*{\alpha+\beta+\gamma+n+2}z_1^{n+1}z_2^{n+1}\mathcal{U}_w\dd z_1\wedge \dd z_2\\=
    w^2\frac{\brac*{\alpha+\gamma+n+1}\brac*{\alpha+n+1}}{\alpha+\beta+2\gamma+n+2}\mathcal{U}_w z_1^n
    z_2^n\dd z_1\wedge \dd z_2+\dd\brac*{\cdots}.
  \end{multline}
  For \(\alpha+\beta+\gamma\notin\ZZ\) this allows one to conclude that any \(2\)-form of the form \eqref{eq:monomial2form1} is cohomologous to a
  scalar multiple of \(\mathcal{U}_w\). If
  \(\alpha+\beta+\gamma=p\in \ZZ\), then \(z_1^{n+1}z_2^{n+1}\mathcal{U}_w\dd
  z_1\wedge \dd z_2\) is cohomologous to a
  scalar multiple of \(\mathcal{U}_w\dd z_1\wedge \dd z_2\) for \(n\ge -p-1\)
  and to \(0\) for \(n\le -p-2\) and so the same conclusion holds.
  Further since
  \(\mathcal{U}_w\) is symmetric under interchanging
  \(\alpha,\beta\) and factors of \(z_i\) with \(z_i-w\) analogous arguments to
  those above imply that any \(2\)-form of the form
  \begin{equation}
    \brac*{\brac*{z_1-w}^n\brac*{z_2-w}^m+\brac*{z_2-w}^n\brac*{z_1-w}^m}\mathcal{U}_w\dd
    z_1\wedge\dd z_2, \qquad m,n\in \ZZ,
    \label{eq:monomial2form2}
  \end{equation}
  is also cohomologous to a scalar multiple of
  \(\mathcal{U}_w\).
  Next we consider more general skew symmetric \(2\)-forms made up of linear combinations of
  \(2\)-forms of the form
  \begin{equation}
    \brac*{\frac{z_1^kz_2^l}{\brac*{z_1-w}^n\brac*{z_2-w}^m}+\frac{z_2^kz_1^l}{\brac*{z_2-w}^n\brac*{z_1-w}^m}}\mathcal{U}_w\dd z_1\wedge\dd z_2, \qquad m,n\in \ZZ_{\ge0},\quad k,l\in \ZZ,
    \label{eq:monomial2form3}
  \end{equation}
  however, using partial fractal decomposition these can be reduced linear
  combinations of \eqref{eq:monomial2form1} and of
  \begin{align}
    &\brac*{\frac{z_1^k}{\brac*{z_2-w}^m}+\frac{z_2^k}{\brac*{z_1-w}^m}}\mathcal{U}_w\dd
    z_1\wedge\dd z_2,
    \nonumber\\
    &\brac*{\frac{1}{z_1^k\brac*{z_2-w}^m}+\frac{1}{z_2^k\brac*{z_1-w}^m}}\mathcal{U}_w\dd
    z_1\wedge\dd z_2,
    \qquad k,m\in \ZZ_{\ge0}.
    \label{eq:monomial2form4}
  \end{align}
  So we consider these two new cases. The first case of
  \eqref{eq:monomial2form4} reduced to the case \eqref{eq:monomial2form2}
  after noting
  \begin{align}
    \frac{z^k_i}{(z_j-w)^m}=\frac{(z_i-w+w)^k}{(z_j-w)^m}=\frac{1}{(z_j-w)^m}\sum_{p=0}^k\binom{k}{p}(z_i-w)^pw^{k-p}.
  \end{align}
  The final case can be reduced to the case \eqref{eq:monomial2form1} after
  noting the relation
  \begin{align}
    &\brac*{\frac{1}{z_1^k\brac*{z_2-w}^m}+\frac{1}{z_2^k\brac*{z_1-w}^m}}\mathcal{U}_w(\alpha,\beta,\gamma,z)\dd
    z_1\wedge\dd z_2\nonumber\\
    &=
    \brac*{\frac{\brac*{z_1-w}^m}{z_1^{k}}+\frac{\brac*{z_2-w}^m}{z_2^{k}}}\mathcal{U}_w(\alpha,\beta-m,\gamma,z)\dd
    z_1\wedge\dd z_2\nonumber\\
    &=\brac*{\sum_{j\ge0}\binom{m}{j}z_1^{j-k}(-w)^{m-j}+\sum_{j\ge0}\binom{m}{j}z_2^{j-k}(-w)^{m-j}}\mathcal{U}_w(\alpha,\beta-m,\gamma,z)\dd
    z_1\wedge\dd z_2
  \end{align}
  and that \((\alpha,\beta,\gamma)\) satisfies the assumptions of the lemma if
  and only if \((\alpha,\beta+m,\gamma)\) does.
  Finally, the most general form that a skew symmetric \(2\)-form can take is
  a linear combination of 
  \(2\)-forms of the form
  \begin{equation}
    \frac{F(z_1,z_2)}{\brac*{z_1-z_2}^{2n}}\mathcal{U}_w\dd z_1\wedge\dd z_2,
    \quad n\in \ZZ_{\ge0},\ 
    F(z_1,z_2)\in \CC\sqbrac*{z_1^\pm,z_2^\pm,\brac*{z_1-w}^{-1},\brac*{z_2-w}^{-1}}.
    \label{eq:monomial2form5}
  \end{equation}
  The factor of \(\brac*{z_1-z_2}^{2n}\) can be absorbed in to the multivalued
  function \(\mathcal{U}_w\) so that, from the reasoning above, we can conclude
  the existence of a rational functions \(c(\alpha,\beta,\gamma)\) satisfying
  \begin{align}
    \frac{F(z_1,z_2)}{\brac*{z_1-z_2}^{2n}}\mathcal{U}_w\dd z_1\wedge\dd z_2
    &=F(z_1,z_2)\mathcal{U}_w(\alpha,\beta,\gamma-2n,z)\dd z_1\wedge\dd z_2\nonumber\\
    &=c(\alpha,\beta,\gamma)\mathcal{U}_w(\alpha,\beta,\gamma-2n,z)\dd z_1\wedge\dd z_2+\dd\brac*{\cdots}.
  \end{align}
  This shift in the \(\gamma\)-argument can be converted to a shift in the
  \(\alpha\) and \(\beta\) arguments by repeatedly using the identity
  \begin{align}
    \dd((z_1-z_2)\mathcal{U}_w(\alpha,\beta,\gamma-2,z_1,z_2)\dd z_2+(z_1-z_2)\mathcal{U}_w(\alpha,\beta,\gamma-2,z_1,z_2)\dd z_1)&\nonumber\\
    =2(\gamma-1)\mathcal{U}_w(\alpha,\beta,\gamma-2,z_1,z_2)\dd z_1\wedge\dd z_2
    -{\alpha}z_1^{-1}z_2^{-1}\mathcal{U}_w(\alpha,\beta,\gamma,z_1,z_2)\dd z_1\wedge\dd z_2& \nonumber\\
    -{\beta}\brac*{z_1-w}^{-1}\brac*{z_2-w}^{-1}\mathcal{U}_w(\alpha,\beta,\gamma,z_1,z_2)\dd
    z_1\wedge\dd z_2&.
    \label{eq:gammashift}
  \end{align}
  which can then be simplified using all the previously discussed
  identities. Hence any 2-form of the form \eqref{eq:monomial2form5} is
  cohomologous to a scalar multiple of \(\mathcal{U}_w\dd z_1\wedge \dd
  z_2\). Thus \(\dim H^2(Z_{w},\nabla_{\omega_w})^{\mathfrak{A}}=1\).
  
  The rationality of the functions \(c(\alpha,\beta,\gamma)\) is immediate
  from the identity \eqref{eq:rec1}, its corresponding version for the case
  \eqref{eq:monomial2form2} and \eqref{eq:gammashift}, because the coefficients that appear in these
  identities have zeros only in \(\mathcal{G}\). The rationality of \(c(\alpha,\beta,\gamma)\) is also a special case of \cite[Cor 1.1.1]{SusInt}.
\end{proof}

For real numbers $a>b$, let $\Delta_{a,b}=\{a>z_1>z_2>b\}$.
Given $f(z_1,z_2)\in \mathbb{C}[z^{\pm 1}_1,z^{\pm 1}_2,(z_1-w)^{-1},(z_2-w)^{-1}]$, we see that for appropriate $(\alpha,\beta,\gamma)$, the integral 
\begin{equation}
\mathcal{I}_w[f](\alpha,\beta,\gamma)=\int_{\Delta_{w,0}}\mathcal{U}_w(\alpha,\beta,\gamma,z)f(z_1,z_2){\rm d}z_1{\rm d}z_2
\end{equation}
converges. 
In the case $f=1$, the explicit formula is given by
\begin{equation}
  \mathcal{I}_w[1](\alpha,\beta,\gamma)=w^{2(\alpha+\beta+1)}\prod_{i=1,2}\frac{\Gamma(1+i\gamma)\Gamma(1+\alpha+(i-1)\gamma)\Gamma(1+\beta+(i-1)\gamma)}{\Gamma(1+\gamma)\Gamma(2+\alpha+\beta+i\gamma)}.
  \label{eq:integrationformula}
\end{equation}  
In \cite[Cor 1.1.1]{SusSel}, it was shown that
$\mathcal{I}_w[f](\alpha,\beta,\gamma)$ admits an analytic continuation to
$\mathbb{C}^3\setminus
\mathcal{H}$, where
\(\mathcal{H}\) is the hyperplane arrangement \eqref{eq:hplanes} in \cref{thm:one-dim}.
In \cite[Def 3.6]{TK}, a non-zero twisted cycle $[{\Delta}_{w,0}(\alpha,\beta,\gamma)]\in H_2(Z_{w},\mathcal{L}^\vee_w(\alpha,\beta,\gamma))$ was constructed from the simplex ${\Delta}_{w,0}$. This cycle $[{\Delta}_{w,0}(\alpha,\beta,\gamma)]$ satisfies
\begin{equation}
  \int_{[\Delta_{w,0}(\alpha,\beta,\gamma)]}\mathcal{U}_w(\alpha,\beta,\gamma,z)f(z_1,z_2){\rm
    d}z_1{\rm d}z_2=\mathcal{I}_w[f](\alpha,\beta,\gamma)
  \label{eq:integralrelation}
\end{equation}
for $f(z_1,z_2)\in \mathbb{C}[z^{\pm 1}_1,z^{\pm
  1}_2,(z_1-w)^{-1},(z_2-w)^{-1}]$. We shall abbreviate this cycle as
\(\sqbrac*{\Delta_{w,0}}\) when the parameters \(\alpha,\beta,\gamma\) are
clear from context.

Let $\lambda,\epsilon\in \RR$ such that $\lambda\neq \pm
\lambda_{u-1,v-1}\pmod{2} = \pm
(\frac{u}{v}-2)\pmod{2}$ and $\lambda+\epsilon\neq \pm
\lambda_{u-1,v-1}\pmod{2}$. We specialise the parameters of
\(\mathcal{U}_w\brac*{\alpha,\beta,\gamma,z}\) to
\(\alpha=\frac{u}{2v}-1-\frac{\lambda+\epsilon}{2},\
\beta=\frac{u}{2v}-1+\frac{\lambda+\epsilon}{2},\ \gamma=1-\frac{u}{v}\) and
define the shorthand
\begin{align}
  \mathcal{U}_w^\epsilon\brac*{z}&=\mathcal{U}_w\brac*{\tfrac{u}{2v}-1-\tfrac{\lambda+\epsilon}{2},\tfrac{u}{2v}-1+\tfrac{\lambda+\epsilon}{2},1-\frac{u}{v},z},
  \nonumber\\
  \mathcal{L}_{w,\epsilon{}}&=\mathcal{L}_w\brac*{\tfrac{u-2v-(\lambda
    +\epsilon)v}{2v},\tfrac{u-2v+(\lambda
    +\epsilon)v}{2v},-\tfrac{u}{v}+1},\nonumber\\
  \mathcal{L}^\vee_{w,\epsilon{}}&=\mathcal{L}^\vee_w\brac*{\tfrac{u-2v-(\lambda
    +\epsilon)v}{2v},\tfrac{u-2v+(\lambda +\epsilon)v}{2v},-\tfrac{u}{v}+1}. 
  \label{eq:U-c}
\end{align}
Now consider a non-zero 
$\lat{\frac{u}{v}}$-intertwining operator
\begin{equation}
  I_{{\lambda+\epsilon,-\lambda}}\in \binom{\mathbb{F}_{\epsilon\frac{v}{u-2v}\brac*{a-b}-2b}}{\mathbb{F}_{\brac*{\lambda+\epsilon}\frac{v}{u-2v}\brac*{a-b}-b},\ \mathbb{F}_{-\lambda\frac{v}{u-2v}\brac*{a-b}-b}},
\end{equation}
and recall that  
\( \slE{\sqbrac*{\lambda}}{1}{1}\cong \mmVM{1}{1}\cten
\mathbb{F}_{\sqbrac*{\lambda\frac{v}{u-2v}\brac*{a-b}-b}}\) and
\(\sigma^{-1} \slE{\sqbrac{\lambda_{u-1,v-1}+\epsilon}}{1}{1}\cong
\mmVM{1}{1}\cten\mathbb{F}_{\sqbrac*{\epsilon\frac{v}{u-2v}\brac*{a-b}-2b}}=\bee^\epsilon\HH^0\). Tensoring
\(I_{\lambda+\epsilon,-\lambda}\) with the action of \(\mmV{u}{v}\) on \(\mmVM{1}{1}\) yields an
\(\mmsl{u}{v}\)-intertwining operator
\begin{equation}
  \mathcal{Y}_{\lambda+\epsilon,-\lambda}\in\binom{\sigma^{-1} \slE{\sqbrac{\lambda_{u-1,v-1}+\epsilon}}{1}{1}}{\slE{\sqbrac*{\lambda+\epsilon}}{1}{1},\ \slE{\sqbrac*{-\lambda}}{1}{1}}
\end{equation}
by restriction.

\begin{cor}
  For any $o\in (\mathbb{H}^1)'$,
  $m\in\mmVM{1}{1}\cten\mathbb{F}_{\sqbrac{\lambda\frac{v}{u-2v}(a-b)-b}}\cong
  \slE{\sqbrac*{\lambda}}{1}{1}$ and
  $n\in
  \mmVM{1}{1}\cten\mathbb{F}_{\sqbrac{-\lambda\frac{v}{u-2v}(a-b)-b}}\cong
  \slE{\sqbrac*{-\lambda}}{1}{1}$
  \begin{equation}
    \pair{ o}{\mathcal{Q}^{-}(z_1)\mathcal{Q}^{+}(z_2)\mathbf{e}^{-\epsilon}\mathcal{Y}_{\lambda+\epsilon,-\lambda}(\mathbf{e}^\epsilon m,w)n}
    \in \mathcal{U}_w^\epsilon\brac*{z}\CC\sqbrac{z^{\pm 1}_1,z^{\pm
        1}_2,(z_1-w)^{-1},(z_2-w)^{-1}}^{\mathfrak{S}}. 
  \end{equation}
  Specifically, \(o_0=\tilde{v}_{1,1}\cten \bra{0}\), \(m_0={v}_{1,1}\cten \ket{\frac{2v\lambda}{u-2v}(a-b)-b}\) and \(n_0={v}_{1,1}\cten \ket{\frac{-2v\lambda}{u-2v}(a-b)-b}\) we have
  \begin{align}
    \pair{o_0 }{\mathcal{Q}^{-}(z_1)\mathcal{Q}^{+}(z_2)\mathbf{e}^{-\epsilon}\mathcal{Y}_{\lambda+\epsilon,-\lambda}\brac*{\mathbf{e}^\epsilon m_0,w}n_0}
    =w^{-\frac{u-2+\epsilon v}{2v}}\mathcal{U}_w^\epsilon\brac*{z}
  \end{align}
  up to a phase.
\end{cor}

Then we can define a bilinear operator 
$\mathcal{Q}^{[2]}*^{(\epsilon)}\mathcal{Y}_{\lambda,-\lambda}(-,w)-  \colon
\slE{\lambda}{1}{1} \times \slE{-\lambda}{1}{1} \to \overline{\HH^1}$. 
as
\begin{equation}
  \mathcal{Q}^{[2]}*^{(\epsilon)}\mathcal{Y}_{\lambda,-\lambda}(m,w)n=\int_{[\Delta_{w,0}]}\mathcal{Q}^{-}(z_1)\mathcal{Q}^{+}(z_2)\mathbf{e}^{-\epsilon}\mathcal{Y}_{\lambda+\epsilon,-\lambda}(\mathbf{e}^\epsilon
  m,w)n\dd z_1\dd z_2,
\end{equation}
where \(m\in\slE{\lambda}{1}{1}\), \(n\in \slE{-\lambda}{1}{1}\).

Let $\iota \colon Z_{w} \to Z_{w^{-1}}$ be a diffeomorphism defined by 
$ \iota(z_1,z_2)=(z^{-1}_1,z^{-1}_2).$
Let $u_1=z^{-1}_1$, $u_2=z^{-1}_2$. Then we have
\begin{equation}
  \iota_*(\mathcal{U}_{w,\epsilon}(z_1,z_2)\dd z_1\wedge \dd z_2)=(u_2-u_1)^{-2\frac{u}{v}+2}\prod_{i=1,2}u^{\frac{u}{v}-2}_i(1-u_iw)^{\frac{u-2v+(\lambda+\epsilon)v}{2v}}\dd u_1\wedge\dd u_2.
  \label{eq:trans}
\end{equation}
Thus we can pair this differential form with the twisted cycle $\Gamma$ defined by (\ref{eq:gammatwist}).
By pulling back $\Gamma$ to the original variables $(z_1,z_2)$, we have a twisted cycle
$ [\iota^*(\Gamma)]\in H_2(Z_{w},\mathcal{L}^\vee_w)$. 
The cycle $\Gamma_{\infty}=\iota^*(\Gamma)$ allows us to define a
  transpose of the operator \(\mathcal{Q}^{[2]}\).
\begin{lemma}
  Let \(\tensor[^t]{\mathcal{Q}}{^{[2]}}:\brac*{\HH^1}^\prime\to
  \brac*{\HH^0}^{\prime}\) be the transpose of \(\mathcal{Q}^{[2]}\).
  That is, the a linear map uniquely characterised by the identity
  \begin{equation}
    \pair{\tensor[^t]{\mathcal{Q}}{^{[2]}}o}{m}=\pair{o}{\mathcal{Q}^{[2]}m}.
    \label{eq:remQ2}
  \end{equation}
  Then an explicity integral formula for
  \(\tensor[^t]{\mathcal{Q}}{^{[2]}}\) is given by
    \begin{equation}
      \int_{[\Gamma_\infty]}o\brac{\mathcal{Q}^{-\,\opp}(z_1)\mathcal{Q}^{+\, \opp}(z_2)-}\dd z_1\dd z_2,
    \end{equation}
    where \(-\) indicates the argument from \(\HH^0\), and where
    \(\mathcal{Q}^{\pm\, \opp}\) denotes the application the opposition formula
    \eqref{eq:oppformula}, that is, \(\mathcal{Q}^{\pm\, \opp}(z)=-z^{-2}\mathcal{Q}^\pm(z^{-1})\).
  \end{lemma}
  \begin{proof}
    The lemma follows immediate from the definition of \([\Gamma_\infty]\) as
    \begin{align}
      \pair{\tensor[^t]{\mathcal{Q}}{^{[2]}}o}{m}
      &=\int_{[\Gamma_\infty]}\pair{o}{\mathcal{Q}^-(z_1^{-1})\mathcal{Q}^+(z_2^{-1})m}\frac{\dd
      z_1 \dd z_2}{z_1z_2}=
      \int_{[\Gamma]}
      \pair{o}{\mathcal{Q}^-(u_1)\mathcal{Q}^+(u_2)m}\dd
        u_1 \dd u_2\nonumber\\
      &=\pair{o}{\mathcal{Q}^{[2]}m}.
    \end{align}
  \end{proof}
Next we consider the operator analogous to
\(\mathcal{Q}^{[2]}{*}^{(\epsilon)}_{\infty}\mathcal{Y}_{\lambda,-\lambda}(-,w)-\),
but integrated over \([\Gamma_\infty]\) instead of \([\Gamma]\).
\begin{lemma}
 Consider linear operator 
$\mathcal{Q}^{[2]}{*}^{(\epsilon)}_{\infty}\mathcal{Y}_{\lambda,-\lambda}(-,w)-\colon \slE{\lambda}{1}{1} \times \slE{-\lambda}{1}{1} \to \overline{\HH^1}$
defined by
\begin{equation}
\mathcal{Q}^{[2]}{*}^{(\epsilon)}_{\infty}\mathcal{Y}_{\lambda,-\lambda}(m,w)n=\int_{[\Gamma_{\infty}]}\mathcal{Q}^{-}(z_1)\mathcal{Q}^{+}(z_2)\bee^{-\epsilon}\mathcal{Y}_{\lambda+\epsilon,-\lambda}(\bee^{\epsilon} m,w)n\dd z_1\dd z_2,
\end{equation}
for $m\in \slE{\lambda}{1}{1}$ and $n\in \slE{-\lambda}{1}{1}$. Then
\begin{equation}
 Q^{[2]}*^{(\epsilon)}\mathcal{Y}_{\lambda,-\lambda}=c_{\lambda,\epsilon}Q^{[2]}*^{(\epsilon)}_{\infty}\mathcal{Y}_{\lambda,-\lambda},
  \label{eq:upto} 
\end{equation}
where
\begin{equation}
  c_{\lambda,\epsilon}
  =\frac{e^{3\pi\ii \frac{u}{v}}}{4\sin(\pi \tfrac{u}{v})\sin
    (2\pi\tfrac{u}{v})}\frac{\Gamma(1-\tfrac{u}{v})}{\Gamma(\tfrac{u}{v}-1)\Gamma(3-2\tfrac{u}{v})}I_1[1]\brac*{\tfrac{u-(\lambda +\epsilon)v}{2v}-1,\tfrac{u+(\lambda +\epsilon)v}{2v}-1,1-\tfrac{u}{v}}.
\end{equation}
Further,
\begin{equation}
\pair{o}{\mathcal{Q}^{[2]}*^{(\epsilon)}\mathcal{Y}_{\lambda,-\lambda}(m,w)n}
=c_{\lambda,\epsilon}\pair{ \mathcal{Q}^{[2]}o}{\mathbf{e}^{-\epsilon}\mathcal{Y}_{\lambda+\epsilon,-\lambda}(\mathbf{e}^\epsilon m,w)n} .
\label{eq:upto2}
\end{equation}
\end{lemma}
\begin{proof}
  The lemma follows by direct computation. By Lemma \ref{thm:one-dim}, we know
  that for any $o\in (\mathbb{H}^1)'$, $m\in \slE{\lambda}{1}{1}$ and
  $n\in \slE{-\lambda}{1}{1}$, there exists a constant $c\in \CC$ such that
\begin{align}
&w^{\frac{u-2v +\epsilon v}{2v}}\pair{ o}{\mathcal{Q}^{-}(z_1)\mathcal{Q}^{+}(z_2)\bee^{-\epsilon}\mathcal{Y}_{\lambda+\epsilon,-\lambda}(\bee^{\epsilon} m,w)n}\dd z_1\wedge\dd z_2\nonumber\\
&\qquad=
c\ \pair{ o_0}{\mathcal{Q}^{-}(z_1)\mathcal{Q}^{+}(z_2)\bee^{-\epsilon}\mathcal{Y}_{\lambda+\epsilon,-\lambda}(\bee^{\epsilon} m_0,w)n_0}\dd z_1\wedge\dd z_2+\dd\brac*{\cdots}\nonumber\\
&\qquad=c\ \mathcal{U}^{\epsilon}_w(z_1,z_2)\dd z_1\wedge\dd z_2+\dd\brac*{\cdots},
\label{eq:matrixelsimp}
\end{align}
where \(\dd\brac*{\cdots}\) is some total derivative, that is, an exact
\(2\)-form, and where $m_0={v}_{1,1}\cten\ket{\frac{2v\lambda}{u-2v}(a-b)-b}$,
$n_0={v}_{1,1}\cten \ket{\frac{-2v\lambda}{u-2v}(a-b)-b}$ and
$o_0=\tilde{v}_{1,1}\cten\bra{0}$.
By integrating \eqref{eq:matrixelsimp} over \([\Gamma]\) and
\([\Gamma_\infty]\) we obtain the identities
\begin{align}
  \pair{o}{Q^{[2]}*^{(\epsilon)}\mathcal{Y}_{\lambda,-\lambda}(m,w)n}&=
  cw^{\frac{2v-u -\epsilon v}{2v}}\int_{\sqbrac*{\Delta_{w,0}}}\ \mathcal{U}^{\epsilon}_w(z_1,z_2)\dd
  z_1\wedge\dd z_2,\nonumber\\
    \pair{o}{Q^{[2]}*^{(\epsilon)}_\infty\mathcal{Y}_{\lambda,-\lambda}(m,w)n}&=
  cw^{\frac{2v-u -\epsilon v}{2v}}\int_{\sqbrac*{{\Gamma_\infty}}}\ \mathcal{U}^{\epsilon}_w(z_1,z_2)\dd z_1\wedge\dd z_2.\
\end{align}
The integral over \(\sqbrac*{\Delta_{w,0}}\) is just \eqref{eq:integralrelation}.
While for the \(\sqbrac*{\Gamma_\infty}\)-integral we use \eqref{eq:trans}
  \begin{align}
  \int_{\sqbrac*{{\Gamma_\infty}}}\ \mathcal{U}^{\epsilon}_w(z_1,z_2)\dd
  z_1\wedge\dd
  z_2&=\int_{\sqbrac{\Gamma}}\brac*{u_2-u_1}^{2-2\frac{u}{v}}\prod_{i=1}^2u_i^{\frac{u}{v}-2}\brac*{1-u_iw}^{\frac{u}{2v}-1+\frac{\lambda+\epsilon}{2}}\dd
  u_1\wedge \dd u_2\nonumber\\
  &=\oint_0\int_{\Gamma_{0,1}}\brac*{xy-x}^{2-2\frac{u}{v}}x^{\frac{u}{v}-2}\brac*{xy}^{\frac{u}{v}-2}\brac*{1-xw}^{\frac{u}{2v}-1+\frac{\lambda+\epsilon}{2}}x\dd
    x\dd y\nonumber\\
    &=(-1)^{-2\frac{u}{v}}\int_{\Gamma_{0,1}}y^{\frac{u}{v}-2}\brac*{1-y}^{2-2\frac{u}{v}}\dd
      y\nonumber\\
    &=e^{-2\pi\ii \frac{u}{v}}\brac*{1-e^{2\pi\ii
      \frac{u}{v}}}\brac*{1-e^{-2\pi\ii
      2\frac{u}{v}}}\frac{\Gamma(\frac{u}{v}-1)\Gamma(3-2\frac{u}{v})}{\Gamma(2-\frac{u}{v})}\nonumber\\
    &=4e^{-3\pi\ii \frac{u}{v}}\sin(\pi \tfrac{u}{v})\sin (2\pi\tfrac{u}{v})\frac{\Gamma(\frac{u}{v}-1)\Gamma(3-2\frac{u}{v})}{\Gamma(2-\frac{u}{v})}.
  \end{align}
  Thus
  \begin{equation}
    c_{\lambda,\epsilon}=4e^{-3\pi\ii \frac{u}{v}}\sin(\pi \tfrac{u}{v})\sin
    (2\pi\tfrac{u}{v})\frac{\Gamma(\frac{u}{v}-1)\Gamma(3-2\frac{u}{v})}{\Gamma(2-\frac{u}{v})}
    I_1[1]\brac*{\tfrac{u-(\lambda +\epsilon)v}{2v}-2,\tfrac{u+(\lambda +\epsilon)v}{2v}-2,1-\frac{u}{v}}.
  \end{equation}
Hence
\begin{align}
  c_{\lambda,\epsilon}&=\frac{\pair{o_0}{\mathcal{Q}^{[2]}*^{(\epsilon)}\mathcal{Y}_{\lambda,-\lambda}(m_0,w)n_0}}
 {\pair{o_0}{\mathcal{Q}^{[2]}*^{(\epsilon)}_{\infty}\mathcal{Y}_{\lambda,-\lambda}(m_0,w)n_0}}\nonumber\\
&=\frac{e^{3\pi\ii \frac{u}{v}}}{4\sin(\pi \tfrac{u}{v})\sin
    (2\pi\tfrac{u}{v})}\frac{\Gamma(1-\tfrac{u}{v})}{\Gamma(\tfrac{u}{v}-1)\Gamma(3-2\tfrac{u}{v})}I_1[1]\brac*{\tfrac{u-(\lambda +\epsilon)v}{2v}-1,\tfrac{u+(\lambda +\epsilon)v}{2v}-1,1-\tfrac{u}{v}}.
\end{align}
To obtain \eqref{eq:upto2} note
\begin{align}
  \pair{o}{\mathcal{Q}^{[2]}*^{(\epsilon)}\mathcal{Y}_{\lambda,-\lambda}(m,w)n}&=
  c_{\lambda,\epsilon}\pair{o}{\mathcal{Q}^{[2]}*^{(\epsilon)}_\infty\mathcal{Y}_{\lambda,-\lambda}(m,w)n}\nonumber\\
&=c_{\lambda,\epsilon}\pair{\tensor[^t]{\mathcal{Q}}{^{[2]}} o}{\mathbf{e}^{-\epsilon}\mathcal{Y}_{\lambda+\epsilon,-\lambda}(\mathbf{e}^\epsilon
  m,w)n} .
  \label{eq:Qadj}
\end{align}
\end{proof}

We set $\mathcal{Q}^{[2]}{*}\mathcal{Y}_{\lambda,-\lambda}=\mathcal{Q}^{[2]}{*}^{(0)}\mathcal{Y}_{\lambda,-\lambda}$.

\begin{prop}
The operator $Q^{[2]}*\mathcal{Y}_{\lambda,-\lambda}$, constructed above, is a
non-zero intertwining operator of type 
\begin{equation}
\begin{pmatrix}
\ \mathbb{H}^1 \\
\slE{\lambda}{1}{1}\ \ \slE{-\lambda}{1}{1}
\end{pmatrix}
.
\end{equation}

Furthermore, the image of this intertwining operator is the submodule
$L_1\subset \mathbb{H}^1$  (recall that \(L_1\) is the \voa{}
\(\mmsl{u}{v}\) as a module over itself and hence the tensor unit).
Denote the intertwining operator \(Q^{[2]}*\mathcal{Y}_{\lambda,-\lambda}\)
with the codomain restricted to \(L_1\) by \(e_E\).
\label{thm:Int-Q2}
\end{prop}
We will show later that \(e_E\) defines an evaluation morphism for \(\slE{\lambda}{1}{1}\).
\begin{proof}
Since $\mathcal{Q}^{[2]}{*}\mathcal{Y}_{\lambda,-\lambda}$ clearly satisfies
the truncation condition, we begin by showing that it also satisfies the convergence condition.
Let $m\in E_{\lambda;1,1}$, $n\in E_{-\lambda;1,1}$, $o \in
(\mathbb{H}^1)' $ and $g\in \mmsl{u}{v}$ and 
consider
\begin{align}
\pair{o}{\mathcal{Q}^{[2]}{*}\mathcal{Y}_{\lambda,-\lambda}(m,w)Y(g,z)n}&=\sum_{k\in \mathbb{Z}}\pair{o}{\mathcal{Q}^{[2]}{*}\mathcal{Y}_{\lambda,-\lambda}(m,w)g_kn} z^{-k-1},\nonumber\\
\pair{o}{\mathcal{Q}^{[2]}{*}\mathcal{Y}_{\lambda,-\lambda}(Y(g,z-w)m,w)n}&=\sum_{k\in \mathbb{Z}}
\pair{o}{\mathcal{Q}^{[2]}{*}\mathcal{Y}_{\lambda,-\lambda}(g_km,w)n}(z-w)^{-k-1},\nonumber\\
\pair{o}{Y(g,z)\mathcal{Q}^{[2]}{*}\mathcal{Y}_{\lambda,-\lambda}(m,w)n}&=\sum_{k\in \mathbb{Z}}\pair{o}{g_k\mathcal{Q}^{[2]}{*}\mathcal{Y}_{\lambda,-\lambda}(m,w)n} z^{-k-1}.
  \label{eq:seri}
\end{align}
The convergence of the first two series in their respective domains follows from
applying the identity \eqref{eq:Qadj} to transpose \(\mathcal{Q}^{[2]}\) and
move it to the left side of the pairing, and further noting that
\(\mathcal{Y}_{\lambda,-\lambda}\) is an intertwining operator and hence
satisfies the convergence condition. The convergence of the third series
follows from
\begin{align}
  \pair{o}{Y(g,z)\mathcal{Q}^{[2]}{*}\mathcal{Y}_{\lambda,-\lambda}(m,w)n}&=
                                                                            \pair{Y(g,z)^{\opp}o}{\mathcal{Q}^{[2]}{*}\mathcal{Y}_{\lambda,-\lambda}(m,w)n}\nonumber\\
&= c_{\lambda,0}\pair{\tensor[^t]{\mathcal{Q}}{^{[2]}}Y(g,z)^{\opp}o}{\mathcal{Y}_{\lambda,-\lambda}(m,w)n}\nonumber\\
&=c_{\lambda,0}\pair{Y(g,z)^{\opp}\tensor[^t]{\mathcal{Q}}{^{[2]}}o}{\mathcal{Y}_{\lambda,-\lambda}(m,w)n}\nonumber\\
&=c_{\lambda,0}\pair{\tensor[^t]{\mathcal{Q}}{^{[2]}}o}{Y(g,z)\mathcal{Y}_{\lambda,-\lambda}(m,w)n},
\end{align}
where the third equality uses that \(\mathcal{Q}^{[2]}\) and hence also its
transpose commutes with the action of the \voa{,} and convergence then follows
from \(\mathcal{Y}_{\lambda,-\lambda}\) satisfying the convergence condition.

The \(P(w)\)-compatibility condition for
$\mathcal{Q}^{[2]}{*}\mathcal{Y}_{\lambda,-\lambda}$ also follows by
transposing \(\mathcal{Q}^{[2]}\) to move it to the left side of the pairing. That
is for $f(z)\in\mathbb{C}[z,z^{-1},(z-w)^{-1}]$, consider 
\begin{align}
  \oint_{0,w}f(z)\pair{ o}{Y(g,z)\mathcal{Q}^{[2]}{*}\mathcal{Y}_{\lambda,-\lambda}(m,w)n}\dd z
  &-
  \oint_wf(z)\pair{o}{\mathcal{Q}^{[2]}{*}\mathcal{Y}_{\lambda,-\lambda}(Y(g,z-w)m,w)n}\dd z\nonumber\\
  &-\oint_wf(z)\pair{o}{\mathcal{Q}^{[2]}{*}\mathcal{Y}_{\lambda,-\lambda}(m,w)Y(g,z)n}\dd
  z\nonumber\\
  =\oint_{0,w}c_{\lambda,\epsilon}f(z)\pair{\tensor[^t]{\mathcal{Q}}{^{[2]}}o}{Y(g,z)\mathcal{Y}_{\lambda,-\lambda}(m,w)n}\dd z
  &-
  \oint_wc_{\lambda,\epsilon}f(z)\pair{\tensor[^t]{\mathcal{Q}}{^{[2]}}o}{\mathcal{Y}_{\lambda,-\lambda}(Y(g,z-w)m,w)n}\dd z\nonumber\\
  &-\oint_wc_{\lambda,\epsilon}f(z)\pair{\tensor[^t]{\mathcal{Q}}{^{[2]}}o}{\mathcal{Y}_{\lambda,-\lambda}(m,w)Y(g,z)n}\dd z.
\end{align}
The \rhs{} of the above equality vanishes due to
\(\mathcal{Y}_{\lambda,-\lambda}\) satisfying the Jacobi identity and thus so does
\(\mathcal{Q}^{[2]}{*}\mathcal{Y}_{\lambda,-\lambda}\). 
\end{proof}

\begin{lemma}
Consider the bilinear operator
${\rm Ad}_{\bee^{\epsilon}}\mathcal{Q}^{[2]}*\mathcal{Y}_{\lambda,-\lambda}(-,w)- \colon \slE{\lambda}{1}{1} \cten \slE{-\lambda}{1}{1} \to \overline{\mathbb{H}}^1$
defined by
\begin{equation}
{\rm Ad}_{\bee^{\epsilon}}\mathcal{Q}^{[2]}*\mathcal{Y}_{\lambda,-\lambda}(m,w)n =w^\epsilon\int_{[\Delta_{w,0}\bigl(\frac{u-2v-(\lambda+\epsilon) v}{2v},\frac{u-2v+\lambda v}{2v},-\frac{u}{v}+1\bigr)]}\bee^{\epsilon}\mathcal{Q}^{-}(z_1)\mathcal{Q}^{+}(z_2)\bee^{-\epsilon}\mathcal{Y}_{\lambda,-\lambda}(m,w)n,
\end{equation}
where \(m\in \slE{\lambda}{1}{1}\), \(n\in \slE{-\lambda}{1}{1}\).  Then
the operator
 \begin{equation}
   \mathcal{Y}^{\Delta}_{\lambda,-\lambda}=\left.\frac{\dd}{\dd \epsilon}{\rm Ad}_{\bee^{\epsilon}}\mathcal{Q}^{[2]}*\mathcal{Y}_{\lambda,-\lambda}\right|_{\epsilon=0}\colon \slE{\lambda}{1}{1} \cten \slE{-\lambda}{1}{1} \to \overline{\mathbb{H}}^1
 \end{equation}
 is well-defined and non-trivial.
\end{lemma}

\begin{proof}
By \cite[Cor 1.1.1]{SusSel}, we see that
\begin{equation}
  \pair{o}{{\rm Ad}_{\bee^{\epsilon}}\mathcal{Q}^{[2]}*\mathcal{Y}_{\lambda,-\lambda}(m,w)n},\qquad m\in \slE{\lambda}{1}{1},\ n\in \slE{-\lambda}{1}{1},\ o \in (\mathbb{H}^1)'
\end{equation}
is holomorphic in the variable \(\epsilon\) in a neighbourhood of
  $\epsilon=0$. This can also be verified by expanding the image of ${\rm
  Ad}_{\bee^{\epsilon}}\mathcal{Q}^{[2]}*\mathcal{Y}_{\lambda,-\lambda}$
in a monomial basis of \(\mathbb{H}^1\) (a basis given by monomials of negative Virasoro and
  Heisenberg modes applied to \(v_{1,1}\cten\ket{\frac{2nv}{u-2v}(a-b)}\)
  with an appropriate choice of ordering), where
we see that each coefficient of this expansion is holomorphic in a
  neighbourhood of   $\epsilon=0$.
  Thus \(\mathcal{Y}^{\Delta}_{\lambda,-\lambda}\) is well-defined.

  To conclude non-triviality consider the 
 vectors $m_0={v}_{1,1}\cten \ket{\frac{2v\lambda}{u-2v}(a-b)-b}$,
 $n_0={v}_{1,1}\cten \ket{\frac{-2v\lambda}{u-2v}(a-b)-b}$ and
 $o_0=\tilde{v}_{1,1}\cten \bra{0}$ and note
\begin{align}
  &w^{\frac{u-2v}{2v}}\pair{ o_0}{\mathcal{Y}^{\Delta}_{\lambda,-\lambda}(m_0,w)n_0}
  =w^{\frac{u-2v}{2v}}\partial_{\epsilon}\pair{ o_0}{{\rm Ad}_{\bee^{\epsilon}}\mathcal{Q}^{[2]}*\mathcal{Y}_{\lambda,-\lambda}(m_0,w)n_0}|_{\epsilon=0}\nonumber\\
  &=\left.\frac{\partial}{\partial \epsilon}\prod_{i=1,2}\frac{\Gamma(1+i(1-\frac{u}{v}))\Gamma(\frac{u-(\lambda+\epsilon) v}{2v}+(i-1)(1-\frac{u}{v}))\Gamma(\frac{u+\lambda v}{2v}+(i-1)(1-\frac{u}{v}))}{\Gamma(1+i(1-\frac{u}{v}))\Gamma(2+\frac{u-2v-\epsilon v}{v}+i(1-\frac{u}{v}))}\right|_{\epsilon=0}.
\end{align}
Then by using the asymptotic expansion of the digamma function 
\begin{equation}
  \frac{\Gamma'(z)}{\Gamma(z)}\sim  \log(z)-\frac{1}{2z}-\sum_{n=1}^\infty\frac{B_{2n}}{2nz^{2n}},\qquad |{\rm arg}z|<\pi,\ z\rightarrow \infty,
\end{equation}
we have an asymptotic behaviour
\begin{equation}
  \frac{w^{\frac{u-2v}{2v}}\pair{ o_0}{\mathcal{Y}^{\Delta}_{\lambda,-\lambda}(m_0,w)n_0}}{\Gamma(\frac{u}{2v}+\frac{\lambda}{2})\Gamma(1-\frac{u}{2v}+\frac{\lambda}{2})}\sim C\Gamma(-2^{-1}\lambda)^2 \log(-\lambda),
\qquad \lambda\rightarrow -\infty,
\end{equation}
where $B_{2n}$ are the Bernoulli numbers and $C$ is a some non-zero constant.
Thus $\mathcal{Y}^{\Delta}_{\lambda,-\lambda}$ is non-trivial.
\end{proof}

With the definition of $\mathcal{Y}^{\Delta}_{\lambda,-\lambda}$ in place,
we are now able to give an explicit construction of an intertwining operator
with non-trivial logarithmic parts.
\begin{theorem}
  Recall the bilinear operators \(\mathcal{Y}_{\lambda,-\lambda}:\slE{\lambda}{1}{1} \cten \slE{-\lambda}{1}{1}\rightarrow\overline{\mathbb{H}^0}\) and
  \(\mathcal{Y}^{\Delta}_{\lambda,-\lambda}:\slE{\lambda}{1}{1} \cten \slE{-\lambda}{1}{1}\rightarrow\overline{\mathbb{H}^1}\) defined above and consider
  $\mathcal{Y}^{\log}_{\lambda,-\lambda}=\mathcal{Y}_{\lambda,-\lambda}+\mathcal{Y}^{\Delta}_{\lambda,-\lambda} :\slE{\lambda}{1}{1} \cten \slE{-\lambda}{1}{1}\rightarrow\overline{\mathbb{H}}$. Then
  $\mathcal{Y}^{\log}_{\lambda,-\lambda}$ is a $P(w)$-intertwining operator of
  type
  \begin{equation}
    \binom{(\mathbb{H},\widetilde{Y}_{c_{\lambda,0}})}{\slE{\lambda}{1}{1},\,\slE{-\lambda}{1}{1}}.
  \end{equation}
  \label{thm:logint}
\end{theorem}

The above theorem implies that \(\sigma^{-1}P_{u-1,v-1}\), the projective
cover of the tensor unit, is a direct summand of the fusion product \(\slE{\lambda}{1}{1}\fuse\slE{-\lambda}{1}{1}\).
To prove this theorem we prepare the following lemmas.

\begin{lemma}
  For all $m\in E_{\lambda;1,1}$, $n\in E_{-\lambda;1,1}$, $o \in
  (\mathbb{H}^1)' $ and $g\in \mmsl{u}{v}$ the series
  \begin{align}
\pair{o}{{\rm Ad}_{\bee^{\epsilon}}\mathcal{Q}^{[2]}{*}\mathcal{Y}_{\lambda,-\lambda}(m,w)Y(g,z)n}&=\sum_{k\in \mathbb{Z}}\pair{o}{{\rm Ad}_{\bee^{\epsilon}}\mathcal{Q}^{[2]}{*}\mathcal{Y}_{\lambda,-\lambda}(m,w)g_kn} z^{-k-1},\nonumber\\
\pair{o}{{\rm Ad}_{\bee^{\epsilon}}\mathcal{Q}^{[2]}{*}\mathcal{Y}_{\lambda,-\lambda}(Y(g,z-w)m,w)n}&=\sum_{k\in \mathbb{Z}}
\pair{o}{{\rm Ad}_{\bee^{\epsilon}}\mathcal{Q}^{[2]}{*}\mathcal{Y}_{\lambda,-\lambda}(g_km,w)n}(z-w)^{-k-1},\nonumber\\
\pair{o}{Y(g,z){\rm Ad}_{\bee^{\epsilon}}\mathcal{Q}^{[2]}{*}\mathcal{Y}_{\lambda,-\lambda}(m,w)n}&=\sum_{k\in \mathbb{Z}}\pair{o}{g_k{\rm Ad}_{\bee^{\epsilon}}\mathcal{Q}^{[2]}{*}\mathcal{Y}_{\lambda,-\lambda}(m,w)n} z^{-k-1},
    \label{eq:seri2}
  \end{align}
are absolutely convergent on $w>|z|>0$, $w>|z-w|>0$ and $|z|>w$, respectively.
\label{thm:conv}
\end{lemma}
\begin{proof}
  We only show the first case, as the other cases can be proved in the same way.
  Consider the correlation function
\begin{multline}
  \pair{o}{\bee^{\epsilon}\mathcal{Q}^-(z_1)\mathcal{Q}^+(z_1)\bee^{-\epsilon}\mathcal{Y}_{\lambda,-\lambda}(m,w)Y(g,z)n}=\\z^{-\frac{\epsilon}{2}}_1z^{-\frac{\epsilon}{2}}_2\pair{o}{\mathcal{Q}^-(z_1)\mathcal{Q}^+(z_1)\mathcal{Y}_{\lambda,-\lambda}(m,w)Y(g,z)n},
\end{multline}
which is absolutely convergent on $w> |z|>0$ for each fixed $(z_1,z_2)$.
By expanding
\begin{multline}
  z^{-\frac{\epsilon}{2}}_1z^{-\frac{\epsilon}{2}}_2\pair{o}{\mathcal{Q}^-(z_1)\mathcal{Q}^+(z_2)\mathcal{Y}_{\lambda,-\lambda}(m,w)Y(g,z)n}=\\
  \sum_{k\in \mathbb{Z}}z^{-\frac{\epsilon}{2}}_1z^{-\frac{\epsilon}{2}}_2\pair{o}{\mathcal{Q}^-(z_1)\mathcal{Q}^+(z_2)\mathcal{Y}_{\lambda,-\lambda}(m,w)g_{k}n} z^{-k-1}
  \label{eq:twof}
\end{multline}
we see that
\begin{align}
\label{eq:mullem}
\begin{aligned}
&w^{\frac{u-2v}{2v}}z^{-\frac{\epsilon}{2}}_1z^{-\frac{\epsilon}{2}}_2\pair{o}{\mathcal{Q}^-(z_1)\mathcal{Q}^+(z_2)\mathcal{Y}_{\lambda,-\lambda}(m,w)g_{k}n}\\
&\in \mathcal{U}_{w}\brac*{\tfrac{u-(\lambda+\epsilon) v}{2v}-1,\tfrac{u+\lambda v}{2v}-1,1-\tfrac{u}{v},z_1,z_2}\mathbb{C}[w^{\pm 1}][z^{\pm 1}_1,z^{\pm 1}_2,(z_1-w)^{-1},(z_2-w)^{-1}]^{\mathfrak{S}_2}.
\end{aligned}
\end{align}
Then by Lemma \ref{thm:one-dim}, there exists a rational function \(c_{k}(\lambda,w)\in \mathbb{C}(\lambda)\cten\mathbb{C}[\epsilon,w^{\pm 1}]\)
such that
\begin{align}
  &w^{\frac{u-2v}{2v}}z^{-\frac{\epsilon}{2}}_1z^{-\frac{\epsilon}{2}}_2\pair{o}{\mathcal{Q}^-(z_1)\mathcal{Q}^+(z_2)\mathcal{Y}_{\lambda,-\lambda}(m,w) g_{k}n}\dd z_1\wedge\dd z_2\nonumber\\
  &= c_{k}(\lambda,w)\mathcal{U}_{w}\brac*{\tfrac{u-(\lambda+\epsilon)
    v}{2v}-1,\tfrac{u+\lambda v}{2v}-1,1-\tfrac{u}{v},z_1,z_2}\dd
  z_1\wedge\dd z_2 +\dd\brac*{\cdots}. 
  \label{eq:oned}
\end{align}
By taking sufficiently small absolute value of $z$, we see that $z$ is contained in the tubular neighbourhood of the twisted cycle $[\Delta_{w,0}\bigl(\frac{u-2v-(\lambda+\epsilon) v}{2v},\frac{u-2v+\lambda v}{2v},-\frac{u}{v}+1\bigr)]$ (see \cite[Sec 3.2.4]{hgeom}).
Thus, the first of \eqref{eq:seri2} converges for small $|z|>0$ and satisfies
\begin{multline}
w^{\frac{u-2v}{2v}}\pair{o}{{\rm Ad}_{\bee^{\epsilon}}\mathcal{Q}^{[2]}{*}\mathcal{Y}_{\lambda,-\lambda}(m,w)Y(g,z)n}\\
=I_1[1]\brac*{\tfrac{u-(\lambda +\epsilon)v}{2v}-1,\tfrac{u+\lambda v}{2v}-1,1-\frac{u}{v}}\sum_{k\in \mathbb{Z}}c_{k}(\lambda,w)z^{-k-1}.
\label{eq:seri-sup}
\end{multline}

After having shown convergence for small $|z|>0$, we next need to extend to the
domain \(|w|>|z|>0\).
Let $\iota:Z_{w}\rightarrow Z_{w^{-1}}$ be a diffeomorphism defined by 
$
\iota(z_1,z_2)=(z^{-1}_1,z^{-1}_2).
$
Let $u_1=z^{-1}_1$, $u_2=z^{-1}_2$. Then we have
\begin{align}
  &\iota_*\brac*{\mathcal{U}_{w}\bigl(\tfrac{u-(\lambda+\epsilon) v}{2v}-1,\tfrac{u+\lambda v}{2v}-1,1-\frac{u}{v},z_1,z_2\bigr)\dd z_1\wedge \dd z_2}\nonumber\\
  &=(u_2-u_1)^{-2\frac{u}{v}+2}\prod_{i=1,2}u^{\frac{u}{v}-2+\frac{\epsilon}{2}}_i(1-u_iw)^{\frac{u-2v+\lambda v}{2v}}\dd u_1\wedge\dd u_2\nonumber\\
  &=w^{\frac{u-2v+\lambda v}{2v}}\mathcal{U}_{w^{-1}}\brac*{\tfrac{u}{v}-2+\tfrac{\epsilon}{2},\tfrac{u-2v+\lambda v}{2v},1-\tfrac{u}{v},u_2,u_1}\dd u_1\wedge \dd u_2,
\label{eq:trans2}
\end{align}
where the multivalued function in the pushforward $\iota_*$ is the multivalued function on the right-hand side of \eqref{eq:mullem}.
Then by pulling back the twisted cycle
\begin{equation}
  \sqbrac*{\Delta_{w^{-1},0}\brac*{\tfrac{u}{v}-2+\tfrac{\epsilon}{2},\tfrac{u-2v+\lambda v}{2v},1-\tfrac{u}{v}}}\in H_2\brac*{Z_{w^{-1}},\mathcal{L}^\vee_{w^{-1}}\brac*{\tfrac{u}{v}-2+\tfrac{\epsilon}{2},\tfrac{u+\lambda v}{2v}-1,1-\tfrac{u}{v}}}
\end{equation}
along $\iota$, we have a twisted cycle
\begin{equation}
  \iota^*\sqbrac*{\Delta_{w^{-1},0}\brac*{\tfrac{u}{v}-2+\tfrac{\epsilon}{2},\tfrac{u+\lambda v}{2v}-1,1-\tfrac{u}{v}}}\in H_2\brac*{Z_{w},\mathcal{L}^\vee_{w}\brac*{\tfrac{u-(\lambda+\epsilon) v}{2v}-1,\tfrac{u+\lambda v}{2v}-1,1-\tfrac{u}{v}}}.
\end{equation}
This twisted cycle is a regularisation of the contour $\{\infty>z_1>z_2>w\}$. Then by \eqref{eq:twof} and by \eqref{eq:oned}, we have
\begin{align}
&w^{\frac{u-2v}{2v}}\int_{\iota^*\sqbrac*{\Delta_{w^{-1},0}\brac*{\frac{u}{v}-2+\frac{\epsilon}{2},\frac{u+\lambda v}{2v}-1,1-\frac{u}{v}}}}
z^{-\frac{\epsilon}{2}}_1z^{-\frac{\epsilon}{2}}_2\pair{o}{\mathcal{Q}^-(z_1)\mathcal{Q}^+(z_2)\mathcal{Y}_{\lambda,-\lambda}(m,w)Y(g,z)n}\dd z_1\wedge\dd z_2\nonumber\\
&=I_1[1]\brac*{\tfrac{u}{v}-2+\tfrac{\epsilon}{2},\tfrac{u+\lambda v}{2v}-1,1-\tfrac{u}{v}}\sum_{k\in \mathbb{Z}}c_{k}(\lambda,w)z^{-k-1},
\label{eq:twof-2}
\end{align}
which is absolutely convergent on $w>|z|>0$. Thus by \eqref{eq:seri-sup}, we see that the first series of \eqref{eq:seri2} is absolutely convergent on $w>|z|>0$.
\end{proof}

\begin{lemma}
 The series
  \begin{align}
    \pair{o}{\mathcal{Q}^{[2]}{*}^{(\epsilon)}\mathcal{Y}_{\lambda,-\lambda}(m,w)Y(g,z)n}&=\sum_{k\in \mathbb{Z}}\pair{o}{\mathcal{Q}^{[2]}{*}^{(\epsilon)}\mathcal{Y}_{\lambda,-\lambda}(m,w)g_kn} z^{-k-1},\nonumber\\
    \pair{o}{\mathcal{Q}^{[2]}{*}^{(\epsilon)}\mathcal{Y}_{\lambda,-\lambda}(\bee^{-\epsilon}Y(g,z-w)\bee^{\epsilon}m,w)n}&=\sum_{k\in \mathbb{Z}}
    \pair{o}{\mathcal{Q}^{[2]}{*}^{(\epsilon)}\mathcal{Y}_{\lambda,-\lambda}(\bee^{-\epsilon}g_k\bee^{\epsilon}m,w)n}(z-w)^{-k-1},\nonumber\\
    \pair{o}{Y(g,z)\mathcal{Q}^{[2]}{*}^{(\epsilon)}\mathcal{Y}_{\lambda,-\lambda}(m,w)n}&=\sum_{k\in \mathbb{Z}}\pair{o}{g_k\mathcal{Q}^{[2]}{*}^{(\epsilon)}\mathcal{Y}_{\lambda,-\lambda}(m,w)n} z^{-k-1},
    \label{eq:seri-a}
  \end{align}
  are absolutely convergent on the domains \(w>|z|>0\), \(w>|z-w|>0\) and
    \(|z|>w\), respectively, 
  for any $m\in E_{\lambda;1,1}$, $n\in E_{-\lambda;1,1}$,
  $o \in (\mathbb{H}^1)' $ and $g\in \mmsl{u}{v}$.
\end{lemma}
\begin{proof}
  The proof follows by a very similar argument to the proof of Lemma \ref{thm:conv}.
\end{proof}

\begin{lemma}
For $m\in E_{\lambda;1,1}$, $n\in E_{-\lambda;1,1}$, $o \in (\mathbb{H}^1)'$, $g\in \mmsl{u}{v}$ and $f(z)\in\mathbb{C}[z,z^{-1},(z-w)^{-1}]$, we have
\begin{align}
&\oint_{w,0}f(z)\pair{o}{Y(g,z)\mathcal{Y}^{\Delta}_{\lambda,-\lambda}( m,w)n}\frac{\dd {z}}{2\pi i}\nonumber\\
&=\oint_{w}f(z)\pair{o}{\mathcal{Y}^{\Delta}_{\lambda,-\lambda}( Y(g,z-w)m,w)n}\frac{\dd {z}}{2\pi i}
+\oint_{0}f(z)\pair{o}{\mathcal{Y}^{\Delta}_{\lambda,-\lambda}( m,w)Y(g,z)n}\frac{\dd {z}}{2\pi i}\nonumber\\
&\qquad +c_{\lambda,-\lambda}\oint_{w,0}f(z)\pair{o}{\Delta_1(g;z)\mathcal{Y}_{\lambda,-\lambda}( m,w)n}\frac{\dd {z}}{2\pi i}.
\label{eq:Deltarel}
\end{align}
\label{thm:Deltarel}
\end{lemma}
\begin{proof}
  From the definition of $\mathcal{Y}^{\Delta}_{\lambda,-\lambda}$, it is enough to show the equality
  \begin{align}
    &\oint_{w}f(z)\pair{o}{{\rm Ad}_{\bee^{\epsilon}}\mathcal{Q}^{[2]}{*}\mathcal{Y}_{\lambda,-\lambda}( Y(g,z-w)m,w)n}\frac{\dd {z}}{2\pi i}
    +\oint_{0}f(z)\pair{o}{{\rm
        Ad}_{\bee^{\epsilon}}\mathcal{Q}^{[2]}{*}\mathcal{Y}_{\lambda,-\lambda}(
      m,w)Y(g,z)n}\frac{\dd {z}}{2\pi i}\nonumber \\
    &-\oint_{w,0}f(z)\pair{o}{Y(g,z){\rm
        Ad}_{\bee^{\epsilon}}\mathcal{Q}^{[2]}{*}\mathcal{Y}_{\lambda,-\lambda}( m,w)n}\frac{\dd {z}}{2\pi i}\nonumber \\
    &= -\epsilon c_{\lambda+\epsilon,-\lambda}\oint_{w,0}f(z)\pair{o}{\Delta_1(g;z;\epsilon)\mathcal{Y}_{\lambda,-\lambda}( m,w)n}\frac{\dd {z}}{2\pi i}+\brac*{\dots},
    \label{eq:Deltarel-e}
  \end{align}
 where \(\brac*{\dots}\) denotes functions which are homomorphic
  near $\epsilon=0$ with
  vanishing derivatives at $\epsilon=0$, since
\eqref{eq:Deltarel} is obtained by differentiating both sides of \eqref{eq:Deltarel-e} by $\epsilon$ and setting $\epsilon=0$.
By the $P(w)$-compatibility condition of $\mathcal{Y}_{\lambda,-\lambda}$ and by the proof of \cref{thm:comq}, we see that there exist $F_1,F_2\in \mathbb{C}[z_1,z^{-1}_1,z_2,z^{-1}_2,(z_1-w)^{-1},(z_2-w)^{-1},w,w^{-1}]$ satisfying
\begin{align}
  &\oint_{w}f(z)\pair{o}{\bee^{\epsilon}\mathcal{Q}^{-}(z_1)\mathcal{Q}^{+}(z_2)\bee^{-\epsilon}\mathcal{Y}_{\lambda,-\lambda}( Y(g,z-w)m,w)n}\frac{\dd {z}}{2\pi i}\nonumber\\
  &\qquad+\oint_{0}f(z)\pair{o}{\bee^{\epsilon}\mathcal{Q}^{-}(z_1)\mathcal{Q}^{+}(z_2)\bee^{-\epsilon}\mathcal{Y}_{\lambda,-\lambda}( m,w)Y(g,z)n}\frac{\dd {z}}{2\pi i}\nonumber\\
  &\qquad-\oint_{w,0}f(z)\pair{o}{Y(g,z)\bee^{\epsilon}\mathcal{Q}^{-}(z_1)\mathcal{Q}^{+}(z_2)\bee^{-\epsilon}\mathcal{Y}_{\lambda,-\lambda}( m,w)n}\frac{\dd {z}}{2\pi i}\nonumber\\
  &=\oint_{w,0}f(z)\pair{o}{\bee^{\epsilon}\mathcal{Q}^{-}(z_1)\mathcal{Q}^{+}(z_2)\bee^{-\epsilon}Y(g,z)\mathcal{Y}_{\lambda,-\lambda}( m,w)n}\frac{\dd {z}}{2\pi i}\nonumber\\
  &\qquad-\oint_{w,0}f(z)\pair{ o}{Y(g,z)\bee^{\epsilon}\mathcal{Q}^{-}(z_1)\mathcal{Q}^{+}(z_2)\bee^{-\epsilon}\mathcal{Y}_{\lambda,-\lambda}( m,w)n}\frac{\dd {z}}{2\pi i}\nonumber\\
  &=w^{-\frac{u-2v }{2v}}z^{-\frac{\epsilon}{2}}_1z^{-\frac{\epsilon}{2}}_2\dd _{z_1,z_2}\brac*{\mathcal{U}_w\brac*{\tfrac{u-\lambda v}{2v}-1,\tfrac{u+\lambda v}{2v}-1,1-\tfrac{u}{v},z_1,z_2}\bigl(F_1\dd z_2+F_2\dd z_1\bigr)}.
  \label{eq:sup1}
\end{align}
Similarly, we see that there exist $F^{(\epsilon)}_1,F^{(\epsilon)}_2\in \mathbb{C}[z_1,z^{-1}_1,z_2,z^{-1}_2,(z_1-w)^{-1},(z_2-w)^{-1},w,w^{-1},\epsilon]$ satisfying
\begin{align}
  &\oint_{w}f(z)\pair{o}{\mathcal{Q}^{-}(z_1)\mathcal{Q}^{+}(z_2)\bee^{-\epsilon}\mathcal{Y}_{\lambda+\epsilon,-\lambda}(Y(g,z-w)\bee^{\epsilon}m,w)n}\frac{\dd {z}}{2\pi i}\nonumber\\
  &\qquad+\oint_{0}f(z)\pair{o}{\mathcal{Q}^{-}(z_1)\mathcal{Q}^{+}(z_2)\bee^{-\epsilon}\mathcal{Y}_{\lambda+\epsilon,-\lambda}(\bee^{\epsilon}m,w)Y(g,z)n}\frac{\dd {z}}{2\pi i}\nonumber\\
  &\qquad-\oint_{w,0}f(z)\pair{o}{\bee^{-\epsilon}Y(g,z)\bee^{\epsilon}\mathcal{Q}^{-}(z_1)\mathcal{Q}^{+}(z_2)\bee^{-\epsilon}\mathcal{Y}_{\lambda+\epsilon,-\lambda}( \bee^{\epsilon} m,w)n}\frac{\dd {z}}{2\pi i}\nonumber\\
  &=\oint_{w,0}f(z)\pair{o}{\mathcal{Q}^{-}(z_1)\mathcal{Q}^{+}(z_2)\bee^{-\epsilon}Y(g,z)\mathcal{Y}_{\lambda+\epsilon,-\lambda}(\bee^{\epsilon}m,w)n}\frac{\dd {z}}{2\pi i}\nonumber\\
  &\qquad-\oint_{w,0}f(z)\pair{o}{\bee^{-\epsilon}Y(g,z)\bee^{\epsilon}\mathcal{Q}^{-}(z_1)\mathcal{Q}^{+}(z_2)\bee^{-\epsilon}\mathcal{Y}_{\lambda+\epsilon,-\lambda}( \bee^{\epsilon} m,w)n}\frac{\dd {z}}{2\pi i}\nonumber\\
  &=z^{-\frac{\epsilon}{2}}_1z^{-\frac{\epsilon}{2}}_2\oint_{w,0}f(z)\pair{o}{\bee^{-\epsilon}\mathcal{Q}^{-}(z_1)\mathcal{Q}^{+}(z_2)Y(g,z)\mathcal{Y}_{\lambda+\epsilon,-\lambda}(\bee^{\epsilon}m,w)n}\frac{\dd {z}}{2\pi i}\nonumber\\
  &\qquad-z^{-\frac{\epsilon}{2}}_1z^{-\frac{\epsilon}{2}}_2\oint_{w,0}f(z)\pair{o}{\bee^{-\epsilon}Y(g,z)\mathcal{Q}^{-}(z_1)\mathcal{Q}^{+}(z_2)\mathcal{Y}_{\lambda+\epsilon,-\lambda}( \bee^{\epsilon} m,w)n}\frac{\dd {z}}{2\pi i}\nonumber\\
  &=w^{-\frac{u-2v+\epsilon v}{2v}}z^{-\frac{\epsilon}{2}}_1z^{-\frac{\epsilon}{2}}_2\dd _{z_1,z_2}\brac*{\mathcal{U}_w\brac*{\tfrac{u-\lambda v}{2v}-1,\tfrac{u+(\lambda+\epsilon) v}{2v}-1,1-\tfrac{u}{v},z_1,z_2}\brac*{F^{(\epsilon)}_1\dd z_2+F^{(\epsilon)}_2\dd z_1}}
  \label{eq:sup2}
\end{align}
and $F^{(0)}_1=F_1$, $F^{(0)}_2=F_2$. 
Thus, by noting \eqref{eq:seri-a}, \eqref{eq:sup1} and \eqref{eq:sup2}, we have
\begin{align}
  &\oint_{w}f(z)\pair{o}{{\rm Ad}_{\bee^{\epsilon}}\mathcal{Q}^{[2]}{*}\mathcal{Y}_{\lambda,-\lambda}( Y(g,z-w)m,w)n}\frac{\dd {z}}{2\pi i}
  +\oint_{0}f(z)\pair{o}{{\rm Ad}_{\bee^{\epsilon}}\mathcal{Q}^{[2]}{*}\mathcal{Y}_{\lambda,-\lambda}( m,w)Y(g,z)n}\frac{\dd {z}}{2\pi i}\nonumber\\
  &\qquad-\oint_{w,0}f(z)\pair{o}{Y(g,z){\rm Ad}_{\bee^{\epsilon}}\mathcal{Q}^{[2]}{*}\mathcal{Y}_{\lambda,-\lambda}( m,w)n}\frac{\dd {z}}{2\pi i}\nonumber\\
  &= w^{\frac{\epsilon}{2}}\Bigl(\oint_{w}f(z)\pair{o}{\mathcal{Q}^{[2]}{*}^{(\epsilon)}\mathcal{Y}_{\lambda,-\lambda}( \bee^{-\epsilon}Y(g,z-w)\bee^{\epsilon}m,w)n}\frac{\dd {z}}{2\pi i}\nonumber\\
  &\qquad+\oint_{0}f(z)\pair{o}{\mathcal{Q}^{[2]}{*}^{(\epsilon)}\mathcal{Y}_{\lambda,-\lambda}( m,w)Y(g,z)n}\frac{\dd {z}}{2\pi i}\nonumber\\
  &\qquad-\oint_{w,0}f(z)\pair{o}{\bee^{-\epsilon}Y(g,z)\bee^{\epsilon}\mathcal{Q}^{[2]}{*}^{(\epsilon)}\mathcal{Y}_{\lambda,-\lambda}( m,w)n}\frac{\dd {z}}{2\pi i}\Bigr)+\brac*{\dots},
  \label{eq:Deltarel-e2}
\end{align}
where \(\brac*{\dots}\) denotes functions which are homomorphic
  near $\epsilon=0$ with
  vanishing derivatives at $\epsilon=0$.
By \eqref{eq:upto2}, we have
\begin{align}
  &\oint_{w}f(z)\pair{o}{\mathcal{Q}^{[2]}{*}^{(\epsilon)}\mathcal{Y}_{\lambda,-\lambda}( \bee^{-\epsilon}Y(g,z-w)\bee^{\epsilon}m,w)n}\frac{\dd {z}}{2\pi i}
  +\oint_{0}f(z)\pair{o}{\mathcal{Q}^{[2]}{*}^{(\epsilon)}\mathcal{Y}_{\lambda,-\lambda}( m,w)Y(g,z)n}\frac{\dd {z}}{2\pi i}\nonumber\\
  &\qquad-\oint_{w,0}f(z)\pair{o}{\bee^{-\epsilon}Y(g,z)\bee^{\epsilon}\mathcal{Q}^{[2]}{*}^{(\epsilon)}\mathcal{Y}_{\lambda,-\lambda}( m,w)n}\frac{\dd {z}}{2\pi i}\nonumber\\
  &=c_{\lambda,\epsilon}\Bigl(\oint_{w}f(z)\pair{\mathcal{Q}^{[2]}o}{\bee^{-\epsilon}\mathcal{Y}_{\lambda+\epsilon,-\lambda}( Y(g,z-w)\bee^{\epsilon}m,w)n}\frac{\dd {z}}{2\pi i}\nonumber\\
  &\qquad+\oint_{0}f(z)\pair{\mathcal{Q}^{[2]}o}{\bee^{-\epsilon}\mathcal{Y}_{\lambda+\epsilon,-\lambda}( \bee^{\epsilon}m,w)Y(g,z)n}\frac{\dd {z}}{2\pi i}\nonumber\\
  &\qquad-\oint_{w,0}f(z)\pair{\mathcal{Q}^{[2]}\bee^{-\epsilon}Y(g,z)^{\rm opp}\bee^{\epsilon}o}{\bee^{-\epsilon}\mathcal{Y}_{\lambda+\epsilon,-\lambda}( \bee^{\epsilon}m,w)n}\frac{\dd {z}}{2\pi i}\Bigr)\nonumber\\
  &=c_{\lambda,\epsilon}\Bigl(\oint_{w,0}f(z)\pair{\mathcal{Q}^{[2]}o}{\bee^{-\epsilon}Y(g,z)\mathcal{Y}_{\lambda+\epsilon,-\lambda}(\bee^{\epsilon}m,w)n}\frac{\dd {z}}{2\pi i}\nonumber\\
  &\qquad-\oint_{w,0}f(z)\pair{\mathcal{Q}^{[2]}\bee^{-\epsilon}Y(g,z)^{\rm opp}\bee^{\epsilon}o}{\bee^{-\epsilon}\mathcal{Y}_{\lambda+\epsilon,-\lambda}( \bee^{\epsilon}m,w)n}\frac{\dd {z}}{2\pi i}\Bigr)\nonumber\\
  &=c_{\lambda,\epsilon}\Bigl(\oint_{w,0}f(z)\pair{\bee^{-\epsilon}Y(g,z)^{\rm opp}\bee^{\epsilon}\mathcal{Q}^{[2]}o}{\bee^{-\epsilon}\mathcal{Y}_{\lambda+\epsilon,-\lambda}(\bee^{\epsilon}m,w)n}\frac{\dd {z}}{2\pi i}\nonumber\\
  &\qquad-\oint_{w,0}f(z)\pair{\mathcal{Q}^{[2]}\bee^{-\epsilon}Y(g,z)^{\rm opp}\bee^{\epsilon}o}{\bee^{-\epsilon}\mathcal{Y}_{\lambda+\epsilon,-\lambda}( \bee^{\epsilon}m,w)n}\frac{\dd {z}}{2\pi i}\Bigr)\nonumber\\
  &=-c_{\lambda,\epsilon}\oint_{w,0}f(z)\res_{y=z}\pair{\mathcal{Q}^{[2]}(y)\bee^{-\epsilon}Y(g,z)^{\rm opp}\bee^{\epsilon}o}{\bee^{-\epsilon}\mathcal{Y}_{\lambda+\epsilon,-\lambda}( \bee^{\epsilon}m,w)n}\frac{\dd {z}}{2\pi i}.
\end{align}
By \eqref{eq:remQ2}, we see that the last term of the above identity becomes
\begin{align}
  &-c_{\lambda,\epsilon}\oint_{w,0}f(z)\res_{y=z}\pair{o}{\mathcal{Q}^{[2]}(y)\bee^{-\epsilon}Y(g,z)\bee^{\epsilon}\bee^{-\epsilon}\mathcal{Y}_{\lambda+\epsilon,-\lambda}( \bee^{\epsilon}m,w)n}\frac{\dd {z}}{2\pi i}\nonumber\\
  &=-\epsilon c_{\lambda,\epsilon}\oint_{w,0}f(z)\pair{o}{\Delta_1(g;z;\epsilon)\bee^{-\epsilon}\mathcal{Y}_{\lambda+\epsilon,-\lambda}( \bee^{\epsilon}m,w)n}\frac{\dd {z}}{2\pi i}.
\end{align}
Thus by \eqref{eq:Deltarel-e2}, we see that the identity \eqref{eq:Deltarel-e} holds.
\end{proof}

\begin{proof}[Proof of \cref{thm:logint}]
  The operator  $\mathcal{Y}^{log}_{\lambda,-\lambda}$ clearly
  satisfies the truncation condition and
  the convergence condition follows from Lemma \ref{thm:conv}.
  Therefore, all that remains to be shown is
  that $\mathcal{Y}^{log}_{\lambda,-\lambda}$ satisfies the Cauchy-Jacobi
  identity. 
  By \cref{thm:Deltarel},
  for $m\in E_{\lambda;1,1}$, $n\in E_{-\lambda;1,1}$, $o \in (\mathbb{H}^1)'$, $g\in \mmsl{u}{v}$ and $f(z)\in\mathbb{C}[z,z^{-1},(z-w)^{-1}]$, $\mathcal{Y}^{log}_{\lambda,-\lambda}$ satisfies
  \begin{align}
    &\oint_{w,0}f(z)\pair{o}{\widetilde{Y}_{c_{\lambda,0}}(g,z)\mathcal{Y}^{log}_{\lambda,-\lambda}( m,w)n}\frac{\dd {z}}{2\pi i}\nonumber\\
    &=\oint_{w}f(z)\pair{o}{\mathcal{Y}^{log}_{\lambda,-\lambda}( Y(g,z-w)m,w)n}\frac{\dd {z}}{2\pi i}
    +\oint_{0}f(z)\pair{o}{\mathcal{Y}^{log}_{\lambda,-\lambda}( m,w)Y(g,z)n}\frac{\dd {z}}{2\pi i},
    \label{eq:CJiden}
  \end{align}
  where $\widetilde{Y}_{c_{\lambda,0}}(g,z)$ is the $\mmsl{u}{v}$-action on $\mathbb{H}$ defined by
\begin{equation}
  \widetilde{Y}_{c_{\lambda,-\lambda}}(g,z)=
  \begin{cases}
    Y(g,z)-c_{\lambda,0}\Delta_1(g;z) &{\rm on}\ \mathbb{H}^0,\\
    Y(g,z) &{\rm on}\ \mathbb{H}^1.
  \end{cases}
\end{equation}
Therefore, from \eqref{eq:CJiden}, $\mathcal{Y}^{log}_{\lambda,-\lambda}$ satisfies the Cauchy-Jacobi identity.
\end{proof}

\begin{remark}
  Because the action of \(L_0\) on \(\HH\) has Jordan blocks of size 2, the
  \(P(w)\) intertwining operator $\mathcal{Y}^{log}_{\lambda,-\lambda}$
  has first order logarithmic terms (that is, \(\log z\) appears but
  \(\brac*{\log z}^n\), \(n\ge2\) does not)
  once the complex variable \(w\) is replaced by a
  formal variable \(z\) by using the identity \cite[Eq. (4.17)]{HLZ3}, that is,
  \begin{equation}
    \mathcal{Y}^{\log}_{\lambda,-\lambda}(v_1,z)v_2=
    y^{L(0)}z^{L(0)}\left.\mathcal{Y}^{\log}_{\lambda,-\lambda}(y^{-L(0)}z^{-L(0)}v_{1}\otimes
      y^{-L(0)}z^{-L(0)}v_{2})\right|_{y=e^{-\log w}},
  \end{equation}
  where \(z^{L_{0}}\) is defined to the formal series
  \begin{equation}
    z^{L(0)}=z^n\sum_{i\in \mathbb{Z}_{\geq 0}}\frac{(L(0)-n)^iv}{i!}(\log z)^i.
  \end{equation}
\end{remark}

\subsection{The rigidity of simple \texorpdfstring{\(\slE{\mu}{1}{1}\)}{E11} modules}
\label{sec:11rigidity}
In this subsection, following the same principles as in \cite[Prop 5.13]{AW},
we will show that the $\slE{\sqbrac*{\lambda}}{1}{1}$, \(\lambda\in \RR\setminus
\set{\pm 
  \lambda_{1,1}+\ZZ}\), \(\sqbrac*{\lambda}=\lambda+2\ZZ\) are rigid in
$\wtmod{\mmsl{u}{v}}$, with rigid duals \(\slE{\sqbrac*{-\lambda}}{1}{1}\) (we
use explicit representatives for the coset \(\sqbrac*{\lambda}\), because it
will be convenient in calculations below).
Note that since the category of weight modules is braided, there is no need to
distinguish left and right duals.
The integration cycles required here are more complicated than in
\cite{AW}, hence the details of the arguments are considerably more subtle.
The methods of \cite[Sec 3.2]{hgeom} on twisted cycles and twisted homology groups will be crucial.

Let $w_1,w_2\in \mathbb{R}_{>0}$ be real numbers satisfying $w_1>w_1-w_2>0$
and let $E=\slE{\sqbrac*{\lambda}}{1}{1}$,
$E^\vee=\slE{\sqbrac*{-\lambda}}{1}{1}$.
By \cite[Lem 4.2.1 and Cor 4.2.2]{CMcY} a sufficient condition for
rigidity is the existence of two morphisms 
$i_{\slE{\lambda}{1}{1}} \colon L_1 \to \slE{\lambda}{1}{1} \fuse_w
\slE{-\lambda}{1}{1}$ and $e_{\slE{\lambda}{1}{1}} \colon \slE{-\lambda}{1}{1}
\fuse_w \slE{\lambda}{1}{1} \to L_1$, where $L_1$ is the unit of the
  category (the \voa{} as a module over itself), such that the composition
\begin{equation}
  R_E = E\xrightarrow{l^{-1}} L_1\fuse E\xrightarrow{i_{E}\fuse {\rm 1}} (E\fuse_{w_1} E^{\vee})\fuse_{w_2} E\xrightarrow{\mathcal{A}^{-1}}E\fuse_{w_1}(E^{\vee}\fuse_{w_2} E)\xrightarrow{{\rm 1}\fuse e_{E}}E\fuse L_1\xrightarrow{r} E
  \label{eq:zigzag}
\end{equation}
is non-zero. Here, $l$ and $r$ are the left and right unitors, respectively,
and $\mathcal{A}$ the associator of the category. Since
\(\slE{\sqbrac*{\lambda}}{1}{1}\) is simple \(R_E\) will be a
scalar multiple of the identity and hence the pair
\((e_{\slE{\lambda}{1}{1}},i_{\slE{\lambda}{1}{1}})\) define evaluation and
coevaluation morphisms after dividing one of them by this scalar.
We begin by constructing candidates for these evaluation and coevaluation maps. For $e_E$, we take,
\begin{align}
  e_E=\mathcal{Q}^{[2]}*\mathcal{Y}_{-\lambda,\lambda},
  \label{eq:ev}
\end{align}
as noted after \cref{thm:Int-Q2}.
To characterise the coevaluation morphism we prepare the notation
\begin{equation}
\overline{\theta_\lambda}=\ket{\frac{v\lambda}{u-2v}(a-b)-b},\qquad
\theta_{\lambda}=\vac\cten\overline{\theta_\lambda}\in
\slE{\sqbrac*{\lambda}}{1}{1},
\end{equation}
where \(\vac\) is the minimal module vacuum vector.
Note that from \cref{thm:logint} we know that $E\fuse E^{\vee}$ has a
direct summand isomorphic to $(\mathbb{H},\widetilde{Y}_{c_{\lambda,0}})$, the
projective cover of \(L_1\). We then characterise the coevaluation map $i_E$
by the image of the vacuum vector \(\Omega\in L_1\).
\begin{multline}
  i_E \colon \Omega
  \rightarrow v_{1,1}\cten\ket{0}\xrightarrow{(\mathcal{Q}^{[2]})^{-1}}
  v_{1,1}\cten\ket{-2b}\\\rightarrow
  \mathcal{Y}_{\lambda,-\lambda}(\ket{\theta_{\lambda}},w)\ket{\theta_{-\lambda}} 
  \xrightarrow{\mathcal{Q}^{[2]}*}
  \mathcal{Q}^{[2]}*\mathcal{Y}_{\lambda,-\lambda}(\ket{\theta_{\lambda}},w)\ket{\theta_{-\lambda}}, 
\label{eq:coev}
\end{multline}
where the first arrow is the inclusion of $L_1$ into $\sigma
\slE{\lambda_{1,1}}{1}{1}\subset (\mathbb{H},\widetilde{Y})$,
$(\mathcal{Q}^{[2]})^{-1}$ denotes picking preimages of
$\mathcal{Q}^{[2]}$. Note that the ambiguity of picking preimages of
$\mathcal{Q}^{[2]}$ is undone by applying $\mathcal{Q}^{[2]}*$
and hence the map is a well-defined morphism (this is a consequence of \eqref{eq:upto2}).
The morphism
\(R_E\) is thus non-zero (and hence \(E\) is rigid) if and only if
\(R_E(\theta_\lambda)\) is non-zero. This in turn is the case if
and only if the matrix element \(\pair{\mu_\lambda}{R_E(\theta_\lambda)}\)
is non-zero, where
\begin{equation}
  \overline{\mu_{\lambda}}=
  \bra{\frac{v\lambda}{u-2v}(a-b)+b},\qquad
  \mu_{\lambda}=\vac^\ast\cten\overline{\mu_{\lambda}},
\end{equation}
and \(\vac^\ast\) is the vector dual to the vacuum vector \(\vac\) of the
Virasoro minimal model. The remainder of this section will be dedicated to
showing that \(\pair{\mu_\lambda}{R_E(\theta_\lambda)}\)
is indeed non-zero.

From the characterisations \eqref{eq:ev} and \eqref{eq:coev} above of \(e_E\)
and \(i_E\), we see that
\begin{align}
  \pair{\mu_\lambda}{R_E(\theta_\lambda)}&=\int\int
  \pair{\mu_\lambda}{\scr^-(z_1)\scr^+(z_2)\mathcal{Y}\brac*{\theta_\lambda;w_1}\scr^-(z_3)\scr^+(z_4)
    \mathcal{Y}\brac*{\theta_{-\lambda};w_2}\theta_{\lambda}}\dd
                      z_1\cdots\dd z_4\nonumber\\
                    &=\int\int 
                      \pair{\vac^\ast}{I^-_{1,2}\brac*{v_{1,2};z_1}I^+_{1,2}\brac*{v_{1,2};z_2}I^-_{1,2}\brac*{v_{1,2};z_1}I^+_{1,2}\brac*{v_{1,2};z_4}\vac}\nonumber\\
                    &\qquad
                      \cdot\pair{\overline{\mu_\lambda}}{\mathcal{Y}\brac*{b;z_1}\mathcal{Y}\brac*{b;z_2}\mathcal{Y}\brac*{b;z_3}\mathcal{Y}\brac*{b;z_4}\mathcal{Y}\brac*{\overline{\theta_\lambda};w_1}\mathcal{Y}\brac*{\overline{\theta_{-\lambda}};w_2}\overline{\theta_\lambda}}\dd
                      z_1\cdots\dd z_4,
  \label{eq:matel}
\end{align}
where, for the momement, we deliberately suppress the integration
cycles to first focus on the formula for the integrand.
The currents in the \(z_1,z_2,w_1\) variables result from the application of the coevaluation morphism \(i_E\) in \eqref{eq:zigzag} to the vacuum vector \(\Omega\), that is, the final term in \eqref{eq:coev}. The evaluation morphism then needs to be applied to the right factor in the image of coevaluation and the right most tensor factor in the fourth entry of the sequence \eqref{eq:zigzag}, which results in the currents in the \(z_3,z_4,w_2\) variables. The morphism $R_E$ is evaluated by first integrating the $z_1,z_2$ variables over the twisted cycle $[\Delta_{w_1,w_2}]$ (the translate of
\(\sqbrac*{\Delta_{w_1-w_2,0}}\) by \(w_2\)) , and then integrating over the $z_3,z_4$ variables over the twisted cycle $[\Delta_{w_2,0}]$.
Since the integrand of $R_E$ has singularities at $z_i=z_j$ $(i=1,2,\; j=3,4)$, these two successive integrations are carried out by expanding the singular terms at $z_i=z_j$ $(i=1,2,\; j=3,4)$ into power series. To make it possible to perform these integrations in the presence of such singularities, we will replace the integration cycles \([\Delta_{w_1,w_2}]\), \([\Delta_{w_2,0}]\) by equivalent locally finite cycles below. The presence of singularities makes convergence of the integral a concern, however, since the evaluation, coevaluation and associativity homomorphisms are all well-defined, we know a priori that the integral must converge. As we shall see, the precise form of the singularities (in particular, the hypergeometric function singularity determined in \cref{thm:corrVir}) will be crucial to being able to conclude that $R_E$ is non-zero in \cref{thm:nonzeromatel}.

The free field part of the matrix element above is easily evaluated to
\begin{align}
  &\pair{\overline{\mu_\lambda}}{\mathcal{Y}\brac*{b;z_1}\mathcal{Y}\brac*{b;z_2}\mathcal{Y}\brac*{\overline{\theta_\lambda};w_1}\mathcal{Y}\brac*{b;z_3}\mathcal{Y}\brac*{b;z_4}\mathcal{Y}\brac*{\overline{\theta_{-\lambda}};w_2}\overline{\theta_\lambda}}\nonumber\\
  &\quad=(-1)^{\frac{u+v\lambda}{v}-2}
  \prod_{1\le i<j\le 4}\brac*{z_i-z_j}^{1-\frac{u}{2v}}
    \prod_{i=1}^4\brac*{z_i-w_1}^{\frac{u+v\lambda}{2v}-1}\prod_{i=1}^4\brac*{z_i-w_2}^{\frac{u-v\lambda}{2v}-1}\prod_{i=1}^4z_i^{\frac{u+v\lambda}{2v}-1}\nonumber\\
  &\qquad\quad\cdot\brac*{w_1-w_2}^{1-\frac{u}{2v}}w_1^{\frac{u+v\lambda}{2v}-1}w_2^{1-\frac{u}{2v}}.
\end{align}
The phase \((-1)^{\frac{u+v\lambda}{v}-2}\) is due to the variable \(w_1\) being written to the right in the differences \(z_3-w_1\), \(z_4-w_1\) appearing in the second product, while the \(w_1\) current in the matrix element is to the left of the \(z_3\) and \(z_4\) currents. We will suppress this phase below, as it has no bearing on whether or not \(R_E\) vanishes.
To compute the \(4\)-point function involving the Virasoro intertwining
operators, we need to solve the corresponding BPZ equation \cite{BPZ}.

\begin{lemma}
The intertwining operators \(I^\pm_{1,2}\) can be normalised such that
\begin{align}
  \Psi(z)&=\pair{\vac^\ast}{I^-_{1,2}\brac*{v_{1,2};z_1}I^+_{1,2}\brac*{v_{1,2};z_2}I^-_{1,2}\brac*{v_{1,2};z_{3}}I^+_{1,2}\brac*{v_{1,2};z_4}\vac}\nonumber\\
  &=\brac*{\frac{z_{1,2}z_{3,4}z_{1,4}z_{2,3}}{z_{1,3} z_{2,4}}}^{1-\frac{3u}{2v}}
  \tensor[_2]{F}{_{1}}\brac*{2-\tfrac{3u}{v},1-\tfrac{u}{v},2-\tfrac{2u}{v};\tfrac{z_{1,4}
      z_{2,3}}{z_{1,3} z_{2,4}}},
  \label{eq:Vir4ptfn}
\end{align}
for $|z_1|>|z_2|>|z_3|>|z_4|$, where 
$ z_{ij}=z_{i}-z_{j}$.
\label{thm:corrVir}
\end{lemma}
\begin{proof}
  Note that each of the four intertwining operators in the \(4\)-point function \eqref{eq:Vir4ptfn} is non-zero and goes from a pair of simple highest weight modules to a simple highest weight module. Their leading term will therefore always be the highest weight vector of their codomain (if not they would be identically equal to \(0\)) and hence the \(4\)-point function is non-zero. So what needs to be determined is what form the correlation function takes. As we recall below, a \(4\)-point function where all four insertions are \(v_{1,2}\) satisfies a degree \(2\) linear differential equation, which hence has a \(2\)-dimensional space of solutions. The solution above corresponds to the intertwining operator types we have chosen. An additional solution can be found by replacing the intertwining operators for the variables \(z_2\) and \(z_3\) by ones of types \(\binom{\mmVM{1}{2}}{\mmVM{1}{2},\ \mmVM{1}{3}}\) and \(\binom{\mmVM{1}{3}}{\mmVM{1}{2},\ \mmVM{1}{2}}\), respectively.

  By global conformal covariance any \(4\)-point function of conformal highest
  weight (also known as primary) vectors can be factorised into a part depending
  on the differences of variables \(z_i-z_j\) and a part depending only on the
  cross ratio
  \begin{align}
    x=\frac{\brac*{z_1-z_2}\brac*{z_3-z_4}}{\brac*{z_1-z_3}\brac*{z_2-z_4}}.
  \end{align}
  For this particular \(4\)-point function this yields
  \begin{align}
    \Psi(z)=z_{1,3}^{-2h_{1,2}}z_{2,4}^{-h_{1,2}}G(x),
  \end{align}
  where \(h_{1,2}=\frac{3u}{4v}-\frac{1}{2}\) is the conformal weight of the
  highest weight vector \(v_{1,2}\). Using the Ward identities to convert
  Virasoro generators into differential operators (see \cite[Sec 8.3]{CFTbook} for an overview), the singular vector
  relation \(\brac*{L_{-1}^2-\frac{u}{v}L_{-2}}v_{1,2}=0\) then implies that \(G(x)\)
  satisfies the differential equation
  \begin{align}
    \sqbrac*{\partial_x^2+\frac{u}{v}\brac*{\frac{1}{x-1}+\frac{1}{x}}\partial_x-h_{1,2}\frac{u}{v}\brac*{\frac{1}{x^2}+\frac{1}{\brac{x-1}^2}}-2h_{1,2}\frac{u}{v}\brac*{\frac{1}{x}-\frac{1}{x-1}}}G(x)=0.
    \label{eq:BPZeq}
  \end{align}
  This is a Fuchsian differential equation (see, for example, \cite[Sec 2.1]{GtoP} for an overview on solving such equations). The characteristic
  exponents at \(x=0\) and at \(x=1\) for this differential equation are the
  roots of the equation
  \begin{align}
    \rho^2+\brac*{\frac{u}{v}-1}\rho -h_{1,2}\frac{u}{v}=0,
    \label{eq:01exp}
  \end{align}
  that is \(\rho\in\set{1-\frac{3u}{2v},\frac{u}{2v}}\). The characteristic
  exponents at \(x=\infty\) are the roots of
  \begin{align}
    \gamma^2+\brac{1-\frac{2u}{v}}\gamma=0
  \end{align}
  and hence \(\gamma\in \set{0,\frac{2u}{v}-1}\). For each choice of root \(\rho\) in
  \eqref{eq:01exp} there is a solution to \eqref{eq:BPZeq} given by a series
  expansion in \(x^{\rho+n}\), \(n\in \ZZ_{\ge0}\). Further, the types of the
  intertwining operators \(I^+_{1,2}\) and \(I^-_{1,2}\) are
  \(\binom{\mmVM{1}{2}}{\mmVM{1}{2},\ \mmVM{1}{1}}\) and
  \(\binom{\mmVM{1}{1}}{\mmVM{1}{2},\ \mmVM{1}{2}}\) respectively, so
  \(I^-_{1,2}\brac*{v_{1,2};z_3}I^+_{1,2}\brac*{v_{1,2};z_4}\) can be expanded as a
  series in \(\brac*{z_3-z_4}^{-2h_{1,2}+n}\), \(n\in \ZZ_{\ge0}\), which
  implies that we need to pick the root \(1-\frac{3u}{2v}=-2h_{1,2}\) as our
  solution to \eqref{eq:01exp} (the root \(\frac{u}{2v}=h_{1,3}-2h_{1,2}\) corresponds to replacing the intertwining operators in the variables \(z_2\) and \(z_3\) by ones of types \(\binom{\mmVM{1}{2}}{\mmVM{1}{2},\ \mmVM{1}{3}}\) and \(\binom{\mmVM{1}{3}}{\mmVM{1}{2},\ \mmVM{1}{2}}\), respectively). Hence
  \begin{align}
    G(x)=x^{1-\frac{3u}{2v}}\brac*{x-1}^{1-\frac{3u}{2v}}\tensor[_2]{F}{_{1}}\brac*{2-3\tfrac{u}{v},1-\tfrac{u}{v},2-2\tfrac{u}{v};x},
  \end{align}
  which proves \eqref{eq:Vir4ptfn}.
\end{proof}

To more easily group the various factors appearing in integrand defining
\(\pair{\mu_\lambda}{R_E(\theta_\lambda)}\), we introduce the notation
\begin{align}
  H(z,w)&=w_1^{-2\lambda}z_1^{\frac{u+v\lambda}{2v}-1}z_2^{\frac{u+v\lambda}{2v}-1}
                              \brac*{z_1-z_{3}}^{\frac{u}{v}}\brac*{z_2-z_{4}}^{\frac{u}{v}}\brac*{z_1-
                          z_{4}}^{2-2\frac{u}{v}}\brac*{z_2-z_{3}}^{2-2\frac{u}{v}}\nonumber\\
 &\qquad\cdot\brac*{w_1- z_3}^{\frac{u+v\lambda}{2v}-1}\brac*{w_1-
   z_4}^{\frac{u+v\lambda}{2v}-1}\nonumber\\
  &=w_1^{-2\lambda}z_1^{1+\frac{v\lambda-u}{2v}}z_2^{1+\frac{v\lambda-u}{2v}}
                              \brac*{1-\frac{z_{3}}{z_1}}^{\frac{u}{v}}\brac*{1-\frac{z_{4}}{z_2}}^{\frac{u}{v}}\brac*{1-
                          \frac{z_{4}}{z_1}}^{2-2\frac{u}{v}}\brac*{1-\frac{z_{3}}{z_2}}^{2-2\frac{u}{v}}\nonumber\\
 &\qquad\cdot\brac*{w_1- z_3}^{\frac{u+v\lambda}{2v}-1}\brac*{w_1-
   z_4}^{\frac{u+v\lambda}{2v}-1},\nonumber\\
  \mathcal{G}^\lambda_{u_1,u_2}(y_1,y_2)&=(y_1-y_2)^{2-2\frac{u}{v}}\prod_{i=1}^2\brac*{y_i-u_1}^{\frac{u+v\lambda}{2v}-1}\prod_{i=1}^2\brac*{y_i-u_2}^{\frac{u-v\lambda}{2v}-1},\nonumber\\
  f(w_1,w_2)&=w_1^{1-\lambda-\frac{u}{2v}}w_2^{1-\frac{u}{2v}}\brac*{w_1-w_2}^{1-\frac{u}{2v}}.
  \label{eq:auxfns}
\end{align}
Then we can write the integrand defining \(\pair{\mu_\lambda}{R_E(\theta_\lambda)}\) as
\begin{multline}
  M(\lambda,w,z)=
  f(w_1,w_2)H(z,w)\\
  \tensor[_2]{F}{_{1}}\brac*{2-3\tfrac{u}{v},1-\tfrac{u}{v},2-2\tfrac{u}{v};\tfrac{\brac*{z_1-z_2}\brac*{z_3-z_4}}{\brac*{z_1-z_3}\brac*{z_2-z_4}}}
  \mathcal{G}^\lambda_{w_1,w_2}(z_1,z_2)
  \mathcal{G}^{-\lambda}_{w_2,0}(z_3,z_4).
\end{multline}
The factor of \(w_1^{-2\lambda}\) in the definition of \(H\)
is there to simplify the
following expansion formula.
\begin{align}
  &H(z,w)\ \tensor[_2]{F}{_{1}}\brac*{2-\tfrac{3u}{v},1-\tfrac{u}{v},2-\tfrac{2u}{v};\tfrac{\brac*{z_1-z_2}\brac*{z_3-z_4}}{\brac*{z_1-z_3}\brac*{z_2-z_4}}} \nonumber\\
   &\qquad
     = w_2^{\frac{v\lambda-u}{v}+2}w_1^{\frac{u+v\lambda}{v}-2}\sum_i H_{1,2}^{(i)}(z_1,z_2,w_2)H_{3,4}^{(i)}(z_3,z_4,w_1),\nonumber\\
&H_{1,2}^{(i)}\in
\CC\sqbrac*{w_2^{\pm1}}\cten \CC\powser{z_i-w_2\st
    i=1,2}^{\mathfrak{S}},\nonumber\\
&H_{3,4}^{(i)}\in\CC\sqbrac*{w_1^{\pm1}}\cten
\CC\sqbrac*{z_3,z_4}^{\mathfrak{S}},
 \label{eq:Hexpansion}
\end{align}
 where the superscript \(\mathfrak{S}\)
indicates invariance under permuting \(z_1\) with \(z_2\) or \(z_3\) with
\(z_4\), respectively, and where the factor \(w_2^{\frac{v\lambda-u}{v}+2}w_1^{\frac{u+v\lambda}{v}-2}\) comes from
\begin{align}
z_1^{1+\frac{v\lambda-u}{2v}}z_2^{1+\frac{v\lambda-u}{2v}}\brac*{w_1- z_3}^{\frac{u+v\lambda}{2v}-1}\brac*{w_1-
   z_4}^{\frac{u+v\lambda}{2v}-1}
                =w_2^{\frac{v\lambda-u}{v}+2}w_1^{\frac{u+v\lambda}{v}-2}\brac*{1+\frac{z_1-w_2}{w_2}}^{1+\frac{v\lambda-u}{2v}}&\nonumber\\
  \cdot\brac*{1+\frac{z_2-w_2}{w_2}}^{1+\frac{v\lambda-u}{2v}}\brac*{1- \frac{z_3}{w_1}}^{\frac{u+v\lambda}{2v}-1}\brac*{1-
   \frac{z_4}{w_1}}^{\frac{u+v\lambda}{2v}-1}.&
\end{align}
The non-integral exponents of the factors mixing the \(z_1,z_2\) variables
with the \(z_3,z_4\) make the integral \eqref{eq:matel}
very difficult to evaluate head on. So instead we will construct regularised
cycles which are homologous (in appropriate homology groups) to
\(\sqbrac*{\Delta_{w_1,w_2}}\), \(\sqbrac*{\Delta_{w_2,0}}\).
These
will allow us to decompose the integral into more manageable
parts. Specifically, they will allow us to pull the sum appearing in
\eqref{eq:Hexpansion} out of the integral.
  \begin{defn}
    Let \(\sqbrac*{\Delta^0_{w_1,w_2}}\) and \(\sqbrac*{\Delta^0_{w_2,0}}\) be
    the surfaces constructed in \cref{fig:cycconstr}.
    \begin{figure}[ht]
      \centering
  \def\svgwidth{5cm} 
\begingroup%
  \makeatletter%
  \providecommand\color[2][]{%
    \errmessage{(Inkscape) Color is used for the text in Inkscape, but the package 'color.sty' is not loaded}%
    \renewcommand\color[2][]{}%
  }%
  \providecommand\transparent[1]{%
    \errmessage{(Inkscape) Transparency is used (non-zero) for the text in Inkscape, but the package 'transparent.sty' is not loaded}%
    \renewcommand\transparent[1]{}%
  }%
  \providecommand\rotatebox[2]{#2}%
  \newcommand*\fsize{\dimexpr\f@size pt\relax}%
  \newcommand*\lineheight[1]{\fontsize{\fsize}{#1\fsize}\selectfont}%
  \ifx\svgwidth\undefined%
    \setlength{\unitlength}{226.77165354bp}%
    \ifx\svgscale\undefined%
      \relax%
    \else%
      \setlength{\unitlength}{\unitlength * \real{\svgscale}}%
    \fi%
  \else%
    \setlength{\unitlength}{\svgwidth}%
  \fi%
  \global\let\svgwidth\undefined%
  \global\let\svgscale\undefined%
  \makeatother%
  \begin{picture}(1,1)%
    \lineheight{1}%
    \setlength\tabcolsep{0pt}%
    \put(0,0){\includegraphics[width=\unitlength,page=1]{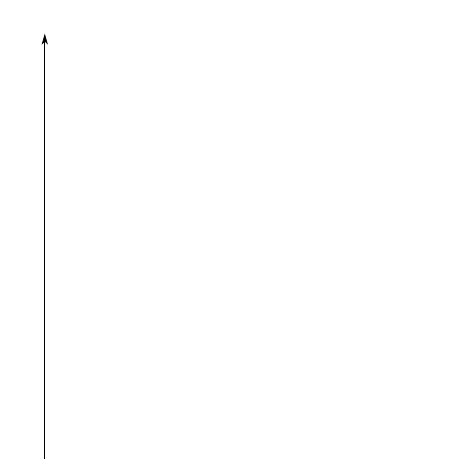}}%
    \put(0.0691626,0.94586346){\color[rgb]{0,0,0}\makebox(0,0)[lt]{\lineheight{1.25}\smash{\begin{tabular}[t]{l}$z_2$\end{tabular}}}}%
    \put(0.91450739,0.09852486){\color[rgb]{0,0,0}\makebox(0,0)[lt]{\lineheight{1.25}\smash{\begin{tabular}[t]{l}$z_1$\end{tabular}}}}%
    \put(0,0){\includegraphics[width=\unitlength,page=2]{simplexSbu.pdf}}%
    \put(0.10443153,0.00342207){\color[rgb]{0,0,0}\makebox(0,0)[lt]{\lineheight{1.25}\smash{\begin{tabular}[t]{l}$(w_2,w_2)$\end{tabular}}}}%
    \put(0.69005065,0.00342207){\color[rgb]{0,0,0}\makebox(0,0)[lt]{\lineheight{1.25}\smash{\begin{tabular}[t]{l}$(w_1,w_2)$\end{tabular}}}}%
    \put(0.68220878,0.83616306){\color[rgb]{0,0,0}\makebox(0,0)[lt]{\lineheight{1.25}\smash{\begin{tabular}[t]{l}$(w_1,w_1)$\end{tabular}}}}%
    \put(0.11968679,0.79591982){\color[rgb]{0,0,0}\makebox(0,0)[lt]{\lineheight{1.25}\smash{\begin{tabular}[t]{l}$(w_2,w_1)$\end{tabular}}}}%
  \end{picture}%
\endgroup%

  \def\svgwidth{5cm} 
\begingroup%
  \makeatletter%
  \providecommand\color[2][]{%
    \errmessage{(Inkscape) Color is used for the text in Inkscape, but the package 'color.sty' is not loaded}%
    \renewcommand\color[2][]{}%
  }%
  \providecommand\transparent[1]{%
    \errmessage{(Inkscape) Transparency is used (non-zero) for the text in Inkscape, but the package 'transparent.sty' is not loaded}%
    \renewcommand\transparent[1]{}%
  }%
  \providecommand\rotatebox[2]{#2}%
  \newcommand*\fsize{\dimexpr\f@size pt\relax}%
  \newcommand*\lineheight[1]{\fontsize{\fsize}{#1\fsize}\selectfont}%
  \ifx\svgwidth\undefined%
    \setlength{\unitlength}{226.77165354bp}%
    \ifx\svgscale\undefined%
      \relax%
    \else%
      \setlength{\unitlength}{\unitlength * \real{\svgscale}}%
    \fi%
  \else%
    \setlength{\unitlength}{\svgwidth}%
  \fi%
  \global\let\svgwidth\undefined%
  \global\let\svgscale\undefined%
  \makeatother%
  \begin{picture}(1,1)%
    \lineheight{1}%
    \setlength\tabcolsep{0pt}%
    \put(0,0){\includegraphics[width=\unitlength,page=1]{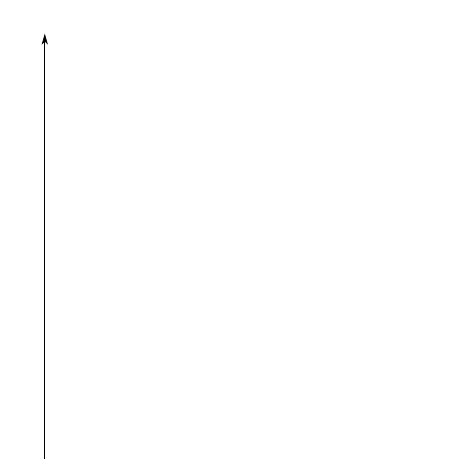}}%
    \put(0.0691626,0.94586346){\color[rgb]{0,0,0}\makebox(0,0)[lt]{\lineheight{1.25}\smash{\begin{tabular}[t]{l}$z_2$\end{tabular}}}}%
    \put(0.91450739,0.09852486){\color[rgb]{0,0,0}\makebox(0,0)[lt]{\lineheight{1.25}\smash{\begin{tabular}[t]{l}$z_1$\end{tabular}}}}%
    \put(0,0){\includegraphics[width=\unitlength,page=2]{simplexSball.pdf}}%
    \put(0.30275166,0.0006855){\color[rgb]{0,0,0}\makebox(0,0)[lt]{\lineheight{1.25}\smash{\begin{tabular}[t]{l}$\partial B(w_2)$\end{tabular}}}}%
    \put(0.11968679,0.79591982){\color[rgb]{0,0,0}\makebox(0,0)[lt]{\lineheight{1.25}\smash{\begin{tabular}[t]{l}$(w_2,w_1)$\end{tabular}}}}%
    \put(0.68053832,0.83616306){\color[rgb]{0,0,0}\makebox(0,0)[lt]{\lineheight{1.25}\smash{\begin{tabular}[t]{l}$(w_1,w_1)$\end{tabular}}}}%
    \put(0.69005065,0.00342207){\color[rgb]{0,0,0}\makebox(0,0)[lt]{\lineheight{1.25}\smash{\begin{tabular}[t]{l}$(w_1,w_2)$\end{tabular}}}}%
    \put(0,0){\includegraphics[width=\unitlength,page=3]{simplexSball.pdf}}%
  \end{picture}%
\endgroup%

  \def\svgwidth{5cm}
  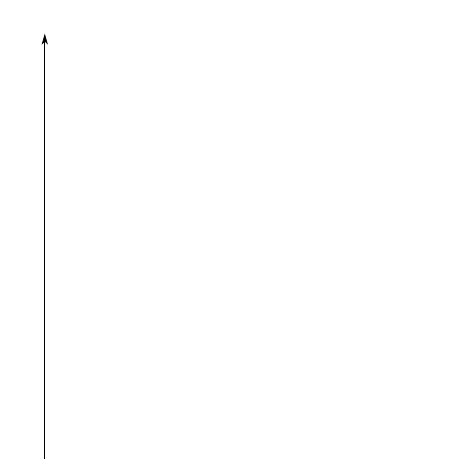\\
  \def\svgwidth{5cm} 
\begingroup%
  \makeatletter%
  \providecommand\color[2][]{%
    \errmessage{(Inkscape) Color is used for the text in Inkscape, but the package 'color.sty' is not loaded}%
    \renewcommand\color[2][]{}%
  }%
  \providecommand\transparent[1]{%
    \errmessage{(Inkscape) Transparency is used (non-zero) for the text in Inkscape, but the package 'transparent.sty' is not loaded}%
    \renewcommand\transparent[1]{}%
  }%
  \providecommand\rotatebox[2]{#2}%
  \newcommand*\fsize{\dimexpr\f@size pt\relax}%
  \newcommand*\lineheight[1]{\fontsize{\fsize}{#1\fsize}\selectfont}%
  \ifx\svgwidth\undefined%
    \setlength{\unitlength}{226.77165354bp}%
    \ifx\svgscale\undefined%
      \relax%
    \else%
      \setlength{\unitlength}{\unitlength * \real{\svgscale}}%
    \fi%
  \else%
    \setlength{\unitlength}{\svgwidth}%
  \fi%
  \global\let\svgwidth\undefined%
  \global\let\svgscale\undefined%
  \makeatother%
  \begin{picture}(1,1)%
    \lineheight{1}%
    \setlength\tabcolsep{0pt}%
    \put(0.0691626,0.94586346){\color[rgb]{0,0,0}\makebox(0,0)[lt]{\lineheight{1.25}\smash{\begin{tabular}[t]{l}$z_2$\end{tabular}}}}%
    \put(0.91450739,0.09852486){\color[rgb]{0,0,0}\makebox(0,0)[lt]{\lineheight{1.25}\smash{\begin{tabular}[t]{l}$z_1$\end{tabular}}}}%
    \put(0,0){\includegraphics[width=\unitlength,page=1]{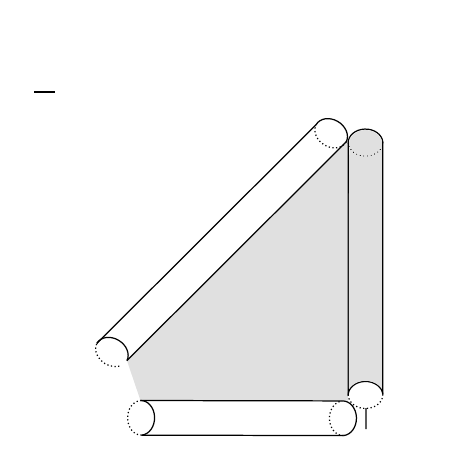}}%
    \put(0.11968679,0.79591982){\color[rgb]{0,0,0}\makebox(0,0)[lt]{\lineheight{1.25}\smash{\begin{tabular}[t]{l}$(w_2,w_1)$\end{tabular}}}}%
    \put(0,0){\includegraphics[width=\unitlength,page=2]{simplexScones.pdf}}%
    \put(0.10443153,0.00342207){\color[rgb]{0,0,0}\makebox(0,0)[lt]{\lineheight{1.25}\smash{\begin{tabular}[t]{l}$(w_2,w_2)$\end{tabular}}}}%
    \put(0.69005065,0.00342207){\color[rgb]{0,0,0}\makebox(0,0)[lt]{\lineheight{1.25}\smash{\begin{tabular}[t]{l}$(w_1,w_2)$\end{tabular}}}}%
    \put(0.68053832,0.83616306){\color[rgb]{0,0,0}\makebox(0,0)[lt]{\lineheight{1.25}\smash{\begin{tabular}[t]{l}$(w_1,w_1)$\end{tabular}}}}%
    \put(0,0){\includegraphics[width=\unitlength,page=3]{simplexScones.pdf}}%
  \end{picture}%
\endgroup%

  \def\svgwidth{5cm} 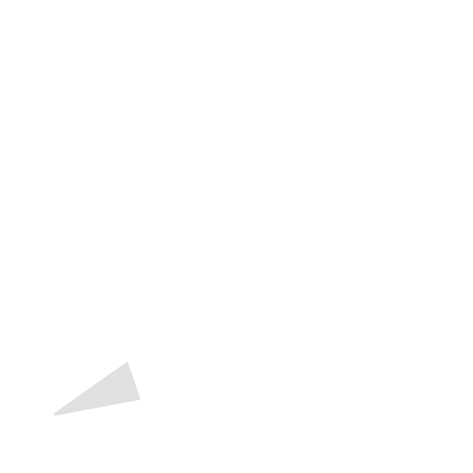
  \def\svgwidth{5cm} 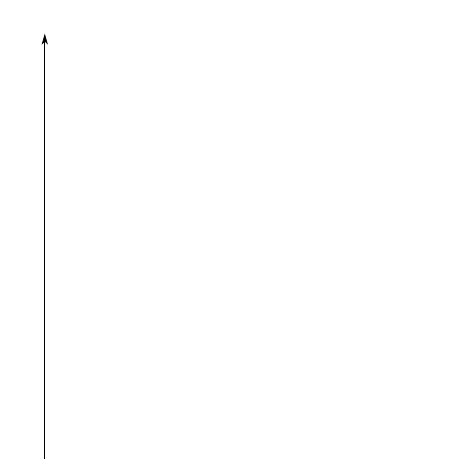
  \caption{Constructing the regularised cycles
    \(\sqbrac*{\Delta^0_{w_1,w_2}}\) and \(\sqbrac*{\Delta^0_{w_2,0}}\):
      Starting from the top left with the standard twisted cycle built from
      the \(2\)-simplex \(\set{w_1>|z_1|>|z_2|>w_2}\) by replacing the edges by tubular neighbourhoods. Note
      that the tubes do not intersect. This is just an artifact of presenting
      the image in two real (as opposed to complex) dimensions. In
      the second image \(B(w_2)=\set{4|z_1-w_2|^2+4|z_2-w_2|^2\le
        \brac*{w_1+w_2}^2}\) denotes a ball centred at
      \((w_2,w_2)\) with radius 
      \(\frac{w_1^2+w_2^2}{4}\). We cut away all of the twisted cycle that is
      within the ball leaving two circles on the tubes connected by a line and
      obtain the third image (see \cite[Sec 5.4]{TK} and
      \cite[Sec 3.2.4]{hgeom}).
      The circles are then connected to \((w_2,w_2)\)
      by a pair of cones, whose orientation matches that of the circle they
      are attached to. This leaves a wedge shaped hole and is depicted in
      the fourth image. The hole is filled in to give
      \(\sqbrac*{\Delta_{w_1,w_2}^0}\) in the fifth image. The same procedure
      but with \((w_1,w_2)\) replaced by \((w_2,0)\) yields
      \(\sqbrac*{\Delta^0_{w_2,0}}\).}
  \label{fig:cycconstr}
\end{figure}
\label{def:twistcyc}
  \end{defn}
\begin{lemma}
Consider the surfaces \(\sqbrac*{\Delta^0_{w_1,w_2}}\) and
\(\sqbrac*{\Delta^0_{w_2,0}}\) constructed in \cref{def:twistcyc}. Then 
\(
[\Delta^0_{w_1,w_2}]\in H^{\lf}_2(\mathcal{Z}_{w_1,w_2},\mathcal{L}^\vee(\mathcal{G}^{\lambda}_{w_1,w_2}))
\) 
and
\(
[\Delta^0_{w_2,0}]\in H^{\lf}_2(\mathcal{Z}_{w_2,0},\mathcal{L}^\vee(\mathcal{G}^{-\lambda}_{w_2,0}))
\) and these cycles satisfy
\begin{align}
\partial[\Delta^0_{w_1,w_2}]&=\{(z_1,z_2)=(w_2,w_2)\},&
  \partial[\Delta^0_{w_2,0}]&=\{(z_3,z_4)=(w_2,w_2)\}.
\label{eq:loccycbdy}
\end{align}
Further, \(\sqbrac*{\Delta^0_{w_1,w_2}}\) and
\(\sqbrac*{\Delta^0_{w_2,0}}\) are, respectively, homologous to
\(\sqbrac*{\Delta_{w_1,w_2}}\) 
and
\(\sqbrac*{\Delta_{w_2,0}}\).
\label{thm:lfcyc}
\end{lemma}
\begin{proof}
  This lemma is a combination of the reasoning and constructions in Sections
  3.2.4 and 3.2.5 of \cite{hgeom}:
  The fact that \(\sqbrac*{\Delta^0_{w_1,w_2}}\) and
  \(\sqbrac*{\Delta^0_{w_2,0}}\) lie in the stated homology groups is \cite[Lem 3.3]{hgeom}
  and these newly constructed cycles being homologous to \(\sqbrac*{\Delta_{w_1,w_2}}\) and
  \(\sqbrac*{\Delta_{w_2,0}}\), respectively, is \cite[Thm 3.1]{hgeom}.
\end{proof}

\begin{prop}
  Let \(\lambda\in\CC\) with \(\Re{\lambda}\le -2 -2\frac{u}{v}\),
  then 
  \begin{align}
    \pair{\mu_\lambda}{R_E(\theta_\lambda)}&=\int_{[\Delta^0_{w_1,w_2}]}\brac*{\int_{[\Delta^0_{w_2,0}]}M(\lambda,w,z)\dd
      z_3\wedge \dd z_4}\dd z_1\wedge \dd z_2,\nonumber\\
    &=f(w_1,w_2)\sum_i\int_{\sqbrac*{\Delta^0_{w_1,w_2}}}
    H_{1,2}^{(i)}(z_1,z_2,w_2)\mathcal{G}^\lambda\brac*{z_1,z_2}\dd z_1\wedge
    \dd z_2
    \nonumber\\
    &\qquad\int_{\sqbrac*{\Delta^0_{w_2,0}}}H_{3,4}^{(i)}(z_3,z_4,w_1)\mathcal{G}^{-\lambda}(z_3,z_4)\dd
    z_3\wedge \dd z_4.
    \label{eq:regint}
  \end{align}
  In particular the \rhs{} is well-defined and the sum converges for \(\Re(\lambda)\le -2 -2\frac{u}{v}\).
  \label{thm:regcyc}
\end{prop}
\begin{proof}
Let \([{\Delta}^\eta_{x_{i},x_{i+1}}]=[\Delta^0_{x_{i},x_{i+1}}]\cap\{|z_i-w_2|^2+|z_{i+1}-w_2|^2\geq \eta^2\}\),
where \(\frac{w_1+w_2}{2}>\eta>0\). 
Let us denote \(\pair{\mu_\lambda}{R_E(\theta_\lambda)}_\eta\) the
integral \eqref{eq:matel} 
integrated over
\([\Delta^\eta_{w_1,w_2}]\) and \([\Delta^\eta_{w_2,0}]\).
From the definition of \([{\Delta}^\eta_{x_{i},x_{i+1}}]\), for any points \((z_1,z_2)\in [{\Delta}^\eta_{x_{1},x_{2}}]\), \((z_3,z_4)\in [{\Delta}^\eta_{x_{3},x_{4}}]\), we see that \(|z_k|>|z_l|\), \(k=1,2,\ l=3,4\).
That is the cycle \([\Delta^\eta_{w_1,w_2}]\times
  [\Delta^{\eta}_{w_2,0}]\) is contained in the domain
\begin{equation}
  \{w_1-w_2>|z_1-w_2|>|z_2-w_2|>0\}\cap\{ |z_2|>|z_3|>|z_4|>0\}.
  \label{eq:domainforsum}
\end{equation}

Further, the right-hand side of the expansion \eqref{eq:Hexpansion} is
  absolutely convergent on the domain \eqref{eq:domainforsum}, hence we can pull the sum out of the integral. That is, we have 
\begin{align}
  \pair{\mu_\lambda}{R_E(\theta_\lambda)}_\eta&=\int_{[\Delta^\eta_{w_1,w_2}]\cten[\Delta^\eta_{w_2,0}]}
  M(\lambda,w,z)\dd z_1\wedge \dd z_2\wedge \dd z_3\wedge \dd z_4\nonumber\\
  &=
  \int_{[\Delta^\eta_{w_1,w_2}]}\brac*{\int_{[\Delta^\eta_{w_2,0}]}M(\lambda,w,z)\dd
    z_3\wedge \dd z_4}\dd z_1\wedge \dd z_2\nonumber\\
  &=f(w_1,w_2)\int_{\sqbrac*{\Delta^\eta_{w_1,w_2}}}\int_{\sqbrac*{\Delta^\eta_{w_2,0}}}\sum_i
  \left(H_{1,2}^{(i)}(z_1,z_2,w_2)H_{3,4}^{(i)}(z_3,z_4,w_1)\mathcal{G}^\lambda\brac*{z_1,z_2}\right.\nonumber\\
  &\qquad\left. \mathcal{G}^{-\lambda}(z_3,z_4)\dd
    z_3\wedge \dd z_4\right) \dd z_1\wedge \dd z_2\nonumber\\
  &=f(w_1,w_2)\sum_i\int_{\sqbrac*{\Delta^\eta_{w_1,w_2}}}
    H_{1,2}^{(i)}(z_1,z_2,w_2)\mathcal{G}^\lambda\brac*{z_1,z_2}\dd z_1\wedge
    \dd z_2
    \nonumber\\
    &\qquad\int_{\sqbrac*{\Delta^\eta_{w_2,0}}}H_{3,4}^{(i)}(z_3,z_4,w_1)\mathcal{G}^{-\lambda}(z_3,z_4)\dd
    z_3\wedge \dd z_4.
\label{eq:splitint}
\end{align}
Since \(\Re\lambda\le -2-2\frac{u}{v}\), we see that \(M(\lambda,w,z)\)
is bounded on the closure of 
\(U_{w_2}\cap ([\Delta^0_{w_1,w_2}]\times[\Delta^0_{w_2,0}])\)
where \(U_{w_2}\) is a small
neighbourhood of \(w_2\). Thus from \eqref{eq:splitint}, we have 
\begin{align}
  \pair{\mu_\lambda}{R_E(\theta_\lambda)}_0&=\int_{[\Delta^0_{w_1,w_2}]}\brac*{\int_{[\Delta^0_{w_2,0}]}M(\lambda,w,z)\dd
    z_3\wedge \dd z_4}\dd z_1\wedge \dd z_2\nonumber\\
      &=f(w_1,w_2)\sum_i\int_{\sqbrac*{\Delta^0_{w_1,w_2}}}
    H_{1,2}^{(i)}(z_1,z_2,w_2)\mathcal{G}^\lambda\brac*{z_1,z_2}\dd z_1\wedge
    \dd z_2
    \nonumber\\
    &\qquad\int_{\sqbrac*{\Delta^0_{w_2,0}}}H_{3,4}^{(i)}(z_3,z_4,w_1)\mathcal{G}^{-\lambda}(z_3,z_4)\dd
    z_3\wedge \dd z_4.
  \label{eq:splitint2}
\end{align}
\end{proof}

Next we show that \(\pair{\mu_\lambda}{R_E(\theta_\lambda)}\) is non-zero
  by showing that the integral \eqref{eq:regint}
is
non-zero.

\begin{lemma}
  Let \(\lambda\in \CC\), then for \(\Re\lambda\le -2-2\frac{u}{v}\) 
  and
  \(\lambda-3\frac{u}{v}\notin 2\ZZ\) 
  we have that \(\pair{\mu_\lambda}{R_E(\theta_\lambda)}\neq0\).
  \label{thm:nonzeromatel}
\end{lemma}
\begin{proof}
We set
\(R(\lambda,w_1,w_2)=(w_1^{\lambda+\frac{u}{v}-2}w_2^{\lambda+2-\frac{u}{v}}f(w_1,w_2))^{-1}\pair{\mu_\lambda}{R_E(\theta_\lambda)}\).
We will see shortly that
\(R(\lambda,w_1,w_2)\)
admits an expansion in
\(\CC\powser{\frac{w_1-w_2}{w_2},\frac{w_2}{w_1}}\) and we will show
that this expansion is non-zero by showing that the constant term is
non-zero. Computing this expansion is most easily done in a new set of
variables \(y_1,\dots,y_4\) related to the \(z_1,\dots,z_4\) as follows.
\begin{equation}
  \label{eq:zytrans}
  z_1=(w_1-w_2)y_1+w_2,\ z_2=(w_1-w_2)y_2+w_2,\ z_3=w_2y_3,\ z_4=w_2y_4.
\end{equation}
In these new variables the \(\mathcal{G}\) functions are expressed as
\begin{equation}
  \mathcal{G}^\lambda_{w_1.w_2}(z_1,z_2)=\brac*{w_1-w_2}^{-2}\mathcal{G}^\lambda_{1,0}(y_1,y_2),\qquad
  \mathcal{G}^{-\lambda}_{w_2,0}(z_3,z_4)=w_2^{-2}\mathcal{G}^{-\lambda}_{1,0}(y_3,y_4).
\end{equation}
In particular, the right-hand side of \eqref{eq:regint} becomes
\begin{align}
  (w_1-w_2)^2w^2_2\int_{[\Delta^0_{1,0}]}\brac*{\int_{[\Delta^0_{1,0}]}\overline{M}(\lambda,w,y)\dd
    y_3\wedge \dd y_4}\dd y_1\wedge \dd y_2,
  \label{eq:yMint}
\end{align}
where \(\overline{M}(\lambda,w,y)\) is 
\({M}(\lambda,w,z)\) after the change of variables \eqref{eq:zytrans}.
We can then integrate \eqref{eq:yMint} term by term as a series in \(\frac{w_1-w_2}{w_2},\frac{w_2}{w_1}\).
Thus we consider how the factors \(H\) and \({}_2F_1\) of
\({M}(\lambda,w,z)\) 
are expand after the change of variables \eqref{eq:zytrans}.
The function \(H\) can then be expressed as
\begin{align}
  &(w_2^{\frac{v\lambda-u}{v}+2}w_1^{\frac{u+v\lambda}{v}-2})^{-1}H(z,w)\nonumber\\
  &=
  \brac*{1+\frac{z_1-w_2}{w_2}}^{1+\frac{v\lambda-u}{2v}}\brac*{1+\frac{z_2-w_2}{w_2}}^{1+\frac{v\lambda-u}{2v}}\brac*{1- \frac{z_3}{w_1}}^{\frac{u+v\lambda}{2v}-1}\brac*{1-
    \frac{z_4}{w_1}}^{\frac{u+v\lambda}{2v}-1}\nonumber\\
  &\qquad\cdot\brac*{1-\frac{z_{3}}{w_2\brac*{1+\frac{z_1-w_2}{w_2}}}}^{\frac{u}{v}}\brac*{1-\frac{z_{4}}{w_2\brac*{1+\frac{z_2-w_2}{w_2}}}}^{\frac{u}{v}}\nonumber\\
  &\qquad\cdot                             \brac*{1-
    \frac{z_{4}}{w_2\brac*{1+\frac{z_1-w_2}{w_2}}}}^{2-2\frac{u}{v}}\brac*{1-\frac{z_{3}}{w_2\brac*{1+\frac{z_2-w_2}{w_2}}}}^{2-2\frac{u}{v}}.\nonumber\\
  &=
  \brac*{1+\frac{w_1-w_2}{w_2}y_1}^{1+\frac{v\lambda-u}{2v}}\brac*{1+\frac{w_1-w_2}{w_2}y_2}^{1+\frac{v\lambda-u}{2v}}\brac*{1- \frac{w_2}{w_1}y_3}^{\frac{u+v\lambda}{2v}-1}\brac*{1-
    \frac{w_2}{w_1}y_4}^{\frac{u+v\lambda}{2v}-1}\nonumber\\
  &\qquad\cdot\brac*{1-\frac{y_{3}}{1+\frac{w_1-w_2}{w_2}y_1}}^{\frac{u}{v}}\brac*{1-\frac{y_{4}}{1+\frac{w_1-w_2}{w_2}y_2}}^{\frac{u}{v}}\nonumber\\
  &\qquad\cdot                             \brac*{1-
    \frac{y_{4}}{1+\frac{w_1-w_2}{w_2}y_1}}^{2-2\frac{u}{v}}\brac*{1-\frac{y_{3}}{1+\frac{w_1-w_2}{w_2}y_2}}^{2-2\frac{u}{v}}\nonumber\\
  &=(1-y_3)^{2-\frac{u}{v}}(1-y_4)^{2-\frac{u}{v}}+o(\tfrac{w_1-w_2}{w_2},\tfrac{w_2}{w_1}).
\end{align}
Before we consider the hypergeometric function note that
\begin{align}
  1-\frac{z_1-z_2}{z_1-z_3}\frac{z_3-z_4}{z_2-z_4}&=\frac{z_1-z_4}{z_1-z_3}\frac{z_2-z_3}{z_2-z_4}
  =\frac{\brac*{\frac{w_1-w_2}{w_2}y_1+1-y_4}}{\brac*{\frac{w_1-w_2}{w_2}y_1+1-y_3}}\frac{\brac*{\frac{w_1-w_2}{w_2}y_2+1-y_3}}{\brac*{\frac{w_1-w_2}{w_2}y_2+1-y_4}}\nonumber\\
  &=1+o(\tfrac{w_1-w_2}{w_2},\tfrac{w_2}{w_1})
\end{align}
and hence we need to expand around \(1\) rather than \(0\). This can be done
using the connection formula
\begin{multline}
  \tensor[_2]{F}{_{1}}(2-3\tfrac{u}{v},1-\tfrac{u}{v},2-2\tfrac{u}{v};x)
  =\frac{\Gamma(2-2\tfrac{u}{v})\Gamma(2\tfrac{u}{v}-1)}{\Gamma(\tfrac{u}{v})\Gamma(1-\tfrac{u}{v})}\tensor[_2]{F}{_{1}}(2-3\tfrac{u}{v},1-\tfrac{u}{v},2-2\tfrac{u}{v};1-x)\\
  +\brac*{1-x}^{2\frac{u}{v}-1}\frac{\Gamma(2-2\tfrac{u}{v})\Gamma(1-2\tfrac{u}{v})}{\Gamma(2-3\tfrac{u}{v})\Gamma(1-\tfrac{u}{v})}\tensor[_2]{F}{_{1}}(\tfrac{u}{v},1-\tfrac{u}{v},2\tfrac{u}{v};1-x). 
\end{multline}
Recalling the well known formulae for hypergeometric functions evaluated at
\(1\) we therefore obtain
\begin{align}
  &\tensor[_2]{F}{_{1}}(2-3\tfrac{u}{v},1-\tfrac{u}{v},2-2\tfrac{u}{v};\tfrac{z_1-z_2}{z_1-z_3}\tfrac{z_3-z_4}{z_2-z_4})\nonumber\\
  &=\frac{\Gamma(2-2\tfrac{u}{v})\Gamma(2\tfrac{u}{v}-1)}{\Gamma(\tfrac{u}{v})\Gamma(1-\tfrac{u}{v})}\frac{\Gamma(2-2\tfrac{u}{v})\Gamma(2\tfrac{u}{v}-1)}{\Gamma(\tfrac{u}{v})\Gamma(1-\tfrac{u}{v})}\nonumber\\
  &\quad    +\frac{\Gamma(2-2\tfrac{u}{v})\Gamma(1-2\tfrac{u}{v})}{\Gamma(2-3\tfrac{u}{v})\Gamma(1-\tfrac{u}{v})}\frac{\Gamma(2\tfrac{u}{v})\Gamma(2\tfrac{u}{v}-1)}{\Gamma(\tfrac{u}{v})\Gamma(3\tfrac{u}{v}-1)}
    +o(\tfrac{w_1-w_2}{w_2},\tfrac{w_2}{w_1})\nonumber\\
  &=\frac{\sin(\pi\tfrac{u}{v})^2+\sin(\pi3\tfrac{u}{v})\sin(\pi\tfrac{u}{v})}{\sin(\pi2\tfrac{u}{v})^2}+o(\tfrac{w_1-w_2}{w_2},\tfrac{w_2}{w_1})=1+o(\tfrac{w_1-w_2}{w_2},\tfrac{w_2}{w_1}),
\end{align}
where in the second equality we have used Euler's reflection formula
\(\Gamma(1-x)\Gamma(x)=\frac{\pi}{\sin(\pi x)}\).
Thus
\begin{align}
  &R(\lambda,w_1,w_2)=\int_{[\Delta_{1,0}]}\mathcal{G}^\lambda_{1,0}(y_1,y_2)\dd
  y_1\wedge\dd y_2\nonumber\\
  &\qquad \qquad\qquad\cdot \int_{[\Delta_{1,0}]}\mathcal{G}^{-\lambda}_{1,0}(y_3,y_4)(1-y_3)^{2-\frac{u}{v}}(1-y_4)^{2-\frac{u}{v}}\dd
    y_3\wedge\dd y_4+o(\tfrac{w_1-w_2}{w_2},\tfrac{w_2}{w_1})\nonumber\\
    &=\mathcal{I}_1\sqbrac*{1}\brac*{\tfrac{u-v\lambda}{2v}-1,\tfrac{u+v\lambda}{2v}-1,1-\tfrac{u}{v}}\cdot
\mathcal{I}_1\sqbrac*{1}\brac*{\tfrac{u+v\lambda}{2v}-1,\tfrac{-u-v\lambda}{2v}+1,1-\tfrac{u}{v}}+o(\tfrac{w_1-w_2}{w_2},\tfrac{w_2}{w_1}),
\end{align}
where we have evaluated the integrals using
\eqref{eq:integrationformula}. From the formula \eqref{eq:integrationformula}
for \(\mathcal{I}\) in terms of \(\Gamma\)-functions, we see that the
\(\Gamma\) functions in the numerator have poles when
\(\lambda-3\tfrac{u}{v}\in\ZZ\), but that
  \(\pair{\mu_\lambda}{R_E(\theta_\lambda)}\neq0\) when \(\lambda-3\tfrac{u}{v}\notin\ZZ\). Hence the lemma follows.
\end{proof}

\begin{theorem}
  The simple projective modules \(\sigma^\ell \slE{\mu}{1}{1}\), \(\mu\in \RR/2\ZZ\),
  \(\mu\neq \pm \lambda_{1,1}+2\ZZ\), \(\lambda\in \ZZ\) are rigid, with rigid
  dual \(\sigma^{-\ell}\slE{-\mu}{1}{1}\).
  \label{thm:11rigid}
\end{theorem}
\begin{proof}
  Without loss of generality we can restrict ourselves to only a single value of the spectral flow parameter
  and hence we choose \(\ell=0\). As mentioned at the beginning of this
  section by \cite[Lem 4.2.1 and Cor 4.2.2]{CMcY} a sufficient
  condition for \(\slE{-\mu}{1}{1}\) being the rigid dual of
  \(\slE{\mu}{1}{1}\) is the non-vanishing of the matrix element \eqref{eq:matel}. So initially
  assume that \(\mu\neq 3\frac{u}{v}+2\ZZ\) 
  and pick a representative \(\lambda\in \mu\) which is sufficiently negative
  for \cref{thm:nonzeromatel} to apply. Hence the matrix element
  \(\pair{\mu_\lambda}{R_E(\theta_\lambda)}\) is non-vanishing and \(\slE{\mu}{1}{1}\) is rigid. To
  extend this conclusion to the previously excluded case, let
  \(\mu_1,\mu_2\in \RR/2\ZZ\), \(\mu_1,\mu_2\notin \set{\pm
    \tfrac{u}{v}+2\ZZ,3\tfrac{u}{v}+2\ZZ}\) satisfying
  \(\mu_1+\mu_2=3\tfrac{u}{v}+2\ZZ\). Then the \(\slE{\mu_i}{1}{1}\) are both
  rigid (by our reasoning so far) and
  projective. From the tensor product formula \eqref{eq:EEfusion}, which we have
  already proved, we see that the tensor product \(\slE{\mu_1}{1}{1}\fuse
  \slE{\mu_2}{1}{1}\) contains \(\sigma\slE{3\tfrac{u}{v}+2\ZZ}{1}{1}\) and
  \(\sigma^{-1}\slE{5\tfrac{u}{v}+2\ZZ}{1}{1}\) as direct summands. Hence
   \(\slE{3\tfrac{u}{v}+2\ZZ}{1}{1}\)is also  rigid and the theorem follows.
\end{proof}

\section{Non-semisimple fusion decompositions}
\label{sec:nonssifusion}

\begin{theorem}
  The non-semisimple fusion product decomposition of \cref{thm:EEfusionconj}
  hold. That is,
  for \(1\le r\le u-1\), \(\mu,\mu'\in \RR/2\ZZ\), \(\pm \lambda_{1,1}\notin \mu\), \(\pm
  \lambda_{r,1}\notin \mu'\) and \(r-1\in \mu+\mu'\)
  \begin{align}
    \slE{\mu}{1}{1}\fuse \slE{\mu'}{r}{1}\cong \sigma^{-1}P_{u-r,v-1}\oplus
    \brac*{1-\delta_{v,2}}\slE{\mu+\mu'}{r}{2},
    \label{eq:r1logfusion2}
  \end{align}
  while for \(1\le r\le u-1\), \(2\le s\le v-2\), \(\mu,\mu'\in \RR/2\ZZ\), \(\pm \lambda_{1,1}\notin \mu\), \(\pm
  \lambda_{r,s}\notin \mu'\),
  \begin{align}
    \slE{\mu}{1}{1}\fuse \slE{\mu'}{r}{s}\cong
    \begin{cases}
      P_{r,s-1}\oplus \sigma^{-1} \slE{\mu+\mu'+\frac{u}{v}}{r}{s}\oplus
      \slE{\mu+\mu'}{r}{s+1},&\text{if }\ \lambda_{r,s-1}\in \mu+\mu',\\
      P_{u-r,v-s-1}\oplus\sigma^{-1}\slE{\mu+\mu'+\frac{u}{v}}{r}{s}\oplus
      \slE{\mu+\mu'}{r}{s-1}&\text{if }\ \lambda_{u-r,v-s-1}\in \mu+\mu',\\
      \sigma^{-1}P_{r,s}\oplus \sigma \slE{\mu+\mu'-\frac{u}{v}}{r}{s}\oplus
      \slE{\mu+\mu'}{r}{s-1},&\text{if }\ \lambda_{r,s+1}\in \mu+\mu',\\
      \sigma^{-1}P_{u-r,v-s}\oplus\sigma \slE{\mu+\mu'-\frac{u}{v}}{r}{s}\oplus
      \slE{\mu+\mu'}{r}{s+ 1}&\text{if }\ \lambda_{u-r,v-s+1}\in \mu+\mu'.
    \end{cases}
   \label{eq:rslogfusion}
  \end{align}
  Moreover, the category \(\wtmod{\mmsl{u}{v}}\) is rigid.
\end{theorem}
\begin{proof}
  Note that a module is rigid if and only if any of its spectral flows are rigid.
  We begin by proving \eqref{eq:r1logfusion2} for \(r=1\). Note that from
  \cref{thm:upperbounds,thm:lowerbounds1,thm:lowerbounds2} we have the equality \(\dim
  \hom_{\mmsl{u}{v}}\brac*{\slE{\mu}{1}{1}\fuse\slE{-\mu}{1}{1},\slE{0}{1}{2}}=1-\delta_{v,2}\),
  hence since \(\slE{0}{1}{2}\) is simple and projective, if it exists (that
  is, if \(v>2\)) it must be a direct summand. From
  \cref{thm:upperbounds,thm:lowerbounds1,thm:lowerbounds2} we also see that there can be no
  simple projective summands other than \(\slE{0}{1}{2}\). Further, from
  \cref{thm:logint} we also know that \(\sigma^{-1}P_{u-1,v-1}\) is a direct
  summand of
  \(\slE{\mu}{1}{1}\fuse\slE{-\mu}{1}{1}\). Thus
  \begin{equation}
    \slE{\mu}{1}{1}\fuse\slE{-\mu}{1}{1}\cong
    \sigma^{-1}P_{u-1,v-1}\oplus
    \brac*{1-\delta_{v,2}}\slE{0}{1}{2}\oplus M,
    \label{eq:1111rulewithM}
  \end{equation}
  where \(M\) is a projective module whose composition factors must all be
  spectral flows of highest weight modules, that is of the form \(\sigma^\ell
  \slD{r}{s}{+}\). We then
  choose \(\nu\in \RR/2\ZZ\) such that \(\nu\notin
  \set{2\ZZ,-\mu,\sqbrac*{\pm 2t},\sqbrac*{\pm t \pm'
      \lambda_{1,2}},\sqbrac*{\pm\lambda_{1,3}}}\), and by associativity we
  can evaluate
  \(\slE{\nu}{1}{1}\fuse\slE{\mu}{1}{1}\fuse\slE{-\mu}{1}{1}\) as
  \begin{align}
    \brac*{\slE{\nu}{1}{1}\fuse\slE{\mu}{1}{1}}\fuse\slE{-\mu}{1}{1}
    &\cong \sigma^{-2}\slE{\nu+2\tfrac{u}{v}}{1}{1}\oplus
    \sigma^{2}\slE{\nu-2\tfrac{u}{v}}{1}{1}
    \oplus
     (3-\delta_{v,2})\slE{\nu}{1}{1}\nonumber\\
     &\quad\oplus (1-\delta_{v,2})
       \brac*{\sigma^{-1}\slE{\nu+\tfrac{u}{v}}{1}{2}\oplus\slE{\nu-\tfrac{u}{V}}{1}{2}\oplus(1-\delta_{v,3})\slE{\nu}{1}{3}}
       \label{eq:leftbracket}
  \end{align}
  and note that the \rhs{} of \eqref{eq:leftbracket} is semisimple with
  \(5-\delta_{v,2}+\brac*{1-\delta_{v,2}}\brac*{3-\delta_{v,3}}\) simple
  projective summands. Evaluating the fusion product with the other bracketing
  yields
  \begin{align}
    \slE{\nu}{1}{1}\fuse\brac*{\slE{\mu}{1}{1}\fuse\slE{-\mu}{1}{1}}&\cong
    \slE{\nu}{1}{1}\fuse\brac*{\sigma^{-1}P_{u-1,v-1}\oplus
      \brac*{1-\delta_{v,2}}\slE{0}{1}{2}\oplus M}\nonumber\\
    &=\slE{\nu}{1}{1}\fuse\brac*{\sigma^{-1}P_{u-1,v-1}\oplus M}\nonumber\\
    &\oplus \brac*{1-\delta_{v,2}}\brac*{\sigma^{-1}\slE{\nu+\tfrac{u}{v}}{1}{2}\oplus\sigma\slE{\nu-\tfrac{u}{v}}{1}{2}\oplus\slE{\nu}{1}{1}\oplus(1-\delta_{v,3})\slE{\nu}{1}{3}}.
  \end{align}
  The projective module \(P_{u-1,v-1}\) has four composition factors and since
  the \rhs{} of  \eqref{eq:leftbracket} is semisimple, \(\slE{\nu}{1}{1}\fuse
  P_{u-1,v-1}\) must contribute at least four simple summands in addition to
  the \((1-\delta_{v,2})(3-\delta_{v,3})\) contributed by
  \(\slE{\nu}{1}{1}\fuse\slE{0}{1}{2}\), which brings the total up to
  \(5-\delta_{v,2}+\brac*{1-\delta_{v,2}}\brac*{3-\delta_{v,3}}\). Hence
  \(\slE{\nu}{1}{1}\fuse M =0\), that is \(M=0\).

  The decomposition formula \eqref{eq:r1logfusion2} for general \(r\) now
follows from the formula for \(r=1\) by repeatedly applying \(L_2\fuse\).

Next we show \eqref{eq:rslogfusion} for \(r=1\) (once this formula has been
establish for \(r=1\), the formula for general \(r\) follows again by applying
\(L_2\fuse\)). We consider first the case
\(\mu+\mu'=\sqbrac{\lambda_{r,s-1}}\), then
\cref{thm:upperbounds,thm:lowerbounds1,thm:lowerbounds2} imply that
\(\slE{\mu}{1}{1}\fuse\slE{\mu'}{1}{s}\) has exactly two direct summands
isomorphic to \(\sigma^{-1}\slE{\sqbrac*{\lambda_{r,s-1}+\tfrac{u}{v}}}{1}{s}\)
and \(\slE{\sqbrac*{\lambda_{r,s-1}+\tfrac{u}{v}}}{1}{s+1}\). From
\cref{thm:lowerbounds1} we also see that there exists a surjective
intertwining operator of type
\(\binom{\sigma\slEr{1}{s}{+}}{\slE{\mu}{1}{1},\,
  \slE{\mu'}{1}{s}}\). Hence \( P_{1,s-1}\) (which by
\eqref{eq:projmchar} is the projective cover of \(\sigma^{-1}\slEr{1}{s}{-}\))
must be a direct summand of \(\slE{\mu}{1}{1}\fuse\slE{\mu'}{1}{s}\). Thus
\begin{equation}
  \slE{\mu}{1}{1}\fuse\slE{\mu'}{1}{s}=\sigma^{-1} P_{1,s-1}\oplus
  \sigma^{-1}\slE{\sqbrac*{\lambda_{r,s-1}+\tfrac{u}{v}}}{1}{s}\oplus
  \slE{\sqbrac*{\lambda_{r,s-1}+\tfrac{u}{v}}}{1}{s+1}\oplus M,
  \label{eq:111srulewithM}
\end{equation}
where again \(M\) is a projective module whose composition factors must all be
  spectral flows of highest weight modules, that is of the form \(\sigma^\ell
 \slD{r}{s}{+}\). We can now use the same reasoning for concluding that
  \(M=0\) as we did for \eqref{eq:1111rulewithM} by applying
  \(\slE{\nu}{1}{1}\) to \eqref{eq:111srulewithM} with \(\nu\in \RR/2\ZZ\)
  chosen such that the fusion product will be semisimple. Counting the number
  of simple summands contributed by each summand in \eqref{eq:111srulewithM}
  to \(\slE{\nu}{1}{1}\fuse\slE{\mu}{1}{1}\fuse\slE{\mu'}{1}{s}\) will again
  imply that \(M=0\).
  Similarly for \(\mu+\mu'=\sqbrac*{ \lambda_{{u-1}, {v-s-1}}}\)
  \cref{thm:lowerbounds1} implies the existence of a surjective intertwining
  operator of type \(\binom{\sigma\slEr{u-1}{v-s}{+}}{\slE{\mu}{1}{1},\,
  \slE{\mu'}{1}{s}}\) with \(P_{u-1,v-s-1}\) as the corresponding projective
summand. While for \(\mu+\mu'=\sqbrac*{\lambda_{1,s+1}}\) and
\(\mu+\mu'=\sqbrac*{\lambda_{u-1,v-s+1}}\) the corresponding intertwining
operator types and projective summands are \(\binom{\sigma^{-1}\slEr{1}{s+1}{-}}{\slE{\mu}{1}{1},\,
  \slE{\mu'}{1}{s}}\) and \(\sigma^{-1}P_{1,s}\), and \(\binom{\sigma^{-1}\slEr{u-1}{v-s+1}{-}}{\slE{\mu}{1}{1},\,
  \slE{\mu'}{1}{s}}\) and \(\sigma^{-1}P_{u-1,v-s}\), respectively.

With the fusion product decomposition formulae of \cref{thm:EEfusionconj} now
proved, we see that all projective modules 
appear as direct summands in repeated products of the simple projective
modules \(\slE{1}{1}{\mu}\) and \(\slE{2}{1}{\mu}\). Hence all projective
modules are rigid. To conclude that all the composition factors of the
  non-simple indecomposable projectives are rigid, we note that in the
  nomenclature of \cite[App A]{KL4} \(\wtmod{\mmsl{u}{v}}\) is weakly rigid
  \cite[Def A.4]{KL4} and semirigid \cite[Def A.6]{KL4} because the \(\mmsl{u}{v}\) is isomorphic to itself and
  is hence a dualising object in the sense of Grothendieck-Verdier categories
  (see \cite[Thm 2.12]{ALSW}). Further, \(\wtmod{\mmsl{u}{v}}\) is Frobenius
  \cite[Def A.9]{KL4} because there are sufficiently many injectives and
  projectives, all projectives are injective and vice versa, and all
  projectives are rigid. Thus we can use \cite[Prop A.2]{KL4}, which states
  that if two terms in an exact sequence in a Frobenius category are rigid,
  then the third is as well. We already know that the simple modules
  \(L_r\cong \sigma^{-1}\slD{u-r}{v-1}{+}\), \(1\le r\le u-1\) are rigid. From
  \cref{thm:sl2modclass} we know that \(P_{r,v-1}\) admits a socle filtration
  \(P_{r,v-1}=M_2\supset M_1\supset M_0\) satisfying
  \begin{equation}
    M_0\cong\slD{r}{v-1}{+}\cong M_2/M_1,\quad
    M_1/M_0\cong\sigma^{-1}\slD{u-r}{v-2}{+}\oplus \sigma^2\slD{u-r}{1}{+}.
  \end{equation}
  Thus \cite[Prop A.2]{KL4} implies that \(M_2/M_0\) is rigid because \(M_0\)
  and \(M_2=P_{r,v-1}\) are rigid and in turn \(M_1/M_0\) is rigid because
  \(M_2/M_1\cong (M_2/M_0)/(M_1/M_0)\) and \(M_2/M_0\) are. Hence
  \(\slD{u-r}{v-2}{+}\) and \(\slD{u-r}{1}{+}\) are rigid because they are
  (spectral flows) of direct summands of rigid modules. By the same argument
  the rigidity of the remaining \(\slD{r}{s}{+}\) can be deduced inductively
  from the socle filtration of \(P_{r,s}\), \(1\le s\le v-2\), which satisfies
  \begin{equation}
    M_0\cong\slD{r}{s}{+}\cong M_2/M_1,\quad
    M_1/M_0\cong\sigma^{-1}\slD{r}{s-1}{+}\oplus \sigma\slD{r}{s+1}{+}.
  \end{equation}
  and hence the rigidity of \(\slD{r}{s}{+}\) implies the rigidity of
  \(\slD{r}{s+1}{+}\).
Thus \(\wtmod{\mmsl{u}{v}}\) is rigid.
 \end{proof}

 \begin{theorem}
   The category \(\wtmod{\mmN{u}{v}}\) is rigid.
 \end{theorem}

 \begin{proof}
   This follows immediately from \cref{thm:rigiditytransport} and the
   rigidity of \(\wtmod{\mmsl{u}{v}}\).
 \end{proof}
 

\begin{thebibliography}{10}

\bibitem{A1}
D.~Adamovi\'{c}.
\newblock Representations of the {$N=2$} superconformal vertex algebra.
\newblock {\em Int. Math. Res. Not.}, 1999:61--79, 1999.
\newblock \textsf{arXiv:math/9809141 \mbox{[math.QA]}}.

\bibitem{A2}
D.~Adamovi\'{c}.
\newblock Vertex algebra approach to fusion rules for {$N=2$} superconformal
  minimal models.
\newblock {\em J. Algebra}, 239:549--572, 2001.

\bibitem{A3}
D.~Adamovi{\'c}.
\newblock Realizations of simple affine vertex algebras and their modules: the
  cases $\widehat{sl(2)}$ and $\widehat{osp (1,2)}$.
\newblock {\em Comm. Math. Phys.}, 366:1025--1067, 2019.
\newblock \textsf{arXiv:1711.11342 \mbox{[math.QA]}},.

\bibitem{bpAKR}
D~Adamovi\'{c}, K~Kawasetsu, and D~Ridout.
\newblock Weight module classifications for {Bershadsky--Polyakov} algebras.
\newblock {\em Commun. Contemp. Math.}, 26:2350063, 2024.
\newblock \textsf{arXiv:2303.03713 \mbox{[math.QA]}}.

\bibitem{AM}
D~Adamovi{\'c} and A.~Milas.
\newblock Vertex operator algebras associated to modular invariant
  representations for {$A_1^{(1)}$}.
\newblock {\em Math. Res. Lett.}, 2:563--575, 1995.
\newblock \textsf {arXiv:q-alg/9509025}.

\bibitem{AM09}
D.~Adamovi\'c and A.~Milas.
\newblock Lattice construction of logarithmic modules for certain vertex
  algebras.
\newblock {\em SM}, 15:535--561, 2009.
\newblock \textsf{arXiv:0902.3417 \mbox{[math.QA]}}.

\bibitem{ALSW}
R.~Allen, S.~Lentner, C.~Schweigert, and S.~Wood.
\newblock Duality structures for module categories of vertex operator algebras
  and the {F}eigin {F}uchs boson.
\newblock {\em Sel. Math. New Ser.}, 31:57, 2025.
\newblock \textsf {arXiv:2107.05718 \mbox{[math.QA]}}.

\bibitem{AW}
R.~Allen and S.~Wood.
\newblock Bosonic ghostbusting -- the bosonic ghost vertex algebra admits a
  logarithmic module category with rigid fusion.
\newblock {\em Comm. Math. Phys.}, 390:959--1015, 2022.
\newblock \textsf{arXiv:2001.05986 \mbox{[math.QA]}}.

\bibitem{hgeom}
K.~Aomoto and M.~Kita.
\newblock {\em Theory of hypergeometric functions}.
\newblock Springer Monographs in Mathematics. Springer-Verlag, Tokyo, 2011.

\bibitem{ACK}
T.~Arakawa, T.~Creutzig, and K.~Kawasetsu.
\newblock Weight representations of affine {K}ac-{M}oody algebras and small
  quantum groups.
\newblock {\em Adv. Math.}, 477:110365, 2025.
\newblock \textsf {arXiv:2311.10233 \mbox{[math.RT]}}.

\bibitem{ACRpar}
J.~Auger, T.~Creutzig, and D.~Ridout.
\newblock Modularity of logarithmic parafermion vertex algebras.
\newblock {\em Lett. Math. Phys.}, 108:2543--–2587, 2018.
\newblock \textsf{arXiv:1704.05168 \mbox{[math.QA]}}.

\bibitem{BPZ}
A.A Belavin, A.M. Polyakov, and A.B. Zamolodchikov.
\newblock Infinite conformal symmetry in two-dimensional quantum field theory.
\newblock {\em Nucl. Phys.}, B241:333--380, 1984.

\bibitem{BerFel90}
D.~Bernard and G.~Felder.
\newblock Fock representations and {BRST} cohomology in {$\mathsf{SL}(2)$}
  current algebra.
\newblock {\em Comm. Math. Phys.}, 127:145--168, 1990.

\bibitem{Blosl279}
R.~Block.
\newblock Classification of the irreducible representations of
  $\mathfrak{sl}(2,\mathbb{C})$.
\newblock {\em Bull. Amer. Math. Soc.}, 1:247--250, 1972.

\bibitem{CRRn11}
M.~Canagasabey, J.~Rasmussen, and D.~Ridout.
\newblock Fusion rules for the logarithmic {$N=1$} superconformal minimal
  models {I}: the {Neveu-Schwarz} sector.
\newblock {\em J. Phys.}, A48:415402, 2015.
\newblock \textsf{arXiv:1504.03155 \mbox{[hep-th]}}.

\bibitem{CRRn12}
M.~Canagasabey and D.~Ridout.
\newblock Fusion rules for the logarithmic {$N=1$} superconformal minimal
  models {II}: including the {Ramond} sector.
\newblock {\em Nucl. Phys.}, B905:132--187, 2016.
\newblock \textsf{arXiv:1512.05837 \mbox{[hep-th]}}.

\bibitem{Cfusion}
T.~Creutzig.
\newblock Resolving {V}erlinde's formula of {L}ogarithmic {CFT}, 2024.
\newblock \textsf{arXiv:2411.11383 \mbox{[math.QA]}}.

\bibitem{Csl2ten}
T.~Creutzig.
\newblock Tensor categories of weight modules of {$\widehat{\mathfrak{sl}}_2$}
  at admissible level.
\newblock {\em J. Lond. Math. Soc.}, 110:e70037, 2024.
\newblock \textsf {arXiv:2311.10240 \mbox{[math.RT]}}.

\bibitem{CGNR}
T.~Creutzig, N.~Genra, S.~Nakatsuka, and R.~Sato.
\newblock Correspondences of categories for subregular {$\mathcal{W}$}-algebras
  and principal {$\mathcal{W}$}-superalgebras.
\newblock {\em Comm. Math. Phys.}, 393:1--60, 2022.
\newblock \textsf{arXiv:2104.00942 \mbox{[math.RT]}}.

\bibitem{CHY}
T.~Creutzig, Y.-Z. Huang, and J.~Yang.
\newblock Braided tensor categories of admissible modules for affine {L}ie
  algebras.
\newblock {\em Comm. Math. Phys.}, 362:827--854, 2018.
\newblock \textsf{arXiv:1709.01865 \mbox{[math.QA]}}.

\bibitem{CJOHRY}
T.~Creutzig, C.~Jiang, F.~Orosz~Hunziker, D.~Ridout, and J.~Yang.
\newblock Tensor categories arising from the {V}irasoro algebra.
\newblock {\em Adv. Math.}, 380:107601, 2021.
\newblock \textsf{arXiv:2002.03180 \mbox{[math.RT]}}.

\bibitem{CKLR}
T.~Creutzig, S.~Kanade, A.R. Linshaw, and D.~Ridout.
\newblock Schur--{W}eyl duality for {H}eisenberg cosets.
\newblock {\em Transf. Groups}, 24:301--354, 2019.
\newblock \textsf {arXiv:1611.00305 \mbox{[math.QA]}}.

\bibitem{ospCKLR}
T.~Creutzig, S.~Kanade, T.~Liu, and D.~Ridout.
\newblock Cosets, characters and fusion for admissible-level
  \(\mathfrak{osp}(1|2)\) minimal models.
\newblock {\em Nucl. Phys.}, B938:22--55, 2019.
\newblock \textsf{arXiv:1806.09146 \mbox{[hep-th]}}.

\bibitem{CKM}
T.~Creutzig, S.~Kanade, and R.~McRae.
\newblock Tensor categories for vertex operator superalgebra extensions.
\newblock {\em Mem. Amer. Math. Soc.}, 2024.
\newblock \textsf{arXiv:1705.05017 \mbox{[math.QA]}}.

\bibitem{CL2}
T.~Creutzig and A.~R Linshaw.
\newblock Cosets of affine vertex algebras inside larger structures.
\newblock {\em J. Algebra}, 517:396--438, 2019.
\newblock \textsf{arXiv:1407.8512 \mbox{[math.RT]}}.

\bibitem{CLRW}
T.~Creutzig, T.~Liu, D.~Ridout, and S.~Wood.
\newblock Unitary and non-unitary {$N=2$} minimal models.
\newblock {\em J. High Energy Phys.}, 2019:1--45, 2019.
\newblock \textsf{arXiv:1902.08370 [math-ph]}.

\bibitem{CMOHYN1}
T.~Creutzig, R.~McRae, F.~Orosz~Hunziker, and J.~Yang.
\newblock ${N} = 1$ super {V}irasoro tensor categories.
\newblock {\em Comm. Math. Phys.}, 407(4):73, 2026.
\newblock \textsf{arXiv:2412.18127 \mbox{[math.QA]}}.

\bibitem{DLcomp}
T.~Creutzig, R.~McRae, and J.~Yang.
\newblock Direct limit completions of vertex tensor categories.
\newblock {\em Commun. Contemp. Math.}, 24:2150033, 2022.
\newblock \textsf{arXiv:2006.09711 \mbox{[math.QA]}}.

\bibitem{CMcY}
T.~Creutzig, R.~McRae, and J.~Yang.
\newblock Tensor structure on the {Kazhdan-Lusztig} category for affine
  $\mathfrak{gl}(1| 1)$.
\newblock {\em Int. Math. Res. Not.}, 2022:12462--2515, 2022.
\newblock \textsf{arXiv:2009.00818 \mbox{[math.QA]}}.

\bibitem{CMY2}
T.~Creutzig, R.~McRae, and J.~Yang.
\newblock Rigid {T}ensor {S}tructure on {B}ig {M}odule {C}ategories for {S}ome
  {W}-(super)algebras in {T}ype {A}.
\newblock {\em Comm. Math. Phys.}, 404:339--400, 2023.
\newblock \textsf{arXiv:2210.04678 \mbox{[math.QA]}}.

\bibitem{CMYrigid}
T.~Creutzig, R.~McRae, and J.~Yang.
\newblock Ribbon categories of weight modules for affine $\mathfrak{sl}_2$ at
  admissible levels.
\newblock {\em Comm. Math. Phys.}, 407(6):123, 2026.
\newblock \textsf{arXiv:2411.11386 \mbox{[math.QA]}}.

\bibitem{CMtheta}
T.~Creutzig and A.~Milas.
\newblock False theta functions and the {Verlinde} formula.
\newblock {\em Adv. Math.}, 262:520--545, 2014.
\newblock \textsf{arXiv:1309.6037 \mbox{[math.QA]}}.

\bibitem{CMSqdim}
T.~Creutzig, A.~Milas, and S.~Wood.
\newblock On regularised quantum dimensions of the singlet vertex operator
  algebra and false theta functions.
\newblock {\em Int. Math. Res. Not.}, 5:1390--1432, 2017.
\newblock \textsf{arXiv:1411.3282 \mbox{[math.QA]}}.

\bibitem{CR1}
T.~Creutzig and D.~Ridout.
\newblock Modular data and {V}erlinde formulae for fractional level {WZW}
  models {I}.
\newblock {\em Nucl. Phys.}, B865:83--114, 2012.
\newblock \textsf{arXiv:1205.6513 \mbox{[hep-th]}}.

\bibitem{CR2}
T.~Creutzig and D.~Ridout.
\newblock Modular data and {V}erlinde formulae for fractional level {WZW}
  models {II}.
\newblock {\em Nucl. Phys.}, B875:423--458, 2013.
\newblock \textsf{arXiv:1205.6513 \mbox{[hep-th]}}.

\bibitem{glCR}
T.~Creutzig and D.~Ridout.
\newblock Relating the archetypes of logarithmic conformal field theory.
\newblock {\em Nucl. Phys.}, B872:348--391, 2013.
\newblock \textsf{arXiv:1107.2135 \mbox{[hep-th]}}.

\bibitem{CFTbook}
P.~Di~Francesco, P.~Mathieu, and D.~S\'en\'echal.
\newblock {\em Conformal Field Theory}.
\newblock Graduate Texts in Contemporary Physics. Springer, New York, 1997.

\bibitem{DoLGVA93}
C.~Dong and J.~Lepowsky.
\newblock {\em Generalized vertex algebras and relative vertex operators}.
\newblock Progress in Mathematics. {Birkh{\"a}user}, {Boston}, 1993.

\bibitem{EG}
W.~Eholzer and M.R. Gaberdiel.
\newblock Unitarity of rational {$N= 2$} superconformal theories.
\newblock {\em Comm. Math. Phys.}, 186:61--85, 1997.
\newblock \textsf{arXiv:hep-th/9601163}.

\bibitem{EGNO}
P.~Etingof, S.~Gelaki, D.~Nikshych, and V.~Ostrik.
\newblock {\em Tensor categories}, volume 205 of {\em Mathematical Surveys and
  Monographs}.
\newblock American Mathematical Society, Providence, RI, 2015.
\newblock \textsf{arXiv:math/0101219 \mbox{[math.QA]}}.

\bibitem{FRRsl3}
J.~Fasquel, C.~Raymond, and D.~Ridout.
\newblock Modularity of admissible-level $\mathfrak{sl}_3$ minimal models with
  denominator 2.
\newblock {\em Comm. Math. Phys.}, 406:60, 2025.
\newblock \textsf{arXiv:2406.10646 \mbox[math.QA]}.

\bibitem{FRbp}
Z.~Fehily and D.~Ridout.
\newblock Modularity of {Bershadsky-Polyakov} minimal models.
\newblock {\em Lett. Math. Phys.}, 112:61, 2022.
\newblock \textsf{arXiv:2110.10336 \mbox{[math.QA]}}.

\bibitem{FST}
B.~L. Feigin, A.~M. Semikhatov, and I.~Yu. Tipunin.
\newblock Equivalence between chain categories of representations of affine
  {$\mathfrak{sl} (2)$} and {$N= 2$} superconformal algebras.
\newblock {\em J. Math. Phys.}, 39:3865--3905, 1998.
\newblock \textsf{arXiv:hep-th/9701043}.

\bibitem{FZ}
I.~Frenkel and Y.~Zhu.
\newblock Vertex operator algebras associated to representations of affine and
  {V}irasoro algebras.
\newblock {\em Duke Math. J.}, 66:123--168, 1992.

\bibitem{GK1}
M.~Gorelik and V.~Kac.
\newblock On simplicity of vacuum modules.
\newblock {\em Adv. Math.}, 211:621--677, 2007.
\newblock \textsf{arXiv:math-ph/0606002}.

\bibitem{HuaVer08}
Y.-Z. Huang.
\newblock Vertex operator algebras and the {Verlinde} conjecture.
\newblock {\em Commun. Contemp. Math.}, 10:103--1054, 2008.
\newblock \textsf{arXiv: \mbox{math}/0406291}.

\bibitem{HL}
Y.-Z. Huang and J.~Lepowsky.
\newblock Tensor categories and the mathematics of rational and logarithmic
  conformal field theory.
\newblock {\em Journal of Physics A}, 46:494009, 2013.

\bibitem{HLZ3}
Y.-Z. Huang, J.~Lepowsky, and L.~Zhang.
\newblock Logarithmic tensor category theory, {III}: Intertwining maps and
  tensor product bifunctors.
\newblock \textsf{arXiv:1012.4197 \mbox{[math.QA]}}.

\bibitem{HLZ}
Y.-Z. Huang, J.~Lepowsky, and L.~Zhang.
\newblock Logarithmic tensor category theory for generalized modules for a
  conformal vertex algebra, {I}: introduction and strongly graded algebras and
  their generalized modules.
\newblock In {\em Conformal Field Theories and Tensor Categories: Proceedings
  of a Workshop Held at Beijing International Center for Mathematical
  Research}, Mathematical Lectures from Peking University, pages 169--248.
  Springer, Heidelberg, 2014.
\newblock \textsf{arXiv:1012.4193 \mbox{[math.QA]}}.

\bibitem{HMII}
Y.-Z. Huang and A.~Milas.
\newblock Intertwining operator superalgebras and vertex tensor categories for
  superconformal algebras, {II}.
\newblock {\em Tans. Amer. Math. Soc.}, 354:363--385, 2002.
\newblock \textsf{arXiv:math/0004039 \mbox{[math.QA]}}.

\bibitem{HYzhudims}
Y.-Z. Huang and J~Yang.
\newblock Logarithmic intertwining operators and associative algebras.
\newblock {\em J. Pure Appl. Algebra}, 216:1467--1492, 2012.
\newblock \textsf{arXiv:1104.4679 \mbox{[math.QA]}}.

\bibitem{GtoP}
K~Iwasaki, H~Kimura, S~Shimomura, and M~Yoshida.
\newblock {\em From Gauss to Painlev\'e: a modern theory of special functions}.
\newblock Aspects of {Mathematics}. {Vieweg+Teubner}, Wiesbaden, 2013.

\bibitem{KacBeg}
V.~Kac.
\newblock {\em Vertex algebras for beginners}, volume~10 of {\em University
  Lecture Series}.
\newblock American Mathematical Society, 1998.

\bibitem{KM2}
V.~Kac and M.~Wakimoto.
\newblock Modular and conformal invariance constraints in representation theory
  of affine algebras.
\newblock {\em Adv. Math.}, 70:156--–236, 1988.

\bibitem{KM1}
V.~Kac and M.~Wakimoto.
\newblock Modular invariant representations of infinite-dimensional {Lie}
  algebras and superalgebras.
\newblock {\em Proc. Nat. Acad. Sci. USA}, 85:4956--4960, 1988.

\bibitem{KM3}
V.~Kac and M.~Wakimoto.
\newblock Classification of modular invariant representations of affine
  algebras.
\newblock In {\em In Infinite-Dimensional Lie Algebras and Groups
  (Luminy-Marseille, 1988)}, volume~7 of {\em Adv. Ser. Math. Phys.}, pages
  138--177, New Jersey, 1989. World Scientific.

\bibitem{KW}
V.~Kac and W.~Wang.
\newblock Vertex operator superalgebras and their representations.
\newblock In {\em Mathematical aspects of conformal and topological field
  theories and quantum groups}, volume 175 of {\em Contemp. Math.}, page 161,
  1994.
\newblock \textsf{arXiv:hep-th/9312065}.

\bibitem{KaRi}
S~Kanade and D~Ridout.
\newblock {\em {{NGK}} and {{HLZ}}: Fusion for Physicists and Mathematicians},
  pages 135--181.
\newblock Springer International Publishing, Cham, 2019.
\newblock \textsf{arXiv:1812.10713 [\mbox{math-ph}]}.

\bibitem{KRWsl3}
K.~Kawasetsu, D.~Ridout, and S.~Wood.
\newblock Admissible-level $\mathfrak{sl}_3$ minimal models.
\newblock {\em Lett. Math. Phys.}, 112:54, 2022.
\newblock \textsf{arXiv:2107.13204 \mbox{[math.QA]}}.

\bibitem{KZ2}
Y.~Kazama and H.~Suzuki.
\newblock Characterization of {N=2} superconformal models generated by the
  coset space method.
\newblock {\em Phys. Lett. B}, 216:112--116, 1989.

\bibitem{KZ1}
Y.~Kazama and H.~Suzuki.
\newblock New {N=2} superconformal field theories and superstring
  compactification.
\newblock {\em Nucl. Phys. B}, 321:232--268, 1989.

\bibitem{KL4}
D.~Kazhdan and G.~Lusztig.
\newblock Tensor structures arising from affine {L}ie algebras. {IV}.
\newblock {\em J. Amer. Math. Soc.}, 7:383--453, 1994.

\bibitem{KoSo}
I.G. Koh and P.~Sorba.
\newblock Fusion rules and (sub-)modular invariant partition functions in
  non-unitary theories.
\newblock {\em Phys. Lett.}, B215:723--729, 1988.

\bibitem{LepRam82}
J.~Lepowsky and R.L. Wilson.
\newblock A {Lie} theoretic interpretation and proof of the {Rogers-Ramanujan}
  identities.
\newblock {\em Adv. Math.}, 45:21--72, 1982.

\bibitem{MCDelPRod}
R.~McRae.
\newblock Deligne tensor products of categories of modules for vertex operator
  algebras, 2023.
\newblock \textsf{arXiv:2304.14023 \mbox{[math.QA]}}.

\bibitem{MS}
R.~McRae and V.~Sopin.
\newblock Fusion and (non)-rigidity of {Virasoro Kac} modules in logarithmic
  minimal models at (p, q)-central charge.
\newblock {\em Physica Scripta}, 99:035233, 2024.
\newblock \textsf{arXiv:2302.08907 \mbox{[math.QA]}}.

\bibitem{MYc=25}
R.~McRae and J.~Yang.
\newblock An ${sl 2}$-type tensor category for the {V}irasoro algebra at
  central charge 25 and applications.
\newblock {\em Mathematische Zeitschrift}, 303:32, 2023.
\newblock \textsf{ arXiv:2202.07351 \mbox{[math.QA]}}.

\bibitem{Myc=1p}
R.~McRae and J.~Yang.
\newblock Structure of {V}irasoro tensor categories at central charge
  {$13-6p-6p^{-1}$} for integers {$ p> 1$}.
\newblock {\em Tans. Amer. Math. Soc.}, 378:7451--7509, 2025.
\newblock \textsf { arXiv:2011.02170 \mbox{[math.QA]}}.

\bibitem{Ridk12}
D.~Ridout.
\newblock $\widehat{\mathfrak{sl}}(2)_{-1/2}$: {A} case study.
\newblock {\em Nucl. Phys.}, B814:485--521, 2009.
\newblock \textsf{arXiv:0810.3532 \mbox{[hep-th]}}.

\bibitem{RSWosp}
D.~Ridout, J.~Snadden, and S.~Wood.
\newblock An admissible level $\widehat{\mathfrak{osp}}\left(1\middle\vert
  2\right)$-model: modular transformations and the {Verlinde} formula.
\newblock {\em Lett. Math. Phys.}, 108:2363--2423, 2018.
\newblock \textsf{arXiv:1705.04006 \mbox{[hep-th]}}.

\bibitem{RWbghosts}
D.~Ridout and S.~Wood.
\newblock Bosonic ghosts at $c = 2$ as a logarithmic {CFT}.
\newblock {\em Lett. Math. Phys.}, 105:279--307, 2014.
\newblock \textsf{arXiv:1408.4185 \mbox{[hep-th]}}.

\bibitem{RWtrip}
D.~Ridout and S.~Wood.
\newblock Modular transformations and {Verlinde} formulae for logarithmic
  $(p_{+}, p_{-})$-models.
\newblock {\em Nucl. Phys.}, B880:175--202, 2014.
\newblock \textsf{arXiv:1310.6479 \mbox{[hep-th]}}.

\bibitem{RW}
D.~Ridout and S.~Wood.
\newblock Relaxed singular vectors, {J}ack symmetric functions and fractional
  level {${\widehat{\mathfrak{sl}}(2)}$} models.
\newblock {\em Nucl. Phys. B}, 894:621--664, 2015.
\newblock \textsf{arXiv:1501.07318 \mbox{[hep-th]}}.

\bibitem{semi}
A.~M. Semikhatov.
\newblock The {MFF} singular vectors in topological conformal theories.
\newblock {\em Modern Phys. Lett. A}, 9:1867--1896, 1994.
\newblock \textsf{arXiv:hep-th/9311180}.

\bibitem{SusInt}
E.~Sussman.
\newblock The regularization of {D}otsenko–{F}ateev integrals.
\newblock {\em Lett. Math. Phys.}, 113:29, 2023.
\newblock \textsf{arXiv:2308.12900 \mbox{[math-ph]}}.

\bibitem{SusSel}
E.~Sussman.
\newblock The singularities of {S}elberg- and {D}otsenko-{F}ateev-like
  integrals.
\newblock {\em Ann. Henri Poincar{\'e}}, 25:3957--4032, 2024.
\newblock \textsf{arXiv:2301.03750 \mbox{[math-ph]}}.

\bibitem{TK}
A.~Tsuchiya and Y.~Kanie.
\newblock Fock space representations of the {V}irasoro algebra. {I}ntertwining
  operators.
\newblock {\em Publ. Res. Inst. Math. Sci.}, 22:259--327, 1986.

\bibitem{TW}
A.~Tsuchiya and S.~Wood.
\newblock On the extended {$W$}-algebra of type {${\mathfrak{sl}}_2$} at
  positive rational level.
\newblock {\em Int. Math. Res. Not.}, 2015:5357--5435, 2014.
\newblock \textsf{arXiv:1302.6435 \mbox{[math.QA]}}.

\bibitem{Z1}
Y.~Zhu.
\newblock Modular invariance of characters of vertex operator algebras.
\newblock {\em J. Amer. Math. Soc.}, 9:237--302, 1996.

\end{thebibliography}

\end{document}